\documentclass[11pt]{article}
\usepackage{amsfonts}

\usepackage{amsthm}
\usepackage{amsmath}
\usepackage{dsfont}
\usepackage{stmaryrd}
\usepackage{graphics}
\usepackage[linktocpage]{hyperref}

\long\def\ig#1{\relax}
\ig{Thanks to Roberto Minio for this def'n.  Compare the def'n of
\comment in AMSTeX.}

\newcount \coefa
\newcount \coefb
\newcount \coefc
\newcount\tempcounta
\newcount\tempcountb
\newcount\tempcountc
\newcount\tempcountd
\newcount\xext
\newcount\yext
\newcount\xoff
\newcount\yoff
\newcount\gap%
\newcount\arrowtypea
\newcount\arrowtypeb
\newcount\arrowtypec
\newcount\arrowtyped
\newcount\arrowtypee
\newcount\height
\newcount\width
\newcount\xpos
\newcount\ypos
\newcount\run
\newcount\rise
\newcount\arrowlength
\newcount\halflength
\newcount\arrowtype
\newdimen\tempdimen
\newdimen\xlen
\newdimen\ylen
\newsavebox{\tempboxa}%
\newsavebox{\tempboxb}%
\newsavebox{\tempboxc}%

\makeatletter
\def\settypes(#1,#2,#3){\arrowtypea#1 \arrowtypeb#2 \arrowtypec#3}
\def\settoheight#1#2{\setbox\@tempboxa\hbox{#2}#1\ht\@tempboxa\relax}%
\def\settodepth#1#2{\setbox\@tempboxa\hbox{#2}#1\dp\@tempboxa\relax}%
\def\settokens[#1`#2`#3`#4]{%
     \def\tokena{#1}\def\tokenb{#2}\def\tokenc{#3}\def\tokend{#4}}
\def\setsqparms[#1`#2`#3`#4;#5`#6]{%
\arrowtypea #1
\arrowtypeb #2
\arrowtypec #3
\arrowtyped #4
\width #5
\height #6
}
\def\setpos(#1,#2){\xpos=#1 \ypos#2}

\def\bfig{\begin{picture}(\xext,\yext)(\xoff,\yoff)}
\def\efig{\end{picture}}

\def\putbox(#1,#2)#3{\put(#1,#2){\makebox(0,0){$#3$}}}

\def\settriparms[#1`#2`#3;#4]{\settripairparms[#1`#2`#3`1`1;#4]}%

\def\settripairparms[#1`#2`#3`#4`#5;#6]{%
\arrowtypea #1
\arrowtypeb #2
\arrowtypec #3
\arrowtyped #4
\arrowtypee #5
\width #6
\height #6
}

\def\resetparms{\settripairparms[1`1`1`1`1;500]\width 500}

\resetparms

\def\mvector(#1,#2)#3{
\put(0,0){\vector(#1,#2){#3}}%
\put(0,0){\vector(#1,#2){30}}%
}
\def\evector(#1,#2)#3{{
\arrowlength #3
\put(0,0){\vector(#1,#2){\arrowlength}}%
\advance \arrowlength by-30
\put(0,0){\vector(#1,#2){\arrowlength}}%
}}

\def\horsize#1#2{%
\settowidth{\tempdimen}{$#2$}%
#1=\tempdimen
\divide #1 by\unitlength
}

\def\vertsize#1#2{%
\settoheight{\tempdimen}{$#2$}%
#1=\tempdimen
\settodepth{\tempdimen}{$#2$}%
\advance #1 by\tempdimen
\divide #1 by\unitlength
}

\def\vertadjust[#1`#2`#3]{%
\vertsize{\tempcounta}{#1}%
\vertsize{\tempcountb}{#2}%
\ifnum \tempcounta<\tempcountb \tempcounta=\tempcountb \fi
\divide\tempcounta by2
\vertsize{\tempcountb}{#3}%
\ifnum \tempcountb>0 \advance \tempcountb by20 \fi
\ifnum \tempcounta<\tempcountb \tempcounta=\tempcountb \fi
}

\def\horadjust[#1`#2`#3]{%
\horsize{\tempcounta}{#1}%
\horsize{\tempcountb}{#2}%
\ifnum \tempcounta<\tempcountb \tempcounta=\tempcountb \fi
\divide\tempcounta by20
\horsize{\tempcountb}{#3}%
\ifnum \tempcountb>0 \advance \tempcountb by60 \fi
\ifnum \tempcounta<\tempcountb \tempcounta=\tempcountb \fi
}

\ig{ In this procedure, #1 is the paramater that sticks out all the way,
#2 sticks out the least and #3 is a label sticking out half way.  #4 is
the amount of the offset.}

\def\sladjust[#1`#2`#3]#4{%
\tempcountc=#4
\horsize{\tempcounta}{#1}%
\divide \tempcounta by2
\horsize{\tempcountb}{#2}%
\divide \tempcountb by2
\advance \tempcountb by-\tempcountc
\ifnum \tempcounta<\tempcountb \tempcounta=\tempcountb\fi
\divide \tempcountc by2
\horsize{\tempcountb}{#3}%
\advance \tempcountb by-\tempcountc
\ifnum \tempcountb>0 \advance \tempcountb by80\fi
\ifnum \tempcounta<\tempcountb \tempcounta=\tempcountb\fi
\advance\tempcounta by20
}

\def\putvector(#1,#2)(#3,#4)#5#6{{%
\xpos=#1
\ypos=#2
\run=#3
\rise=#4
\arrowlength=#5
\arrowtype=#6
\ifnum \arrowtype<0
    \ifnum \run=0
        \advance \ypos by-\arrowlength
    \else
        \tempcounta \arrowlength
        \multiply \tempcounta by\rise
        \divide \tempcounta by\run
        \ifnum\run>0
            \advance \xpos by\arrowlength
            \advance \ypos by\tempcounta
        \else
            \advance \xpos by-\arrowlength
            \advance \ypos by-\tempcounta
        \fi
    \fi
    \multiply \arrowtype by-1
    \multiply \rise by-1
    \multiply \run by-1
\fi
\ifnum \arrowtype=1
    \put(\xpos,\ypos){\vector(\run,\rise){\arrowlength}}%
\else\ifnum \arrowtype=2
    \put(\xpos,\ypos){\mvector(\run,\rise)\arrowlength}%
\else\ifnum\arrowtype=3
    \put(\xpos,\ypos){\evector(\run,\rise){\arrowlength}}%
\fi\fi\fi
}}

\def\putsplitvector(#1,#2)#3#4{
\xpos #1
\ypos #2
\arrowtype #4
\halflength #3
\arrowlength #3
\gap 140
\advance \halflength by-\gap
\divide \halflength by2
\ifnum \arrowtype=1
    \put(\xpos,\ypos){\line(0,-1){\halflength}}%
    \advance\ypos by-\halflength
    \advance\ypos by-\gap
    \put(\xpos,\ypos){\vector(0,-1){\halflength}}%
\else\ifnum \arrowtype=2
    \put(\xpos,\ypos){\line(0,-1)\halflength}%
    \put(\xpos,\ypos){\vector(0,-1)3}%
    \advance\ypos by-\halflength
    \advance\ypos by-\gap
    \put(\xpos,\ypos){\vector(0,-1){\halflength}}%
\else\ifnum\arrowtype=3
    \put(\xpos,\ypos){\line(0,-1)\halflength}%
    \advance\ypos by-\halflength
    \advance\ypos by-\gap
    \put(\xpos,\ypos){\evector(0,-1){\halflength}}%
\else\ifnum \arrowtype=-1
    \advance \ypos by-\arrowlength
    \put(\xpos,\ypos){\line(0,1){\halflength}}%
    \advance\ypos by\halflength
    \advance\ypos by\gap
    \put(\xpos,\ypos){\vector(0,1){\halflength}}%
\else\ifnum \arrowtype=-2
    \advance \ypos by-\arrowlength
    \put(\xpos,\ypos){\line(0,1)\halflength}%
    \put(\xpos,\ypos){\vector(0,1)3}%
    \advance\ypos by\halflength
    \advance\ypos by\gap
    \put(\xpos,\ypos){\vector(0,1){\halflength}}%
\else\ifnum\arrowtype=-3
    \advance \ypos by-\arrowlength
    \put(\xpos,\ypos){\line(0,1)\halflength}%
    \advance\ypos by\halflength
    \advance\ypos by\gap
    \put(\xpos,\ypos){\evector(0,1){\halflength}}%
\fi\fi\fi\fi\fi\fi
}

\def\putmorphism(#1)(#2,#3)[#4`#5`#6]#7#8#9{{%
\run #2
\rise #3
\ifnum\rise=0
  \puthmorphism(#1)[#4`#5`#6]{#7}{#8}{#9}%
\else\ifnum\run=0
  \putvmorphism(#1)[#4`#5`#6]{#7}{#8}{#9}%
\else
\setpos(#1)%
\arrowlength #7
\arrowtype #8
\ifnum\run=0
\else\ifnum\rise=0
\else
\ifnum\run>0
    \coefa=1
\else
   \coefa=-1
\fi
\ifnum\arrowtype>0
   \coefb=0
   \coefc=-1
\else
   \coefb=\coefa
   \coefc=1
   \arrowtype=-\arrowtype
\fi
\width=2
\multiply \width by\run
\divide \width by\rise
\ifnum \width<0  \width=-\width\fi
\advance\width by60
\if l#9 \width=-\width\fi
\putbox(\xpos,\ypos){#4}
{\multiply \coefa by\arrowlength
\advance\xpos by\coefa
\multiply \coefa by\rise
\divide \coefa by\run
\advance \ypos by\coefa
\putbox(\xpos,\ypos){#5} }%
{\multiply \coefa by\arrowlength
\divide \coefa by2
\advance \xpos by\coefa
\advance \xpos by\width
\multiply \coefa by\rise
\divide \coefa by\run
\advance \ypos by\coefa
\if l#9%
   \put(\xpos,\ypos){\makebox(0,0)[r]{$#6$}}%
\else\if r#9%
   \put(\xpos,\ypos){\makebox(0,0)[l]{$#6$}}%
\fi\fi }%
{\multiply \rise by-\coefc
\multiply \run by-\coefc
\multiply \coefb by\arrowlength
\advance \xpos by\coefb
\multiply \coefb by\rise
\divide \coefb by\run
\advance \ypos by\coefb
\multiply \coefc by70
\advance \ypos by\coefc
\multiply \coefc by\run
\divide \coefc by\rise
\advance \xpos by\coefc
\multiply \coefa by140
\multiply \coefa by\run
\divide \coefa by\rise
\advance \arrowlength by\coefa
\ifnum \arrowtype=1
   \put(\xpos,\ypos){\vector(\run,\rise){\arrowlength}}%
\else\ifnum\arrowtype=2
   \put(\xpos,\ypos){\mvector(\run,\rise){\arrowlength}}%
\else\ifnum\arrowtype=3
   \put(\xpos,\ypos){\evector(\run,\rise){\arrowlength}}%
\fi\fi\fi}\fi\fi\fi\fi}}

\def\puthmorphism(#1,#2)[#3`#4`#5]#6#7#8{{%
\xpos #1
\ypos #2
\width #6
\arrowlength #6
\putbox(\xpos,\ypos){#3\vphantom{#4}}%
{\advance \xpos by\arrowlength
\putbox(\xpos,\ypos){\vphantom{#3}#4}}%
\horsize{\tempcounta}{#3}%
\horsize{\tempcountb}{#4}%
\divide \tempcounta by2
\divide \tempcountb by2
\advance \tempcounta by30
\advance \tempcountb by30
\advance \xpos by\tempcounta
\advance \arrowlength by-\tempcounta
\advance \arrowlength by-\tempcountb
\putvector(\xpos,\ypos)(1,0){\arrowlength}{#7}%
\divide \arrowlength by2
\advance \xpos by\arrowlength
\vertsize{\tempcounta}{#5}%
\divide\tempcounta by2
\advance \tempcounta by20
\if a#8 %
   \advance \ypos by\tempcounta
   \putbox(\xpos,\ypos){#5}%
\else
   \advance \ypos by-\tempcounta
   \putbox(\xpos,\ypos){#5}%
\fi}}

\def\putvmorphism(#1,#2)[#3`#4`#5]#6#7#8{{%
\xpos #1
\ypos #2
\arrowlength #6
\arrowtype #7
\settowidth{\xlen}{$#5$}%
\putbox(\xpos,\ypos){#3}%
{\advance \ypos by-\arrowlength
\putbox(\xpos,\ypos){#4}}%
{\advance\arrowlength by-140
\advance \ypos by-70
\ifdim\xlen>0pt
   \if m#8%
      \putsplitvector(\xpos,\ypos){\arrowlength}{\arrowtype}%
   \else
      \putvector(\xpos,\ypos)(0,-1){\arrowlength}{\arrowtype}%
   \fi
\else
   \putvector(\xpos,\ypos)(0,-1){\arrowlength}{\arrowtype}%
\fi}%
\ifdim\xlen>0pt
   \divide \arrowlength by2
   \advance\ypos by-\arrowlength
   \if l#8%
      \advance \xpos by-40
      \put(\xpos,\ypos){\makebox(0,0)[r]{$#5$}}%
   \else\if r#8%
      \advance \xpos by40
      \put(\xpos,\ypos){\makebox(0,0)[l]{$#5$}}%
   \else
      \putbox(\xpos,\ypos){#5}%
   \fi\fi
\fi
}}

\def\topadjust[#1`#2`#3]{%
\yoff=10
\vertadjust[#1`#2`{#3}]%
\advance \yext by\tempcounta
\advance \yext by 10
}
\def\botadjust[#1`#2`#3]{%
\vertadjust[#1`#2`{#3}]%
\advance \yext by\tempcounta
\advance \yoff by-\tempcounta
}
\def\leftadjust[#1`#2`#3]{%
\xoff=0
\horadjust[#1`#2`{#3}]%
\advance \xext by\tempcounta
\advance \xoff by-\tempcounta
}
\def\rightadjust[#1`#2`#3]{%
\horadjust[#1`#2`{#3}]%
\advance \xext by\tempcounta
}
\def\rightsladjust[#1`#2`#3]{%
\sladjust[#1`#2`{#3}]{\width}%
\advance \xext by\tempcounta
}
\def\leftsladjust[#1`#2`#3]{%
\xoff=0
\sladjust[#1`#2`{#3}]{\width}%
\advance \xext by\tempcounta
\advance \xoff by-\tempcounta
}
\def\adjust[#1`#2;#3`#4;#5`#6;#7`#8]{%
\topadjust[#1``{#2}]
\leftadjust[#3``{#4}]
\rightadjust[#5``{#6}]
\botadjust[#7``{#8}]}

\def\putsquarep<#1>(#2)[#3;#4`#5`#6`#7]{{%
\setsqparms[#1]%
\setpos(#2)%
\settokens[#3]%
\puthmorphism(\xpos,\ypos)[\tokenc`\tokend`{#7}]{\width}{\arrowtyped}b%
\advance\ypos by \height
\puthmorphism(\xpos,\ypos)[\tokena`\tokenb`{#4}]{\width}{\arrowtypea}a%
\putvmorphism(\xpos,\ypos)[``{#5}]{\height}{\arrowtypeb}l%
\advance\xpos by \width
\putvmorphism(\xpos,\ypos)[``{#6}]{\height}{\arrowtypec}r%
}}

\def\putsquare{\@ifnextchar <{\putsquarep}{\putsquarep%
   <\arrowtypea`\arrowtypeb`\arrowtypec`\arrowtyped;\width`\height>}}
\def\square{\@ifnextchar< {\squarep}{\squarep
   <\arrowtypea`\arrowtypeb`\arrowtypec`\arrowtyped;\width`\height>}}
\def\squarep<#1>[#2`#3`#4`#5;#6`#7`#8`#9]{{
\setsqparms[#1]
\xext=\width                                          
\yext=\height                                         
\topadjust[#2`#3`{#6}]
\botadjust[#4`#5`{#9}]
\leftadjust[#2`#4`{#7}]
\rightadjust[#3`#5`{#8}]
\begin{picture}(\xext,\yext)(\xoff,\yoff)
\putsquarep<\arrowtypea`\arrowtypeb`\arrowtypec`\arrowtyped;\width`\height>%
(0,0)[#2`#3`#4`#5;#6`#7`#8`{#9}]%
\end{picture}%
}}

\def\putptrianglep<#1>(#2,#3)[#4`#5`#6;#7`#8`#9]{{%
\settriparms[#1]%
\xpos=#2 \ypos=#3
\advance\ypos by \height
\puthmorphism(\xpos,\ypos)[#4`#5`{#7}]{\height}{\arrowtypea}a%
\putvmorphism(\xpos,\ypos)[`#6`{#8}]{\height}{\arrowtypeb}l%
\advance\xpos by\height
\putmorphism(\xpos,\ypos)(-1,-1)[``{#9}]{\height}{\arrowtypec}r%
}}

\def\putptriangle{\@ifnextchar <{\putptrianglep}{\putptrianglep
   <\arrowtypea`\arrowtypeb`\arrowtypec;\height>}}
\def\ptriangle{\@ifnextchar <{\ptrianglep}{\ptrianglep
   <\arrowtypea`\arrowtypeb`\arrowtypec;\height>}}

\def\ptrianglep<#1>[#2`#3`#4;#5`#6`#7]{{
\settriparms[#1]%
\width=\height                         
\xext=\width                           
\yext=\width                           
\topadjust[#2`#3`{#5}]
\botadjust[#3``]
\leftadjust[#2`#4`{#6}]
\rightsladjust[#3`#4`{#7}]
\begin{picture}(\xext,\yext)(\xoff,\yoff)
\putptrianglep<\arrowtypea`\arrowtypeb`\arrowtypec;\height>%
(0,0)[#2`#3`#4;#5`#6`{#7}]%
\end{picture}%
}}

\def\putqtrianglep<#1>(#2,#3)[#4`#5`#6;#7`#8`#9]{{%
\settriparms[#1]%
\xpos=#2 \ypos=#3
\advance\ypos by\height
\puthmorphism(\xpos,\ypos)[#4`#5`{#7}]{\height}{\arrowtypea}a%
\putmorphism(\xpos,\ypos)(1,-1)[``{#8}]{\height}{\arrowtypeb}l%
\advance\xpos by\height
\putvmorphism(\xpos,\ypos)[`#6`{#9}]{\height}{\arrowtypec}r%
}}

\def\putqtriangle{\@ifnextchar <{\putqtrianglep}{\putqtrianglep
   <\arrowtypea`\arrowtypeb`\arrowtypec;\height>}}
\def\qtriangle{\@ifnextchar <{\qtrianglep}{\qtrianglep
   <\arrowtypea`\arrowtypeb`\arrowtypec;\height>}}

\def\qtrianglep<#1>[#2`#3`#4;#5`#6`#7]{{
\settriparms[#1]
\width=\height                         
\xext=\width                           
\yext=\height                          
\topadjust[#2`#3`{#5}]
\botadjust[#4``]
\leftsladjust[#2`#4`{#6}]
\rightadjust[#3`#4`{#7}]
\begin{picture}(\xext,\yext)(\xoff,\yoff)
\putqtrianglep<\arrowtypea`\arrowtypeb`\arrowtypec;\height>%
(0,0)[#2`#3`#4;#5`#6`{#7}]%
\end{picture}%
}}

\def\putdtrianglep<#1>(#2,#3)[#4`#5`#6;#7`#8`#9]{{%
\settriparms[#1]%
\xpos=#2 \ypos=#3
\puthmorphism(\xpos,\ypos)[#5`#6`{#9}]{\height}{\arrowtypec}b%
\advance\xpos by \height \advance\ypos by\height
\putmorphism(\xpos,\ypos)(-1,-1)[``{#7}]{\height}{\arrowtypea}l%
\putvmorphism(\xpos,\ypos)[#4``{#8}]{\height}{\arrowtypeb}r%
}}

\def\putdtriangle{\@ifnextchar <{\putdtrianglep}{\putdtrianglep
   <\arrowtypea`\arrowtypeb`\arrowtypec;\height>}}
\def\dtriangle{\@ifnextchar <{\dtrianglep}{\dtrianglep
   <\arrowtypea`\arrowtypeb`\arrowtypec;\height>}}

\def\dtrianglep<#1>[#2`#3`#4;#5`#6`#7]{{
\settriparms[#1]
\width=\height                         
\xext=\width                           
\yext=\height                          
\topadjust[#2``]
\botadjust[#3`#4`{#7}]
\leftsladjust[#3`#2`{#5}]
\rightadjust[#2`#4`{#6}]
\begin{picture}(\xext,\yext)(\xoff,\yoff)
\putdtrianglep<\arrowtypea`\arrowtypeb`\arrowtypec;\height>%
(0,0)[#2`#3`#4;#5`#6`{#7}]%
\end{picture}%
}}

\def\putbtrianglep<#1>(#2,#3)[#4`#5`#6;#7`#8`#9]{{%
\settriparms[#1]%
\xpos=#2 \ypos=#3
\puthmorphism(\xpos,\ypos)[#5`#6`{#9}]{\height}{\arrowtypec}b%
\advance\ypos by\height
\putmorphism(\xpos,\ypos)(1,-1)[``{#8}]{\height}{\arrowtypeb}r%
\putvmorphism(\xpos,\ypos)[#4``{#7}]{\height}{\arrowtypea}l%
}}

\def\putbtriangle{\@ifnextchar <{\putbtrianglep}{\putbtrianglep
   <\arrowtypea`\arrowtypeb`\arrowtypec;\height>}}
\def\btriangle{\@ifnextchar <{\btrianglep}{\btrianglep
   <\arrowtypea`\arrowtypeb`\arrowtypec;\height>}}

\def\btrianglep<#1>[#2`#3`#4;#5`#6`#7]{{
\settriparms[#1]
\width=\height                         
\xext=\width                           
\yext=\height                          
\topadjust[#2``]
\botadjust[#3`#4`{#7}]
\leftadjust[#2`#3`{#5}]
\rightsladjust[#4`#2`{#6}]
\begin{picture}(\xext,\yext)(\xoff,\yoff)
\putbtrianglep<\arrowtypea`\arrowtypeb`\arrowtypec;\height>%
(0,0)[#2`#3`#4;#5`#6`{#7}]%
\end{picture}%
}}

\def\putAtrianglep<#1>(#2,#3)[#4`#5`#6;#7`#8`#9]{{%
\settriparms[#1]%
\xpos=#2 \ypos=#3
{\multiply \height by2
\puthmorphism(\xpos,\ypos)[#5`#6`{#9}]{\height}{\arrowtypec}b}%
\advance\xpos by\height \advance\ypos by\height
\putmorphism(\xpos,\ypos)(-1,-1)[#4``{#7}]{\height}{\arrowtypea}l%
\putmorphism(\xpos,\ypos)(1,-1)[``{#8}]{\height}{\arrowtypeb}r%
}}

\def\putAtriangle{\@ifnextchar <{\putAtrianglep}{\putAtrianglep
   <\arrowtypea`\arrowtypeb`\arrowtypec;\height>}}
\def\Atriangle{\@ifnextchar <{\Atrianglep}{\Atrianglep
   <\arrowtypea`\arrowtypeb`\arrowtypec;\height>}}

\def\Atrianglep<#1>[#2`#3`#4;#5`#6`#7]{{
\settriparms[#1]
\width=\height                         
\xext=\width                           
\yext=\height                          
\topadjust[#2``]
\botadjust[#3`#4`{#7}]
\multiply \xext by2 
\leftsladjust[#3`#2`{#5}]
\rightsladjust[#4`#2`{#6}]
\begin{picture}(\xext,\yext)(\xoff,\yoff)%
\putAtrianglep<\arrowtypea`\arrowtypeb`\arrowtypec;\height>%
(0,0)[#2`#3`#4;#5`#6`{#7}]%
\end{picture}%
}}

\def\putAtrianglepairp<#1>(#2)[#3;#4`#5`#6`#7`#8]{{
\settripairparms[#1]%
\setpos(#2)%
\settokens[#3]%
\puthmorphism(\xpos,\ypos)[\tokenb`\tokenc`{#7}]{\height}{\arrowtyped}b%
\advance\xpos by\height
\advance\ypos by\height
\putmorphism(\xpos,\ypos)(-1,-1)[\tokena``{#4}]{\height}{\arrowtypea}l%
\putvmorphism(\xpos,\ypos)[``{#5}]{\height}{\arrowtypeb}m%
\putmorphism(\xpos,\ypos)(1,-1)[``{#6}]{\height}{\arrowtypec}r%
}}

\def\putAtrianglepair{\@ifnextchar <{\putAtrianglepairp}{\putAtrianglepairp%
   <\arrowtypea`\arrowtypeb`\arrowtypec`\arrowtyped`\arrowtypee;\height>}}
\def\Atrianglepair{\@ifnextchar <{\Atrianglepairp}{\Atrianglepairp%
   <\arrowtypea`\arrowtypeb`\arrowtypec`\arrowtyped`\arrowtypee;\height>}}

\def\Atrianglepairp<#1>[#2;#3`#4`#5`#6`#7]{{%
\settripairparms[#1]%
\settokens[#2]%
\width=\height
\xext=\width
\yext=\height
\topadjust[\tokena``]%
\vertadjust[\tokenb`\tokenc`{#6}]
\tempcountd=\tempcounta                       
\vertadjust[\tokenc`\tokend`{#7}]
\ifnum\tempcounta<\tempcountd                 
\tempcounta=\tempcountd\fi                    
\advance \yext by\tempcounta                  
\advance \yoff by-\tempcounta                 %
\multiply \xext by2 
\leftsladjust[\tokenb`\tokena`{#3}]
\rightsladjust[\tokend`\tokena`{#5}]%
\begin{picture}(\xext,\yext)(\xoff,\yoff)%
\putAtrianglepairp
<\arrowtypea`\arrowtypeb`\arrowtypec`\arrowtyped`\arrowtypee;\height>%
(0,0)[#2;#3`#4`#5`#6`{#7}]%
\end{picture}%
}}

\def\putVtrianglep<#1>(#2,#3)[#4`#5`#6;#7`#8`#9]{{%
\settriparms[#1]%
\xpos=#2 \ypos=#3
\advance\ypos by\height
{\multiply\height by2
\puthmorphism(\xpos,\ypos)[#4`#5`{#7}]{\height}{\arrowtypea}a}%
\putmorphism(\xpos,\ypos)(1,-1)[`#6`{#8}]{\height}{\arrowtypeb}l%
\advance\xpos by\height
\advance\xpos by\height
\putmorphism(\xpos,\ypos)(-1,-1)[``{#9}]{\height}{\arrowtypec}r%
}}

\def\putVtriangle{\@ifnextchar <{\putVtrianglep}{\putVtrianglep
   <\arrowtypea`\arrowtypeb`\arrowtypec;\height>}}
\def\Vtriangle{\@ifnextchar <{\Vtrianglep}{\Vtrianglep
   <\arrowtypea`\arrowtypeb`\arrowtypec;\height>}}

\def\Vtrianglep<#1>[#2`#3`#4;#5`#6`#7]{{
\settriparms[#1]
\width=\height                         
\xext=\width                           
\yext=\height                          
\topadjust[#2`#3`{#5}]
\botadjust[#4``]
\multiply \xext by2 
\leftsladjust[#2`#3`{#6}]
\rightsladjust[#3`#4`{#7}]
\begin{picture}(\xext,\yext)(\xoff,\yoff)%
\putVtrianglep<\arrowtypea`\arrowtypeb`\arrowtypec;\height>%
(0,0)[#2`#3`#4;#5`#6`{#7}]%
\end{picture}%
}}

\def\putVtrianglepairp<#1>(#2)[#3;#4`#5`#6`#7`#8]{{
\settripairparms[#1]%
\setpos(#2)%
\settokens[#3]%
\advance\ypos by\height
\putmorphism(\xpos,\ypos)(1,-1)[`\tokend`{#6}]{\height}{\arrowtypec}l%
\puthmorphism(\xpos,\ypos)[\tokena`\tokenb`{#4}]{\height}{\arrowtypea}a%
\advance\xpos by\height
\putvmorphism(\xpos,\ypos)[``{#7}]{\height}{\arrowtyped}m%
\advance\xpos by\height
\putmorphism(\xpos,\ypos)(-1,-1)[``{#8}]{\height}{\arrowtypee}r%
}}

\def\putVtrianglepair{\@ifnextchar <{\putVtrianglepairp}{\putVtrianglepairp%
    <\arrowtypea`\arrowtypeb`\arrowtypec`\arrowtyped`\arrowtypee;\height>}}
\def\Vtrianglepair{\@ifnextchar <{\Vtrianglepairp}{\Vtrianglepairp%
    <\arrowtypea`\arrowtypeb`\arrowtypec`\arrowtyped`\arrowtypee;\height>}}

\def\Vtrianglepairp<#1>[#2;#3`#4`#5`#6`#7]{{%
\settripairparms[#1]%
\settokens[#2]
\xext=\height                  
\width=\height                 
\yext=\height                  
\vertadjust[\tokena`\tokenb`{#4}]
\tempcountd=\tempcounta        
\vertadjust[\tokenb`\tokenc`{#5}]
\ifnum\tempcounta<\tempcountd%
\tempcounta=\tempcountd\fi
\advance \yext by\tempcounta
\botadjust[\tokend``]%
\multiply \xext by2
\leftsladjust[\tokena`\tokend`{#6}]%
\rightsladjust[\tokenc`\tokend`{#7}]%
\begin{picture}(\xext,\yext)(\xoff,\yoff)%
\putVtrianglepairp
<\arrowtypea`\arrowtypeb`\arrowtypec`\arrowtyped`\arrowtypee;\height>%
(0,0)[#2;#3`#4`#5`#6`{#7}]%
\end{picture}%
}}

\def\putCtrianglep<#1>(#2,#3)[#4`#5`#6;#7`#8`#9]{{%
\settriparms[#1]%
\xpos=#2 \ypos=#3
\advance\ypos by\height
\putmorphism(\xpos,\ypos)(1,-1)[``{#9}]{\height}{\arrowtypec}l%
\advance\xpos by\height
\advance\ypos by\height
\putmorphism(\xpos,\ypos)(-1,-1)[#4`#5`{#7}]{\height}{\arrowtypea}l%
{\multiply\height by 2
\putvmorphism(\xpos,\ypos)[`#6`{#8}]{\height}{\arrowtypeb}r}%
}}

\def\putCtriangle{\@ifnextchar <{\putCtrianglep}{\putCtrianglep
    <\arrowtypea`\arrowtypeb`\arrowtypec;\height>}}
\def\Ctriangle{\@ifnextchar <{\Ctrianglep}{\Ctrianglep
    <\arrowtypea`\arrowtypeb`\arrowtypec;\height>}}

\def\Ctrianglep<#1>[#2`#3`#4;#5`#6`#7]{{
\settriparms[#1]
\width=\height                          
\xext=\width                            
\yext=\height                           
\multiply \yext by2 
\topadjust[#2``]
\botadjust[#4``]
\sladjust[#3`#2`{#5}]{\width}
\tempcountd=\tempcounta                 
\sladjust[#3`#4`{#7}]{\width}
\ifnum \tempcounta<\tempcountd          
\tempcounta=\tempcountd\fi              
\advance \xext by\tempcounta            
\advance \xoff by-\tempcounta           %
\rightadjust[#2`#4`{#6}]
\begin{picture}(\xext,\yext)(\xoff,\yoff)%
\putCtrianglep<\arrowtypea`\arrowtypeb`\arrowtypec;\height>%
(0,0)[#2`#3`#4;#5`#6`{#7}]%
\end{picture}%
}}

\def\putDtrianglep<#1>(#2,#3)[#4`#5`#6;#7`#8`#9]{{%
\settriparms[#1]%
\xpos=#2 \ypos=#3
\advance\xpos by\height \advance\ypos by\height
\putmorphism(\xpos,\ypos)(-1,-1)[``{#9}]{\height}{\arrowtypec}r%
\advance\xpos by-\height \advance\ypos by\height
\putmorphism(\xpos,\ypos)(1,-1)[`#5`{#8}]{\height}{\arrowtypeb}r%
{\multiply\height by 2
\putvmorphism(\xpos,\ypos)[#4`#6`{#7}]{\height}{\arrowtypea}l}%
}}

\def\putDtriangle{\@ifnextchar <{\putDtrianglep}{\putDtrianglep
    <\arrowtypea`\arrowtypeb`\arrowtypec;\height>}}
\def\Dtriangle{\@ifnextchar <{\Dtrianglep}{\Dtrianglep
   <\arrowtypea`\arrowtypeb`\arrowtypec;\height>}}

\def\Dtrianglep<#1>[#2`#3`#4;#5`#6`#7]{{
\settriparms[#1]
\width=\height                         
\xext=\height                          
\yext=\height                          
\multiply \yext by2 
\topadjust[#2``]
\botadjust[#4``]
\leftadjust[#2`#4`{#5}]
\sladjust[#3`#2`{#5}]{\height}
\tempcountd=\tempcountd                
\sladjust[#3`#4`{#7}]{\height}
\ifnum \tempcounta<\tempcountd         
\tempcounta=\tempcountd\fi             
\advance \xext by\tempcounta           %
\begin{picture}(\xext,\yext)(\xoff,\yoff)
\putDtrianglep<\arrowtypea`\arrowtypeb`\arrowtypec;\height>%
(0,0)[#2`#3`#4;#5`#6`{#7}]%
\end{picture}%
}}

\def\setrecparms[#1`#2]{\width=#1 \height=#2}%
\def\recursep<#1`#2>[#3;#4`#5`#6`#7`#8]{{%
\width=#1 \height=#2
\settokens[#3]
\settowidth{\tempdimen}{$\tokena$}
\ifdim\tempdimen=0pt
  \savebox{\tempboxa}{\hbox{$\tokenb$}}%
  \savebox{\tempboxb}{\hbox{$\tokend$}}%
  \savebox{\tempboxc}{\hbox{$#6$}}%
\else
  \savebox{\tempboxa}{\hbox{$\hbox{$\tokena$}\times\hbox{$\tokenb$}$}}%
  \savebox{\tempboxb}{\hbox{$\hbox{$\tokena$}\times\hbox{$\tokend$}$}}%
  \savebox{\tempboxc}{\hbox{$\hbox{$\tokena$}\times\hbox{$#6$}$}}%
\fi
\ypos=\height
\divide\ypos by 2
\xpos=\ypos
\advance\xpos by \width
\xext=\xpos \yext=\height
\topadjust[#3`\usebox{\tempboxa}`{#4}]%
\botadjust[#5`\usebox{\tempboxb}`{#8}]%
\sladjust[\tokenc`\tokenb`{#5}]{\ypos}%
\tempcountd=\tempcounta
\sladjust[\tokenc`\tokend`{#5}]{\ypos}%
\ifnum \tempcounta<\tempcountd
\tempcounta=\tempcountd\fi
\advance \xext by\tempcounta
\advance \xoff by-\tempcounta
\rightadjust[\usebox{\tempboxa}`\usebox{\tempboxb}`\usebox{\tempboxc}]%
\bfig
\putCtrianglep<-1`1`1;\ypos>(0,0)[`\tokenc`;#5`#6`{#7}]%
\puthmorphism(\ypos,0)[\tokend`\usebox{\tempboxb}`{#8}]{\width}{-1}b%
\puthmorphism(\ypos,\height)[\tokenb`\usebox{\tempboxa}`{#4}]{\width}{-1}a%
\advance\ypos by \width
\putvmorphism(\ypos,\height)[``\usebox{\tempboxc}]{\height}1r%
\efig
}}

\def\recurse{\@ifnextchar <{\recursep}{\recursep<\width`\height>}}

\def\puttwohmorphisms(#1,#2)[#3`#4;#5`#6]#7#8#9{{%
\puthmorphism(#1,#2)[#3`#4`]{#7}0a
\ypos=#2
\advance\ypos by 20
\puthmorphism(#1,\ypos)[\phantom{#3}`\phantom{#4}`#5]{#7}{#8}a
\advance\ypos by -40
\puthmorphism(#1,\ypos)[\phantom{#3}`\phantom{#4}`#6]{#7}{#9}b
}}

\def\puttwovmorphisms(#1,#2)[#3`#4;#5`#6]#7#8#9{{%
\putvmorphism(#1,#2)[#3`#4`]{#7}0a
\xpos=#1
\advance\xpos by -20
\putvmorphism(\xpos,#2)[\phantom{#3}`\phantom{#4}`#5]{#7}{#8}l
\advance\xpos by 40
\putvmorphism(\xpos,#2)[\phantom{#3}`\phantom{#4}`#6]{#7}{#9}r
}}

\def\puthcoequalizer(#1)[#2`#3`#4;#5`#6`#7]#8#9{{%
\setpos(#1)%
\puttwohmorphisms(\xpos,\ypos)[#2`#3;#5`#6]{#8}11%
\advance\xpos by #8
\puthmorphism(\xpos,\ypos)[\phantom{#3}`#4`#7]{#8}1{#9}
}}

\def\putvcoequalizer(#1)[#2`#3`#4;#5`#6`#7]#8#9{{%
\setpos(#1)%
\puttwovmorphisms(\xpos,\ypos)[#2`#3;#5`#6]{#8}11%
\advance\ypos by -#8
\putvmorphism(\xpos,\ypos)[\phantom{#3}`#4`#7]{#8}1{#9}
}}

\def\putthreehmorphisms(#1)[#2`#3;#4`#5`#6]#7(#8)#9{{%
\setpos(#1) \settypes(#8)
\if a#9 %
     \vertsize{\tempcounta}{#5}%
     \vertsize{\tempcountb}{#6}%
     \ifnum \tempcounta<\tempcountb \tempcounta=\tempcountb \fi
\else
     \vertsize{\tempcounta}{#4}%
     \vertsize{\tempcountb}{#5}%
     \ifnum \tempcounta<\tempcountb \tempcounta=\tempcountb \fi
\fi
\advance \tempcounta by 60
\puthmorphism(\xpos,\ypos)[#2`#3`#5]{#7}{\arrowtypeb}{#9}
\advance\ypos by \tempcounta
\puthmorphism(\xpos,\ypos)[\phantom{#2}`\phantom{#3}`#4]{#7}{\arrowtypea}{#9}
\advance\ypos by -\tempcounta \advance\ypos by -\tempcounta
\puthmorphism(\xpos,\ypos)[\phantom{#2}`\phantom{#3}`#6]{#7}{\arrowtypec}{#9}
}}

\def\putarc(#1,#2)[#3`#4`#5]#6#7#8{{%
\xpos #1
\ypos #2
\width #6
\arrowlength #6
\putbox(\xpos,\ypos){#3\vphantom{#4}}%
{\advance \xpos by\arrowlength
\putbox(\xpos,\ypos){\vphantom{#3}#4}}%
\horsize{\tempcounta}{#3}%
\horsize{\tempcountb}{#4}%
\divide \tempcounta by2
\divide \tempcountb by2
\advance \tempcounta by30
\advance \tempcountb by30
\advance \xpos by\tempcounta
\advance \arrowlength by-\tempcounta
\advance \arrowlength by-\tempcountb
\halflength=\arrowlength \divide\halflength by 2
\divide\arrowlength by 5
\put(\xpos,\ypos){\bezier{\arrowlength}(0,0)(50,50)(\halflength,50)}
\ifnum #7=-1 \put(\xpos,\ypos){\vector(-3,-2)0} \fi
\advance\xpos by \halflength
\put(\xpos,\ypos){\xpos=\halflength \advance\xpos by -50
   \bezier{\arrowlength}(0,50)(\xpos,50)(\halflength,0)}
\ifnum #7=1 {\advance \xpos by
   \halflength \put(\xpos,\ypos){\vector(3,-2)0}} \fi
\advance\ypos by 50
\vertsize{\tempcounta}{#5}%
\divide\tempcounta by2
\advance \tempcounta by20
\if a#8 %
   \advance \ypos by\tempcounta
   \putbox(\xpos,\ypos){#5}%
\else
   \advance \ypos by-\tempcounta
   \putbox(\xpos,\ypos){#5}%
\fi
}}

\makeatother
\newtheorem{theorem}{Theorem}[section]
\newtheorem{lemma}[theorem]{Lemma}
\newtheorem{corollary}[theorem]{Corollary}

\newtheorem{proposition}[theorem]{Proposition}

\makeindex \makeglossary
\begin{document}

\sloppy

\newcommand{\nl}{\hspace{2cm}\\ }

\def\nec{\Box}
\def\pos{\Diamond}
\def\diam{{\tiny\Diamond}}

\def\lc{\lceil}
\def\rc{\rceil}
\def\lf{\lfloor}
\def\rf{\rfloor}
\def\lk{\langle}
\def\rk{\rangle}
\def\blk{\dot{\langle\!\!\langle}}
\def\brk{\dot{\rangle\!\!\rangle}}

\newcommand{\pa}{\parallel}
\newcommand{\lra}{\longrightarrow}
\newcommand{\hra}{\hookrightarrow}
\newcommand{\hla}{\hookleftarrow}
\newcommand{\ra}{\rightarrow}
\newcommand{\la}{\leftarrow}
\newcommand{\lla}{\longleftarrow}
\newcommand{\da}{\downarrow}
\newcommand{\ua}{\uparrow}
\newcommand{\dA}{\downarrow\!\!\!^\bullet}
\newcommand{\uA}{\uparrow\!\!\!_\bullet}
\newcommand{\Da}{\Downarrow}
\newcommand{\DA}{\Downarrow\!\!\!^\bullet}
\newcommand{\UA}{\Uparrow\!\!\!_\bullet}
\newcommand{\Ua}{\Uparrow}
\newcommand{\Lra}{\Longrightarrow}
\newcommand{\Ra}{\Rightarrow}
\newcommand{\Lla}{\Longleftarrow}
\newcommand{\La}{\Leftarrow}
\newcommand{\nperp}{\perp\!\!\!\!\!\setminus\;\;}
\newcommand{\pq}{\preceq}

\newcommand{\lms}{\longmapsto}
\newcommand{\ms}{\mapsto}
\newcommand{\subseteqnot}{\subseteq\hskip-4 mm_\not\hskip3 mm}

\def\o{{\omega}}
\def\sM{{\bf sM}}
\def\bA{{\bf A}}
\def\bEM{{\bf EM}}
\def\bM{{\bf M}}
\def\bN{{\bf N}}
\def\bC{{\bf C}}
\def\bI{{\bf I}}
\def\bK{{\bf K}}
\def\bL{{\bf L}}
\def\bT{{\bf T}}
\def\bS{{\bf S}}
\def\bD{{\bf D}}
\def\bB{{\bf B}}
\def\bW{{\bf W}}
\def\bP{{\bf P}}
\def\bX{{\bf X}}
\def\bU{{\bf U}}
\def\bY{{\bf Y}}
\def\ba{{\bf a}}
\def\bb{{\bf b}}
\def\bc{{\bf c}}
\def\bd{{\bf d}}
\def\bh{{\bf h}}
\def\bi{{\bf i}}
\def\bj{{\bf j}}
\def\bk{{\bf k}}
\def\bm{{\bf m}}
\def\bn{{\bf n}}
\def\bp{{\bf p}}
\def\bq{{\bf q}}
\def\be{{\bf e}}
\def\br{{\bf r}}
\def\bi{{\bf i}}
\def\bs{{\bf s}}
\def\bt{{\bf t}}
\def\bu{{\bf u}}
\def\bv{{\bf v}}
\def\bw{{\bf w}}
\def\bz{{\bf z}}

\def\jeden{{\bf 1}}
\def\dwa{{\bf 2}}
\def\trzy{{\bf 3}}
\def\Lam{{\bf \Lambda}}

\def\cBL{{\cal BL}}
\def\cB{{\cal B}}
\def\cA{{\cal A}}
\def\cC{{\cal C}}
\def\cD{{\cal D}}
\def\cE{{\cal E}}
\def\cEM{{\cal EM}}
\def\cF{{\cal F}}
\def\cG{{\cal G}}
\def\cI{{\cal I}}
\def\cJ{{\cal J}}
\def\cK{{\cal K}}
\def\cL{{\cal L}}
\def\cN{{\cal N}}
\def\cM{{\cal M}}
\def\cO{{\cal O}}
\def\cP{{\cal P}}
\def\cQ{{\cal Q}}
\def\cR{{\cal R}}
\def\cS{{\cal S}}
\def\cT{{\cal T}}
\def\cU{{\cal U}}
\def\cV{{\cal V}}
\def\cW{{\cal W}}
\def\cX{{\cal X}}
\def\cY{{\cal Y}}

\def\LMF{{\bf LMF}}
\def\Mon{{\bf Mon}}
\def\LAdj{{\bf LAdj}}
\def\Adj{{\bf Adj}}
\def\RAdj{{\bf RAdj}}
\def\Act{{\bf Act}}
\def\Clsd{{\bf Closed}}
\def\Mlcv{{\bf Malcev}}
\def\cMlcv{{\bf coMalcev}}
\def\Mnd{{\bf Mnd}}
\def\bCat{{{\bf Cat}}}
\def\Cat{{{\bf Cat}}}
\def\Gpd{{{\bf Gpd}}}
\def\Gph{{{\bf Gph}}}
\def\Fib{{{\bf Fib}}}
\def\DFib{{{\bf DFib}}}
\def\BiFib{{{\bf BiFib}}}
\def\Catrc{{{\bf Cat}_{rc}}}
\def\Monrc{{{\bf Mon}_{rc}}}
\def\Mod{{\bf Mod}}
\def\Modm{{\bf cMod}}
\def\cMod{{\bf cMod}}
\def\mon{{{\bf mon}}}
\def\act{{{\bf act}}}
\def\bem{{{\bf em}}}
\def\bkem{{{\bf kem}}}
\def\bKl{{{\bf Kl}}}

\def\oG{{{\omega}Gr}}
\def\mts{{MltSet}}


\def\oCat{{\bf \o Cat}}
\def\oGph{{\bf \o Gph}}
\def\AMon{{\bf AnMnd}}
\def\An{{\bf An}}
\def\Poly{{\bf Poly}}
\def\San{{\bf San}}
\def\Taut{{\bf Taut}}
\def\PMnd{{\bf PolyMnd}}
\def\SanMnd{{\bf SanMnd}}
\def\RiMnd{{\bf RiMnd}}
\def\End{{\bf End}}

\def\ET{\bf ET}
\def\RegET{\bf RegET}
\def\RET{\bf RegET}
\def\LrET{\bf LrET}
\def\RiET{\bf RiET}
\def\SregET{\bf SregET}
\def\Cart{\bf Cart}
\def\wCart{\bf wCart}
\def\CartMnd{\bf CartMnd}
\def\wCartMnd{\bf wCartMnd}

\def\LT{\bf LT}
\def\RegLT{\bf RegLT}
\def\ALT{\bf AnLT}
\def\RiLT{\bf RiLT}

\def\FOp{\bf FOp}
\def\RegOp{\bf RegOp}
\def\SOp{\bf SOp}
\def\RiOp{\bf RiOp}

\def\MonCat{{{\bf MonCat}}}
\def\ActMonCat{{{\bf ActMonCat}}}

\def\F{\mathds{F}}
\def\E{\mathds{E}}
\def\S{\mathds{S}}
\def\I{\mathds{I}}
\def\B{\mathds{B}}

\def\V{\mathds{V}}
\def\W{\mathds{W}}
\def\M{\mathds{M}}
\def\N{\mathds{N}}
\def\R{\mathds{R}}

\pagenumbering{arabic} \setcounter{page}{1}

\title{\bf\Large Polynomial and Analytic Functors and Monads, revisited}

\author{ Stanis\l aw Szawiel,\\ Marek Zawadowski\\
Instytut Matematyki, Uniwersytet Warszawski\\
}

\maketitle

\begin{abstract} We describe an abstract 2-categorical setting to study various notions of polynomial and analytic functors and monads.  \end{abstract}
\section{Introduction}

The analytic and polynomial functors and monads constitute a useful tool in combinatorics, geometry, topology, and logic, to name some areas of their application, c.f. \cite{J1}, \cite{J2}, \cite{Ke}, \cite{CJ}, \cite{BD}, \cite{HMP}, \cite{H}, \cite{AV}, \cite{FGHW}, \cite{Fi2}, \cite{KJBM}, \cite{Z2}, \cite{SZ3}. The main reason we became interested in using them and eventually to studying them on their own was the development of a convenient algebraic definition of the category of opetopic of sets. They proved very useful in the development of the intricate combinatorics of the opetopic sets. However, upon reflection it turned out that the theory of these tools could be itself further developed to explain better how they are related and what can be expected from them. This was done in \cite{Sz} to compare various algebraic definitions of opetopic sets, cf. \cite{BD}, \cite{HMP}, \cite{Z2}, \cite{KJBM}, \cite{SZ2}.  The polynomial functors are related to Kleisli algebras, whereas analytic functors are related to Eilenberg-Moore algebras. This is why, even if polynomial functors are simpler and easier to handle in many contexts, they should not be expected to be closed under limits and colimits. By contrast, even if analytic functors are usually defined in a less transparent way, they have many desirable closure properties.

This paper can be considered as an extension of a part of \cite{Z2}, done in an abstract way. Our objective is to study the tools that already proved  useful.

To explain the idea on a `toy model' of endofunctors on the category of sets $Set$, consider the following. The category of (untyped) signatures $Sig=Set_{/\o}$ is the slice category of $Set$ over the set of natural numbers. It has a (non-symmetric) substitution tensor such that the monoids for this tensor are non-$\Sigma$-operads on $Set$. As the category $Set$ can be identified with the subcategory of $Sig$ of signatures of constants, $Sig$ acts on $Set$. The action
$$\star : Sig\times Set \lra Set$$
\[ \lk \{ A_n\}_n,X\rk\;\mapsto \sum_n\; A_n\times X^n\]
has an exponential adjoint representation
$$\br: Sig \ra End(Set)$$
which is a strong monoidal functor ($End(Set)$, the category of endofunctors on $Set$ is a strict monoidal category with the tensor given by composition). The functor $\br$ has a (lax monoidal) right adjoint, say $U$.  Note that $End(Set)$ has and  $U$ preserves reflexive coequalizers. Thus the resulting monad $\cF=U\br$ is lax monoidal and preserves reflexive coequalizers. The Kleisli category $Sig_\cF$ and the Eilenberg-Moore category $Sig^\cF$ fully embed into  $End(Set)$
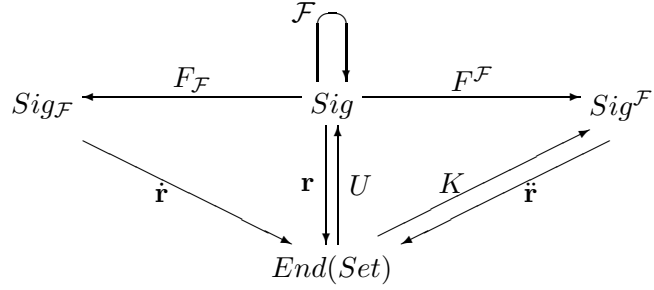
\begin{center} \xext=2000 \yext=950
\begin{picture}(\xext,\yext)(\xoff,\yoff)

 \putmorphism(0,600)(1,0)[Sig_\cF`Sig`]{1000}{0}a
 \putmorphism(0,630)(1,0)[\phantom{Sig_\cF}`\phantom{Sig}`F_\cF]{1000}{-1}a

 \putmorphism(1000,600)(1,0)[\phantom{Sig}`Sig^\cF`]{1000}{0}a
 \putmorphism(1000,630)(1,0)[\phantom{Sig}`\phantom{Sig^\cF}`F^\cF]{1000}{1}a

\putmorphism(1000,530)(0,-1)[\phantom{Sig}`End(Set)`]{500}{0}l
\putmorphism(980,600)(0,-1)[\phantom{Sig}`\phantom{End(Set)}`\br]{550}{1}l
\putmorphism(1020,600)(0,-1)[\phantom{Sig}`\phantom{End(Set)}`U]{550}{-1}r
\putmorphism(0,550)(2,-1)[`\phantom{End(Set)}`\dot{\br}]{1000}{1}l
\putmorphism(1300,220)(2,1)[``K]{450}{1}l
\putmorphism(2100,550)(-2,-1)[`\phantom{End(Set)}`\ddot{\br}]{1000}{1}r

  \put(1000,870){\oval(100,100)[t]}
  \put(950,870){\line(0,-1){185}}
  \put(1050,870){\vector(0,-1){185}}
   \put(860,890){${\cF}$}

 \end{picture}
\end{center}
via $\dot{\br}$ and $\ddot{\br}\dashv K$, respectively. The essential image of $Sig_\cF$ in $End(Set)$ is the category of polynomial functors and arbitrary natural transformations, whereas the essential image of the image of $Sig^\cF$ in $End(Set)$ is the category of all finitary functors and arbitrary natural transformations. The categories $Sig_\cF$ and $Sig^\cF$  are not only Kleisli and Eilenberg-Moore objects in the 2-category $\Cat$ but they are also  Kleisli and Eilenberg-Moore objects in the 2-category $\Mon_l(\Cat)$ of monoidal categories, lax monoidal functors, and monoidal transformations. The Kleisli part is a consequence of Theorem 4.1 from \cite{Z3}, see also \cite{McC}, and the Eilenberg-Moore part is a consequence of Theorem \ref{thm-em-Linton}.

We can also extend the representation of $Sig$ in some finer ways by taking submonads of $\cF$. Consider the factorization of the representation $\br$ via subcategory $WPb(Set)\ra End(Set)$ of endofunctors that preserve weak pullbacks and natural transformations that are weakly cartesian, i.e. with naturality squares being weak pullbacks. The category $WPb(Set)$ has reflexive coequalizers\footnote{If we were to insist on preservation of pullbacks instead of weak pullbacks, the category would not have reflexive coequalizers.} that are preserved by the right adjoint $U'$ to the restricted monoidal representation $\br:Sig\ra WPb(Set)$. The resulting lax monoidal monad $\cS=U'\br$ is the symmetrization monad and images of monoidal representations $Sig_\cS$ and $Sig^\cS$ in $End(Set)$ are polynomial functors with cartesian natural transformations and analytic functors with weakly cartesian natural transformations, respectively, c.f. \cite{Z2}. The story can be lifted to rigid and symmetric operads, on one hand, and polynomial and analytic monads, on the other.
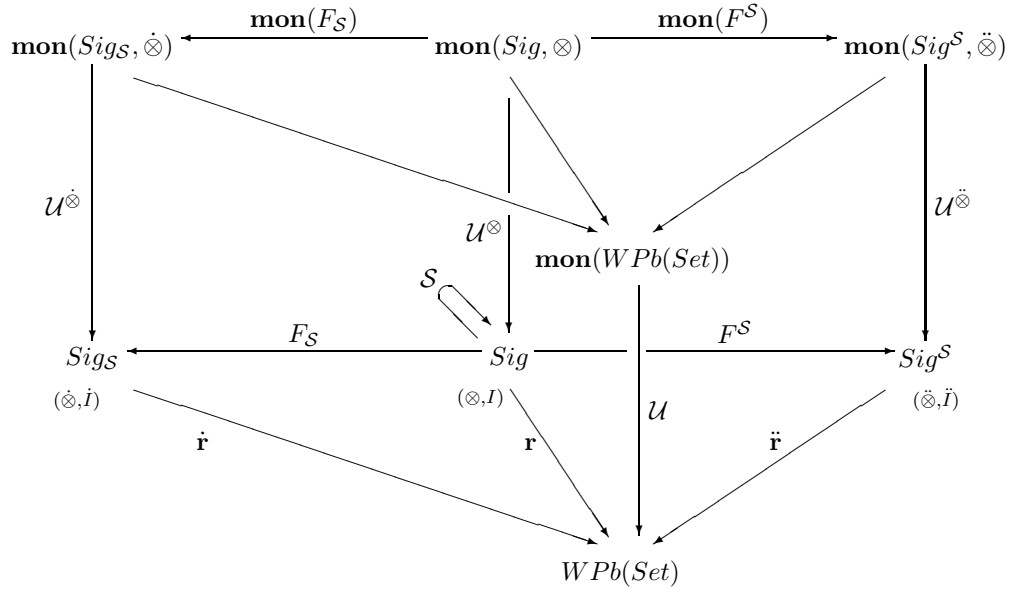
\begin{center} \xext=3200 \yext=2350
\scalebox{0.9}{
\begin{picture}(\xext,\yext)(\xoff,\yoff)
 \setsqparms[0`1`0`0;1600`1200]
 \putsquare(0,850)[\mon(Sig_\cS,\dot{\otimes})`\mon(Sig,\otimes)`Sig_\cS`Sig;`\cU^{\dot{\otimes}}``]
 \put(1600,1850){\line(0,-1){360}}
  \put(1600,1400){\vector(0,-1){450}}
 \put(1440,1300){$\cU^{\otimes}$}

 \putmorphism(0,2080)(1,0)[\phantom{\mon(Sig_\cS,\dot{\otimes})}`\phantom{\mon(Sig,\otimes)}`\mon(F_\cS)]{1600}{-1}a

 \putmorphism(0,880)(1,0)[\phantom{Sig_\cS}`\phantom{Sig}`F_\cS]{1600}{-1}a

 \setsqparms[0`0`1`0;1600`1200]
 \putsquare(1600,850)[\phantom{\mon(Sig,{\otimes})}`\mon(Sig^{\cS},\ddot{\otimes})`\phantom{Sig}`Sig^\cS;``\cU^{\ddot{\otimes}}`]

 \putmorphism(1600,2080)(1,0)[\phantom{\mon(Sig,\otimes)}`\phantom{\mon(Sig,\otimes)^{\widetilde{\cS}}}`{\mon(F^\cS)}]{1600}{1}a


  \put(1700,880){\line(1,0){350}} \put(2130,880){\vector(1,0){950}} \put(2400,910){$F^\cS$}

\put(1400,700){$_{(\otimes, I )}$}
\put(-150,700){$_{(\dot{\otimes}, \dot{I} )}$}
\put(3150,700){$_{(\ddot{\otimes}, \ddot{I} )}$}

\put(1800,0){$WPb(Set)$}
 \putmorphism(-50,800)(3,-1)[\phantom{Sig}`\phantom{Sig^\cS}`]{2200}{1}b
 \put(400,500){$\dot{\br}$}

 \putmorphism(1560,800)(2,-3)[\phantom{Sig}`\phantom{Sig^\cS}`]{470}{1}b
 \put(1660,500){$\br$}

  \putmorphism(3150,800)(-3,-2)[\phantom{Sig}`\phantom{Sig^\cS}`]{1100}{1}b
  \put(2600,500){$\ddot{\br}$}

\put(1700,1200){$\mon(WPb(Set))$}
 \putmorphism(-50,2000)(3,-1)[\phantom{Sig}`\phantom{Sig^\cS}`]{2200}{1}b

 \putmorphism(1560,2000)(2,-3)[\phantom{Sig}`\phantom{Sig^\cS}`]{470}{1}b

  \putmorphism(3150,2000)(-3,-2)[\phantom{Sig}`\phantom{Sig^\cS}`]{1100}{1}b

 \putmorphism(2100,1200)(0,-1)[\phantom{\mon(Sig,\otimes)}`\phantom{Sig}`\cU]{1100}{1}r

  \put(1380,1080){\oval(100,100)[tl]}
  \put(1330,1080){\line(1,-1){150}}
  \put(1380,1130){\vector(1,-1){150}}
   \put(1260,1125){${\cS}$}
 \end{picture}
 }
\end{center}
The category of monoids $\mon(Sig^{\cS},\ddot{\otimes})$ is the category of symmetric operads in $Set$, whereas $\mon(Sig_{\cS},\dot{\otimes})$ is the category of rigid operads, i.e. those symmetric operads whose symmetric actions are free. Their images in the category $\mon(End(Set))$  of monads on $Set$  are categories of (finitary) analytic and polynomial monads, respectively. The symmetrization monad $\cS$ can be defined explicitly, c.f. \cite{SZ4}. For a signature $\{B_n\}_{n\in \o}$ $\cS$ acts as
\[ \cS(\{B_n\}_n) = \{ B_n\times S_n\}_n \]
where $S_n$ is the permutation group of $\{1,\ldots, n\}$, i.e. n-ary operations are taken with all their `versions' obtained by permuting entries. As we said, the lower part of the above diagram is in fact a diagram in 2-category of monoidal categories $\Mon_l(\Cat)$. But it can be even further lifted to the 2-category $\Act_l\Mon_l(\Cat, Set)$ of actions of monoidal categories on $Set$. As a consequence of Theorem \ref{thm-kem-diag-in-actions}, the lower part of the diagram
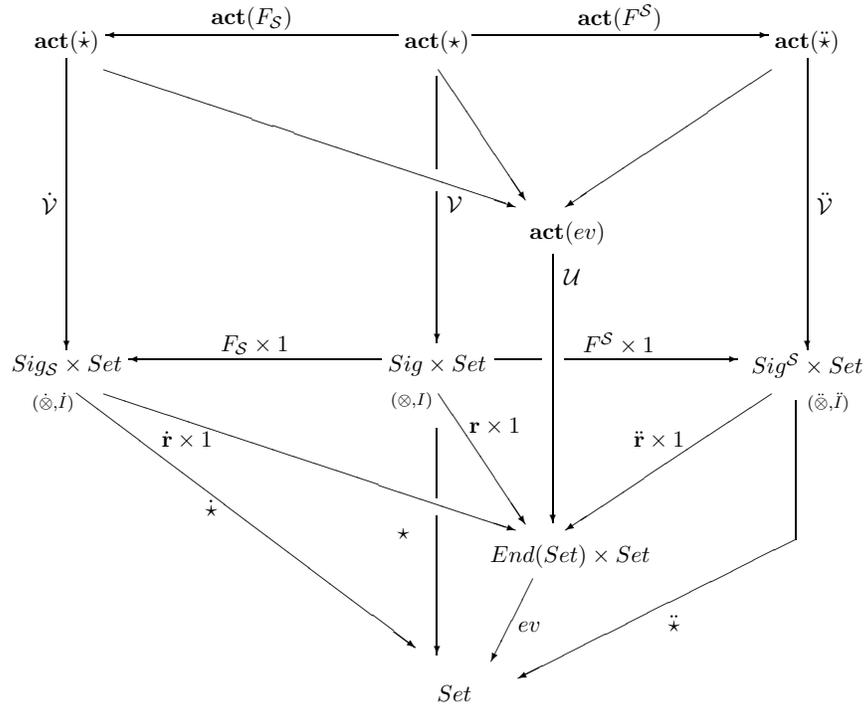
\begin{center} \xext=3200 \yext=3050
\scalebox{0.8}{
\begin{picture}(\xext,\yext)(\xoff,\yoff)
 \setsqparms[0`1`0`0;1600`1400]
 \putsquare(0,1450)[\act(\dot{\star})`\act(\star)`Sig_\cS\times Set`Sig\times Set;`\dot{\cV}`\cV`]
 \put(1600,2700){\line(0,-1){400}}
  \put(1600,2200){\vector(0,-1){650}}

 \putmorphism(0,2880)(1,0)[\phantom{\act(\dot{\star})}`\phantom{\act(\star)}`\act(F_\cS)]{1600}{-1}a

 \putmorphism(0,1480)(1,0)[\phantom{Sig_\cS\times Set}`\phantom{Sig\times Set}`F_\cS\times 1]{1600}{-1}a

 \setsqparms[0`0`1`0;1600`1400]
 \putsquare(1600,1450)[\phantom{\act(\star)}`\act(\ddot{\star})`\phantom{Sig}`Sig^\cS\times Set;``\ddot{\cV}`]

 \putmorphism(1600,2880)(1,0)[\phantom{\act(\star}`\phantom{\act(\ddot{\star})}`{\act(F^\cS)}]{1600}{1}a

 \putmorphism(1600,1480)(1,0)[\phantom{Sig\times Set}`\phantom{Sig^\cS\times Set}`F^\cS\times 1]{1600}{0}a
  \put(1850,1480){\line(1,0){200}}
  \put(2150,1480){\vector(1,0){750}}

\put(1400,1300){$_{(\otimes, I )}$}
\put(-150,1300){$_{(\dot{\otimes}, \dot{I} )}$}
\put(3200,1300){$_{(\ddot{\otimes}, \ddot{I} )}$}

\put(2000,2000){$\act(ev)$}
 \putmorphism(-50,2800)(3,-1)[\phantom{Sig}`\phantom{Sig^\cS}`]{2200}{1}b

 \putmorphism(1560,2800)(2,-3)[\phantom{Sig}`\phantom{Sig^\cS}`]{470}{1}b

  \putmorphism(3150,2800)(-3,-2)[\phantom{Sig}`\phantom{Sig^\cS}`]{1100}{1}b

 \putmorphism(2100,2000)(0,-1)[\phantom{\mon(Sig,\otimes)}`\phantom{Sig}`\phantom{\cU}]{1300}{1}r
  \put(2150,1800){$\cU$}

\put(1830,600){$End(Set)\times Set$}
 \putmorphism(-50,1400)(3,-1)[\phantom{Sig}`\phantom{Sig^\cS}`]{2200}{1}b
 \put(410,1100){$\dot{\br}\times 1$}

 \putmorphism(1560,1400)(2,-3)[\phantom{Sig}`\phantom{Sig^\cS}`]{470}{1}b
 \put(1740,1160){$\br\times 1$}

  \putmorphism(3150,1400)(-3,-2)[\phantom{Sig}`\phantom{Sig^\cS}`]{1100}{1}b
  \put(2450,1100){$\ddot{\br}\times 1$}

  \put(1600,0){$ Set$}

    \put(1600,1180){\line(0,-1){300}}
  \put(1600,800){\vector(0,-1){600}}
  \put(1430,700){$\star$}

  \putmorphism(-50,1400)(4,-3)[\phantom{Sig}`\phantom{Sig^\cS}`]{1650}{1}b
  \put(600,800){$\dot{\star}$}

  \put(3150,700){\vector(-2,-1){1200}}
  \put(3150,1300){\line(0,-1){600}}
  \put(2600,300){$\ddot{\star}$}

  \putmorphism(2050,600)(-1,-2)[\phantom{Sig}`\phantom{Sig^\cS}`]{250}{1}b
  \put(1950,300){$ev$}
  \end{picture}
  }
\end{center}
is again a diagram of Kleisli and Eilenberg-Moore objects but this time for the monad $\cS\times 1$ in the action $\star: Sig\times Set\ra Set$ in the 2-category $\Act_l\Mon_l(\Cat, Set)$. The upper part of the above diagram arises by applying the 2-functor $\act$ (sending actions to actions along actions) to the bottom part of the diagram. This produces the objects of actions and the representations of these objects in $\act(ev)$, which is the category of the actions of monads on $Set$, i.e. the category of algebras for monads on $Set$.

As the title suggests, the main motivation for this work is to study polynomial and analytic functors and monads on their own. However, to exhibit the abstract 2-categorical pattern behind the above story, the paper develops a substantial amount of category theory done in 2-categories with finite limits. These results are of independent interest. They can be particularly useful in the contexts where the operation of substitution plays a role. The main contributions of this kind are Theorems \ref{thm-em-Linton} and \ref{thm-kem-diag-in-actions}. Theorem \ref{thm-em-Linton} says that internally to a 2-category with finite products the structure described by F. Linton gives rise to an Eilenberg-Moore object in the 2-category of monoidal objects $\Mon_l(\cA)$. This problem was suggested by F. Linton in \cite{Li} and its solution has already some history, see \cite{G}, \cite{SZ1}, \cite{Se}. Theorem \ref{thm-kem-diag-in-actions} extends the statement of Theorem \ref{thm-em-Linton} from (lax monoidal) monads on monoidal categories to monads on their actions. Its proof, apart of Theorem \ref{thm-em-Linton}, is based on two facts. The first, Theorem \ref{thm-lift-kem}, says that the lax slice 2-fibration  creates certain Kleisli and Eilenberg-Moore objects and the second, Proposition \ref{prop-Mon-slice-eq-to-Act}, says that certain lax slice 2-category is isomorphic to the 2-category of some actions.

The abstract pattern behind this story considered in a 2-category $\cA$ with finite products is developed in Sections \ref{sec-Alg-2-Cats}, \ref{sec-Mon-prelim}, \ref{sec-Mon}. The applications it produces - among other things, analytic and polynomial functors and monads on slices of $Set$, c.f. \cite{J2}, \cite{Z2}, and on presheaf categories, cf. \cite{FGHW} - are presented in Section \ref{sec-case-study}.

In section \ref{sec-Alg-2-Cats} we describe some general 2-categorical preliminaries. Apart from the last subsection it does not contain anything essentially new but it points out to some new aspects of the very well-known stories. For example, it explains how the Kleisli and Eilenberg-Moore algebras can be used to get better representations of some interesting categories. In subsection \ref{subsec-lax-slice} we study a lax slice as a 2-fibration over an arbitrary 2-category. We show that in some cases this 2-fibration creates both Kleisli and Eilenberg-Moore objects. In section \ref{sec-Mon-prelim} we discuss the monoidal preliminaries. The main part of this long section is the proof of Proposition \ref{prop-Mon-slice-eq-to-Act} saying that the lax slice 2-category over a monoidal object of endomorphisms $\cX^\cX$ is isomorphic to the 2-category of actions on $\cX$. This proof can be safely skipped on first reading. In Section \ref{sec-Mon} we explain under what assumptions on a lax monoidal monad on a 0-cell in a 2-category we get various 0-cells of algebras, monoids and various 1-cells between them. In Section \ref{subsec-monoidal-em-obj} we complete the proof of the result suggested by F. Linton in \cite{Li} concerning the monoidal structure on a category of Eilenberg-Moore algebras for a lax monoidal monad in arbitrary 2-category with finite products. We show (Theorem \ref{thm-em-Linton}) that, under some mild conditions concerning the existence and preservation of reflexive coequalizers, the structure defined by F. Linton in \cite{Li}, does indeed give rise to an Eilenberg-Moore object in any 2-category $\cA$ with finite products. This result extends \cite{G}, \cite{Z1}, \cite{Se}. In subsection \ref{subsec-Act}  we study actions of monoidal categories together with lax monoidal monads. We show (Theorem \ref{thm-kem-diag-in-actions}) that in some circumstances Kleisli and Eilenberg-Moore objects from the category of monoidal categories lift to the categories of actions. This gives a refinement of the previous  result.

Then in Section \ref{sec-case-study} we give some applications. First we discuss the 2-categories we would like to consider in the examples. We are mainly interested in the 2-categories that have fibrations with a fixed base as 0-cells. The bases we consider here are $Set$, $Cat$, $\o Gph$ ($\o$-graphs, globular sets). But there are many other interesting possibilities. However, as substitution is NOT cartesian i.e. the substitution tensor on a fibration of typed signatures {\tiny (almost)} NEVER preserves prone\footnote{Also called cartesian.} morphisms, we allow as 1-cells functors that do not preserve prone morphisms. Then we show how some particular notions of polynomial and analytic functors and monads - including those - on slices of $Set$ c.f. \cite{J2}, and on presheaf categories, c.f. \cite{FGHW}, fit into the general scheme. In the 2-category $\Fib_{/Set}$ of fibrations over $Set$ we repeat the corresponding story from \cite{Z2} using the abstract setup introduced here and then extend it. As the fibrations of rigid, symmetric signatures, rigid and symmetric multicategories and then polynomial and analytic functors and monads arise in these contexts via a very precise abstract procedure (described at the beginning of Section \ref{subsec-general scheme}), it gives a much better understanding what is to be expected from these structures. Similar considerations in the 2-category $\Fib_{/Cat}$, with suitably refined notion of Burroni fibration, give rise to the notion of polynomial and analytic (endo)functors and monads considered in \cite{FGHW}. Similar general considerations in $\Fib_{/\o Gph}$ are less interesting and the resulting representation is from the beginning full on isomorphisms. However, even in this case one can get something non-trivial.

{\bf Notation}. There are two warnings concerning notation that we want to make already in the introduction. When considering various entities the single dot $\dot{(-)}$ will always indicate that the entity is related to Kleisli objects, whereas two dots $\ddot{(-)}$ will always indicate that the entity is related to Eilenberg-Moore objects.

The other convention is of a more serious nature. We call it {\em diagram chasing in 0-cells}. It concerns the way the proofs in a 2-categories are described. As all of the abstract arguments will be made in a 2-category with finite limits, we can use equational logic to simplify the arguments. Such a logic was partially developed in \cite{LiHa} and yet to use it with the full precision would be quite a challenge, not only to write proofs but possibly even more to read them. Instead we decided to use logic but on the `semantic side'. Thus usually when we consider a product of two objects $A\otimes B$ in a monoidal category $(\cC,\otimes, \ldots,)$, we will insist that we are given two projections $\mathbf{A},\mathbf{B}:\cC\times\cC\ra \cC$. When composed with $\otimes: \cC\times\cC\ra \cC$, it will produce a 1-cell $\mathbf{A}\otimes\mathbf{B}: \cC\times\cC\ra \cC$, which is incidentally equal to $\otimes$. This allows manipulating variables while using infix notation for tensor and significantly simplifies many arguments. The details of these conventions the reader will find in Subsections \ref{subsec-diagram-chasing-in-0} and \ref{subsec-monoidal-obj-and-actions}.

\section{2-categorical preliminaries}\label{sec-Alg-2-Cats}

\subsection{(Co)completeness of $0$-cells and properties of $1$-cells}
It makes sense to talk about both completeness and cocompleteness of $0$-cells in any 2-category $\cA$.
A $0$-cell $\cE$ in $\cA$ has (co)limits of type $J$ if for any $0$-cell $\cX$ in $\cA$, the category $\cA(\cX,\cE)$ has (co)limits of type $J$, and any pre-composition functor
\[ \cA(G,\cE): \cA(\cX,\cE)\lra \cA(\cY,\cE) \]
induced by any $1$-cell $G:\cY\ra \cX$ preserves the (co)limit of type $J$. Thus we can talk of a $0$-cell which is finitely complete, has reflexive coequalizers etc.  A 1-cell $F:\cE \ra \cE'$ preserves (co)limits of type $J$ if the post-composition
\[ \cA(\cX,F): \cA(\cX,\cE)\lra \cA(\cX,\cE')\]
preserves the (co)limits of type $J$ for every $0$-cell $\cX$.

We can also mimic in any 2-category $\cA$ the properties that usually apply to functors like faithfulness, fullness and faithfulness, fullness on isos and faithfulness. We say that a $1$-cell $F:\cE \ra \cE'$ has such a property iff for any $0$-cell $\cX$, the post-composition functor
\[ \cA(\cX,F): \cA(\cX,\cE)\lra \cA(\cX,\cE')\]
has this property.

The $0$-cell that has reflexive coequalizers and $1$-cells that preserve them will be called {\em rc $0$-cells} and {\em rc 1-cells}, respectively. A monad is rc iff it is defined on an rc-0-cell and its 1-cell part is rc. A monoidal object (see Section \ref{sec-Mon-prelim}) $(\cC,\otimes,I,\alpha,\lambda,\varrho)$ is rc iff $0$-cell $\cC$ and tensor $1$-cell $\otimes$ are rc.

\subsection{(Co)completeness and some exactness properties of $2$-categories}
All limits, colimits (possibly weighted) and exponentials in 2-categories occurring in this paper will be considered in the strictest possible 2-categorical sense, i.e. they are determined uniquely up to an isomorphism. The results clearly extend to some other more `flexible' `weaker' contexts but we don't need it here for our applications and it would only add unnecessary complications.

We will be working mainly in a 2-category $\cA$ that has finite products. This allows us to talk about monoidal categories inside $\cA$, which we call {\em monoidal objects}. As monoidal categories can act on categories, monoidal objects in $\cA$ can act on 0-cells in $\cA$. The categories of monoids and the categories of actions of monoids along an action of a monoidal category are weighted limits, cf. \cite{SiZ}, so these constructions also make sense in any 2-category $\cA$ with finite limits. As in any 2-category, the notion of a monad, Kleisli-object (or $\bk$-object for short) and Eilenberg-Moore object (or $\bem$-object for short) make sense, c.f. \cite{St}. We have also the notion of a monoidal monad in $\cA$ on a monoidal object.

Even if we usually do not assume that our $2$-category $\cA$ has all the mentioned limits, we will be assuming two exactness properties of $k$-objects in $\cA$ (whenever it makes sense):
\begin{enumerate}
  \item $\bk$-objects commute with finite products;
  \item if $F:\cC\ra \cD$ is a left adjoint $1$-cell ($F\dashv G$), then in the factorization in $\cA$ via Kleisli object for the induced monad $\cR=GF$
  \begin{center} \xext=800 \yext=520
\begin{picture}(\xext,\yext)(\xoff,\yoff)
 \settriparms[1`1`-1;400]
  \putVtriangle(0,50)[\cC`\cD`\cC_\cR ;F`F_\cR`K]
\end{picture}
\end{center}
$1$-cell $K$ is full and faithful.
\end{enumerate}

Clearly, both exactness properties hold trivially in the $2$-category $\Cat$.

We will be also considering exponential $0$-cells in $\cA$ in the clear (strict) 2-categorical sense.

\subsection{Monads and algebras}\label{sec-mnds-algs}
We recall below the definitions of monads, Kleisli and Eilenberg-Moore objects in arbitrary 2-category $\cA$, c.f. \cite{St}, and introduce the related notation.

A tuple $(\cC,\cR,\eta,\mu)$ is a {\em monad} in $\cA$ if $\cC$ is a 0-cell, $\cR:\cC\ra\cC$ is a 1-cell, and $\eta: 1_\cC\ra \cR$, $\mu: \cR^2\ra \cR$ are 2-cells such that
\[ \mu\circ \eta_\cR = 1_\cR = \mu\circ \cR(\eta),\hskip 5mm \mu\circ \mu_\cR =\mu\circ \cR(\mu). \]
As usual, we often say that $\cR$ is a monad when $\cC$, $\eta$, and $\mu$ are understood. A {\em lax morphism of monads} $(T,\tau): (\cC,\cR,\eta,\mu)\lra (\cC',\cR',\eta',\mu')$ is a 1-cell $T:\cC\ra \cC'$ together with a 2-cell $\tau:\cR'T\ra T\cR$ in $\cA$ such that
\[ \tau \circ \eta'_T = T(\eta), \hskip 5mm \tau \circ \mu' = T(\mu)\circ \tau_\cR\circ \cR'(\tau).\]
A {\em transformation of monad morphisms} $\sigma :  (T,\tau)\ra (T',\tau'): (\cC,\cR,\eta,\mu)\ra (\cC',\cR',\eta',\mu')$ is a 2-cell $\sigma:T\ra T'$ such that
\[ \sigma_\cR\circ \tau = \tau'\circ \cR'(\sigma). \]
In this way we have defined a 2-category $\Mnd_l(\cA)$ of monads in $\cA$ with lax morphisms and transformations of lax morphisms.

We say that a triple $(\cX,u,\nu)$ is a {\em subequalizing of} $\cR$ {\em (on $\cX$)} or that it {\em subequalizes the monad} $\cR$,  if $\cX$ is a 0-cell, $u:\cX\ra \cC$ is a 1-cell, $\nu:\cR u \ra u$ is a 2-cell such that
\[ \nu\circ \eta_u=1_u,\hskip 5mm \nu\circ \cR(\nu)= \nu \circ \mu_u. \]
A {\em morphism of subequalizings} $\sigma:  (\cX,u,\nu)\ra(\cX,u',\nu')$ is a 2-cell $\sigma:u\ra u'$ making the square
\begin{center} \xext=800 \yext=500
\begin{picture}(\xext,\yext)(\xoff,\yoff)
 \setsqparms[1`1`1`1;800`400]
 \putsquare(0,50)[\cR u`\cR u'`u`u';\cR(\sigma)`\nu`\nu'`\sigma]
\end{picture}
\end{center}
commute. Let $Subeq_\cR(\cX,\cR)$ denote the category of subequalizings of the monad $\cR$ on a 0-cell $\cX$.
Precomposing a subequalizing $(\cX,u,\nu)$ with 1-cell $L:\cX'\ra \cX$, we get again a subequalizing $(\cX',uL,\nu_L)$ and in fact we have
a 2-functor
\[  Subeq_\cR(-,\cR) : \cA \ra \Cat \]
of subequalizings of $\cR$. If it is representable, i.e. if there is a 0-cell $\cC^\cR$ (also denoted $\bem(\cR)$) and a natural 2-isomorphism
\[ \bEM_\cR : \cA(-,\cC^\cR)\lra Subeq_\cR(-,\cR) \]
then we say that $\cR$ admits the {\em Eilenberg-Moore object}. The value of $\bEM_\cR$ on the identity on $\cC^\cR$ denoted, $(\cC^\cR,U^\cR,\beta)$, is called the Eilenberg-Moore object for $\cR$ or {\em $\bem$-object for} $\cR$, for short. The $\cX$-component of this 2-isomorphism
\[ \bEM_\cR : \cA(\cX,\cC^\cR)\lra Subeq_\cR(\cX,\cR) \]
is given then by
\[ n:L\ra L':\cX\ra \cC^\cR \;\longmapsto U^\cR(n):(\cX,U^\cR L,\beta_L)\ra  (\cX,U^\cR L',\beta_{L'}).\]
Thus, we have the following principles identifying 1- and 2-cells into $\cC^\cR$:
\begin{enumerate}
  \item Two 1-cells $L,L':\cX\ra \cC^\cR$ are equal iff $U^\cR L= U^\cR L'$ and $\beta_L=\beta_{L'}$.
  \item Two parallel 2-cells $n,n': L\ra L':\cX\ra \cC^\cR$ are equal iff $U^\cR(n)=U^\cR(n')$.
\end{enumerate}
In particular, 1-cell $U^\cR$ is faithful and reflects isomorphisms (i.e. precomposing with $U^\cR$ is such a functor).

Assume that the $\bem$-object exists for $\cR$. Then, as $(C,\cR,\mu)$ subequalizes $\cR$, there is a 1-cell $F^\cR: \cC \ra \cC^\cR$ such that
\[ U^\cR\circ F^\cR = \cR, \hskip 5mm \beta_{F^\cR} =\mu. \]
Moreover, $\beta:(\cC^\cR,\cR U^\cR,\mu_{U^\cR}) \ra (\cC^\cR,U^\cR,\beta)$ is a morphism of subequalizings of $\cR$. Hence there is a unique 2-cell $\varepsilon^\cR: F^\cR U^\cR \ra 1_{\cC^\cR}$ such that $U^\cR(\varepsilon^\cR)=\beta$. As we have
\[ U^\cR(\varepsilon^\cR)\circ \eta_{U^\cR}=\beta \circ \eta_{U^\cR}=1_{U^\cR},\hskip 5mm U^\cR(\varepsilon^\cR_{F^\cR}\circ F^\cR(\eta)) = \mu\circ \cR(\eta)=1_{F^\cR} \]
and $U^\cR$ is faithful, it follows that $(F^\cR\dashv U^\cR,\eta,\varepsilon^\cR)$ is an adjunction giving rise to the monad $\cR$.

Dually, we say that a triple $(\cX,u,\nu)$ is a {\em subcoequalizing of} $\cR$ {\em on $\cX$} or that it {\em subcoequalizes the monad} $\cR$,  if $\cX$ is a 0-cell, $u:\cC\ra \cX$ is a 1-cell, $\nu:u\cR  \ra u$ is a 2-cell such that
\[ \nu\circ u(\eta)=1_u,\hskip 5mm \nu\circ \nu_\cR= \nu \circ u(\mu). \]
A {\em morphism of subcoequalizings} $\sigma:  (\cX,u,\nu)\ra(\cX,u',\nu')$ is a 2-cell $\sigma:u\ra u'$ making the square
\begin{center} \xext=800 \yext=500
\begin{picture}(\xext,\yext)(\xoff,\yoff)
 \setsqparms[1`1`1`1;800`400]
 \putsquare(0,50)[u\cR `u'\cR `u`u';\sigma_\cR`\nu`\nu'`\sigma]
\end{picture}
\end{center}
commute. Let $Subcoeq_\cR(\cR,\cX)$ denote the category of subequalizings of the monad $\cR$ on a 0-cell $\cX$.
Postcomposing a subcoequalizing $(\cX,u,\nu)$ with 1-cell $L:\cX\ra \cX'$, we get again a subcoequalizing $(\cX',Lu,L(\nu))$ and in fact we have
a 2-functor
\[  Subcoeq_\cR(\cR,-) : \cA \ra \Cat \]
of subcoequalizings of $\cR$.  If it is representable, i.e. if there is a 0-cell $\cC_\cR$ (also denoted $\bk(\cR)$) and a natural 2-isomorphism
\[ \bKl_\cR : \cA(\cC^\cR,-)\lra Sucobeq_\cR(\cR,-) \]
then we say that $\cR$ admits the {\em Kleisli object}. The value of $\bKl_\cR$ on the identity on $\cC_\cR$, denoted $(\cC_\cR,F_\cR,\kappa)$, is called the Kleisli object for $\cR$, or {\em $\bk$-object for} $\cR$, for short. The $\cX$-component of this 2-isomorphism
\[ \bKl_\cR : \cA(\cC_\cR,\cX)\lra Subcoeq_\cR(\cR,\cX) \]
is given then by
\[ n:L\ra L':\cC_\cR\ra \cX \;\longmapsto n_{F_\cR}:(\cX,L F_\cR ,L(\kappa))\ra  (\cX,L' F_\cR ,L'(\kappa)).\]
Thus we have the analogous principles identifying 1- and 2-cells from $\cC_\cR$:
\begin{enumerate}
  \item Two 1-cells $L,L':\cC_\cR\ra \cX$ are equal iff $L F_\cR = L' F_\cR$ and $L(\kappa)=L'(\kappa)$.
  \item Two parallel 2-cells $n,n': L\ra L':\cC_\cR\ra \cX$ are equal iff $n_{F^\cR}=n'_{F^\cR}$.
\end{enumerate}
In particular, the functor of postcomposing with  $F_\cR$ is faithful and reflects both isomorphisms and identities. In $\Cat$, a functor $F:\cC\ra \cD$ is surjective on objects iff the functor $F\circ(-):\Cat(\cX,\cC)\ra \Cat(\cX,\cD)$ of postcomopsing with $F$ reflects identities, for any category $\cX$.
We will adopt this convention in all 2-categories. Thus, in particular, we shall say that $F_\cR$ is a {\em surjective on objects} 1-cell.

{\em Remarks.}
\begin{enumerate}
  \item Note that even in $\Cat$ $F_\cR$ does not need to be injective on objects. For the closure operator on the powerset $C:\cP(X)\ra\cP(X)$,  the $\bk$-object (and $\bem$-object as well) is the poset of closed sets, and $F^C:\cP(X)\ra \cP(X)^C$ sends subsets of $X$ to their closures. However, if the monad $\cR$ is injective on objects, so is $F_\cR$.
  \item In $\Cat$  $F_\cR$ is always surjective\footnote{In general we expect just essential surjectivity of $F_\cR$. But since we use the strictest possible definition of $\bk$-objects, $F_\cR$ has to reflect identities. For this surjectivity on objects is needed.} and $U^\cR$ is faithful and conservative. If $\bk$-and $\bem$-objects exist for a monad $\cR$ in any 2-category $\cA$, then the functors of  precomposing with $F_\cR$ and of postcomposing with, $U^\cR$ are faithful and conservative.
\end{enumerate}

Assume that  $\bk$-object exists for $\cR$. Then, as $(C,\cR,\mu)$ subcoequalizes $\cR$, there is a 1-cell $U_\cR: \cC_\cR\ra \cC$ such that
\[ U_\cR\circ F_\cR = \cR, \hskip 5mm U_\cR(\kappa)=\mu. \]
Moreover, $\kappa:(\cC_\cR,F_\cR\cR  ,F_\cR(\mu)) \ra (\cC_\cR,F_\cR,\kappa)$ is a morphism of subcoequalizings of $\cR$. Hence there is a unique 2-cell $\varepsilon_\cR: F_\cR U_\cR \ra 1_{\cC_\cR}$ such that $(\varepsilon_\cR)_{F^\cR}=\kappa$. As we have
\[ (U_\cR(\varepsilon_\cR)\circ \eta_{U_\cR})_{F_\cR}=\mu \circ \eta_{\cR}=1_{\cR},\hskip 5mm (\varepsilon_\cR)_{F_\cR}\circ F_\cR(\eta) = \kappa\circ F_\cR(\eta)=1_{F_\cR} \]
and $F_\cR$ is surjective on objects, it follows that $(F_\cR\dashv U_\cR,\eta,\varepsilon_\cR)$ is an adjunction, again giving rise to the monad $\cR$.

To illustrate how we can identify 1-cells and 2-cells from $\cC_\cR$ and to $\cC^\cR$, we shall define the comparison 1-cell from $\cC_\cR$ to $\cC^\cR$   using the universal and the couniversal property of these objects and then we shall prove that they are equal.

As $(\cC^\cR,F^\cR,(\varepsilon^\cR)_{F^\cR}:F^\cR\cR\ra F^\cR)$ is a subcoequalizing of $\cR$, by the couniversal property of $\cC_\cR$ we have a 1-cell $\Phi_\cR:\cC_\cR\ra \cC^\cR$ such that
\[ \Phi_\cR \circ F_\cR= F^\cR, \hskip 5mm \Phi_\cR(\kappa)=(\varepsilon^\cR)_{F^\cR}. \]
As $(\cC_\cR,U_\cR,U_\cR(\varepsilon^\cR):\cR U_\cR\ra U_\cR)$ is a subequalizing of $\cR$, by the universal property of $\cC^\cR$ we have a 1-cell $\Phi^\cR:\cC_\cR\ra \cC^\cR$ such that
\[ U^\cR\circ \Phi^\cR = U_\cR, \hskip 5mm \beta_{\Phi^\cR}=U_\cR(\varepsilon_\cR). \]

First we show that $\Phi^\cR$ satisfies the properties uniquely determining $\Phi_\cR$.
We have
\[ U^\cR\circ \Phi^\cR \circ F_\cR=  U_\cR \circ F_\cR= \cR = U^\cR \circ F^\cR, \]
and
\[ \beta_{\Phi^\cR\circ F_\cR} = U_\cR(\varepsilon_\cR)_{F_\cR} =\mu= U^\cR(\varepsilon^\cR)_{F^\cR}=\beta_{F^\cR}. \]
It follows from  the principle of equality for the 1-cells with codomain $\cC^\cR$ that
\[ \Phi^\cR\circ F_\cR = F^\cR. \]
and hence the 2-cells
\[ \Phi_\cR(\kappa):\Phi^\cR F_\cR\cR\lra \Phi^\cR F^\cR \]
and
\[ (\varepsilon^\cR)_{F^\cR}:F^\cR\cR \lra F^\cR\]
are parallel. Moreover
\[ U^\cR\circ \Phi^\cR(\kappa)=U_\cR(\varepsilon_\cR)_{F_\cR} = \mu =U^\cR(\varepsilon^\cR)_{F^\cR} \]
and then by the principle of equality for the 2-cells with codomain $\cC^\cR$ we have that
\[ \Phi_\cR(\kappa)=(\varepsilon^\cR)_{F^\cR}.\]
Thus $\Phi^\cR=\Phi_\cR$ indeed.

The second proof of the above equality will use the principles of equality of 1- and 2-cells with domain $\cC_\cR$ to show that $\Phi_\cR$ satisfies the properties uniquely determining $\Phi^\cR$. We have
\[ U^\cR\circ \Phi_\cR \circ F_\cR=  U^\cR \circ F^\cR= \cR = U_\cR \circ F_\cR, \]
and
\[ U^\cR\circ \Phi_\cR (\kappa) = U^\cR(\varepsilon^\cR)_{F^\cR} =\mu= U_\cR(\kappa). \]
It follows from  the principle of equality for the 1-cells with domain $\cC_\cR$ that
\[ U^\cR\circ \Phi_\cR  = U_\cR. \]
and hence the 2-cells
\[ \beta_{\Phi_\cR}:\cR U^\cR \Phi_\cR \lra U^\cR \Phi_\cR \]
and
\[ U^\cR(\varepsilon_\cR):\cR U_\cR \lra U_\cR\]
are parallel. Moreover
\[ \beta_{\Phi_\cR F_\cR} = \beta_{F^\cR}= \mu =U_\cR(\varepsilon_\cR)_{F_\cR} \]
and then by the principle of equality for the 2-cells with domain $\cC_\cR$ we have that
\[ \beta_{\Phi_\cR }=U_\cR(\varepsilon_\cR).\]
This ends the second proof of $\Phi^\cR=\Phi_\cR$.

By a {\em $\bkem$-diagram for a monad} $(\cC,\cR,\eta,\varepsilon)$ we mean the diagram in a 2-category $\cA$ that includes the monad $\cR$, the $\bk$- and $\bem$-objects for $\cR$ and all the other data described in this subsection, i.e. a diagram with the data as displayed below
\begin{center} \xext=3200 \yext=1050
\begin{picture}(\xext,\yext)(\xoff,\yoff)
 \putmorphism(0,550)(1,0)[\cC_\cR`\cC`]{1600}{0}a
 \putmorphism(0,580)(1,0)[\phantom{\cC_\cR}`\phantom{\cC}`F_\cR]{1600}{-1}a
 \putmorphism(0,520)(1,0)[\phantom{\cC_\cR}`\phantom{\cC}`U_\cR]{1600}{1}b
  \putmorphism(1600,580)(1,0)[\phantom{\cC}`\cC^\cR`]{1600}{0}a
 \putmorphism(1600,580)(1,0)[\phantom{\cC}`\phantom{\cC^\cR}`F^\cR]{1600}{1}a
 \putmorphism(1600,520)(1,0)[\phantom{\cC}`\phantom{\cC^\cR}`U^\cR]{1600}{-1}b
  \put(1590,300){\oval(100,100)[b]}
  \put(1540,300){\line(0,1){175}}
  \put(1640,300){\vector(0,1){175}}
  \put(1400,100){${(\cR,\eta,\varepsilon)}$}
   \put(100,200){${(F_\cR\dashv U_\cR,\eta,\varepsilon_\cR)}$}
   \put(2300,200){${(F^\cR\dashv U^\cR,\eta,\varepsilon^\cR)}$}
  \put(0,650){\line(1,1){200}}
  \put(200,850){\line(1,0){2800}}
  \put(3000,850){\vector(1,-1){200}}
  \put(1500,900){${\Phi^\cR}$}
 \end{picture}
\end{center}
We often omit the comparison morphism ${\Phi^\cR}$ explicitly.

We end this subsection by proving internally the main part of the Beck's monadicity theorem, i.e. the existence of the left adjoint to the comparison functor.

\begin{proposition}\label{prop-split-coeq} If $(\cC,\cR,\eta,\mu)$ is a monad in a 2-category $\cR$ admitting an $\bem$-object $(\cC^\cR,U^\cR,\beta)$, then $U^\cR:\cC^\cR\ra \cC$ creates coequalizers of $U^\cR$-split pairs.  \end{proposition}
{\it Proof.}~ Let $\cX$ be a 0-cell in $\cA$,
\begin{center} \xext=1600 \yext=350
\begin{picture}(\xext,\yext)(\xoff,\yoff)
 \putmorphism(0,150)(1,0)[\mathds{A}`\mathds{B}`]{600}{0}a
 \putmorphism(0,180)(1,0)[\phantom{\mathds{A}}`\phantom{\mathds{B}}`f]{600}{1}a
  \putmorphism(0,120)(1,0)[\phantom{\mathds{A}}`\phantom{\mathds{B}}`g]{600}{1}b
 \putmorphism(600,150)(1,0)[\phantom{\mathds{B}}`\mathds{C}`q]{600}{1}a
\end{picture}
\end{center}
a diagram in $\cA(\cX,\cC^\cR)$ with $q\circ f=q \circ g$ and
\begin{center} \xext=1600 \yext=350
\begin{picture}(\xext,\yext)(\xoff,\yoff)
 \putmorphism(0,200)(1,0)[U^\cR(\mathds{A})`U^\cR(\mathds{B})`]{800}{0}a
 \putmorphism(0,230)(1,0)[\phantom{U^\cR(\mathds{A}))}`\phantom{U^\cR(\mathds{B})}`U^\cR(f)]{800}{1}a
  \putmorphism(0,170)(1,0)[\phantom{U^\cR(\mathds{A})}`\phantom{U^\cR(\mathds{B})}`U^\cR(g)]{800}{1}b
   \putmorphism(0,20)(1,0)[\phantom{U^\cR(\mathds{A})}`\phantom{U^\cR(\mathds{B})}`t]{800}{-1}b

 \putmorphism(800,200)(1,0)[\phantom{U^\cR(\mathds{B})}`U^\cR(\mathds{C})`U^\cR(q)]{800}{1}a
 \putmorphism(800,130)(1,0)[\phantom{U^\cR(\mathds{B})}`\phantom{U^\cR(\mathds{C})}`s]{800}{-1}b
\end{picture}
\end{center}
\[ U^\cR(q)\circ s = 1_{U^\cR(\mathds{C})},\hskip 5mm 1_{U^\cR(\mathds{B})}= U^\cR(f)\circ t,\hskip 5mm s\circ U^\cR(q)= U^\cR(g)\circ t \]
its $U^\cR$-splitting in $\cA(\cX,\cC)$. We shall show that $q: \mathds{B}\ra \mathds{C}$ is a coequalizer in  $\cA(\cX,\cC^\cR)$.

First we show that
\[ U^\cR(p)\circ s : (U^\cR(\mathds{C},\beta_{\mathds{C}}))\lra (U^\cR(\mathds{C},\beta_{\mathds{C}}))\]
is a morphism of subequalizings of $\cR$. We have
\[ U^\cR(p)\circ s \circ \beta_\mathds{C} = \]
\[ = U^\cR(p)\circ s \circ \beta_\mathds{C}\circ \cR U^\cR(q)\circ \cR(s) = \]
\[ = U^\cR(p)\circ s  \circ U^\cR(q)\circ \beta_\mathds{B}\circ \cR(s) = \]
\[ = U^\cR(p)\circ U^\cR(g)\circ t  \circ \beta_\mathds{B}\circ \cR(s) = \]
\[ = U^\cR(p)\circ U^\cR(f)\circ t  \circ \beta_\mathds{B}\circ \cR(s) = \]
\[ = U^\cR(p) \circ \beta_\mathds{B}\circ \cR(s) = \]
\[ = \beta_\mathds{D}\circ \cR U^\cR(p) \circ  \cR(s) = \]
\[ = \beta_\mathds{D}\circ \cR (U^\cR(p) \circ  U^\cR(f)\circ t \circ  s) = \]
\[ = \beta_\mathds{D}\circ \cR (U^\cR(p) \circ  U^\cR(g)\circ t \circ  s) = \]
\[ = \beta_\mathds{D}\circ \cR (U^\cR(p) \circ s \circ   U^\cR(q)\circ s) = \]
\[ = \beta_\mathds{D}\circ \cR( U^\cR(p) \circ s) \]
Thus $U^\cR(p)\circ s$ is a subequalizing indeed and there is a unique $w:\mathds{C}\ra \mathds{D}$ in $\cA(\cX,\cC^\cR)$ such that $U^\cR(w)= U^\cR(p)\circ s$.
As 2-cells $w\circ q$ and $p$ are parallel and
\[ U^\cR(w\circ q) = U^\cR(p)\circ s \circ U^\cR(q) = \]
\[ = U^\cR(p) \circ U^\cR(g)\circ t = U^\cR(p) \circ U^\cR(f)\circ t = U^\cR(p)  \]
we have $w\circ q=p$ as required.
$\Box$

\begin{proposition}\label{prop-left-adj-ddot-r}
Let  $r:\cC\ra \cM$ be a 1-cell with a right adjoint $U$, so that $(r\dashv U,\eta,\varepsilon)$ is an adjunction. Let $(\cC,\cR=Ur,\eta,\mu=U\varepsilon_r)$
be the induced monad admitting the $\bem$-object $(\cC^\cR,U^\cR,\beta)$. If $\cM$ is an rc-$0$-cell, then the comparison 1-cell $K:\cM\ra \cC^\cR$ has a left adjoint $\ddot{r}$
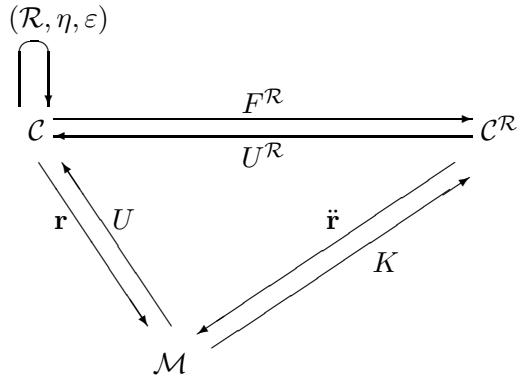
\begin{center} \xext=1700 \yext=1350
\begin{picture}(\xext,\yext)(\xoff,\yoff)
\putmorphism(100,850)(1,0)[\cC`\cC^\cR`]{1600}{0}a
 \putmorphism(100,880)(1,0)[\phantom{\cC}`\phantom{\cC^\cR}`F^\cR]{1600}{1}a
 \putmorphism(100,820)(1,0)[\phantom{\cC}`\phantom{\cC^\cR}`U^\cR]{1600}{-1}b

  \put(90,1100){\oval(100,100)[t]}
  \put(40,1100){\line(0,-1){175}}
  \put(140,1100){\vector(0,-1){175}}
  \put(0,1200){${(\cR,\eta,\varepsilon)}$}
\put(500,0){$\cM$}

 \putmorphism(60,800)(2,-3)[\phantom{\cC}`\phantom{\cC^\cR}`]{470}{1}b
 \put(160,500){$\br$}
  \putmorphism(140,800)(2,-3)[\phantom{\cC}`\phantom{\cC^\cR}`]{470}{-1}b
 \put(360,500){$U$}

  \putmorphism(1650,800)(-3,-2)[\phantom{\cC}`\phantom{\cC^\cR}`]{1100}{1}b
  \put(1100,500){$\ddot{\br}$}
   \putmorphism(1700,750)(-3,-2)[\phantom{\cC}`\phantom{\cC^\cR}`]{1100}{-1}b
   \put(1260,350){$K$}
 \end{picture}
\end{center}
\end{proposition}
{\it Proof.}~ The comparison 1-cell $K:\cM\ra \cC^\cR$ is determined by the subequalizing $(\cM,U,U(\varepsilon))$ of $\cR$ so that
\[ U^\cR K = U, \hskip 5mm \beta_K =U(\varepsilon). \]
As we have
\[ U^\cR K r = U r = \cR U^\cR F^\cR, \hskip 5mm \beta_{Kr} = U\varepsilon_r = \mu =\beta_{F^\cR} \]
it follows that $Kr=F^\cR$.

The left adjoint $\ddot{r}:\cC^\cR\ra \cM$ is given by the (reflexive) coequalizer
\begin{center} \xext=1600 \yext=350
\begin{picture}(\xext,\yext)(\xoff,\yoff)
 \putmorphism(0,150)(1,0)[r \cR U^\cR`r U^\cR`]{800}{0}a
 \putmorphism(0,200)(1,0)[\phantom{r \cR U^\cR}`\phantom{r U^\cR}`r U^\cR(\varepsilon^\cR)]{800}{1}a
 \putmorphism(0,100)(1,0)[\phantom{r \cR U^\cR}`\phantom{r U^\cR}`\varepsilon_{r U^\cR}]{800}{1}b
 \putmorphism(800,150)(1,0)[\phantom{r U^\cR}`\ddot{r}`q]{800}{1}a
\end{picture}
\end{center}
in $\cA(\cC^\cR,\cM)$ with the common inverse $r(\eta)_{U^\cR}$.

The diagram in $\cA(\cC^\cR,\cC^\cR)$
\begin{center} \xext=1800 \yext=350
\begin{picture}(\xext,\yext)(\xoff,\yoff)
 \putmorphism(0,150)(1,0)[F^\cR \cR U^\cR`F^\cR U^\cR`]{1000}{0}a
 \putmorphism(0,200)(1,0)[\phantom{F^\cR \cR U^\cR}`\phantom{F^\cR U^\cR}`F^\cR U^\cR(\varepsilon^\cR)]{1000}{1}a
 \putmorphism(0,100)(1,0)[\phantom{F^\cR \cR U^\cR}`\phantom{F^\cR U^\cR}`K(\varepsilon_{r U^\cR})]{1000}{1}b
 \putmorphism(1000,150)(1,0)[\phantom{F^\cR U^\cR}`1_{\cC^\cR}`\varepsilon^\cR]{800}{1}a
\end{picture}
\end{center}
is a $U^\cR$-split coequalizer as its image in $\cA(\cC^\cR,\cC)$ obtained by postcomposing with $U^\cR$
\begin{center} \xext=1800 \yext=350
\begin{picture}(\xext,\yext)(\xoff,\yoff)
 \putmorphism(0,150)(1,0)[\cR^2 U^\cR`\cR U^\cR`]{1000}{0}a
 \putmorphism(0,200)(1,0)[\phantom{\cR^2 U^\cR}`\phantom{\cR U^\cR}`\cR U^\cR(\varepsilon^\cR)]{1000}{1}a
 \putmorphism(0,100)(1,0)[\phantom{\cR^2 U^\cR}`\phantom{\cR U^\cR}`\mu_{r U^\cR}]{1000}{1}b
 \putmorphism(1000,150)(1,0)[\phantom{\cR U^\cR}`U^\cR`U^\cR(\varepsilon^\cR)]{800}{1}a
\end{picture}
\end{center}
is split by
\begin{center} \xext=1800 \yext=350
\begin{picture}(\xext,\yext)(\xoff,\yoff)
 \putmorphism(0,150)(1,0)[\cR^2 U^\cR`\cR U^\cR`\eta_{\cR U^\cR}]{1000}{-1}a
 \putmorphism(1000,150)(1,0)[\phantom{\cR U^\cR}`U^\cR`\eta_{U^\cR}]{800}{-1}a
\end{picture}
\end{center}
By Proposition \ref{prop-split-coeq}, $\varepsilon^\cR$ is a coequalizer of $F^\cR U^\cR(\varepsilon^\cR)$ and $K(\varepsilon_{r U^\cR})$ in $\cA(\cC^\cR,\cC^\cR)$. Thus by its couniversal property we get the unit 2-cell $\ddot{\eta}$ from the diagram

\begin{center} \xext=2000 \yext=650
\begin{picture}(\xext,\yext)(\xoff,\yoff)
 \putmorphism(0,350)(1,0)[F^\cR \cR U^\cR`F^\cR U^\cR`]{1000}{0}a
 \putmorphism(0,400)(1,0)[\phantom{F^\cR \cR U^\cR}`\phantom{F^\cR U^\cR}`F^\cR U^\cR(\varepsilon^\cR)]{1000}{1}a
 \putmorphism(0,300)(1,0)[\phantom{F^\cR \cR U^\cR}`\phantom{F^\cR U^\cR}`K(\varepsilon_{r U^\cR})]{1000}{1}b

  \put(1200,360){\vector(3,1){600}}
  \put(1360,460){$\varepsilon^\cR$}
   \put(1200,300){\vector(3,-1){600}}
   \put(1360,90){$K(q)$}
\putmorphism(1930,600)(0,-1)[1_{\cC^\cR}`K\ddot{r}`\ddot{\eta}]{500}{1}r

\put(1560,320){$(*)$}
\end{picture}
\end{center}
We get the counit $\ddot{\varepsilon}$ by the couniversal property of $q_K$ from the diagram
\begin{center} \xext=2000 \yext=650
\begin{picture}(\xext,\yext)(\xoff,\yoff)
 \putmorphism(0,350)(1,0)[r \cR U^\cR K`r U^\cR K`]{1000}{0}a
 \putmorphism(0,400)(1,0)[\phantom{r \cR U^\cR K}`\phantom{r U^\cR K}`r U^\cR(\varepsilon_K)]{1000}{1}a
 \putmorphism(0,300)(1,0)[\phantom{r \cR U^\cR K}`\phantom{r U^\cR K}`\varepsilon_{r U^\cR K})]{1000}{1}b

  \put(1200,360){\vector(3,1){600}}
  \put(1460,530){$q_{K}$}
   \put(1200,300){\vector(3,-1){600}}
   \put(1460,120){$\varepsilon$}
\putmorphism(1930,600)(0,-1)[\ddot{r}K`1_\cM`\ddot{\varepsilon}]{500}{1}r

\put(1560,320){$(**)$}
\end{picture}
\end{center}
Precomposing the triangle $(*)$ with $K$ and postcomposing triangle $(**)$ with $K$, we get a pair of triangles
\begin{center} \xext=1000 \yext=1050
\begin{picture}(\xext,\yext)(\xoff,\yoff)
 \putmorphism(0,500)(1,0)[F^\cR U^\cR K`K\ddot{r}K`K(q)_K]{1000}{1}a
 \putmorphism(1000,1000)(0,-1)[K`\phantom{K\ddot{r}K}`\ddot{\eta}_K]{500}{1}r
 \putmorphism(1000,500)(0,-1)[\phantom{K\ddot{r}K}`K`K(\ddot{\varepsilon})]{500}{1}r

  \put(200,600){\vector(2,1){700}}
  \put(400,860){$(\varepsilon^\cR)_{K}$}
   \put(200,400){\vector(2,-1){700}}
   \put(400,120){$K(\varepsilon)$}
\end{picture}
\end{center}
As
\[  U^\cR (\varepsilon^\cR)_K = \beta_K = U(\varepsilon)=U^\cR K(\varepsilon)  \]
we have $(\varepsilon^\cR)_K=K(\varepsilon)$. Thus
\[ K(\ddot{\varepsilon})\circ \ddot{\eta}_K= 1_K \]
as $1_K$ is the unique 2-cell such that $1_K\circ (\varepsilon^\cR)_K=K(\varepsilon)$.

In order to see the other triangular equality, consider the diagram defining both $\ddot{r}(\ddot{\eta})$ and $\ddot{\varepsilon}_{\ddot{r}}$
\begin{center} \xext=2000 \yext=1550
\begin{picture}(\xext,\yext)(\xoff,\yoff)

 \putmorphism(0,1230)(1,0)[\phantom{r\cR U^\cR}`\phantom{r U^\cR}`r U^\cR(\varepsilon^\cR)]{1200}{1}a
 \putmorphism(0,1170)(1,0)[\phantom{r\cR U^\cR}`\phantom{r U^\cR}`\varepsilon_{r U^\cR}]{1200}{1}b

 \putmorphism(0,630)(1,0)[\phantom{r\cR U^\cR K\ddot{r}}`\phantom{r U^\cR K\ddot{r}}`r U^\cR(\varepsilon^\cR) K\ddot{r}]{1200}{1}a
 \putmorphism(0,570)(1,0)[\phantom{r\cR U^\cR K\ddot{r}}`\phantom{r U^\cR K\ddot{r}}`\varepsilon_{r U^\cR K\ddot{r}}]{1200}{1}b

 \setsqparms[0`1`1`0;1200`600]
 \putsquare(0,600)[r R U^\cR`r U^\cR`r\cR U^\cR K \ddot{r}`r U^\cR K \ddot{r};`r \cR U^\cR (\ddot{\eta})`r U^\cR(\ddot{\eta})`]

 \setsqparms[1`0`1`1;800`600]
 \putsquare(1200,600)[\phantom{r U^\cR}`\ddot{r}`\phantom{r U^\cR K \ddot{r}}`\ddot{r}K\ddot{r};q``\ddot{r}(\ddot{\eta})`q_{K\ddot{r}}]

 \putmorphism(2000,600)(0,-1)[\phantom{K\ddot{r}K}`\ddot{r}`\ddot{\varepsilon}_{\ddot{r}}]{600}{1}r

  \put(1300,500){\vector(4,-3){650}}
  \put(1450,260){$\varepsilon_{\ddot{r}}$}

  \put(-500,1200){\line(1,0){300}}
   \put(-500,1200){\line(0,-1){900}}
   \put(-500,300){\line(1,0){1300}}
   \put(800,300){\vector(1,1){200}}
   \put(150,200){$r U(q)$}

\put(1200,1500){\line(0,-1){200}}
\put(0,1500){\line(1,0){1200}}
 \put(0,1500){\vector(0,-1){200}}
   \put(500,1550){$r(\eta_{U^\cR})$}

 \put(2100,1200){\line(1,0){200}}
\put(2300,1200){\line(0,-1){1200}}
 \put(2300,00){\vector(-1,0){200}}
   \put(2320,600){$1_{\ddot{r}}$}
\end{picture}
\end{center}
We have
\begin{align*}
 q &= q\circ r U^\cR(\varepsilon^\cR)\circ r(\eta)_{U^\cR} = \;\; \rm{(q-coequalizer)}  \\
 & = q\circ \varepsilon^\cR_{r U^\cR}\circ r(\eta)_{U^\cR} = \;\; {\rm (MEL\; on\; \varepsilon\; and\; q)} \\
 & = \varepsilon_{\ddot{r}}\circ r U(q)\circ  r(\eta)_{U^\cR} =  \\
 & = \varepsilon_{\ddot{r}}\circ r U^\cR K(q)\circ  r(\eta)_{U^\cR} = \;\; {\rm (def\; of\; \ddot{\eta})} \\
 & = \varepsilon_{\ddot{r}}\circ r U^\cR(\ddot{\eta}) \circ r U^\cR(\varepsilon^\cR) \circ  r(\eta)_{U^\cR} = \;\; {\rm (F^\cR\dashv U^\cR)} \\
 & = \varepsilon_{\ddot{r}}\circ r U^\cR(\ddot{\eta}).
\end{align*}
Thus
\[ \ddot{\varepsilon}_{\ddot{r}}\circ \ddot{r}(\ddot{\eta}) = 1_{\ddot{r}} \]
as $1_{\ddot{r}}$ is the unique 2-cell such that $1_{\ddot{r}}\circ q= \varepsilon_{\ddot{r}}\circ r U^\cR(\ddot{\eta})$.
$\Box$

{\em Remark.} We use the notation from the proof above.  As we have equality of the right adjoints $U=U^\cR \circ K$, it follows that the left adjoints must be isomorphic, i.e. $r\cong \ddot{r}\circ F^\cR$. But we can also see this isomrphism directly as follows. Precomposing the coequalizer $q$ with $F^\cR$, we get a coequalizer $q_{F^\cR}$ which has also a split coequalizer $\varepsilon_r$ with splitting $r(\eta)$ and $r\cR(\eta)$. Thus we have a comparison isomorphism $\tau$ between these coequalizers of the same parallel pair:
\begin{center} \xext=1600 \yext=550
\begin{picture}(\xext,\yext)(\xoff,\yoff)
 \putmorphism(0,350)(1,0)[r \cR^2`r \cR`]{800}{0}a
 \putmorphism(0,400)(1,0)[\phantom{r \cR^2 }`\phantom{r \cR}`r \mu]{800}{1}a
 \putmorphism(0,300)(1,0)[\phantom{r \cR^2}`\phantom{r \cR}`\varepsilon_{r \cR}]{800}{1}b

  \put(900,360){\vector(3,1){600}}
  \put(1160,530){$q_{F^\cR}$}
   \put(900,300){\vector(3,-1){600}}
   \put(1160,120){$\varepsilon_r$}
\putmorphism(1630,600)(0,-1)[\ddot{r} F^\cR`r`\cong]{500}{1}r
\end{picture}
\end{center}

\subsection{Extensions of representations}\label{subsec-extensions-of-rep}

This section rephrases the well-known story of Kleisli and Eilenberg-Moore objects with the emphasis on its representational potential.

Suppose we are given a representation $1$-cell
\[ \br: \cC \lra \cM\]
in an algebraic 2-category $\cA$, so that $\cC$ is (considered as) an `abstract' 0-cell and $\cM$ is (considered as) a `concrete' $0$-cell.
Clearly, it is natural to ask $\br$ to be
\begin{enumerate}
  \item faithful,
  \item even better also full,
  \item or at least faithful, and full on isomorphisms,
  \item or, in addition to faithfulness, to have some intermediate form of fullness that will allow to identify `objects' (0- and 1-cells) of $\cM$ as those that came from $\cC$. This would allow to perform constructions in $\cM$ and come to $\cC$ back with their results in an essentially unique way.
\end{enumerate}
Thus, we may want to add some morphisms to $\cC$ in a `reasonable way' to improve/extend the representation $r$ to one satisfying some of the above properties.

If $\br$ has a right adjoint, say $U$ ($\br\dashv U$ with unit and counit of the adjunction $\eta$ and $\varepsilon$, respectively), this can be done by passing to a category of algebras for the induced monad $\cR=U\circ \br$. In the minimalist way, we can pass from $\cC$ to the Kleisli object  $(\cC_\cR,F_\cR,\bk)$  for $\cR$ and consider the representation $\dot{\br}: \cC_\cR\lra \cM$ such that
\[ \br=\dot{\br}\circ F_\cR,\hskip 2cm \dot{r}(\bk)=\varepsilon_r. \]
In case $\cA$ is $\Cat$, $\dot{\br}$ is necessarily full and faithful. In general, we need to assume separately that such an exactness property holds in the 2-category $\cA$.

However, if $\cC$ has good properties that we may want to keep, it might be better to extend $\br$ by passing to Eilenberg-Moore object $(\cC^\cR,U^\cR,\ba)$ for $\cR$ and consider the representation of $\cC^\cR$, as $\cC^\cR$ is expected to retain many properties of $\cC$. As usual, we have a comparison functor $K: \cM\ra \cC^\cR$ such that
\[ U=U^\cR\circ K, \hskip 2cm \ba_K=U(\varepsilon)  \]
If $\cM$ is rc-$0$-cell, then $K$ has a left adjoint $\ddot{r}:\cC^\cR\lra \cM$, ($\ddot{\br}\dashv K$ with unit and counit $\ddot{\eta}$ and $\ddot{\varepsilon}$, respectively), so that
\[ \br\cong\ddot{\br}\circ F^\cR. \]
If $U$ preserves suitable coequalizers (i.e. a class of coequalizers that contains those that were used to define $\ddot{\br}$, e.g. reflexive coequalizers), then the unit $\ddot{\eta}$ of adjunction $\ddot{\br}\dashv K$ is an iso, i.e. $\ddot{\br}$ is full and faithful. Now the whole diagram in $\cA$ described above looks as follows
\begin{center} \xext=2000 \yext=750
\begin{picture}(\xext,\yext)(\xoff,\yoff)

 \putmorphism(0,600)(1,0)[\cC_\cR`\cC`]{1000}{0}a
 \putmorphism(0,630)(1,0)[\phantom{\cC_\cR}`\phantom{\cC}`F_\cR]{1000}{-1}a
 \putmorphism(0,570)(1,0)[\phantom{\cC_\cR}`\phantom{\cC}`U_\cR]{1000}{1}b

 \putmorphism(1000,600)(1,0)[\phantom{\cC}`\cC^\cR`]{1000}{0}a
 \putmorphism(1000,630)(1,0)[\phantom{\cC}`\phantom{\cC^\cR}`F^\cR]{1000}{1}a
 \putmorphism(1000,570)(1,0)[\phantom{\cC}`\phantom{\cC^\cR}`U^\cR]{1000}{-1}b

\putmorphism(1000,530)(0,-1)[\phantom{\cC}`\cM`]{500}{0}l
\putmorphism(980,530)(0,-1)[\phantom{\cC}`\phantom{\cM}`\br]{500}{1}l
\putmorphism(1020,530)(0,-1)[\phantom{\cC}`\phantom{\cM}`U]{500}{-1}r
\putmorphism(-50,500)(2,-1)[`\phantom{\cM}`\dot{\br}]{1000}{1}l
\putmorphism(1250,170)(2,1)[``K]{450}{1}l
\putmorphism(2000,450)(-2,-1)[`\phantom{\cM}`\ddot{\br}]{1000}{1}r

  \put(790,370){\oval(100,100)[bl]}
  \put(740,370){\line(1,1){185}}
  \put(790,320){\vector(1,1){185}}
   \put(700,215){${\cR}$}

  \put(2800,400){$(1)$}
 \end{picture}
\end{center}
We spell the notation for the adjunctions and sub(co)equalizing 2-cells:
\[ (\br\dashv U, \eta, \varepsilon)\hskip 5 mm (F_\cR\dashv U_\cR, \eta, \varepsilon_\cR),\hskip 5 mm (F^\cR\dashv U^\cR, \eta, \varepsilon^\cR), \hskip 5 mm (\ddot{\br}\dashv K, \ddot{\eta}, \ddot{\varepsilon}) \]
\[ \bk : F_\cR\circ \cR\lra F_\cR, \hskip 5 mm \ba :\cR \circ U^\cR\lra U^\cR\]

In order to get a completion of
\[ \br: \cC\lra \cM \]
that are less full (drastic) we can do two things:
\begin{enumerate}
  \item either consider (non-full) sub-0-cell $\cM'$ of $\cM$ through which $\br$ factors
  \begin{center} \xext=600 \yext=320
\begin{picture}(\xext,\yext)(\xoff,\yoff)
 \settriparms[1`1`1;300]
  \putAtriangle(0,0)[\cC`\cM'`\cM;\br'`\br`]
\end{picture}
\end{center} so that $\cC$ is coreflective in $\cM'$
  and repeat the above for $\br'$ in place of $\br$
  \item or consider a submonad $\cR'$ of $\cR$ and extend $r$ as above but using $\cR'$ instead of $\cR$.
\end{enumerate}
In this paper, we will see examples of both strategies.

\subsection{Representations vs monads}\label{subsec-rep-vs-monads}

In this subsection, we shall show that there is an adjoint correspondence between the factorizations of a given `representation' 1-cell $\br$ with rc right adjoint $U$ and the rc monads over the rc monad $\cR=Ur$. The material in this subsection is mostly rephrasing a part of the content of \cite{Jo} in a way suitable for our context.

Let $\cA$ be a 2-category and $\cC$ a 0-cell in $\cA$. Let $\Mnd_l(\cA,\cC)$ be the category of monads on $\cC$ and lax morphisms of monads.
By $\Adj(\cA,\cC)$ we denote a category whose objects are adjunctions $(\br\dashv U:\cC\ra \cM,\eta,\varepsilon)$ in $\cA$ such that the domain of the left adjoint is the 0-cell $\cC$. Two adjunctions are considered equal iff they have equal right adjoints\footnote{We could say that this category is in fact a category of right adjoints, but we want the left adjoint to be always named.}. The morphisms in $\Adj(\cA,\cC)$ are triangles of adjunctions, such that the triangles of the right adjoints commute (on the nose).

Let us fix an rc adjunction $(\br\dashv U,\eta,\varepsilon)$ in $\cA$ with $\br:\cC\ra\cM$, i.e. $\cC$, $\cM$, $U$ (and hence also $\br$) are rc. Then we have a functor between slice categories
\[ \widehat{(-)} : (\Mnd_l(\cA_{rc},\cC)^{op})_{/\cR} \ra \Adj(\cA_{rc},\cC)_{(r\dashv U)} \]
such that the 2-cell $\tau : \cR'\ra \cR$ which is a 0-cell in $(\Mnd_l(\cA_{rc},\cC)^{op})_{/\cR}$, is sent to the morphism of adjunctions
\begin{center} \xext=1700 \yext=1000
\begin{picture}(\xext,\yext)(\xoff,\yoff)

 \putmorphism(100,850)(1,0)[\cC`\cC^{\cR'}`]{1600}{0}a
 \putmorphism(100,880)(1,0)[\phantom{\cC}`\phantom{\cC^{\cR'}}`F^{\cR'}]{1600}{1}a
 \putmorphism(100,820)(1,0)[\phantom{\cC}`\phantom{\cC^{\cR'}}`U^{\cR'}]{1600}{-1}b

\put(500,0){$\cM$}

 \putmorphism(60,800)(2,-3)[\phantom{\cC}`\phantom{\cC^{\cR'}}`]{470}{1}b
 \put(160,500){$\br$}
  \putmorphism(140,800)(2,-3)[\phantom{\cC}`\phantom{\cC^{\cR'}}`]{470}{-1}b
 \put(360,500){$U$}

  \putmorphism(1650,800)(-3,-2)[\phantom{\cC}`\phantom{\cC^\cR}`]{1100}{1}b
  \put(1100,500){$F^\tau$}
   \putmorphism(1700,750)(-3,-2)[\phantom{\cC}`\phantom{\cC^\cR}`]{1100}{-1}b
   \put(1260,350){$U^\tau$}
 \end{picture}
\end{center}
where $U^\tau$ is the obvious morphism induced by the morphism of monads and $F^\tau$ is its left adjoint existing by Theorem 2 of \cite{Jo}, see also \cite{Li}, \cite{Di}. We also have a functor
\[  \widetilde{(-)}:  \Adj(\cA_{rc},\cC)_{(r\dashv U)} \ra (\Mnd_l(\cA_{rc},\cC)^{op})_{/\cR} \]
such that an adjunction $(r'\dashv U')$ over $(r\dashv U)$
\begin{center} \xext=1700 \yext=1000
\begin{picture}(\xext,\yext)(\xoff,\yoff)

 \putmorphism(100,850)(1,0)[\cC`\cM'`]{1600}{0}a
 \putmorphism(100,880)(1,0)[\phantom{\cC}`\phantom{\cM'}`r']{1600}{1}a
 \putmorphism(100,820)(1,0)[\phantom{\cC}`\phantom{\cM'}`U']{1600}{-1}b

\put(500,0){$\cM$}

 \putmorphism(60,800)(2,-3)[\phantom{\cC}`\phantom{\cM'}`]{470}{1}b
 \put(160,500){$\br$}
  \putmorphism(140,800)(2,-3)[\phantom{\cC}`\phantom{\cM'}`]{470}{-1}b
 \put(360,500){$U$}

  \putmorphism(1650,800)(-3,-2)[\phantom{\cC}`\phantom{\cM'}`]{1100}{1}b
  \put(1100,500){$t$}
   \putmorphism(1700,750)(-3,-2)[\phantom{\cC}`\phantom{\cM'}`]{1100}{-1}b
   \put(1260,350){$V$}
 \end{picture}
\end{center}
is sent to
\begin{center} \xext=1400 \yext=150
\begin{picture}(\xext,\yext)(\xoff,\yoff)
 \putmorphism(0,50)(1,0)[\cR'=U'\,r'`U'\,V\,t\,r'=\cR`U'(\tilde{\eta})_{r'}]{1800}{1}a
 \end{picture}
\end{center}
where with $\tilde{\eta}$ is the unit of the adjunction $t\dashv V$.

We have
\begin{proposition}\label{prop-corefl} If $\cA$ admits $\bem$-objects, then the above functors are well defined and are adjoint  $\widehat{(-)}\dashv\widetilde{(-)}$.
\end{proposition}
{\it Proof.}~ Exercise. $\Box$

{\em Remarks.}
\begin{enumerate}
  \item As usual, the above adjunction restricts to an equivalence between fixed objects. As the unit of this adjunction in an iso, all the monads over $\cR$ are fixed. The fixed adjunctions are monadic adjunctions in the strong sense: the comparison 1-cell to the $\bem$-object needs to be an isomorphism.
  \item Submonads of $\cR$ correspond to coreflective sub-0-cells of $\cM$.
\end{enumerate}

\subsection{Lax slices}\label{subsec-lax-slice}

If $\cA$ is a 2-category and $\cX$ is a 0-cell in $\cA$, {\em the lax slice $\cA_{/_l\cX}$ of $\cA$ over $\cX$} is a 2-category that has as its 0-cells 1-cells $b:B\ra \cX$ in $\cA$ with codomain $\cX$. A 1-cell $(F,\theta): (B,b)\ra (C,c)$ in $\cA_{/_l\cX}$ is a 1-cell $F:C\ra C$ in $\cA$ together with a 2-cell $\theta : c\circ F\ra b$
  \begin{center} \xext=800 \yext=520
\begin{picture}(\xext,\yext)(\xoff,\yoff)
 \settriparms[1`1`1;400]
  \putVtriangle(0,50)[B`C`\cX ;F`b`c]
  \put(370,270){${\theta}$}
   \put(350,180){${\Leftarrow}$}
\end{picture}
\end{center}
A 2-cell $\tau : (F,\theta)\ra (F',\theta')$ in $\cA_{/_l\cX}$ is a 2-cell $\tau:F\ra F'$ in $\cA$ such that the triangle of 2-cells
  \begin{center} \xext=700 \yext=520
\begin{picture}(\xext,\yext)(\xoff,\yoff)
 \settriparms[1`1`1;350]
  \putVtriangle(0,50)[c F`c F'`b ;c(\tau)`\theta`\theta']
     \put(2000,180){$({\bf LST})$}
\end{picture}
\end{center}
commutes.

Throughout the paper by a 2-fibration we mean 2-fibration in the sense of Hermida, cf. \cite{He}.

\begin{proposition}\label{prop-lax-slice-2-fib}
Let $\cM$ be a $0$-cell in a 2-category $\cA$, $p:\cA_{/_l\cM}\lra \cA$ the domain 2-functor. Then
\begin{enumerate}
  \item $p$ is a 2-fibration;
  \item $p$ creates right adjoints of strong morphisms, i.e. if $(F\dashv U,\eta,\varepsilon)$ is an adjunction in $\cA$, $b:B\ra \cM$, $c:C\ra \cM$  1-cells, and  $\theta :c F=b$ an invertible 2-cell, then this adjunction lifts uniquely to an adjunction
\[ ((B,b),(C,c),(F,\theta)\dashv (U,c(\varepsilon)\circ\theta^{-1}_U,\eta,\varepsilon)\]
in the lax slice $\cA_{/_l \cM}$ as
\begin{center} \xext=1200 \yext=720
\begin{picture}(\xext,\yext)(\xoff,\yoff)
 \settriparms[0`1`1;600]
  \putVtriangle(0,50)[B`C`\cM ;`b`c]
  \putmorphism(50,680)(1,0)[\phantom{B}`\phantom{C}`(F,\theta)]{1100}{1}a
  \putmorphism(50,620)(1,0)[\phantom{B}`\phantom{C}`(U,q(\varepsilon)\circ\theta^{-1}_U)]{1100}{-1}b
\end{picture}
\end{center}
\end{enumerate}
\end{proposition}
{\it Proof.}~
Ad 1. To see that $p$ is a 2-fibration, we need to verify that,  cf. \cite{He} Thm. 2.8,
\begin{description}
  \item[(i)]  for any $f:A\ra B$ in $\cA$ and a 0-cell $(A,a:A\ra \cM)$ in $\cA_{/_l\cM}$, there is a prone (cartesian) 2-morphism $(f,\theta):(B,b)\ra (A,a)$
  \item[(ii)] for any two 0-cells $(A,a)$, $(B,b)$, the functor $p_{B,A}:\cA_{/_l\cM}((B,b),(A,a)) \ra \cA(B,A)$ induced by $p$ is a fibration
  \item[(iii)] and finally for any 1-cell $(k,\kappa):(C,c)\ra (B,b)$ in  $\cA_{/_l\cM}$, the precomposition functor
  \[  \cA_{/_l\cM}((k,\kappa),(A,a)): \cA_{/_l\cM}((B,b),(A,a))\lra \cA_{/_l\cM}((C,c),(A,a))\]
  is a morphism of fibrations.
\end{description}

To verify {\bf (i)}, we shall show that a morphism in $\cA_{/_l\cM}$ is 2-prone iff it is strong. Then we can take $(f,id_{a\circ f}): (B,a\circ f)\ra (A,a)$ as a 2-prone morphism in $\cA_{/_l\cM}$ with codomain $(A,a)$ over $f:B\ra A$.

Let $(f,\theta): (B,b)\ra (A,a)$ be a 1-cell in $\cA_{/_l\cM}$ with $\theta:a\circ f\ra b$ an invertible 2-cell. Let $(g,\gamma):(C,c)\ra (A,a)$ be a 1-cell in $\cA_{/_l\cM}$ and $\bar{g}:B\ra A$ a 1-cell in $\cA$ so that we have $f\circ \bar{g}=g$ in $\cA$. Then $(\bar{g},\bar{\gamma})= (\bar{g},\gamma\circ (\theta\circ \bar{g})^{-1}): (C,c)\ra (B,b) $ is the unique 1-cell in $\cA_{/_l\cM}$
\begin{center} \xext=1400 \yext=800
\begin{picture}(\xext,\yext)(\xoff,\yoff)
 \settriparms[1`1`1;400]
  \putVtriangle(600,50)[B`A`\cM ;f`b`a]

  \put(100,700){\vector(2,-1){420}}
  \put(0,700){$C$}
  \put(350,620){$\bar{g}$}

  \put(120,750){\line(1,0){800}}
  \put(920,750){\vector(2,-1){420}}
  \put(1250,620){$g$}

  \put(50,650){\line(1,-2){200}}
  \put(250,250){\vector(3,-1){570}}
  \put(130,250){$c$}

  \put(970,280){$\theta$}
  \put(950,180){$\Leftarrow$}

  \put(350,380){$\bar{\gamma}$}
  \put(350,280){$\Leftarrow$}
\end{picture}
\end{center}
such that
\[ (f,\theta)\circ (\bar{g},\bar{\gamma})= (f\circ \bar{g},\bar{\gamma}\circ (\theta\circ \bar{g})) = (g,\gamma).\]
Thus $(f,\theta)$ is 1-prone.

To see that $(f,\theta)$ is 2-prone, let $\alpha:(g,\gamma)\ra (h,\chi)$ be a 2-cell in $\cA_{/_l\cM}$, i.e.
\begin{center} \xext=600 \yext=400
\begin{picture}(\xext,\yext)(\xoff,\yoff)
 \settriparms[1`1`1;300]
  \putVtriangle(0,0)[a\circ\bar{g}`a\circ h`c;a(\alpha)`\gamma`\chi]
\end{picture}
\end{center}
commutes and $\sigma:\bar{g}\ra\bar{h}:C\ra B$ a 2-cell in $\cA$ such that $p(\alpha)=\alpha= f\circ \sigma$.
Then, putting $\bar{\chi}=\chi\circ(\theta\circ \bar{h})^{-1})$, we have
\[ \bar{\chi}\circ b(\sigma)=  \chi\circ(\theta\circ \bar{h})^{-1})\circ b(\sigma)= \]
\[ = \chi\circ (a\circ f)(\sigma)\circ(\theta\circ \bar{g})^{-1} =\]
\[ = \chi\circ a(\alpha)\circ(\theta\circ \bar{g})^{-1} =\]
\[ = \gamma\circ(\theta\circ \bar{g})^{-1} = \bar{\gamma}\]
and hence  $\sigma:(\bar{g},\bar{\gamma})\ra (\bar{h},\bar{\chi})$ is a morphism in $\cA_{/_l\cM}$ with $f\circ \sigma=\alpha$, as this was true already in $\cA$.

For {\bf (ii)} and {\bf (iii)}, it is enough to show that for any 0-cells $(A,a)$, $(B,b)$ in $\cA_{/_l\cM}$ the functor
\[ p_{(A,a), (B,b)}: \cA_{/_l\cM}((A,a), (B,b)) \lra \cA(A,B) \]
is a fibration with all morphisms being prone. The remaining details are left for the reader.

Ad 2. Let $(F,\theta): (B,b)\ra (C,c)$ be a strong 1-cell in $\cA_{/_l\cM}$ such that $(F\dashv U,\eta,\varepsilon)$ is an adjunction in $\cA$. Then $(U,u)=(U,c(\varepsilon)\circ\theta_U^{-1}):(C,c)\ra (B,b)$ is a 1-cell in $\cA_{/_l\cM}$. As triangular equalities hold true for $\eta$ and $\varepsilon$,  it remains to show that both 2-cells are in $\cA_{/_l\cM}$. To show that the triangle
\begin{center} \xext=600 \yext=400
\begin{picture}(\xext,\yext)(\xoff,\yoff)
 \settriparms[1`1`1;300]
  \putVtriangle(0,0)[b`bUF`b;b(\eta)`1_b`\theta\circ u_F]
\end{picture}
\end{center}
commutes, we have
\[ \theta\circ u_F \circ b(\eta)= \]
\[ = \theta\circ c(\varepsilon)_F\circ\theta_{UF}^{-1} \circ b(\eta)=  \]
\[ =  \theta\circ c(\varepsilon)_F \circ cF(\eta)\circ\theta^{-1}=  \]
\[ =  \theta\circ \theta^{-1}= 1_b \]
where the third equality follows from the middle exchange law for $\eta$ and $\theta^{-1}$.

The commutation of the triangle
\begin{center} \xext=600 \yext=400
\begin{picture}(\xext,\yext)(\xoff,\yoff)
 \settriparms[1`1`1;300]
  \putVtriangle(0,0)[cFU`c`c;c(\varepsilon)` u\circ\theta_U`1_c]
\end{picture}
\end{center}
is even simpler
\[  u\circ\theta_U = c(\varepsilon)\circ\theta_{U}^{-1}\circ \theta_U = c(\varepsilon)= 1_c\circ c(\varepsilon). \]
$\Box$

\begin{theorem}\label{thm-lift-kem}
Let $\cM$ be an rc-$0$-cell in a 2-category $\cA$, $p:\cA_{/_l\cM}\lra \cA$ the domain 2-functor, $(\cR,\eta,\mu)$ a monad given by the adjunction $(r\dashv U,\eta,\varepsilon): \cC\ra \cM$ in $\cA$. Then $p$ creates the $\bkem$-diagram for $\cR$, i.e. the (unique) lift of the $\bkem$-diagram for the monad $(\cR,\eta,\varepsilon)$ in $\cA$
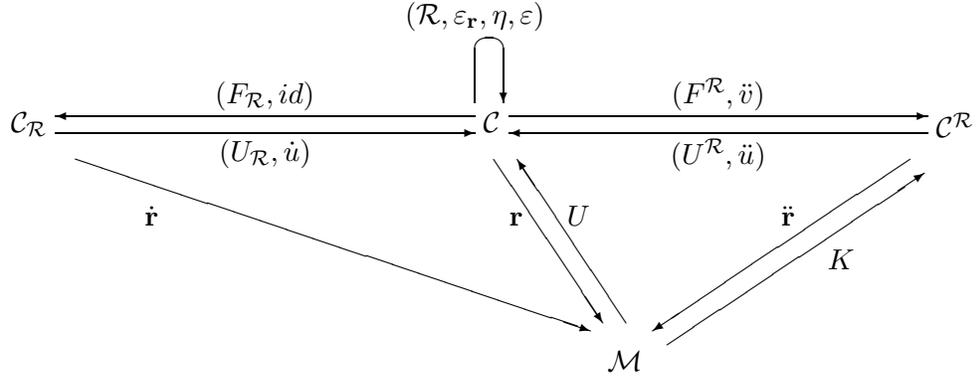
\begin{center} \xext=3200 \yext=1350
\begin{picture}(\xext,\yext)(\xoff,\yoff)
 \setsqparms[0`0`0`0;1600`1200]
 \putsquare(0,850)[\phantom{\mon(\cC_\cR,\dot{\otimes})}`\phantom{\mon(\cC,\otimes)}`\cC_\cR`\cC;``\phantom{\cU^{\otimes}}`]

 \putmorphism(0,880)(1,0)[\phantom{\cC_\cR}`\phantom{\cC}`(F_\cR,id)]{1600}{-1}a
 \putmorphism(0,820)(1,0)[\phantom{\cC_\cR}`\phantom{\cC}`(U_\cR,\dot{u})]{1600}{1}b

 \setsqparms[0`0`0`0;1600`1200]
 \putsquare(1600,850)[\phantom{\mon(\cC,{\otimes})}`\phantom{\mon(\cC^{\cR},\ddot{\otimes})}`\phantom{\cC}`\cC^\cR;``\phantom{\cU^{\ddot{\otimes}}}`]

 \putmorphism(1600,880)(1,0)[\phantom{\cC}`\phantom{\cC^\cR}`(F^\cR,\ddot{v})]{1600}{1}a
 \putmorphism(1600,820)(1,0)[\phantom{\cC}`\phantom{\cC^\cR}`(U^\cR,\ddot{u})]{1600}{-1}b

  \put(1590,1100){\oval(100,100)[t]}
  \put(1540,1100){\line(0,-1){175}}
  \put(1640,1100){\vector(0,-1){175}}
  \put(1300,1200){${(\cR,\varepsilon_\br,\eta,\varepsilon)}$}
\put(2000,0){$\cM$}
 \putmorphism(-50,800)(3,-1)[\phantom{\cC}`\phantom{\cC^\cR}`]{2200}{1}b
 \put(400,500){$\dot{\br}$}

 \putmorphism(1560,800)(2,-3)[\phantom{\cC}`\phantom{\cC^\cR}`]{470}{1}b
 \put(1660,500){$\br$}
  \putmorphism(1640,800)(2,-3)[\phantom{\cC}`\phantom{\cC^\cR}`]{470}{-1}b
 \put(1860,500){$U$}

  \putmorphism(3150,800)(-3,-2)[\phantom{\cC}`\phantom{\cC^\cR}`]{1100}{1}b
  \put(2600,500){$\ddot{\br}$}
   \putmorphism(3200,750)(-3,-2)[\phantom{\cC}`\phantom{\cC^\cR}`]{1100}{-1}b
   \put(2760,350){$K$}
 \end{picture}
\end{center}
is a $\bkem$-diagram for $(\cR,\varepsilon_\br,\eta,\varepsilon)$ in $\cA_{/_l\cM}$, where
  \begin{enumerate}
    \item $\dot{r}:\cC_\cR\ra \cM$ is a morphism from the $\bk$-object $(\cC_\cR,F_\cR,U_\cR)$ existing by its couniversal property;
    \item $\dot{u}=\dot{r}(\varepsilon_\cR):r\,U_\cR=\dot{r}\,F_\cR\, U_\cR \lra \dot{r}$;
    \item $\ddot{r}$ is a left adjoint to $K$ (see Proposition \ref{prop-left-adj-ddot-r}) such that $\ddot{r}\, F^\cR=r$;
    \item $\ddot{u}=\ddot{r}(\varepsilon^\cR):r\,U^\cR=\ddot{r}\,F^\cR\, U^\cR \lra \ddot{r}$.
  \end{enumerate}
\end{theorem}
{\it Proof.}~ Creation of the Kleisli object for $(\cR,\varepsilon_r,\eta,\mu)$. Let $(\cE, \bar{r}:\cE\ra \cM)$ be a 0-cell in $\cA_{/_l\cM}$. We need to show that we have an isomorphism of categories
\[ \bKl_{\cR,\varepsilon_r} : \cA_{/_l\cM}((\cC_\cR,\dot{r}),(\cE,\bar{r}))\lra Subcoeq_{\cA_{/_l\cM}}((\cR,\varepsilon_r,\eta,\mu),(\cE,\bar{r}))  \]
natural in $(\cE, \bar{r})$.

To see that $\bKl_{\cR,\varepsilon_r}$ is bijective on objects, let $L:\cC\ra\cE$ be a 1-cell in $\cA$,  $\theta: \bar{r}L\ra r$ and $\zeta: L\cR\ra L$ 2-cells in $\cA$ so that
$(L,\theta,\zeta)$ be a subcoequalizing of $(\cR,\varepsilon_r)$ in $\cA_{/_l\cM}$.  By the couniversal property of $(\cC_\cR,F_\cR,(\varepsilon_\cR)_{F_\cR})$ in $\cA$ we have a 1-cell $\bar{L}:\cC_\cR\ra \cE$ such that
\[ L= \bar{L}\circ F_\cR, \hskip 1cm \bar{L}((\varepsilon_\cR)_{F_\cR})=\zeta. \]
By the (2-dimensional) couniversal property of $\cC_\cR$, as $\theta: (\bar{r}L,\bar{r}(\zeta)) \ra (r,\varepsilon_r)$ is a morphism of coequalizings, i.e. the square of 2-cells

\begin{center} \xext=800 \yext=500
\begin{picture}(\xext,\yext)(\xoff,\yoff)
 \setsqparms[1`1`1`1;800`400]
 \putsquare(0,50)[\bar{r}L\cR`r\cR`\bar{r}L`r;\theta_\cR`\bar{r}(\zeta)`\varepsilon_r`\theta]
\end{picture}
\end{center}
commutes, there is a (unique) 2-cell $\bar{\theta}:\bar{r}\bar{L}\ra\dot{r}$ such that $\bar{\theta}_{F_\cR}=\theta$. Hence
\[ (\bar{L},\bar{\theta})\circ (F_\cR,id_{\bar{r}})=(\bar{L}\circ F_\cR,\bar{\theta}_{F_\cR})=(L,\theta).\]
Thus $\bKl_{\cR,\varepsilon_r}$ is bijective on objects.

To see that $\bKl_{\cR,\varepsilon_r}$ is full and faithful, consider two subcoequalizings $(L,\theta,\zeta)$ and $(L,\theta,\zeta)$ of $(\cR,\varepsilon_r,\eta,\mu)$ in $\cA_{/_l\cM}$.  Then $(L,\zeta)$ and $(L',\zeta')$ are subcoequalizings of $(\cR,\eta,\mu)$ in $\cA$ that, by the couniversal property of $\cC_\cR$, correspond bijectively to 1-cells in $\bar{L},\bar{L}':\cC_\cR\ra \cE$ in $\cA$. Again by the couniversal property of $\cC_\cR$, the morphisms of subcoequalizings $\sigma : (L,\zeta) \ra (L',\zeta')$ correspond bijectively to 2-cells $\bar{\sigma} : \bar{L}\ra \bar{L}'$ in $\cA$.
As above, we have 2-cells $\bar{\theta}: \bar{r}\bar{L}\ra \dot{r}$ and $\bar{\theta}': \bar{r}\bar{L}'\ra \dot{r}$ such that $\bar{\theta}_{F_\cR}=\theta$ and $\bar{\theta}'_{F_\cR}=\theta'$. Thus we have two parallel 1-cells $(\bar{L},\bar{\theta}), (\bar{L}',\bar{\theta}'): \dot{r}\ra \bar{r}$ in $\cA_{/_l\cM}$.  It remains to show that
\[ \sigma : (L,\theta,\zeta) \ra (L',\theta',\zeta') \]
is a morphism of subcoequalizings of $(\cR,\varepsilon_r,\eta,\mu)$ iff
\[ \bar{\sigma} : (\bar{L},\bar{\theta}\ra (\bar{L}',\bar{\theta}')\]
is a 2-cell in $\cA_{/_l\cM}$, i.e. $\sigma$ is a morphism in $\cA_{/_l\cM}$ iff  $\bar{\sigma}$ is. The former means that the triangle
\begin{center} \xext=700 \yext=450
\begin{picture}(\xext,\yext)(\xoff,\yoff)
 \settriparms[1`1`1;350]
  \putVtriangle(0,0)[\bar{r}L`\bar{r}L'`r;\bar{r}(\sigma)`\theta`\theta']
\end{picture}
\end{center}
commutes, and the latter means that the triangle
\begin{center} \xext=700 \yext=450
\begin{picture}(\xext,\yext)(\xoff,\yoff)
 \settriparms[1`1`1;350]
  \putVtriangle(0,0)[\bar{r}\bar{L}`\bar{r}\bar{L}'`\dot{r};\bar{r}(\bar{\sigma})`\bar{\theta}`\bar{\theta}']
\end{picture}
\end{center}
commutes. The first is obtained from the second by composing with $F_\cR$. By the couniversal property, composing with $F_\cR$ is a full and faithful functor. Thus one of the above triangles commutes iff the other does.

Creation of Eilenberg-Moore object for $(\cR,\varepsilon_r,\eta,\mu)$.

Let $(L:\cE\ra \cC,\theta:r L\ra \bar{r},\xi: \cR L\ra L)$ be a subequalizing of the monad $(\cR,\varepsilon_r,\eta,\mu)$ in $\cA_{/_l\cM}$, i.e. $(L,\xi)$ is a  subequalizing of $(\cR,\eta,\mu)$ in $\cA$, and
\begin{center}\xext=800 \yext=500
\begin{picture}(\xext,\yext)(\xoff,\yoff)
 \settriparms[1`1`-1;400]
  \putVtriangle(0,50)[r \cR L`r L`\bar{r};r(\xi)`\theta \varepsilon_{rL}`\vartheta]
\end{picture}
\end{center}
commutes. By the universal property of $\cC^\cR$ in $\cA$ we have $\bar{L}:\cE\ra \cC^\cR$ such that
\[ L= U^\cR \bar{L}, \hskip 5mm U^\cR(\ddot{\varepsilon})_{\bar{L}} =\xi. \]
We need to show that there is a unique $\bar{\theta}:\ddot{r}\bar{L}\ra \bar{r}$ such that $q_{\bar{L}}\circ \bar{\theta}=\theta$.

Precomposing the coequalizer $q$ with $\bar{L}$, we get a coequalizer
\begin{center} \xext=2000 \yext=650
\begin{picture}(\xext,\yext)(\xoff,\yoff)
 \putmorphism(0,350)(1,0)[r \cR U^\cR \bar{L}`r U^\cR \bar{L}`]{1000}{0}a
 \putmorphism(0,400)(1,0)[\phantom{r \cR U^\cR \bar{L}}`\phantom{r U^\cR \bar{L}}`r U^\cR(\varepsilon_{\bar{L}})]{1000}{1}a
 \putmorphism(0,300)(1,0)[\phantom{r \cR U^\cR \bar{L}}`\phantom{r U^\cR \bar{L}}`\varepsilon_{r U^\cR \bar{L}})]{1000}{1}b

  \put(1200,360){\vector(3,1){600}}
  \put(1460,530){$q_{\bar{L}}$}
   \put(1200,300){\vector(3,-1){600}}
   \put(1460,120){$\theta$}
\putmorphism(1930,600)(0,-1)[\ddot{r}\bar{L}`\bar{r}`\bar{\theta}]{500}{1}r
\end{picture}
\end{center}
whose parallel pair is coequalized by $\theta$. Hence we have a unique $\bar{\theta}$ such that $\bar{\theta}\circ q_{\bar{L}}=\theta$.
This means that
\[ (U^\cR,q) \circ (\bar{L},\theta) = (L,\theta) \]
Thus $\bEM_{(\cR,\varepsilon_r)}$ is bijective on objects.

To see that $\bEM_{(\cR,\varepsilon_r)}$ is also full and faithful, consider two subequalizings $(L,\theta,\xi)$ and  $(L',\theta',\xi')$ of $(\cR,\varepsilon_r)$ in $\cA_{/_l\cM}$ and a morphism of subequalizings  $\tau : (L,\xi)\ra (L',\xi')$ in $\cA$, i.e. the diagram
\begin{center} \xext=600 \yext=500
\begin{picture}(\xext,\yext)(\xoff,\yoff)
 \setsqparms[1`1`1`1;600`400]
 \putsquare(0,50)[\cR L`\cR L'`L`L';\cR(\tau)`\xi`\xi'`\tau]
\end{picture}
\end{center}
commutes. By the universal property of $\cC^\cR$ in $\cA$ there is a unique $\tau: \bar{L}\ra \bar{L}'$ such that $U^\cR(\bar{\tau})=\tau$.
Now it is enough to show that $\tau: (L,\theta,\xi) \ra (L',\theta',\xi')$ is a morphism of subequalizings of  $(\cR,\varepsilon_r)$ in $\cA_{/_l\cM}$  i.e.
the triangle
\begin{center}\xext=800 \yext=500
\begin{picture}(\xext,\yext)(\xoff,\yoff)
 \settriparms[1`1`-1;350]
  \putVtriangle(0,50)[r L`r L'`\bar{r};r(\tau)`\theta `\theta']
  \put(280,220){$^{[\tau]}$}
\end{picture}
\end{center}
commutes iff $\bar{\tau}: (\bar{L},\bar{\theta)} \ra (\bar{L}',\bar{\theta}')$ is a morphism in $\cA_{/_l\cM}$, i.e. the triangle commutes.
 \begin{center}\xext=800 \yext=500
\begin{picture}(\xext,\yext)(\xoff,\yoff)
 \settriparms[1`1`-1;350]
  \putVtriangle(0,50)[\ddot{r} L`\ddot{r} L'`\bar{r};\ddot{r}(\bar{\tau})`\bar{\theta} `\bar{\theta}']
  \put(300,200){$^{[\bar{\tau}]}$}
\end{picture}
\end{center}
In the diagram below
\begin{center} \xext=1000 \yext=900
\begin{picture}(\xext,\yext)(\xoff,\yoff)
 \setsqparms[1`-1`-1`1;1000`700]
 \putsquare(0,100)[\ddot{r} \bar{L}`\ddot{r} \bar{L}'`r U^\cR \bar{L}`r U^\cR \bar{L}';\ddot{r}(\bar{\tau})`q_{\bar{L}}`q_{\bar{L}'}`rU^\cR(\bar{\tau})]

 \put(100,200){\vector(3,2){320}}
 \put(150,520){$\bar{\theta}$}
 \put(900,200){\vector(-3,2){320}}
 \put(800,520){$\bar{\theta}'$}

 \put(100,700){\vector(3,-2){320}}
 \put(150,280){$\theta$}
 \put(900,700){\vector(-3,-2){320}}
  \put(800,280){$\theta'$}

 \put(480,420){$\bar{r}$}
 \put(450,600){$^{[\bar{\tau}]}$}
 \put(480,200){$^{[\tau]}$}
\end{picture}
\end{center}
the outer square commutes by MEL on $q$ and $\bar{\tau}$. The left and right triangles commute by the definitions of $\bar{\theta}$ and $\bar{\theta}'$.
As $q_{\bar{L}}$ is a coequalizer and hence an epi, it follows that the upper triangle commutes iff the lower does.  $\Box$

\section{Monoidal preliminaries in 2-categories}\label{sec-Mon-prelim}
\subsection{Diagram chasing in 0-cells}\label{subsec-diagram-chasing-in-0}

The definitions of a monoidal category, an action of a monoidal category on a category, a category of monoids and, the category of actions of monoids along an action of a monoidal category can be formulated in any 2-category with suitable limits. However, when we want to manipulate such objects, the usual 2-categorical language quickly gets cumbersome. The reason is that the convenient infix notation used for the tensor needs variables (or something that behaves like them). One way to deal with this problem is to develop a suitable language to talk about 2-categorical structures and then interpret it in 2-categories with suitable structure. Then with the help of a sound (and preferably complete) logical system we could prove facts about such 2-categories. A possible base for such a system has been developed in \cite{LiHa}. Even though we believe that the development of such a system, suitable for the considerations of this paper, needs to be done, it is a separate issue (going beyond the scope of this paper). It will need to include what we call {\em diagram chasing in 0-cells} that are not necessarily categories of any sort. Instead of describing such a system we shall develop a notation that, even if being still on the side of semantics, will allow us to present the considerations `inside' 0-cells of 2-categories in the notation familiar from usual categories where we use variables (eventually interpreted as projections) at will. In this notation the composition is, not surprisingly, a substitution of terms, cf. \cite{Law}.  The notation will be explained in the next Subsection  when defining a monoidal object in a 2-category.

\subsection{Monoidal objects and their actions}\label{subsec-monoidal-obj-and-actions}

In this subsection $\cA$ is a 2-category with finite products.

A {\em monoidal object (0-cell)} $(\cC,\otimes,I,\alpha,\lambda,\rho)$ in $\cA$ consists of
\begin{enumerate}
  \item one $0$-cell $\cC$,
  \item two 1-cells $\mathds{A}\otimes\mathds{B} :\cC\times\cC\ra \cC$, $I : 1\ra \cC$,
  \item three invertible 2-cells $\alpha_{\mathds{A},\mathds{B},\mathds{C}}:\mathds{A}\otimes(\mathds{B}\otimes\mathds{C})\lra (\mathds{A}\otimes\mathds{B})\otimes\mathds{C}$, $\lambda_{\mathds{A}}:\mathds{I}\otimes \mathds{A}\ra \mathds{A}$, $\rho_{\mathds{B}}:\mathds{A}\otimes \mathds{I}\ra \mathds{A}$
\end{enumerate}
such that the following two diagrams $\bf MC1$ and $\bf MC2$ of 1- and 2-cells commute. The diagram
\begin{center} \xext=2300 \yext=1000
\begin{picture}(\xext,\yext)(\xoff,\yoff)
\put(800,1020){$\mathds{A}\otimes(\mathds{B}\otimes (\mathds{C}\otimes \mathds{D}))$}
      \put(900,970){\vector(-3,-2){400}}
      \put(100,820){$1_\mathds{A}\otimes \alpha_{\mathds{B},\mathds{C},\mathds{D}}$}

      \put(1300,970){\vector(3,-2){400}}
      \put(1560,820){$\alpha_{\mathds{A},\mathds{B},\mathds{C}\otimes \mathds{D}}$}

\put(0,620){$\mathds{A}\otimes((\mathds{B}\otimes \mathds{C})\otimes \mathds{D})$}
\put(1600,620){$(\mathds{A}\otimes \mathds{B})\otimes (\mathds{C}\otimes \mathds{D})$}
      \put(500,570){\vector(1,-2){180}}
      \put(180,320){$\alpha_{\mathds{A},\mathds{B}\otimes \mathds{C},\mathds{D}}$}

      \put(1700,570){\vector(-1,-2){180}}
      \put(1650,320){$\alpha_{\mathds{A}\otimes \mathds{B},\mathds{C},\mathds{D}}$}

\put(200,120){$(\mathds{A}\otimes(\mathds{B}\otimes \mathds{C}))\otimes \mathds{D}$}
\put(1400,120){$((\mathds{A}\otimes \mathds{B})\otimes \mathds{C})\otimes \mathds{D}$}

  \put(1050,120){\vector(1,0){300}}
      \put(950,0){$\alpha_{\mathds{A},\mathds{B},\mathds{C}}\otimes 1_\mathds{D}$}
       \put(-800,400){$\bf MC1$}
\end{picture}
\end{center}
commutes in the category $\cA(\cC\times\cC\times\cC\times\cC,\cC)$ and where $\mathds{A},\mathds{B}, \mathds{C}, \mathds{D}: \cC\times\cC\times\cC\ra\cC$ are the first, the second, the third, and the forth projections, respectively. For example, $\mathds{B} \otimes\mathds{C} : \cC\times\cC\times\cC\times\cC\ra\cC$ is a 1-cell in $\cA$ which is a composition of a projection to the 2nd and 3rd component followed by $\otimes$, i.e.
\begin{center} \xext=2000 \yext=150
\begin{picture}(\xext,\yext)(\xoff,\yoff)
 \putmorphism(0,50)(1,0)[\cC\times\cC\times\cC\times\cC`\cC\times\cC`\lk \mathds{B},\mathds{C}\rk]{1000}{1}a
 \putmorphism(1000,50)(1,0)[\phantom{\cC\times\cC}`\cC`\otimes]{600}{1}a
\end{picture}
\end{center}
The diagram
  \begin{center}\xext=1000 \yext=700
\begin{picture}(\xext,\yext)(\xoff,\yoff)
 \settriparms[1`1`1;500]
  \putVtriangle(0,50)[\mathds{A}\otimes(\mathds{I}\otimes \mathds{B})`(\mathds{A}\otimes \mathds{I})\otimes \mathds{B}`\mathds{A}\otimes\mathds{B};\alpha_{\mathds{A}, \mathds{I}, \mathds{B}}`\mathds{A}\otimes \lambda_{\mathds{B}}`\rho_{\mathds{A}} \otimes \mathds{B}]

   \put(-1400,280){$\bf MC2$}
\end{picture}
\end{center}
commutes in $\cA(\cC\times\cC,\cC)$. Here, according to the convention $\mathds{A},\mathds{B}: \cC\times\cC\ra\cC$ are the first, and the second projections, respectively, and, for example, $\mathds{A}\otimes(\mathds{I}\otimes \mathds{B}):\cC\times\cC\ra \cC$ is a 1-cell that is a composition displayed below
\begin{center} \xext=2200 \yext=150
\begin{picture}(\xext,\yext)(\xoff,\yoff)
 \putmorphism(0,50)(1,0)[\cC\times\cC`\cC\times\cC\times \cC`\lk \mathds{A},\mathds{I},\mathds{B}\rk]{800}{1}a
 \putmorphism(800,50)(1,0)[\phantom{\cC\times\cC\times \cC}`\cC\times\cC`1\times\otimes]{800}{1}a
  \putmorphism(1600,50)(1,0)[\phantom{\cC\times\cC}`\cC`\otimes]{600}{1}a
\end{picture}
\end{center}
and $\alpha_{\mathds{A}, \mathds{I}, \mathds{B}}$ is the 2-cell $\alpha$ wiskered along the 1-cell $\lk\mathds{A}, \mathds{I}, \mathds{B}\rk$. $\mathds{I}$ is a constant 1-cell `equal' to $I$, i.e. it is a 1-cell of the form $\cE\lra 1 \stackrel{I} {\lra}\cC$ for a suitable domain $\cE$ that is understood in the context (in this case $\cE=\cC\times \cC$).

We can say, for short, that $(C,\otimes)$ is a monoidal object iff the other parts of the data are understood.

Note that, according to the convention, when we have a 1-cell $F:\cC\ra\cD$, then it can be denoted by $F(\mathds{A})$ and the latter notation has the advantage that can be easily extended to, for example,  1-cell $F((\mathds{A}\otimes \mathds{I})\otimes \mathds{B}):\cC\times\cC\ra\cC$. We will often name 1-cell $(\mathds{A}\otimes \mathds{I})$,  $(\mathds{A}\otimes \mathds{I})\otimes \mathds{B})$ etc. as objects of $\cC$, reluctantly adding that they they are `at stage' $\cC$ or $\cC\times\cC$. In fact, the stage of an object can change without notice, as  $F(\mathds{A})$ is at stage $\cC$ when it stands alone but in a formula $(F(\mathds{A})\otimes F(\mathds{A}))\otimes F(\mathds{B})$ it is at stage $\cC\times\cC$. This process of `changing stage' is nothing but the semantical operation corresponding to the weakening in logic.

A {\em lax monoidal morphism (1-cell)}
\[(F,\phi,\bar{\phi}):(\cC,\otimes,I,\alpha,\lambda,\rho)\lra (\cC',\otimes',I',\alpha',\lambda',\rho')\]
is a 1-cell $F:\cC\ra \cC'$ together with 2-cells $\phi : F(\mathds{A}\otimes'\mathds{B})\ra F(\mathds{A}\otimes\mathds{B})$
and $\bar{\phi}: I'\ra F(I)$ such that the three diagrams {\bf MF1}, {\bf MF2}, {\bf MF3} of 2-cells commute.   The diagram
  \begin{center} \xext=2000 \yext=1150
\begin{picture}(\xext,\yext)(\xoff,\yoff)
 \setsqparms[1`1`1`0;2000`450]
  \putsquare(0,550)[F(\mathds{A})\otimes' (F(\mathds{B})\otimes' F(\mathds{C}))`(F(\mathds{A})\otimes' F(\mathds{B}))\otimes' F(\mathds{C})`F(\mathds{A})\otimes' F(\mathds{B}\otimes \mathds{C})`F(\mathds{A}\otimes \mathds{B})\otimes' F(\mathds{C});
  \alpha'_{F(\mathds{A}),F(\mathds{B}),F(\mathds{C})}`1_{F(\mathds{A})}\otimes'\phi`\phi\otimes' 1_{F(\mathds{C})}`]
  \setsqparms[0`1`1`1;2000`450]
  \putsquare(0,100)[\phantom{F(\mathds{A})\otimes' F(\mathds{B}\otimes \mathds{C})}`\phantom{F(\mathds{A}\otimes \mathds{B})\otimes' F(\mathds{C})}`F(\mathds{A}\otimes (\mathds{B}\otimes \mathds{C}))`F((\mathds{A}\otimes \mathds{B})\otimes \mathds{C});
  `\phi`\phi`F(\alpha_{\mathds{A},\mathds{B},\mathds{C}})]
  \put(-1000,280){$\bf MF1$}
\end{picture}
\end{center}
commutes in $\cA(\cC\times\cC\times\cC,\cC')$.
The diagrams
\begin{center} \xext=1500 \yext=600
\begin{picture}(\xext,\yext)(\xoff,\yoff)
 \setsqparms[1`1`-1`1;1500`400]
  \putsquare(0,50)[\mathds{I}'\otimes' F(\mathds{\mathds{A}})`F(\mathds{\mathds{A}})`F(\mathds{I})\otimes F(\mathds{\mathds{A}})`F(\mathds{I}\otimes \mathds{\mathds{A}});
  \lambda'_{F(\mathds{\mathds{A}})}`\bar{\phi}\otimes' 1_{F(\mathds{A})}`F(\lambda_{\mathds{\mathds{A}}})`\phi]
   \put(-1000,280){$\bf MF2$}
\end{picture}
\end{center}
and
\begin{center} \xext=1500 \yext=600
\begin{picture}(\xext,\yext)(\xoff,\yoff)
 \setsqparms[1`1`-1`1;1500`400]
  \putsquare(0,50)[F(\mathds{A})\otimes \mathds{I}'`F(\mathds{A})`F(\mathds{A})\otimes' F(\mathds{I})`F( \mathds{A}\otimes \mathds{I});
\rho'_{F(\mathds{A})}`1_{F(A)}\otimes  \bar{\phi}`F(\rho_{\mathds{A}})`\phi]
\put(-1000,280){$\bf MF3$}
\end{picture}
\end{center}
commute in $\cA(\cC\times\cC,\cC')$. A lax monoidal 1-cell is called {\em strong} ({\em strict}) if the 2-cells $\phi$ and $\bar{\phi}$ are isomorphisms (identities). We obtain the notion of an {\em oplax monoidal 1-cell} by reversing the directions of the coherence 2-cells in the definition  of a lax monoidal 1-cell. As the coherence morphism for the units $\bar{\phi}$ will be denoted, when possible, as the coherence morphism for the tensor $\phi$ with the bar on top, we usually denote the monoidal 1-cells as $(F,\phi)$ rather than $(F,\phi,\bar{\phi})$, for short.

A {\em monoidal transformation between two lax monoidal 1-cells}
\[\tau:(F,\phi)\ra (F',\phi'):(\cC,\otimes,I,\alpha,\lambda,\rho)\lra (\cC',\otimes',I',\alpha',\lambda',\rho')\]
is a 2-cell $\tau:F\ra F'$ in $\cA$ such that the diagrams
\begin{center}
\xext=1600 \yext=650
\begin{picture}(\xext,\yext)(\xoff,\yoff)
 \setsqparms[1`1`1`1;1600`500]
 \putsquare(0,0)[F(\mathds{A})\otimes F(\mathds{B})`F'(\mathds{A})\otimes F'(\mathds{B})`F(\mathds{A}\otimes \mathds{B})`F'(\mathds{A}\otimes \mathds{B});
 \;\;\;\;\tau_\mathds{A}\otimes\tau_{\mathds{B}}`\phi_{\mathds{A},\mathds{B}}`\phi'_{\mathds{A},\mathds{B}}`\tau_{\mathds{A}\otimes \mathds{B}}]
  \put(-1000,280){$\bf MT1$}
\end{picture}
\end{center}
and
  \begin{center} \xext=600 \yext=450
\begin{picture}(\xext,\yext)(\xoff,\yoff)
 \settriparms[1`1`1;350]
  \putAtriangle(0,50)[I'`F(I)`F'(I);\bar{\phi}`\bar{\phi}'`\tau_I]
  \put(-1400,180){$\bf MT2$}
\end{picture}
\end{center}
commute. The convention is as above, but we just point out that the 1-cell $I'$, $F(I)$ etc., are objects of $\cC'$ and have domain $1$ which stands for the empty context.

In this way we have defined the 2-category $\Mon_l(\cA)$ of monoidal objects, lax monoidal 1-cells, and monoidal 2-cells.

\vskip 2mm

Now fix a 0-cell $\cX$ in $\cA$. We shall define the 2-category $\Act_l\Mon_l(\cA,\cX)$  of lax actions of monoidal objects in $\cA$ on $\cX$.

A {\em monoidal action $(\cC,\otimes,I,\alpha,\lambda,\rho,\cX,\star, \psi,\bar{\psi})$ of a monoidal object $(\cC,\otimes,I,\alpha,\lambda,\rho)$  on a 0-cell} $\cX$ consists of
\begin{enumerate}
  \item a monoidal object $(\cC,\otimes,I,\alpha,\lambda,\rho)$,
  \item a 0-cell $\cX$,
  \item a 1-cell $\mathds{A}\star\mathds{X}: \cC\times \cX \ra \cX$,
  \item two 2-cells $\psi: \mathds{A}\star(\mathds{B}\star\mathds{X})\ra (\mathds{A}\otimes(\mathds{B})\star\mathds{X}$, $\bar{\psi}:\mathds{X}\ra \mathds{I}\star\mathds{X}$,
\end{enumerate}
such that the diagrams {\bf MA1}, {\bf MA2}, {\bf MA3} commute. The diagram
 \begin{center} \xext=2300 \yext=1150
\begin{picture}(\xext,\yext)(\xoff,\yoff)
\put(800,1000){$\mathds{A}\star(\mathds{B}\star (\mathds{C}\star \mathds{X}))$}
      \put(900,950){\vector(-3,-2){400}}
      \put(1520,840){$\psi_{\mathds{A}, \mathds{B}, \mathds{C}\star \mathds{X}}$} 

      \put(1300,950){\vector(3,-2){400}}
      \put(180,840){$1_\mathds{A}\star\psi_{\mathds{B},\mathds{C},\mathds{X}}$} 

\put(1600,600){$(\mathds{A}\otimes \mathds{B})\star (\mathds{C}\star \mathds{X})$} 
\put(0,600){$\mathds{A}\star ((\mathds{B}\otimes \mathds{C})\star \mathds{X})$}
      \put(500,550){\vector(1,-2){180}}
      \put(1640,330){$\psi_{\mathds{A}\otimes \mathds{B}, \mathds{C},\mathds{X}}$} 

      \put(1700,550){\vector(-1,-2){180}}
      \put(140,330){$\psi_{\mathds{A},\mathds{B}\otimes \mathds{C}, \mathds{X}}$} 

\put(1400,100){$((\mathds{A}\otimes \mathds{B})\otimes \mathds{C})\star \mathds{X}$} 
\put(200,100){$(\mathds{A}\otimes (\mathds{B}\otimes \mathds{C}))\star \mathds{X}$}

  \put(1080,120){\vector(1,0){300}}
      \put(950,0){$\alpha_{\mathds{A},\mathds{B},\mathds{C}}\star 1_{\mathds{X}}$} 
\put(-700,280){$\bf MA1$}
\end{picture}
\end{center}
commutes in $\cA(\cC\times\cC\times\cC\times\cX,\cX)$. The two squares
\begin{center} \xext=900 \yext=600
\begin{picture}(\xext,\yext)(\xoff,\yoff)
\setsqparms[1`1`-1`1;900`450]
 \putsquare(0,50)[\mathds{A}\star \mathds{X}`\mathds{A}\star \mathds{X}`\mathds{I}\star(\mathds{A}\star \mathds{X})`(\mathds{I}\otimes \mathds{A})\star \mathds{X} ;1_{\mathds{A}\star \mathds{X}}`\bar{\psi}_{\mathds{A}\star \mathds{X}}`\lambda_\mathds{A}\star 1_{\mathds{X}}`\psi_{\mathds{I},\mathds{A},\mathds{X}}]
  \put(-1400,280){$\bf MA2$}
  \end{picture}
\end{center}
and
\begin{center} \xext=900 \yext=600
\begin{picture}(\xext,\yext)(\xoff,\yoff)
 \put(-1400,280){$\bf MA3$}
 \setsqparms[1`1`-1`1;900`450]
 \putsquare(0,50)[\mathds{A}\star(\mathds{I}\star \mathds{X})`\mathds{A}\star \mathds{X}`\mathds{A}\star \mathds{X} `(\mathds{A}\otimes \mathds{I})\star \mathds{X};1_{\mathds{A}\star \mathds{X}}`1_{\mathds{A}}\star \bar{\psi}_{\mathds{X}}`\rho_\mathds{A}\star 1_\mathds{X}`\psi_{\mathds{A},\mathds{I},\mathds{X}}]
 \end{picture}
\end{center}
commute in $\cA(\cC\times\cX,\cX)$.  As with monoidal 1-cells, when possible, we will denote the coherence morphism for the unit $\bar{\psi}$ as coherence morphism for the action $\psi$ with bar on the top, and we usually drop $\bar{\psi}$ from the notation. An action $(\star,\psi,\bar{\psi})$ is called {\em strong} ({\em strict}) if both $\psi$ and $\bar{\psi}$ are isomorphisms (identities).

A {\em morphism (1-cell) of monoidal actions}
\[ (F,\phi,\xi): (\cC,\otimes,I,\alpha,\lambda,\rho,\cX,\star, \psi)\ra (\cC',\otimes',I',\alpha',\lambda',\rho',\cX,\star', \psi')\]
consists of
\begin{enumerate}
  \item  a  monoidal 1-cell  $(F,\phi): (\cC,\otimes,I,\alpha,\lambda,\rho)\lra  (\cC',\otimes',I',\alpha',\lambda',\rho')$,
  \item a 2-cell  $\xi:F(\mathds{A})\star'\mathds{X}\ra \mathds{A}\star\mathds{X}$
\end{enumerate}
such that the diagrams {\bf MAF1}, {\bf MAF2} commute. The diagram
  \begin{center} \xext=2000 \yext=1500
\begin{picture}(\xext,\yext)(\xoff,\yoff)
\put(600,1350){$F(\mathds{A})\star'(F(\mathds{B})\star' \mathds{X})$}
      \put(900,1300){\vector(-3,-2){400}}
      \put(100,1180){$\psi'_{F(\mathds{A}),F(\mathds{B}),\mathds{X}}$}

     \put(1300,1300){\vector(3,-2){400}}
    \put(1520,1180){$1_{F(\mathds{A})}\star'\xi_{\mathds{B},\mathds{X}}$}

  \put(0,900){$(F(\mathds{A})\otimes'F(\mathds{B}))\star' \mathds{X}$}
  \put(1600,900){$F(\mathds{A})\star' \mathds{B}\star \mathds{X}$}
     \put(450,850){\vector(0,-1){350}}
     \put(-50,650){$\phi_{\mathds{A},\mathds{B}}\star' 1_{\mathds{X}}$}

    \put(1750,850){\vector(0,-1){350}}
     \put(1800,650){$\xi_{\mathds{A}, \mathds{B}\star \mathds{X}}$}

\put(150,400){$F(\mathds{A}\otimes \mathds{B})\star'\mathds{X}$}
\put(380,130){$\xi_{\mathds{A}\otimes \mathds{B},\mathds{X}}$}

\put(1500,400){$\mathds{A}\star (\mathds{B}\star \mathds{X})$}
\put(800,0){$(\mathds{A}\otimes \mathds{B})\star \mathds{X}$}
\put(1660,350){\vector(-3,-2){400}}

     \put(600,350){\vector(3,-2){400}}

\put(1500,130){$\psi_{\mathds{A},\mathds{B},\mathds{X}}$}
\put(-800,680){$\bf MAF1$}
\end{picture}
\end{center}
commutes in $\cA(\cC\times\cC\times\cX,\cX)$ and the diagram
 \begin{center} \xext=1000 \yext=650
\begin{picture}(\xext,\yext)(\xoff,\yoff)
 \setsqparms[1`1`-1`1;1000`500]
\putsquare(0,50)[\mathds{X}`\mathds{I}\star \mathds{X}`\mathds{I}'\star' \mathds{X}`F(\mathds{I})\star' \mathds{X};\bar{\psi}_{\mathds{X}}`{\bar{\psi}'}_\mathds{X}`\xi_{\mathds{I},\mathds{X}}`\bar{\phi}\star'1_{\mathds{X}}]
\put(-1400,280){$\bf MAF2$}
\end{picture}
\end{center}
commutes in $\cA(\cX,\cX)$.

A {\em transformation (2-cell) of morphisms of actions}
\[ \tau: (F,\phi,\xi)\ra(F',\phi',\xi'): (\cC,\otimes,I,\alpha,\lambda,\rho,\cX,\star, \psi)\ra (\cC',\otimes',I',\alpha',\lambda',\rho',\cX',\star', \psi')\]
consists of a 2-cell $\tau:F\ra F'$ such that the triangle
 \begin{center}\xext=1000 \yext=700
\begin{picture}(\xext,\yext)(\xoff,\yoff)
 \settriparms[1`1`1;500]
  \putVtriangle(0,50)[F(\mathds{A})\star' \mathds{X}`F'(\mathds{A})\star' \mathds{X}`\mathds{A}\star \mathds{X};\tau\star'1_{\mathds{X}}` \xi`\xi']
\put(-1400,280){$\bf MAT$}
\end{picture}
\end{center}
commutes.

In this way we have defined a 2-category $\Act_l\Mon_l(\cA,\cX)$ of lax action of monoidal objects in 2-category $\cA$ on a 0-cell $\cX$ in $\cA$,
with lax monoidal morphisms and monoidal transformations in a 2-category $\cA$ with finite products. When the data is understood from the context, we shall use as compact notation as possible. We can say that $\tau: (F,\phi,\xi)\ra (F',\phi',\xi'):(\star,\psi)\ra (\star',\psi')$ is a transformation of actions if it is understood what monoidal categories $(C,\otimes)$ and $(\cC',\otimes')$ are involved.

\begin{proposition}{\bf Very special (strict) monoidal objects.} Let $\cX$ be exponentiable 0-cell in 2-category $\cA$ with finite products. Then $\cX^\cX$ is a monoidal object in $\cA$. Moreover, $\cX^\cX\times \cX \lra \cX$ is an action. If $\cX$ is rc-0-cell, so is the monoidal object $\cX^\cX$.
\end{proposition}

{\it Proof.}~ The first two statements are routine. Assume that $\cX$ is exponentiable 0-cell in $\cA$. Fix a 0-cell $\cY$ in $\cA$ and a pair of 2-cells $\sigma,\tau:F\ra G: \cY\ra \cX^\cX$ with common inverse $\iota:G\ra F$. Taking the exponential adjunctions of  $\sigma,\tau, \iota$, we get that $\tilde{\sigma}, \tilde{\tau}:\tilde{F}\ra \tilde{G} : \cY\times \cX\ra \cX$ are 2-cells with a common inverse $\tilde{\iota}:\tilde{G}\ra \tilde{F} : \cY\times \cX\ra \cX$. Thus, as $\cX$ is rc, the pair $\tilde{\sigma}, \tilde{\tau}$ has a coequalizer $q:\tilde{G}\ra Q$ in $\cA(\cY\times \cX,\cX)$. Its adjoint $\tilde{q}: G\ra \tilde{Q}: \cY\ra \cX^\cX$ is a coequalizer of $\tau$ and $\sigma$. It is preserved by precomposition, as $q$ is.  $\Box$

\begin{proposition}
Let $(\cM,\otimes)$ be a monoidal object in a 2-category $\cA$ with finite products (or monoidal 2-category), then
the lax slice $\cA_{/_l\cM}$ is again a monoidal 2-category with the domain projection $\cA_{/_l\cM}\ra \cA$  being a strict monoidal functor.
\end{proposition}

{\it Proof.}~ This is a 2-dimensional analog of the fact that a slice of a monoidal category over a monoid is a monoidal 2-category in a canonical way, with the domain projection being a strict monoidal functor. $\Box$

\begin{proposition}\label{prop-Mon-slice-eq-to-Act}
Let $\cX$ be 0-cell in 2-category $\cA$ with finite products such that $\cX^\cX$ exists. Then the 2-categories  $\Mon_l(\cA)_{/_l\cX^\cX}$ and $\Act_l\Mon_l(\cA,\cX)$ are naturally isomorphic. Under this correspondence strong representations are sent to strong actions.
\end{proposition}

{\it Proof.}~ By assumption for any 0-cell $\cE$ we have a (2-natural) isomorphism of categories
\[ \widehat{(-)} : \cA(\cE\times\cX,\cX) \lra \cA(\cE,\cX^\cX)  \]
\[ \tau : F\ra F': \cE\times\cX\ra \cX \longmapsto \widehat{\tau} : \widehat{F}\ra \widehat{F}': \cE\ra \cX^\cX \]
so that we have
\[ ev\circ (\widehat{F}\times 1_\cX)= F, \hskip 5mm ev \circ (\widehat{\tau},1_{1_\cX}) =\tau \]
The 2-cell $ev(=\mathds{H}(\mathds{X})):\cX^\cX\times\cX \ra \cX$ is the usual evaluation morphism, i.e. $\widehat{ev}=1_{\cX^\cX}$. By
\[ \diamond:\cX^\cX\times\cX^\cX\lra \cX^\cX\]
we denote the composition 1-cell, i.e. the adjoint to
\[ ev\circ (1_{\cX^\cX}\times ev): \cX^\cX\times\cX^\cX\times\cX \lra\cX. \]

We shall show that the isomorphisms $\widehat{(-)}$ induce bijective correspondences between 0- 1-, and 2-cells of $\Act_l\Mon_l(\cA,\cX)$ and $\Mon_l(\cA)_{/_l\cX^\cX}$  that respect compositions and identities.

{\em Correspondence of 0-cells.} Below we describe the data of the 0-cells in 2-categories $\Act_l\Mon_l(\cA,\cX)$ and $\Mon_l(\cA)_{/_l\cX^\cX}$ and how they are related. Let us fix a monoidal object $(\cC,\otimes,I,\alpha,\lambda,\rho)$ in $\cA$.
\begin{enumerate}
  \item Let \[ \star=\mathds{A}\star\mathds{X}: \cC\times \cX \ra \cX\]
  be a 1-cell in $\cA$ and
  \[ r=r(\mathds{A}):\cC\lra \cX^\cX\]
  be its exponential adjoint 1-cell, i.e. $r=\widehat{\star}$.

  \item Let \[ \psi: \mathds{A}\star (\mathds{B}\star\mathds{X})\ra (\mathds{A}\otimes \mathds{B})\star\mathds{X}\]
  or in a diagram
\begin{center} \xext=1000 \yext=650
\begin{picture}(\xext,\yext)(\xoff,\yoff)
 \setsqparms[1`1`1`1;1000`500]
\putsquare(0,50)[\cC\times\cC\times\cX`\cC\times\cX`\cC\times\cX`\cX;\otimes\times 1`1\times\star`\star`\star]
\put(700,180){$\Ra$}
\put(700,280){$\psi$}
\end{picture}
\end{center}
be a 2-cell in $\cA$, and
\[ \varphi: r(\mathds{A})\diamond r(\mathds{B})\ra r(\mathds{A}\otimes \mathds{B})\]
or in a diagram
\begin{center} \xext=1000 \yext=650
\begin{picture}(\xext,\yext)(\xoff,\yoff)
 \setsqparms[1`1`1`1;1000`500]
\putsquare(0,50)[\cC\times\cC`\cC`\cX^\cX\times\cX^\cX`\cX^\cX;\otimes`r\times r`r`\diamond]
\put(700,180){$\Ra$}
\put(700,280){$\varphi$}
\end{picture}
\end{center}
be its exponential adjoint 2-cell in $\cA$, i.e. $\varphi=\widehat{\psi}$.

\item Let \[ \bar{\psi}:\mathds{X}\ra \mathds{I}\star\mathds{X}\]
or in a diagram
   \begin{center}\xext=1000 \yext=600
\begin{picture}(\xext,\yext)(\xoff,\yoff)
 \settriparms[1`1`1;600]
  \putqtriangle(0,50)[\cX`\cC\times\cX`\cX;\lk I,1_{\cX}\rk`1_{\cX}`\star]
\put(400,300){$\Ra$}
\put(400,400){$\bar{\psi}$}
\end{picture}
\end{center}
be a 2-cell in $\cA$ and
\[ \bar{\varphi}:\widehat{1_\cX}\ra r(\mathds{I}) \]
or in a diagram
   \begin{center}\xext=1000 \yext=600
\begin{picture}(\xext,\yext)(\xoff,\yoff)
 \settriparms[1`1`1;600]
  \putqtriangle(0,50)[1`\cC`\cX^\cX;I`\widehat{1_{\cX}}`r]
\put(400,300){$\Ra$}
\put(400,400){$\bar{\varphi}$}
\end{picture}
\end{center}
be its exponential adjoint 2-cell in $\cA$, i.e. $\bar{\varphi}=\widehat{\bar{\psi}}$.
\end{enumerate}

We shall show that  \[ (\cC,\otimes,I,\alpha,\lambda,\rho,\cX,\star,\psi)\]
is an action in $\cA$ on $\cX$ i.e. the diagrams {\bf MA1}, {\bf MA2}, {\bf MA3} commute iff
\[ (r,\varphi): (\cC,\otimes,I,\alpha,\lambda,\rho) \lra (\cX^\cX,\diamond,\widehat{1_\cX})\]
is a monoidal 1-cell, i.e. the diagrams {\bf MF1}, {\bf MF2}, {\bf MF3} commute for $(r,\varphi)$, i.e. the diagrams

 \begin{center} \xext=2300 \yext=1150
\begin{picture}(\xext,\yext)(\xoff,\yoff)
\put(-700,280){$\bf MF1'$}
\put(800,1000){$r(\mathds{A})\diamond r(\mathds{B})\diamond r(\mathds{C})$}
      \put(900,950){\vector(-3,-2){400}}
      \put(1520,840){$\psi\diamond 1_{r(\mathds{C})}$} 

      \put(1300,950){\vector(3,-2){400}}
      \put(340,840){$1_{r(\mathds{A})}\diamond\psi$} 

\put(1450,600){$r(\mathds{A}\otimes \mathds{B})\diamond r(\mathds{C})$} 
\put(0,600){$r(\mathds{A})\diamond r(\mathds{B}\otimes \mathds{C}))$}
      \put(500,550){\vector(1,-2){180}}
      \put(1660,330){$\psi$} 

      \put(1700,550){\vector(-1,-2){180}}
      \put(470,330){$\psi$} 

\put(1400,100){$r((\mathds{A}\otimes \mathds{B})\otimes \mathds{C})$} 
\put(200,100){$r(\mathds{A}\otimes (\mathds{B}\otimes \mathds{C}))$}

  \put(1080,120){\vector(1,0){300}}
      \put(950,0){$r(\alpha_{\mathds{A},\mathds{B},\mathds{C}})$} 

\end{picture}
\end{center}
and
\begin{center} \xext=1000 \yext=600
\begin{picture}(\xext,\yext)(\xoff,\yoff)
 \setsqparms[1`1`-1`1;1000`400]
  \putsquare(0,50)[r(\mathds{\mathds{A}})`r(\mathds{\mathds{A}})`r(\mathds{I})\diamond r(\mathds{\mathds{A}})`r(\mathds{I}\otimes \mathds{A});
  1`\bar{\varphi}\diamond 1_{r(\mathds{A})}`r(\lambda_{\mathds{\mathds{A}}})`\varphi]
   \put(-1000,280){$\bf MF2'$}
\end{picture}
\end{center}
and
\begin{center} \xext=1000 \yext=600
\begin{picture}(\xext,\yext)(\xoff,\yoff)
 \setsqparms[1`1`-1`1;1000`400]
  \putsquare(0,50)[r(\mathds{A})`r(\mathds{A})`r(\mathds{A})\diamond r(\mathds{I})`r(\mathds{A}\otimes \mathds{I});1`1_{r(\mathds{A})}\diamond  \bar{\varphi}`r(\rho_{\mathds{A}})`\varphi]
\put(-1000,280){$\bf MF3'$}
\end{picture}
\end{center}
commute.

First we shall show that {\bf MA1} commutes iff ${\bf MF1'}$ does. To this end we shall show that the compositions of three and of two morphisms in these diagrams are adjoint to one another and hence, as $(\widehat{-)}$ preserves and reflects compositions, the commutation of {\bf MA1} will be equivalent to the commutation of {\bf MF1}.
The composition of the three morphisms
\[ (\alpha\star 1_{\mathds{X}})\circ \psi \circ (1_{\mathds{A}}\star \psi) \]
in {\bf MA1} is the composition of the following pasting diagram
  \begin{center} \xext=1600 \yext=2500
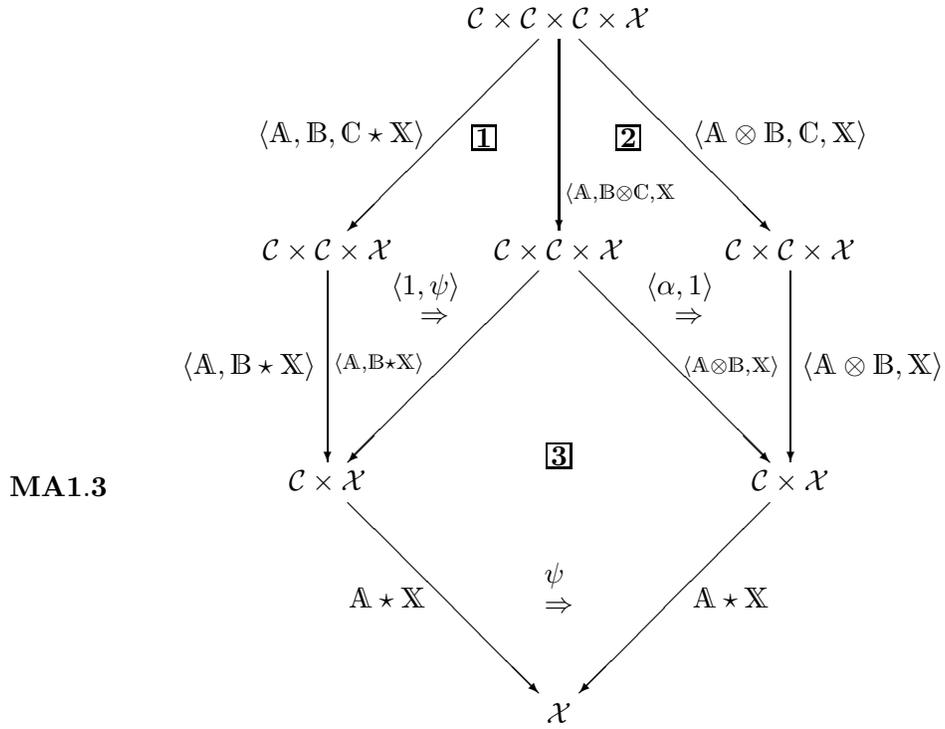
\begin{picture}(\xext,\yext)(\xoff,\yoff)
\put(-1100,800){$\bf \bf MA1.3$}
\put(500,2000){\framebox(80,80){$\bf 1$}}
\put(1000,2000){\framebox(80,80){$\bf 2$}}
\put(760,900){\framebox(80,80){$\bf 3$}}

 \settriparms[1`1`0;800]
  \putAtriangle(0,1650)[\cC\times\cC\times\cC\times\cX`\cC\times\cC\times\cX`\cC\times\cC\times\cX;\lk \mathds{A},\mathds{B},\mathds{C}\star\mathds{X}\rk`\lk \mathds{A}\otimes\mathds{B},\mathds{C},\mathds{X} \rk`]
  \settriparms[1`1`0;800]
  \putAtriangle(0,850)[\cC\times\cC\times\cX`\cC\times\cX`\cC\times\cX;\hskip 5mm^{\lk \mathds{A},\mathds{B}\star\mathds{X}\rk}`\!\!_{\lk\mathds{A}\otimes\mathds{B},\mathds{X}\rk}`]
  \settriparms[0`1`1;800]
  \putVtriangle(0,50)[\phantom{\cC\times\cX}`\phantom{\cC\times\cX}`\cX;`\mathds{A}\star\mathds{X}`\mathds{A}\star\mathds{X}]

  \putmorphism(0,1650)(0,-1)[\phantom{\cC\times\cC\times\cX}`\phantom{\cC\times\cX}`\lk \mathds{A},\mathds{B}\star\mathds{X}\rk]{800}{1}l
  \putmorphism(1600,1650)(0,-1)[\phantom{\cC\times\cC\times\cX}`\phantom{\cC\times\cX}`\lk \mathds{A}\otimes\mathds{B},\mathds{X}\rk]{800}{1}r
  \putmorphism(800,2450)(0,-1)[\phantom{\cC\times\cC\times\cC\times\cX}`\phantom{\cC\times\cC\times\cX}`]{800}{1}l
\put(820,1850){$_{\lk \mathds{A},\mathds{B}\otimes\mathds{C},\mathds{X}}$}

\put(220,1500){$\lk 1,\psi\rk$}
\put(320,1400){$\Ra$}

\put(750,500){$\psi$}
\put(750,400){$\Ra$}

\put(1100,1500){$\lk \alpha,1\rk$}
\put(1200,1400){$\Ra$}
\end{picture}
\end{center}
or in traditional notation
  \begin{center} \xext=1600 \yext=2500
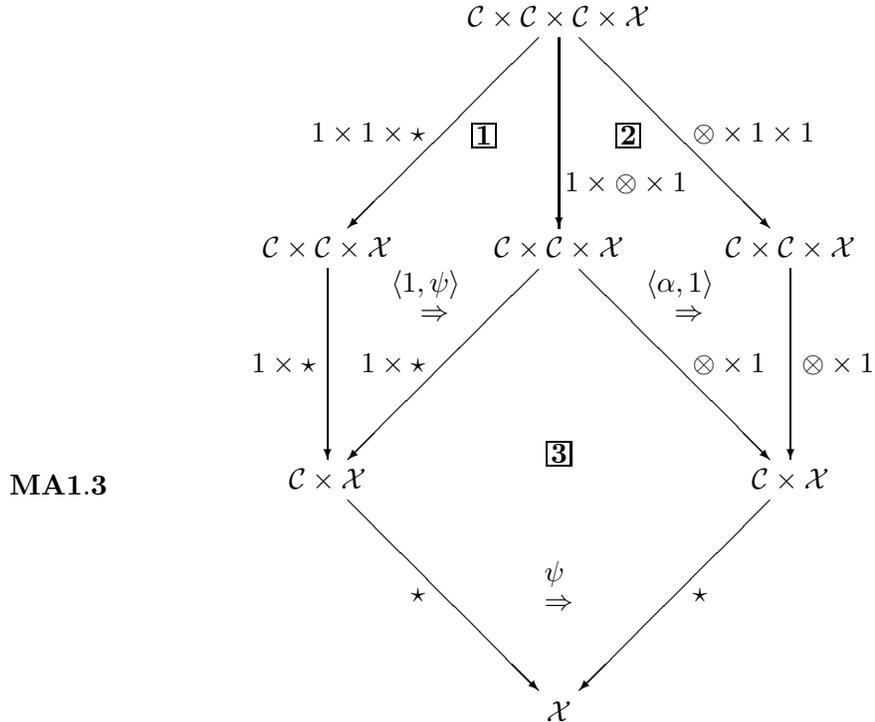
\begin{picture}(\xext,\yext)(\xoff,\yoff)
\put(-1100,800){$\bf \bf MA1.3$}
\put(500,2000){\framebox(80,80){$\bf 1$}}
\put(1000,2000){\framebox(80,80){$\bf 2$}}
\put(760,900){\framebox(80,80){$\bf 3$}}

 \settriparms[1`1`0;800]
  \putAtriangle(0,1650)[\cC\times\cC\times\cC\times\cX`\cC\times\cC\times\cX`\cC\times\cC\times\cX;1\times1\times\star`\otimes\times 1 \times 1`]
  \settriparms[1`1`0;800]
  \putAtriangle(0,850)[\cC\times\cC\times\cX`\cC\times\cX`\cC\times\cX;1\times\star`\otimes \times 1`]
  \settriparms[0`1`1;800]
  \putVtriangle(0,50)[\phantom{\cC\times\cX}`\phantom{\cC\times\cX}`\cX;`\star`\star]

  \putmorphism(0,1650)(0,-1)[\phantom{\cC\times\cC\times\cX}`\phantom{\cC\times\cX}`1\times\star]{800}{1}l
  \putmorphism(1600,1650)(0,-1)[\phantom{\cC\times\cC\times\cX}`\phantom{\cC\times\cX}`\otimes\times 1]{800}{1}r
  \putmorphism(800,2450)(0,-1)[\phantom{\cC\times\cC\times\cC\times\cX}`\phantom{\cC\times\cC\times\cX}`]{800}{1}l
\put(820,1850){$1 \times\otimes\times 1$}

\put(220,1500){$\lk 1,\psi\rk$}
\put(320,1400){$\Ra$}

\put(750,500){$\psi$}
\put(750,400){$\Ra$}

\put(1100,1500){$\lk \alpha,1\rk$}
\put(1200,1400){$\Ra$}

\end{picture}
\end{center}
From now on we shall only use the term notation, leaving the traditional one for the reader.

The composition of the two morphisms
\[ \psi\circ \psi \]
in {\bf MA1} is the composition of the following pasting diagram
  \begin{center} \xext=1600 \yext=2500
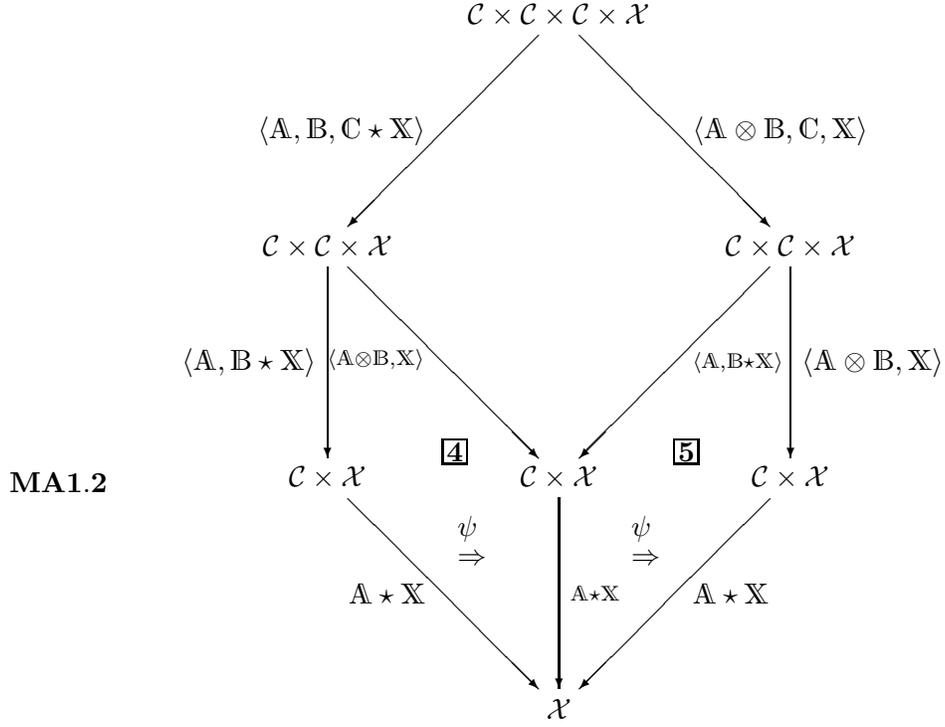
\begin{picture}(\xext,\yext)(\xoff,\yoff)
\put(-1100,800){$\bf MA1.2$}
\put(400,900){\framebox(80,80){$\bf 4$}}
\put(1200,900){\framebox(80,80){$\bf 5$}}

 \settriparms[1`1`0;800]
  \putAtriangle(0,1650)[\cC\times\cC\times\cC\times\cX`\cC\times\cC\times\cX`\cC\times\cC\times\cX;\lk \mathds{A},\mathds{B},\mathds{C}\star\mathds{X}\rk`\lk \mathds{A}\otimes\mathds{B},\mathds{C},\mathds{X}\rk`]
  \settriparms[0`1`1;800]
  \putVtriangle(0,850)[\phantom{\cC\times\cC\times\cX}`\phantom{\cC\times\cC\times\cX}`\cC\times\cX;`^{\lk \mathds{A}\otimes\mathds{B},\mathds{X}\rk}`_{\lk \mathds{A},\mathds{B}\star\mathds{X}\rk}]
  \settriparms[0`1`1;800]
  \putVtriangle(0,50)[\cC\times\cX`\cC\times\cX`\cX;`\mathds{A}\star\mathds{X}`\mathds{A}\star\mathds{X}]

  \putmorphism(0,1650)(0,-1)[\phantom{\cC\times\cC\times\cX}`\phantom{\cC\times\cX}`\lk \mathds{A},\mathds{B}\star\mathds{X}\rk]{800}{1}l
  \putmorphism(1600,1650)(0,-1)[\phantom{\cC\times\cC\times\cX}`\phantom{\cC\times\cX}`\lk \mathds{A}\otimes\mathds{B},\mathds{X}\rk]{800}{1}r
  \putmorphism(800,850)(0,-1)[\phantom{\cC\times\cX}`\phantom{\cX}`_{\mathds{A}\star\mathds{X}}]{800}{1}r

\put(450,650){$\psi$}
\put(450,550){$\Ra$}

\put(1050,650){$\psi$}
\put(1050,550){$\Ra$}
\end{picture}
\end{center}
The composition of the three morphisms
\[ r(\alpha)\circ \varphi \circ(1\diamond\varphi)\]
in  ${\bf MF1'}$ is the composition of the following pasting diagram
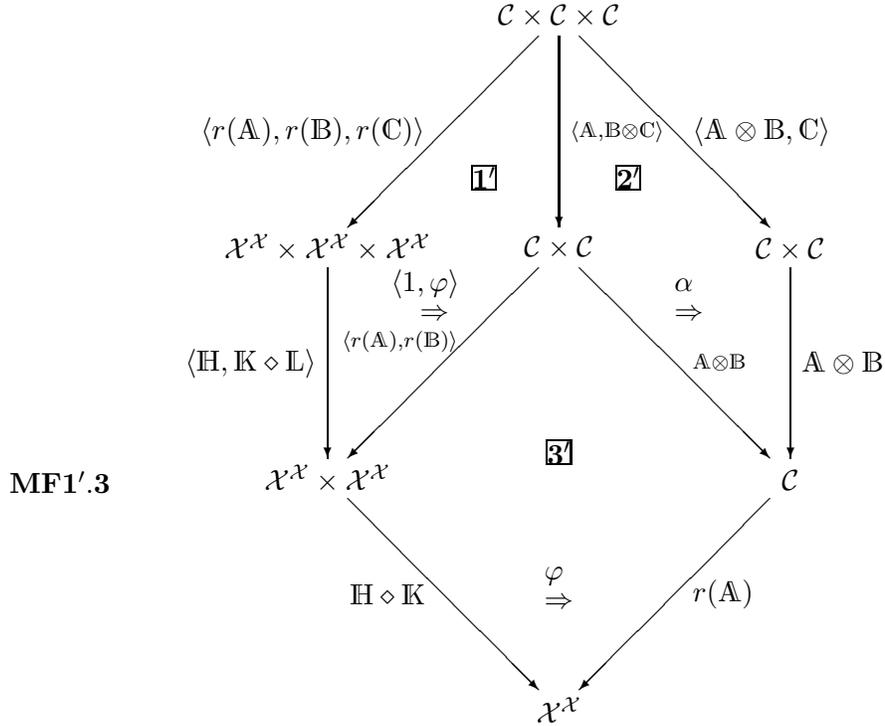
\begin{center} \xext=1600 \yext=2500
\begin{picture}(\xext,\yext)(\xoff,\yoff)
\put(-1100,800){$\bf MF1'.3$}
\put(500,1850){\framebox(80,80){$\bf 1'$}}
\put(1000,1850){\framebox(80,80){$\bf 2'$}}
\put(760,900){\framebox(80,80){$\bf 3'$}}

 \settriparms[1`1`0;800]
  \putAtriangle(0,1650)[\cC\times\cC\times\cC`\cX^\cX\times\cX^\cX\times\cX^\cX`\cC\times\cC;\lk r(\mathds{A}),r(\mathds{B}),r(\mathds{C})\rk`\lk \mathds{A}\otimes\mathds{B},\mathds{C}\rk`]
  \settriparms[1`1`0;800]
  \putAtriangle(0,850)[\cC\times\cC`\cX^\cX\times\cX^\cX`\cC;`_{\mathds{A}\otimes\mathds{B}}`]
  \put(50,1270){$^{\lk r(\mathds{A}),r(\mathds{B})\rk}$}
  \settriparms[0`1`1;800]
  \putVtriangle(0,50)[\phantom{\cX^\cX\times\cX^\cX}`\phantom{\cC}`\cX^\cX;`\mathds{H}\diamond\mathds{K}`r(\mathds{A})]

  \putmorphism(0,1650)(0,-1)[\phantom{\cC\times\cC\times\cX}`\phantom{\cC\times\cX}`\lk \mathds{\mathds{H}},\mathds{K}\diamond\mathds{L}\rk]{800}{1}l
  \putmorphism(800,2450)(0,-1)[\phantom{\cC\times\cC\times\cC\times\cX}`\phantom{\cC\times\cC\times\cX}`_{\lk \mathds{A},\mathds{B}\otimes\mathds{C}\rk}]{800}{1}r
  \putmorphism(1600,1650)(0,-1)[\phantom{\cC\times\cC\times\cX}`\phantom{\cC\times\cX}`\mathds{A}\otimes\mathds{B}]{800}{1}r

\put(220,1500){$\lk 1,\varphi\rk$}
\put(320,1400){$\Ra$}

\put(750,500){$\varphi$}
\put(750,400){$\Ra$}

\put(1200,1500){$\alpha$}
\put(1200,1400){$\Ra$}
\end{picture}
\end{center}
Finally,  the composition of the two morphisms
\[ \varphi \circ (\varphi\diamond 1_r)\]
in  ${\bf MF1'}$ is the composition of the following pasting diagram
  \begin{center} \xext=1600 \yext=2500
 \scalebox{0.9}
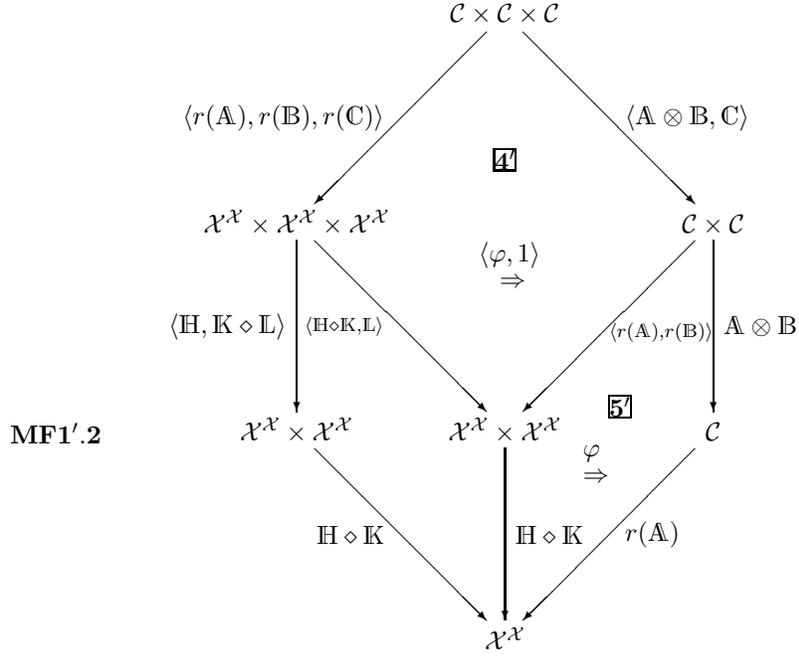
{
\begin{picture}(\xext,\yext)(\xoff,\yoff)
\put(-1100,800){$\bf  MF1'.2$}
\put(760,1850){\framebox(80,80){$\bf 4'$}}
\put(1200,900){\framebox(80,80){$\bf 5'$}}

 \settriparms[1`1`0;800]
  \putAtriangle(0,1650)[\cC\times\cC\times\cC`\cX^\cX\times\cX^\cX\times\cX^\cX`\cC\times\cC;\lk r(\mathds{A}),r(\mathds{B}),r(\mathds{C})\rk`\lk\mathds{A}\otimes\mathds{B},\mathds{C}\rk`]
  \settriparms[0`1`1;800]
  \putVtriangle(0,850)[\phantom{\cX^\cX\times\cX^\cX\times\cX^\cX}`\phantom{\cC\times\cC}`\cX^\cX\times\cX^\cX;`_{\lk \mathds{H}\diamond\mathds{K},\mathds{L}\rk}`]
  \put(1200,1170){$^{\lk r(\mathds{A}),r(\mathds{B})\rk}$}
  \settriparms[0`1`1;800]
  \putVtriangle(0,50)[\cX^\cX\times\cX^\cX`\cC`\cX^\cX;`\mathds{H}\diamond\mathds{K}`r(\mathds{A})]

  \putmorphism(0,1650)(0,-1)[\phantom{\cX^\cX\times\cX^\cX\times\cX^\cX}`\phantom{\cC\times\cX}`\lk \mathds{H},\mathds{K}\diamond\mathds{L}\rk]{800}{1}l
  \putmorphism(800,850)(0,-1)[\phantom{\cX^\cX\times\cX^\cX}`\phantom{\cC\times\cC\times\cX}`\mathds{H}\diamond\mathds{K}]{800}{1}r
  \putmorphism(1600,1650)(0,-1)[\phantom{\cC\times\cC}`\phantom{\cC\times\cX}`\mathds{A}\otimes\mathds{B}]{800}{1}r

\put(700,1500){$\lk \varphi,1\rk$}
\put(780,1400){$\Ra$}

\put(1100,750){$\varphi$}
\put(1100,650){$\Ra$}
\end{picture}
}
\end{center}
Thus we need to show that the composition of the pasting diagrams of three 2-cells  {\bf MA1.3} and ${\bf MF1'.3}$ are adjoint, and moreover that
the composition of the pasting diagrams of two 2-cells  ${\bf MF1'.3}$ and ${\bf MF1'.2}$ are adjoint as well.
We shall show that the (whiskerings of the) corresponding 2-cells are adjoint.

The 2-cell
\begin{center} \xext=1900 \yext=300
\begin{picture}(\xext,\yext)(\xoff,\yoff)
 \putmorphism(0,150)(1,0)[\cC\times\cC\times\cC\times\cX`\cC\times\cX`]{1200}{0}a
 \putmorphism(0,250)(1,0)[\phantom{\cC\times\cC\times\cC\times\cX}`\phantom{\cC\times\cX}`]{1200}{1}a
  \putmorphism(0,50)(1,0)[\phantom{\cC\times\cC\times\cC\times\cX}`\phantom{\cC\times\cX}`]{1200}{1}b
 \putmorphism(1200,150)(1,0)[\phantom{\cC\times\cX}`\cX`\mathds{A}\star\mathds{X}]{700}{1}a
  \put(600,130){$\lk \alpha,1\rk$}

\put(500,130){$\Da$}
\end{picture}
\end{center}
is adjoint to the 2-cell
\begin{center} \xext=1900 \yext=350
\begin{picture}(\xext,\yext)(\xoff,\yoff)
 \putmorphism(0,150)(1,0)[\cC\times\cC\times\cC`\cC`]{1200}{0}a
 \putmorphism(0,250)(1,0)[\phantom{\cC\times\cC\times\cC}`\phantom{\cC}`]{1200}{1}a
  \putmorphism(0,50)(1,0)[\phantom{\cC\times\cC\times\cC}`\phantom{\cC}`]{1200}{1}b
 \putmorphism(1200,150)(1,0)[\phantom{\cC}`\cX^\cX`r(\mathds{A})]{700}{1}a
  \put(600,130){$\alpha$}
\put(500,130){$\Da$}
\end{picture}
\end{center}
This shows the correspondence $\bf 2$ and $\bf 2'$. Similarly, the 2-cell
\begin{center} \xext=2000 \yext=300
\begin{picture}(\xext,\yext)(\xoff,\yoff)
 \putmorphism(0,150)(1,0)[\cC\times\cC\times\cC\times\cX`\cC\times\cX`\lk \mathds{A},\mathds{B}\otimes\mathds{C},\mathds{X}\rk]{1200}{1}a
 \putmorphism(1200,250)(1,0)[\phantom{\cC\times\cX}`\phantom{\cX}`]{800}{1}a
 \putmorphism(1200,150)(1,0)[\phantom{\cC\times\cX}`\cX`]{800}{0}a
 \putmorphism(1200,50)(1,0)[\phantom{\cC\times\cX}`\phantom{\cX}`]{800}{1}b
 \put(1600,130){$\psi$}
\put(1500,130){$\Da$}
\end{picture}
\end{center}
is adjoint to the 2-cell
\begin{center} \xext=1800 \yext=300
\begin{picture}(\xext,\yext)(\xoff,\yoff)
 \putmorphism(0,150)(1,0)[\cC\times\cC\times\cC`\cC`\lk \mathds{A},\mathds{B}\otimes\mathds{C}\rk]{1000}{1}a
 \putmorphism(1000,250)(1,0)[\phantom{\cC\times\cX}`\phantom{\cX^\cX}`]{800}{1}a
 \putmorphism(1000,150)(1,0)[\phantom{\cC\times\cX}`\cX^\cX`]{800}{0}a
 \putmorphism(1000,50)(1,0)[\phantom{\cC\times\cX}`\phantom{\cX^\cX}`]{800}{1}b
 \put(1400,130){$\varphi$}
\put(1300,130){$\Da$}
\end{picture}
\end{center}
This shows the correspondence $\bf 3$ and $\bf 3'$. To see the correspondence of $\bf 1$ and $\bf 1'$, first note that the 2-cell
\begin{center} \xext=1950 \yext=350
\begin{picture}(\xext,\yext)(\xoff,\yoff)
 \putmorphism(0,150)(1,0)[\cC\times\cC\times\cC`\cX^\cX\times\cX^\cX`]{1200}{0}a
 \putmorphism(0,250)(1,0)[\phantom{\cC\times\cC\times\cC}`\phantom{\cX^\cX\times\cX^\cX}`\lk r(\mathds{A}), r(\mathds{B})\diamond r(\mathds{C}) \rk]{1200}{1}a
  \putmorphism(0,50)(1,0)[\phantom{\cC\times\cC\times\cC}`\phantom{\cX^\cX\times\cX^\cX}`\lk r(\mathds{A}), r(\mathds{B}\otimes\mathds{C})\rk]{1200}{1}b
 \putmorphism(1200,150)(1,0)[\phantom{\cX^\cX\times\cX^\cX}`\cX^\cX`\mathds{H}\diamond\mathds{K}]{750}{1}a
  \put(600,130){$\lk 1,\varphi\rk$}
\put(500,130){$\Da$}
\end{picture}
\end{center}
is adjoint to the 2-cells
\begin{center} \xext=3300 \yext=350
\begin{picture}(\xext,\yext)(\xoff,\yoff)
 \putmorphism(0,150)(1,0)[\cC\times\cC\times\cC\times\cX`\cX^\cX\times\cX^\cX\times\cX`]{1800}{0}a
 \putmorphism(0,250)(1,0)[\phantom{\cC\times\cC\times\cC\times\cX}`\phantom{\cX^\cX\times\cX^\cX\times\cX}`\lk r(\mathds{A}), r(\mathds{B})\diamond r(\mathds{C})\rk \times 1_{\mathds{X}}]{1800}{1}a
  \putmorphism(0,50)(1,0)[\phantom{\cC\times\cC\times\cC\times\cX}`\phantom{\cX^\cX\times\cX^\cX\times\cX}`\lk r(\mathds{A}), r(\mathds{B}\otimes\mathds{C})\rk \times 1_{\mathds{X}}]{1800}{1}b
   \putmorphism(1800,150)(1,0)[\phantom{\cX^\cX\times\cX^\cX\times\cX}`\cX^\cX\times\cX`\lk \mathds{H},\mathds{K}(\mathds{X}) \rk]{1000}{1}a
 \putmorphism(2800,150)(1,0)[\phantom{\cX^\cX\times\cX}`\cX`\mathds{H}(\mathds{X})]{500}{1}a
  \put(700,130){$\lk 1,\varphi\rk$}
\put(600,130){$\Da$}
\end{picture}
\end{center}
This 2-cell is equal to
\begin{center} \xext=2600 \yext=350
\begin{picture}(\xext,\yext)(\xoff,\yoff)
 \putmorphism(0,150)(1,0)[\cC\times\cC\times\cC\times\cX`\cX^\cX\times\cX`]{2000}{0}a
 \putmorphism(0,250)(1,0)[\phantom{\cC\times\cC\times\cC\times\cX}`\phantom{\cX^\cX\times\cX}`
 \lk r(\mathds{A}), (r(\mathds{B})\diamond r(\mathds{C}))(\mathds{X})\rk]{2000}{1}a
  \putmorphism(0,50)(1,0)[\phantom{\cC\times\cC\times\cC\times\cX}`\phantom{\cX^\cX\times\cX}`
  \lk r(\mathds{A}), r(\mathds{B}\otimes\mathds{C})(\mathds{X})\rk]{2000}{1}b
 \putmorphism(2000,150)(1,0)[\phantom{\cX^\cX\times\cX}`\cX`\mathds{H}(\mathds{X})]{600}{1}a
  \put(800,130){$\lk 1,\psi\rk$}
\put(700,130){$\Da$}
\end{picture}
\end{center}
and then to the 2-cell
%
%
\begin{center} \xext=3000 \yext=500
\begin{picture}(\xext,\yext)(\xoff,\yoff)
 \putmorphism(0,150)(1,0)[\cC\times\cC\times\cC\times\cX`\cC\times\cX`]{1500}{0}a
 \putmorphism(0,250)(1,0)[\phantom{\cC\times\cC\times\cC\times\cX}`\phantom{\cC\times\cX}`
 \lk \mathds{A}, (r(\mathds{B})\diamond r(\mathds{C}))(\mathds{X})\rk]{1500}{1}a
  \putmorphism(0,50)(1,0)[\phantom{\cC\times\cC\times\cC\times\cX}`\phantom{\cC\times\cX}`
  \lk \mathds{A}, r(\mathds{B}\otimes\mathds{C})(\mathds{X})\rk]{1500}{1}b
   \putmorphism(1500,150)(1,0)[\phantom{\cC\times\cX}`\cX^\cX\times\cX`\lk r(\mathds{A}),\mathds{X}\rk]{900}{1}a
 \putmorphism(2400,150)(1,0)[\phantom{\cX^\cX\times\cX}`\cX`\mathds{H}(\mathds{X})]{600}{1}a
  \put(800,130){$\lk 1,\psi\rk$}
\put(700,130){$\Da$}

 \put(1400,250){\line(0,1){150}}
 \put(1400,400){\line(1,0){1600}}

  \put(3000,400){\vector(0,-1){150}}
   \put(2300,420){$\mathds{A}\star\mathds{X}$}
\end{picture}
\end{center}
This shows the adjointness of cells corresponding to $\bf 1$ and $\bf 1'$. Thus the compositions of pasting diagrams  {\bf MA1.3}  and ${\bf MF1'.3}$  are adjoint.

Now we will show that the compositions of pasting diagrams  {\bf MA1.2}  and ${\bf MF1'.2}$  are adjoint as well. The adjunction of 2-cells corresponding to $\bf 5$ and $\bf 5'$ is easy. We shall show that 2-cells corresponding to $\bf 4$ and $\bf 4'$ are adjoint.

The 2-cell $\diamond\circ \lk \varphi,1\rk$ is equal to
\begin{center} \xext=3200 \yext=350
\begin{picture}(\xext,\yext)(\xoff,\yoff)
 \putmorphism(0,150)(1,0)[\cC\times\cC\times\cC`\cX^\cX\times\cC`]{1400}{0}a
 \putmorphism(0,250)(1,0)[\phantom{\cC\times\cC\times\cC}`\phantom{\cX^\cX\times\cX^\cX}`\lk r(\mathds{A}\diamond r(\mathds{B}), \mathds{C}\rk]{1400}{1}a
  \putmorphism(0,50)(1,0)[\phantom{\cC\times\cC\times\cC}`\phantom{\cX^\cX\times\cX^\cX}`\lk r(\mathds{A}\otimes\mathds{B}),\mathds{C})\rk ]{1400}{1}b
   \putmorphism(1400,150)(1,0)[\phantom{\cX^\cX\times\cC}`\cX^\cX\times\cX^\cX`\lk \mathds{H},r(\mathds{A}')\rk]{900}{1}a
 \putmorphism(2300,150)(1,0)[\phantom{\cX^\cX\times\cX^\cX}`\cX^\cX`\mathds{H}'\diamond\mathds{K}]{600}{1}a
  \put(700,130){$\lk \varphi,1\rk$}
\put(600,130){$\Da$}
\end{picture}
\end{center}
Its adjoint 2-cell is
%
%
  \begin{center} \xext=1600 \yext=3400
  \scalebox{0.9}
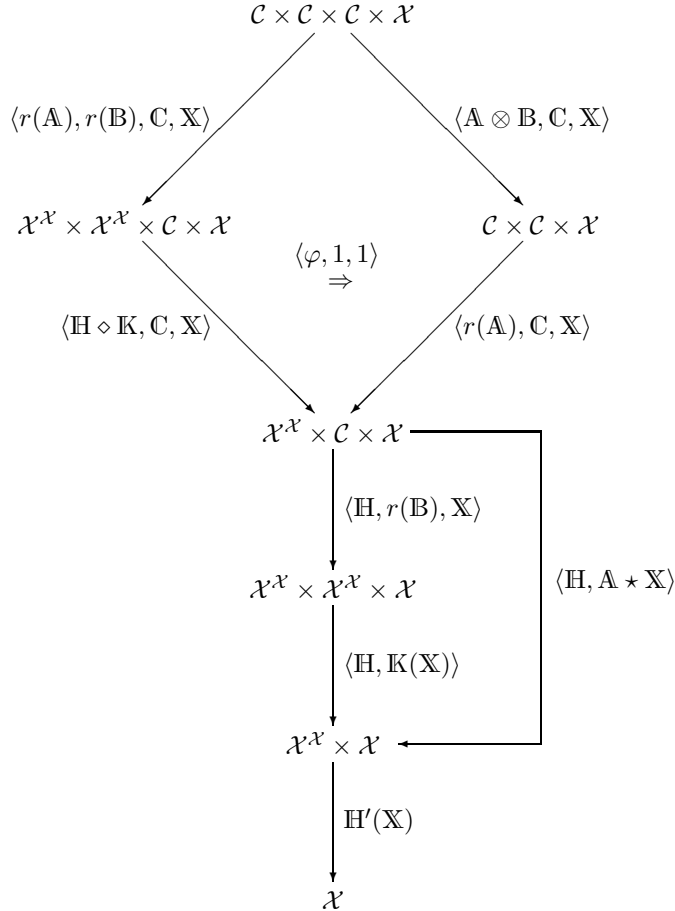
{
\begin{picture}(\xext,\yext)(\xoff,\yoff)
 \settriparms[1`1`0;800]
  \putAtriangle(0,2600)[\cC\times\cC\times\cC\times\cX`\cX^\cX\times\cX^\cX\times\cC\times\cX`\cC\times\cC\times\cX;\lk r(\mathds{A}), r(\mathds{B}),\mathds{C},\mathds{X\rk}`\lk \mathds{A}\otimes\mathds{B},\mathds{C},\mathds{X}\rk`]
  \settriparms[0`1`1;800]
  \putVtriangle(0,1800)[\phantom{\cX^\cX\times\cX^\cX\times\cC\times\cX}`\phantom{\cC\times\cC\times\cX}`\cX^\cX\times\cC\times\cX;`\lk \mathds{H}\diamond\mathds{K},\mathds{C},\mathds{X}\rk`\lk r(\mathds{A}), \mathds{C}, \mathds{X}\rk]

\put(650,2450){$\lk \varphi,1,1\rk$}
\put(780,2350){$\Ra$}

\putmorphism(800,1800)(0,-1)[\phantom{\cX^\cX\times\cC\times\cX}`\cX^\cX\times\cX^\cX\times\cX`\lk \mathds{H},r(\mathds{B}),\mathds{X}\rk]{600}{1}r
\putmorphism(800,1200)(0,-1)[\phantom{\cX^\cX\times\cX^\cX\times\cX}`\cX^\cX\times\cX`\lk \mathds{H},\mathds{K}(\mathds{X})\rk]{600}{1}r
\putmorphism(800,600)(0,-1)[\phantom{\cX^\cX\times\cX}`\cX`\mathds{H}'(\mathds{X})]{600}{1}r

 \put(1100,1800){\line(1,0){500}}
 \put(1600,1800){\line(0,-1){1200}}

  \put(1600,600){\vector(-1,0){550}}
   \put(1650,1200){$\lk \mathds{H},\mathds{A}\star\mathds{X}\rk$}
\end{picture}
}
\end{center}
The above 2-cell is equal to
  \begin{center} \xext=3000 \yext=2900
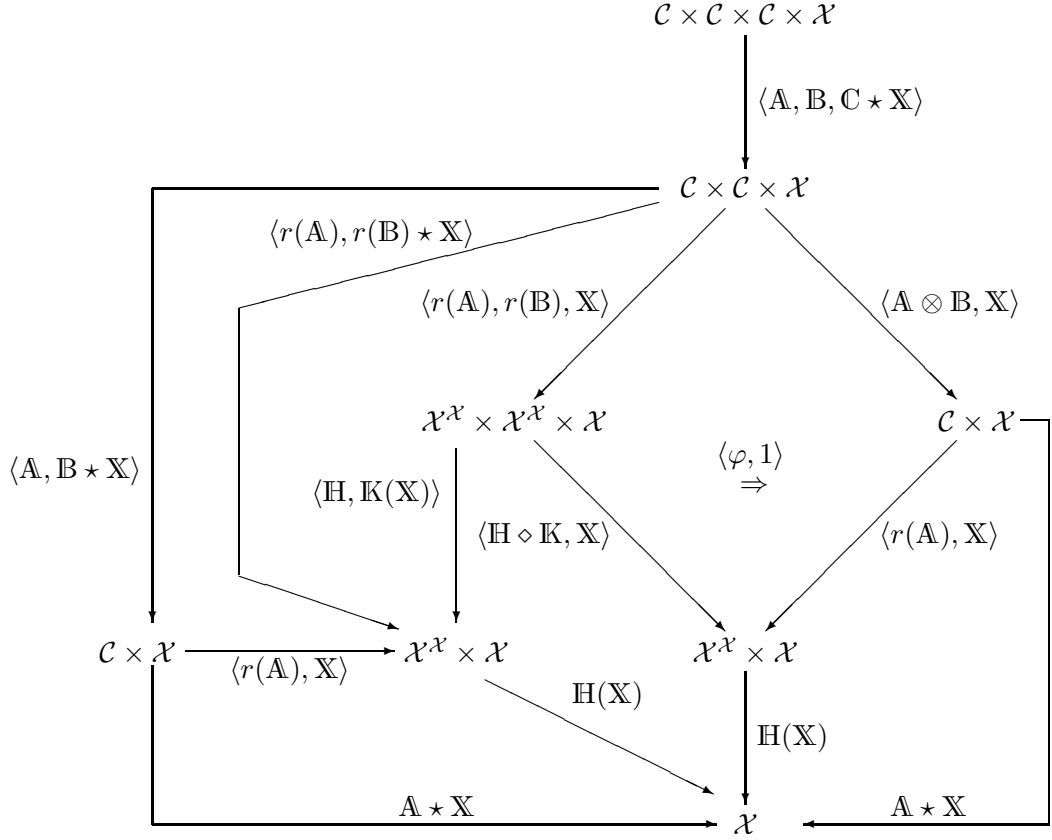
\begin{picture}(\xext,\yext)(\xoff,\yoff)

\putmorphism(2100,2800)(0,-1)[\cC\times\cC\times\cC\times\cX`\phantom{\cC\times\cC\times\cX}`\lk \mathds{A},\mathds{B},\mathds{C}\star\mathds{X}\rk]{600}{1}r

 \settriparms[1`1`0;800]
  \putAtriangle(1300,1400)[\cC\times\cC\times\cX`\cX^\cX\times\cX^\cX\times\cX`\cC\times\cX;\lk r(\mathds{A}), r(\mathds{B}),\mathds{X}\rk`\lk \mathds{A}\otimes\mathds{B},\mathds{X}\rk`]
  \settriparms[0`1`1;800]
  \putVtriangle(1300,600)[\phantom{\cX^\cX\times\cX^\cX\times\cC\times\cX}`\phantom{\cC\times\cC\times\cX}`\cX^\cX\times\cX;`\lk \mathds{H}\diamond\mathds{K},\mathds{X}\rk`\lk r(\mathds{A}),\mathds{X}\rk]

\put(2000,1250){$\lk \varphi,1\rk$}
\put(2080,1150){$\Ra$}

\putmorphism(0,600)(1,0)[\cC\times\cX`\cX^\cX\times\cX`\lk r(\mathds{A}),\mathds{X}\rk]{1100}{1}r
\putmorphism(2100,600)(0,-1)[\phantom{\cX^\cX\times\cX}`\cX`\mathds{H}(\mathds{X})]{600}{1}r

 \put(1800,2200){\line(-1,0){1750}}
  \put(50,2200){\vector(0,-1){1500}}
   \put(-450,1200){$\lk \mathds{A},\mathds{B}\star\mathds{X}\rk$}
 \put(50,550){\line(0,-1){550}}
 \put(50,0){\vector(1,0){1950}}
 \put(900,30){$\mathds{A}\star\mathds{X}$}
 \put(1800,2150){\line(-4,-1){1450}}
 \put(350,1785){\line(0,-1){920}}
 \put(350,860){\vector(3,-1){550}}
 \put(450,2020){$\lk r(\mathds{A}), r(\mathds{B})\star\mathds{X}\rk$}
\put(1100,1300){\vector(0,-1){600}}
\put(600,1120){$\lk \mathds{H},\mathds{K}(\mathds{X})\rk$}
\put(1200,500){\vector(2,-1){790}}
\put(1500,420){$\mathds{H}(\mathds{X})$}
 \put(3050,1400){\line(1,0){100}}
 \put(3150,1400){\line(0,-1){1400}}
 \put(3150,0){\vector(-1,0){850}}
 \put(2600,30){$\mathds{A}\star\mathds{X}$}
\end{picture}
\end{center}
with $ev\circ \lk \varphi,1\rk=\psi$. Thus the 2-cells corresponding to $\bf 4$ and $\bf 4'$ are adjoint and hence $\bf MF1'$ commutes iff $\bf MA1$ does.

Now we shall show that {\bf MA2} commutes iff ${\bf MF2'}$ does. The proof that {\bf MA3} commutes iff ${\bf MF3'}$ does, is similar.
Clearly, $1_{\mathds{A}\star\mathds{X}}$ is adjoint to $1_{r{\mathds{A}}}$.

The composition 2-cell
\[ (\lambda_{\mathds{A}}\star 1_{\mathds{X}})\circ  \psi_{\mathds{I},\mathds{A},\mathds{X}} \circ \bar{\psi}_{\mathds{A}\star\mathds{X}} \]
is the composition of the following pasting diagram
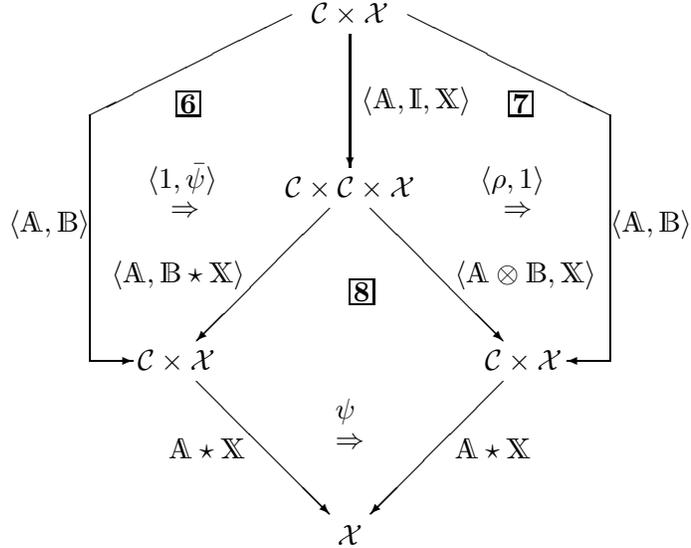
\begin{center} \xext=1800 \yext=1900
\begin{picture}(\xext,\yext)(\xoff,\yoff)

\putmorphism(900,1850)(0,-1)[\cC\times\cX`\phantom{\cC\times\cC\times\cX}`\lk \mathds{A}, \mathds{I}, \mathds{X}\rk]{600}{1}r

 \settriparms[1`1`0;600]
  \putAtriangle(300,650)[\cC\times\cC\times\cX`\cC\times\cX`\cC\times\cX;
  \lk \mathds{A}, \mathds{B}\star\mathds{X}\rk`\lk\mathds{A}\otimes\mathds{B},\mathds{X}\rk` ]
  \settriparms[0`1`1;600]
  \putVtriangle(300,50)[\phantom{\cX^\cX\times\cX^\cX\times\cC\times\cX}`\phantom{\cC\times\cC\times\cX}`\cX;
  `\mathds{A}\star\mathds{X}`\mathds{A}\star\mathds{X}]

\put(200,1250){$\lk 1,\bar{\psi}\rk$}
\put(280,1150){$\Ra$}

\put(850,450){$\psi$}
\put(850,350){$\Ra$}

\put(1350,1250){$\lk \rho,1\rk$}
\put(1430,1150){$\Ra$}

 \put(700,1850){\line(-2,-1){700}}
 \put(0,1500){\line(0,-1){850}}
  \put(0,650){\vector(1,0){150}}
   \put(-280,1100){$\lk \mathds{A},\mathds{B}\rk$}
  \put(1100,1850){\line(2,-1){700}}
 \put(1800,1500){\line(0,-1){850}}
  \put(1800,650){\vector(-1,0){150}}

 \put(1800,1100){$\lk \mathds{A},\mathds{B}\rk$}

 \put(300,1500){\framebox(80,80){$\bf 6$}}
  \put(1450,1500){\framebox(80,80){$\bf 7$}}
   \put(900,850){\framebox(80,80){$\bf 8$}}
\end{picture}
\end{center}
The composition 2-cell
\[ r(\lambda_{\mathds{A}})\circ  \varphi_{\mathds{I},\mathds{A}} \circ (\bar{\varphi}\diamond 1_{r(\mathds{A})}) \]
is the composition of the following pasting diagram
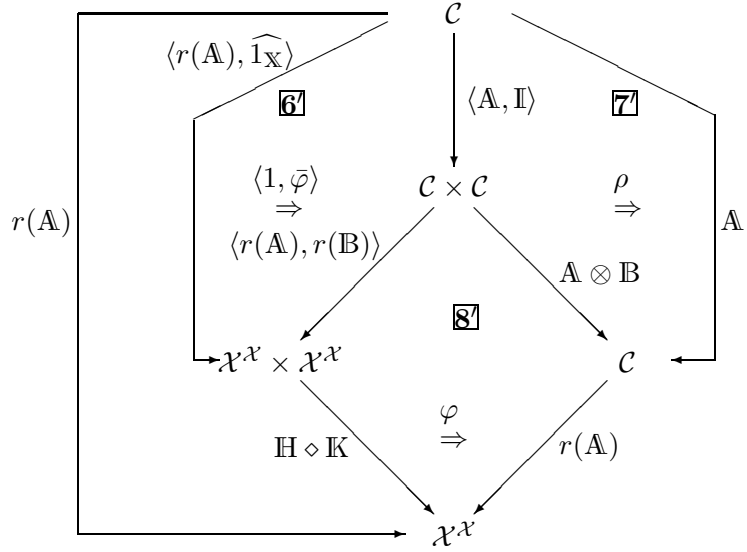
\begin{center} \xext=1800 \yext=1900
\begin{picture}(\xext,\yext)(\xoff,\yoff)

\putmorphism(1350,1850)(0,-1)[\cC`\phantom{\cC\times\cC}`\lk \mathds{A}, \mathds{I}\rk]{600}{1}r

 \settriparms[1`1`0;600]
  \putAtriangle(750,650)[\cC\times\cC`\cX^\cX\times\cX^\cX`\cC; `\mathds{A}\otimes\mathds{B}` ]
   \put(570,1020){$\lk r(\mathds{A}), r(\mathds{B})\rk$}
  \settriparms[0`1`1;600]
  \putVtriangle(750,50)[\phantom{\cX^\cX\times\cX^\cX}`\phantom{\cC\times\cC}`\cX^\cX;  `\mathds{H}\diamond \mathds{K}`r(\mathds{A})]

\put(1300,450){$\varphi$}
\put(1300,350){$\Ra$}

\put(650,1250){$\lk 1,\bar{\varphi}\rk$}
\put(730,1150){$\Ra$}

\put(1900,1250){$\rho$}
\put(1900,1150){$\Ra$}

 \put(1120,1850){\line(-1,0){1070}}
 \put(50,1850){\line(0,-1){1800}}
  \put(50,50){\vector(1,0){1130}}
   \put(-180,1100){$r(\mathds{A})$}
 \put(1120,1820){\line(-2,-1){670}}
 \put(450,1480){\line(0,-1){830}}
  \put(450,650){\vector(1,0){90}}
   \put(350,1680){$\lk r(\mathds{A}), \widehat{1_{\mathds{X}}}\rk$}
  \put(1550,1850){\line(2,-1){700}}
 \put(2250,1500){\line(0,-1){850}}
  \put(2250,650){\vector(-1,0){150}}

 \put(2270,1100){$\mathds{A}$}

 \put(750,1500){\framebox(80,80){$\bf 6'$}}
  \put(1900,1500){\framebox(80,80){$\bf 7'$}}
   \put(1350,760){\framebox(80,80){$\bf 8'$}}
\end{picture}
\end{center}
The 2-cells corresponding to ${\bf 7}$ and ${\bf 8}$ are clearly adjoint to the 2-cells corresponding to ${\bf 7'}$ and ${\bf 8'}$, respectively.

The 2-cell $1\diamond\bar{\varphi}: r(\mathds{A})\ra r(\mathds{A})\diamond r(\mathds{I}):\cC\ra \cX^\cX$ corresponding to ${\bf 6'}$ is adjoint to the  the composite 2-cell of the following pasting diagram
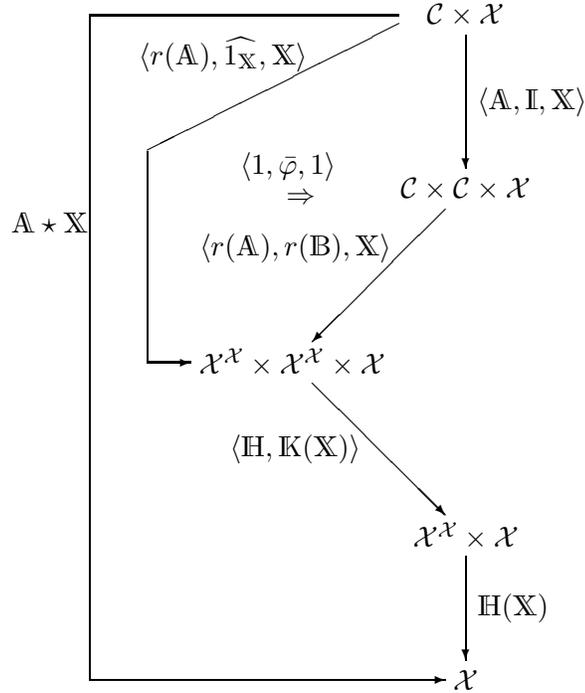
\begin{center} \xext=1900 \yext=2400
\begin{picture}(\xext,\yext)(\xoff,\yoff)

\putmorphism(1350,2350)(0,-1)[\cC\times\cX`\phantom{\cC\times\cC\times\cX}`\lk \mathds{A}, \mathds{I},\mathds{X}\rk]{600}{1}r

 \settriparms[1`0`0;600]
  \putAtriangle(750,1150)[\cC\times\cC\times\cX`\cX^\cX\times\cX^\cX\times\cX`; `` ]
   \put(430,1520){$\lk r(\mathds{A}), r(\mathds{B}),\mathds{X}\rk$}
  \settriparms[0`1`0;600]
  \putVtriangle(750,550)[\phantom{\cX^\cX\times\cX^\cX\times\cC}`\phantom{\cC\times\cC}`\cX^\cX\times\cX;  `\lk \mathds{H},\mathds{K}(\mathds{X})\rk`]

\put(570,1800){$\lk 1,\bar{\varphi},1\rk$}
\put(730,1700){$\Ra$}

 \put(1120,2350){\line(-1,0){1070}}
 \put(50,2350){\line(0,-1){2300}}
  \put(50,50){\vector(1,0){1230}}
   \put(-220,1600){$\mathds{A}\star\mathds{X}$}
 \put(1120,2320){\line(-2,-1){870}}
 \put(250,1880){\line(0,-1){730}}
  \put(250,1150){\vector(1,0){150}}
   \put(220,2180){$\lk r(\mathds{A}), \widehat{1_{\mathds{X}}},\mathds{X}\rk$}

\putmorphism(1350,550)(0,-1)[\phantom{\cX^\cX\times\cX}`\cX`\mathds{H}(\mathds{X})]{500}{1}r
\end{picture}
\end{center}
The above 2-cell is equal to the 2-cell
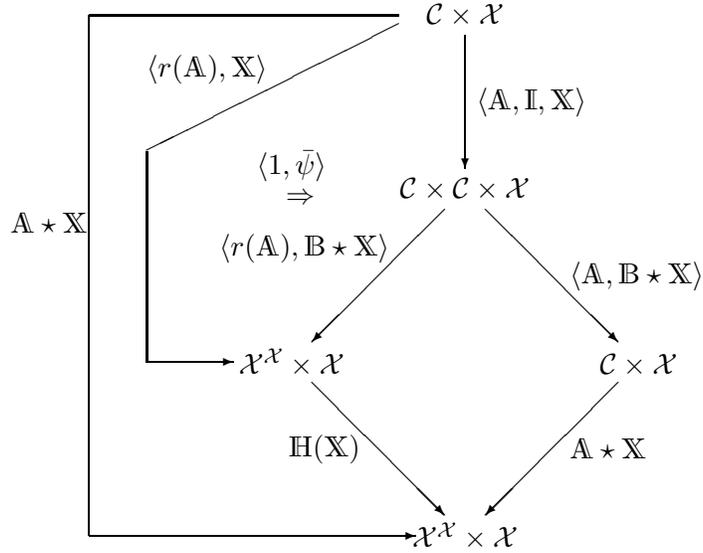
\begin{center} \xext=1800 \yext=1900
\begin{picture}(\xext,\yext)(\xoff,\yoff)

\putmorphism(1350,1850)(0,-1)[\cC\times\cX`\phantom{\cC\times\cC\times\cX}`\lk \mathds{A}, \mathds{I},\mathds{X}\rk]{600}{1}r

 \settriparms[1`1`0;600]
  \putAtriangle(750,650)[\cC\times\cC\times\cX`\cX^\cX\times\cX`\cC\times\cX; `\lk \mathds{A},\mathds{B}\star \mathds{X}\rk` ]
   \put(500,1020){$\lk r(\mathds{A}), \mathds{B}\star\mathds{X}\rk$}
  \settriparms[0`1`1;600]
  \putVtriangle(750,50)[\phantom{\cX^\cX\times\cX^\cX\times\cC}`\phantom{\cC\times\cC}`\cX^\cX\times\cX;  `\mathds{H}(\mathds{X})`\mathds{A}\star\mathds{X}]

\put(630,1300){$\lk 1,\bar{\psi}\rk$}
\put(730,1200){$\Ra$}

 \put(1120,1850){\line(-1,0){1070}}
 \put(50,1850){\line(0,-1){1800}}
  \put(50,50){\vector(1,0){1130}}
   \put(-220,1100){$\mathds{A}\star\mathds{X}$}
 \put(1120,1820){\line(-2,-1){870}}
 \put(250,1380){\line(0,-1){730}}
  \put(250,650){\vector(1,0){300}}
   \put(250,1640){$\lk r(\mathds{A}),\mathds{X}\rk$}
\end{picture}
\end{center}
This ends the proof of bijective correspondence of 0-cells.

{\em Correspondence of 1-cells.} Let $(\star,\psi)$, $(\star',\psi')$  be actions of monoidal objects $(\cC,\otimes,I,\alpha,\lambda,\rho)$, $(\cC',\otimes',I',\alpha',\lambda',\rho')$, respectively, on $\cX$, i.e. 0-cells in $\Mon_l(\cA)_{/_l\cX^\cX}$ and let $(r,\varphi): (\cC,\otimes)\ra (\cX^\cX,\diamond)$, $(r',\varphi'): (\cC',\otimes')\ra (\cX^\cX,\diamond)$ be corresponding monoidal 1-cells, i.e. 0-cells in $\Act_l\Mon_l(\cA,\cX)$.
Let $(F,\phi): (\cC,\otimes,I,\alpha,\lambda,\rho)\ra (\cC,\otimes,I,\alpha,\lambda,\rho)$ be monoidal 1-cells. Let the 2-cells
\[ \xi : F(\mathds{A})\star'\mathds{\mathds{X}}\ra \mathds{A}\star\mathds{X}:\cC\times\cX\ra \cX,\hskip 5mm
\theta : r'F(\mathds{A})\ra r(\mathds{A}):\cC\ra \cX^\cX \]
i.e. 2-cells in the diagrams
  \begin{center} \xext=2200 \yext=550
\begin{picture}(\xext,\yext)(\xoff,\yoff)
 \settriparms[1`1`-1;400]
  \putVtriangle(0,50)[\cC\times\cX`\cC'\times\cX`\cX ;F\times 1`\star`\star']
  \put(360,280){$\xi$}
\put(350,200){$\Ra$}

 \settriparms[1`1`-1;400]
  \putVtriangle(1400,50)[\cC`\cC'`\cX^\cX ;F`r`r']
  \put(1760,280){$\theta$}
\put(1750,200){$\Ra$}
\end{picture}
\end{center}
be adjoint. We need to show that the 1-cell
\[ \xi: (\star,\psi)\ra (\star',\psi')\]
in $\Act_l\Mon_l(\cA,\cX)$, i.e. satisfies $\bf MAF1$  and $\bf MAF2$  iff the 1-cell
\[ \theta: (r'F,r'(\phi)\circ\varphi'_{F\times F},r'(\bar{\phi})\circ \bar{\varphi}')\ra (r,\varphi)\]
in $\Mon_l(\cA)_{/_l\cX^\cX}$ i.e. satisfies $\bf MT1$ and $\bf MT2$.

The commutation of $\bf MF1$ means that the pasting diagrams (corresponding to the two composable sequences of morphisms)
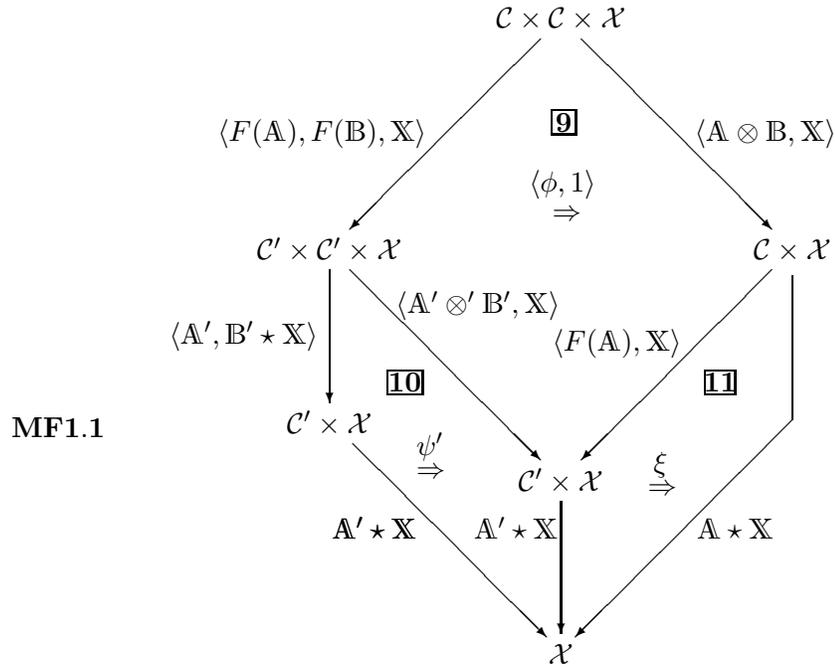
\begin{center} \xext=1600 \yext=2500
\begin{picture}(\xext,\yext)(\xoff,\yoff)
\put(-1100,800){$\bf MF1.1$}
\put(770,1850){\framebox(80,80){$\bf 9$}}
\put(200,950){\framebox(120,80){$\bf 10$}}
\put(1300,950){\framebox(100,80){$\bf 11$}}

\settriparms[1`1`0;800]
\putAtriangle(0,1450)[\cC\times\cC\times\cX`\cC'\times\cC'\times\cX`\cC\times\cX;\lk F(\mathds{A}),F(\mathds{B}),\mathds{X}\rk`\lk \mathds{A}\otimes\mathds{B},\mathds{X}\rk`]

\settriparms[0`1`1;800]
\putVtriangle(0,650)[\phantom{\cC'\times\cC'\times\cX}`\phantom{\cC\times\cX}`\cC'\times\cX;`` ]

  \put(230,1230){$\lk \mathds{A}'\otimes'\mathds{B}',\mathds{X}\rk$}
  \put(770,1100){$\lk F(\mathds{A}),\mathds{X}\rk$}

  \putmorphism(800,650)(0,-1)[\phantom{\cC'\times\cX}`\cX`]{600}{1}l
  \put(500,450){$\mathds{A}'\star\mathds{X}$}

  \putmorphism(0,1450)(0,-1)[\phantom{\cC'\times\cC'\times\cX}`\cC'\times\cX`\lk \mathds{A}',\mathds{B}'\star\mathds{X}\rk]{600}{1}l

\put(690,1650){$\lk \phi,1\rk$}
\put(770,1550){$\Ra$}

\put(300,740){$\psi'$}
\put(300,660){$\Ra$}

\put(1120,680){$\xi$}
\put(1100,600){$\Ra$}

  \put(80,780){\vector(1,-1){670}}
   \put(10,450){$\mathds{A}'\star\mathds{X}$}

   \put(1600,860){\line(0,1){500}}
    \put(1600,860){\vector(-1,-1){750}}
   \put(10,450){$\mathds{A}'\star\mathds{X}$}
    \put(1270,450){$\mathds{A}\star\mathds{X}$}
\end{picture}
\end{center}
and
\begin{center} \xext=1600 \yext=2500
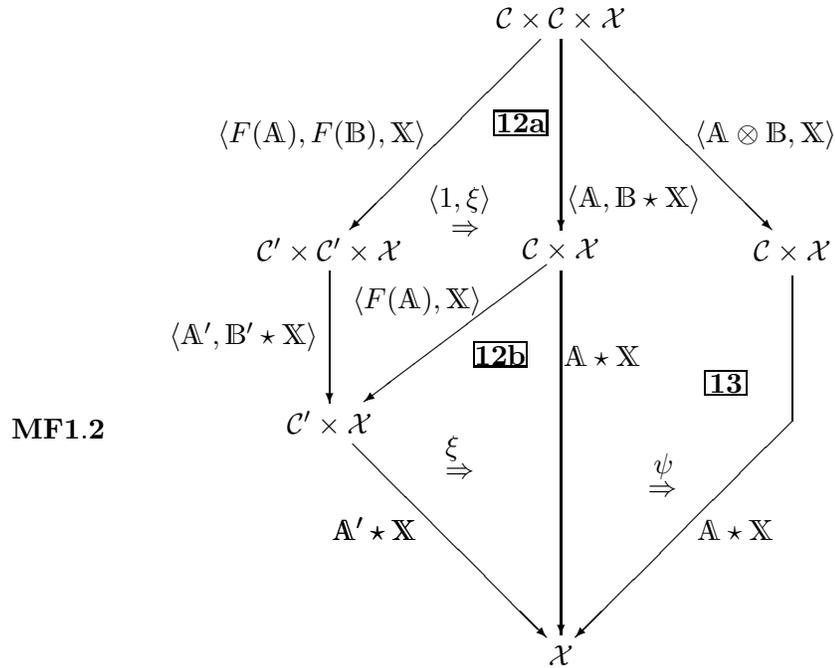
\begin{picture}(\xext,\yext)(\xoff,\yoff)
\put(-1100,800){$\bf MF1.2$}
\put(570,1850){\framebox(180,80){$\bf 12a$}}
\put(500,1050){\framebox(180,80){$\bf 12b$}}
\put(1300,950){\framebox(140,80){$\bf 13$}}

\settriparms[1`1`0;800]
\putAtriangle(0,1450)[\cC\times\cC\times\cX`\cC'\times\cC'\times\cX`\cC\times\cX;\lk F(\mathds{A}),F(\mathds{B}),\mathds{X}\rk`\lk \mathds{A}\otimes\mathds{B},\mathds{X}\rk`]

\putmorphism(800,2250)(0,-1)[\phantom{\cC\times\cC\times\cX}`\cC\times\cX`]{800}{1}r
\put(820,1600){$\lk\mathds{A},\mathds{B}\star\mathds{X}\rk$}

\putmorphism(800,1450)(0,-1)[\phantom{\cC\times\cX}`\phantom{\cX}`]{1400}{1}r
\put(810,1050){$\mathds{A}\star\mathds{X}$}

  \putmorphism(800,650)(0,-1)[\phantom{\cC'\times\cX}`\cX`]{600}{0}l

  \putmorphism(0,1450)(0,-1)[\phantom{\cC'\times\cC'\times\cX}`\cC'\times\cX`\lk \mathds{A}',\mathds{B}'\star\mathds{X}\rk]{600}{1}l

\put(340,1600){$\lk 1,\xi\rk$}
\put(420,1500){$\Ra$}

\put(400,740){$\xi$}
\put(400,660){$\Ra$}

\put(1120,680){$\psi$}
\put(1100,600){$\Ra$}

  \put(80,780){\vector(1,-1){670}}
   \put(10,450){$\mathds{A}'\star\mathds{X}$}

   \put(1600,860){\line(0,1){500}}
    \put(1600,860){\vector(-1,-1){750}}
   \put(10,450){$\mathds{A}'\star\mathds{X}$}
    \put(1270,450){$\mathds{A}\star\mathds{X}$}

     \put(750,1400){\vector(-4,-3){630}}
   \put(80,1250){$\lk F(\mathds{A}),\mathds{X}\rk$}
\end{picture}
\end{center}
compose to the same 2-cell. The commutation of $\bf MT1$ means that the pasting diagrams
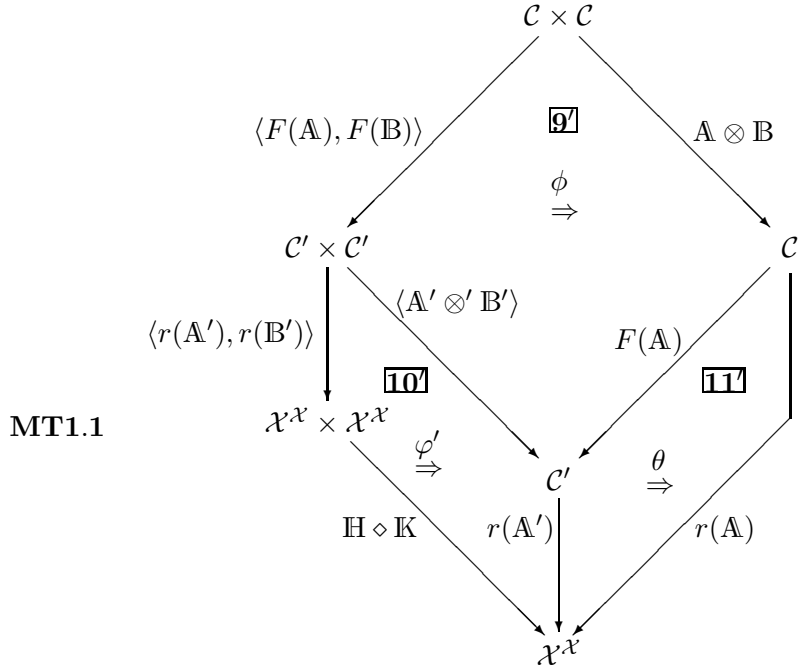
\begin{center} \xext=1600 \yext=2500
\begin{picture}(\xext,\yext)(\xoff,\yoff)
\put(-1100,800){$\bf MT1.1$}
\put(770,1850){\framebox(90,80){$\bf 9'$}}
\put(200,950){\framebox(140,80){$\bf 10'$}}
\put(1300,950){\framebox(140,80){$\bf 11'$}}

\settriparms[1`1`0;800]
\putAtriangle(0,1450)[\cC\times\cC`\cC'\times\cC'`\cC;\lk F(\mathds{A}),F(\mathds{B})\rk`\mathds{A}\otimes\mathds{B}`]

\settriparms[0`1`1;800]
\putVtriangle(0,650)[\phantom{\cC'\times\cC'}`\phantom{\cC}`\cC';`` ]

  \put(230,1230){$\lk \mathds{A}'\otimes'\mathds{B}'\rk$}
  \put(990,1100){$F(\mathds{A})$}

  \putmorphism(800,650)(0,-1)[\phantom{\cC'\times\cX}`\cX^\cX`]{600}{1}l
  \put(550,450){$r(\mathds{A}')$}

  \putmorphism(0,1450)(0,-1)[\phantom{\cC'\times\cC'}`\cX^\cX\times\cX^\cX`\lk r(\mathds{A}'),r(\mathds{B}')\rk]{600}{1}l

\put(770,1650){$\phi$}
\put(770,1550){$\Ra$}

\put(300,740){$\varphi'$}
\put(300,660){$\Ra$}

\put(1120,680){$\theta$}
\put(1100,600){$\Ra$}

  \put(80,780){\vector(1,-1){670}}
   \put(50,450){$\mathds{H}\diamond\mathds{K}$}

   \put(1600,860){\line(0,1){500}}
    \put(1600,860){\vector(-1,-1){750}}
    \put(1270,450){$r(\mathds{A})$}
\end{picture}
\end{center}
and
\begin{center} \xext=1600 \yext=2500
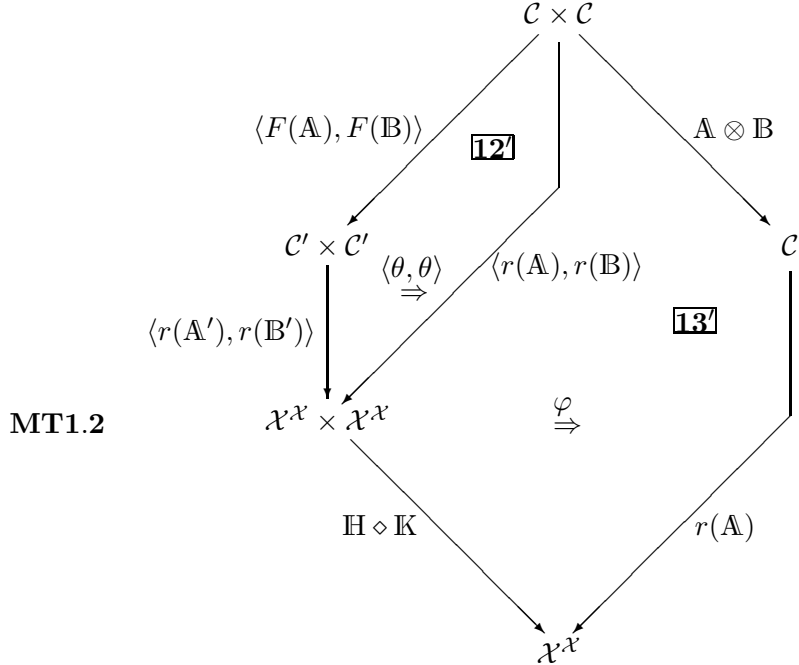
\begin{picture}(\xext,\yext)(\xoff,\yoff)
\put(-1100,800){$\bf MT1.2$}
\put(500,1750){\framebox(140,80){$\bf 12'$}}
\put(1200,1150){\framebox(140,80){$\bf 13'$}}

\settriparms[1`1`0;800]
\putAtriangle(0,1450)[\cC\times\cC`\cC'\times\cC'`\cC;\lk F(\mathds{A}),F(\mathds{B})\rk`\mathds{A}\otimes\mathds{B}`]

  \putmorphism(800,650)(0,-1)[\phantom{\cC'\times\cX}`\cX^\cX`]{600}{0}l

  \putmorphism(0,1450)(0,-1)[\phantom{\cC'\times\cC'}`\cX^\cX\times\cX^\cX`\lk r(\mathds{A}'),r(\mathds{B}')\rk]{600}{1}l

\put(180,1340){$\lk \theta,\theta\rk$}
\put(250,1260){$\Ra$}

\put(780,880){$\varphi$}
\put(780,800){$\Ra$}

  \put(80,780){\vector(1,-1){670}}
   \put(50,450){$\mathds{H}\diamond\mathds{K}$}

   \put(1600,860){\line(0,1){500}}
    \put(1600,860){\vector(-1,-1){750}}
    \put(1270,450){$r(\mathds{A})$}

    \put(800,1650){\line(0,1){500}}
    \put(800,1650){\vector(-1,-1){750}}
    \put(560,1350){$\lk r(\mathds{A}),r(\mathds{B})\rk$}
\end{picture}
\end{center}
compose to the same 2-cell. It is easy to see that the pasting diagrams $\bf MF1.1$ and $\bf MT1.1$ compose to the adjoint cells.
The cells $\bf 13$ and $\bf 13'$ in pasting diagrams $\bf MF1.2$ and  $\bf MT1.2$, respectively, are adjoint as well.
To see that compositions of the pasting diagrams are adjoint, it is enough to show that the composite 2-cell of $\bf 12a$ and $\bf 12b$ is adjoint to
the 2-cell corresponding to $\bf 12'$, i.e. $\theta\diamond\theta$. The exponential adjoint to $\theta\diamond\theta$ is the composition of the
following pasting diagram

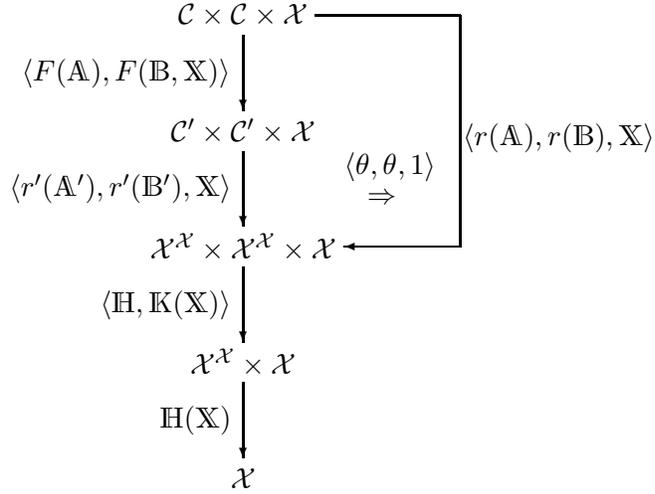
\begin{center} \xext=900 \yext=1700
\begin{picture}(\xext,\yext)(\xoff,\yoff)
\putmorphism(250,1650)(0,-1)[\cC\times\cC\times\cX`\cC'\times\cC'\times\cX`\lk F(\mathds{A}),F(\mathds{B},\mathds{X})\rk]{400}{1}l
\putmorphism(250,1250)(0,-1)[\phantom{\cC'\times\cC'\times\cX}`\cX^\cX\times\cX^\cX\times\cX`\lk r'(\mathds{A}'),r'(\mathds{B}'),\mathds{X}\rk]{400}{1}l
\putmorphism(250,850)(0,-1)[\phantom{\cX^\cX\times\cX^\cX\times\cX}`\cX^\cX\times\cX`\lk \mathds{H},\mathds{K}(\mathds{X})\rk]{400}{1}l
\putmorphism(250,450)(0,-1)[\phantom{\cX^\cX\times\cX}`\cX`\mathds{H}(\mathds{X})]{400}{1}l

\put(600,1100){$\lk \theta,\theta,1\rk$}
\put(680,1000){$\Ra$}

 \put(500,1650){\line(1,0){500}}
 \put(1000,1650){\line(0,-1){800}}
  \put(1000,850){\vector(-1,0){400}}
   \put(1010,1200){$\lk r(\mathds{A}),r(\mathds{B}),\mathds{X} \rk$}
\end{picture}
\end{center}
which is equal to the composition of
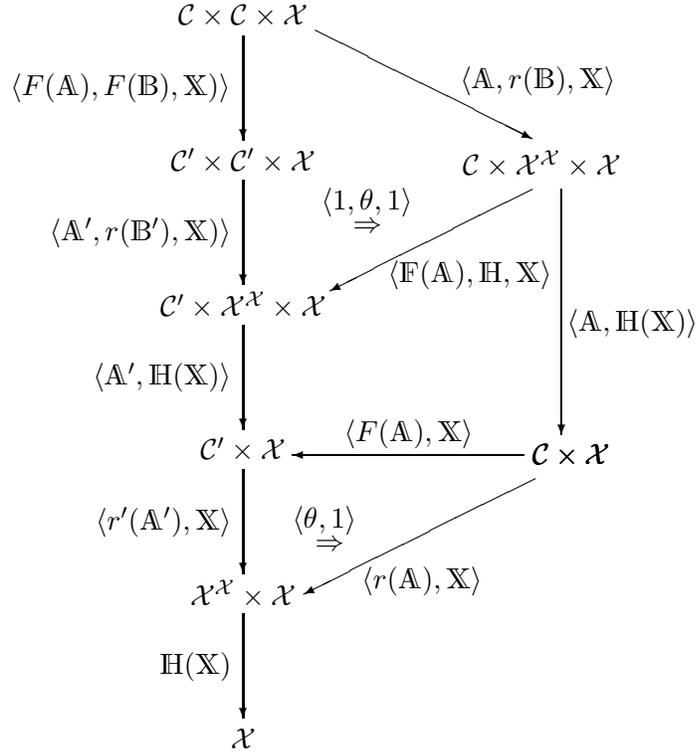
\begin{center} \xext=900 \yext=2600
\begin{picture}(\xext,\yext)(\xoff,\yoff)
\putmorphism(250,2550)(0,-1)[\cC\times\cC\times\cX`\phantom{\cC'\times\cC'\times\cX}`\lk F(\mathds{A}),F(\mathds{B}),\mathds{X})\rk]{500}{1}l
\putmorphism(250,2050)(0,-1)[\cC'\times\cC'\times\cX`\cC'\times\cX^\cX\times\cX`\lk \mathds{A}',r(\mathds{B}'),\mathds{X})\rk]{500}{1}l
\putmorphism(250,1550)(0,-1)[\phantom{\cC'\times\cX^\cX\times\cX}`\cC'\times\cX`\lk\mathds{A}',\mathds{H}(\mathds{X})\rk]{500}{1}l
\putmorphism(250,1050)(0,-1)[\phantom{\cX^\cX\times\cX^\cX\times\cX}`\cX^\cX\times\cX`\lk r'(\mathds{A}'),\mathds{X}\rk]{500}{1}l
\putmorphism(250,550)(0,-1)[\phantom{\cX^\cX\times\cX}`\cX`\mathds{H}(\mathds{X})]{500}{1}l

\put(520,1880){$\lk 1,\theta,1\rk$}
\put(630,1800){$\Ra$}

\put(420,780){$\lk \theta,1\rk$}
\put(500,700){$\Ra$}

  \put(500,2500){\vector(2,-1){750}}
     \put(1010,2000){$\cC\times\cX^\cX\times\cX$}
   \put(1000,2300){$\lk \mathds{A},r(\mathds{B}),\mathds{X} \rk$}
 \put(1250,1950){\vector(-2,-1){700}}
   \put(750,1630){$\lk \F(\mathds{A}),\mathds{H},\mathds{X} \rk$}
 \put(1350,1950){\vector(0,-1){850}}
     \put(1250,1000){$\cC\times\cX$}
   \put(1370,1460){$\lk \mathds{A},\mathds{H}(\mathds{X}) \rk$}
 \put(1220,1030){\vector(-1,0){800}}
     \put(1250,1000){$\cC\times\cX$}
   \put(600,1080){$\lk F(\mathds{A}),\mathds{X} \rk$}
 \put(1260,950){\vector(-2,-1){800}}
     \put(1250,1000){$\cC\times\cX$}
   \put(660,570){$\lk r(\mathds{A}),\mathds{X} \rk$}
\end{picture}
\end{center}
The last pasting diagram can be easily seen to have the same composition as $\bf 12a$ and $\bf 12b$. This ends the proof that the pasting diagram $\bf MF1$ commutes iff $\bf MT1$ does.

The diagram $\bf MF2$ says that the pasting diagram
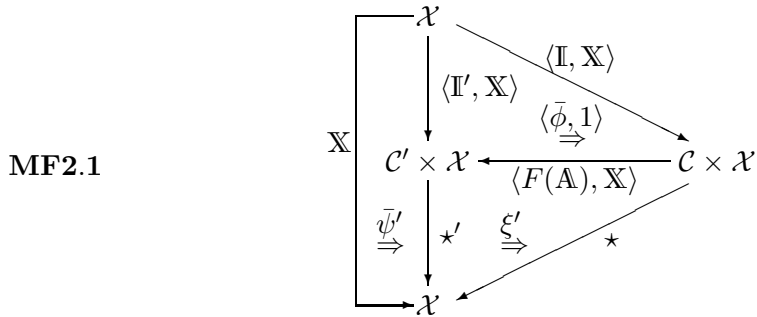
\begin{center} \xext=900 \yext=1100
\begin{picture}(\xext,\yext)(\xoff,\yoff)
\putmorphism(250,1050)(0,-1)[\cX`\cC'\times\cX`\lk\mathds{I}',\mathds{X}\rk]{500}{1}r
\putmorphism(250,550)(0,-1)[\phantom{\cC'\times\cX}`\cX`\star']{500}{1}r
\putmorphism(250,550)(1,0)[\phantom{\cC'\times\cX}`\cC\times\cX`\lk F(\mathds{A}),\mathds{X}\rk]{1000}{-1}l

\put(70,300){$\bar{\psi}'$}
\put(70,230){$\Ra$}

\put(500,300){$\xi'$}
\put(500,230){$\Ra$}

\put(630,670){$\lk \bar{\phi},1\rk $}
\put(700,600){$\Ra$}

\put(200,1050){\line(-1,0){200}}
\put(0,1050){\line(0,-1){1000}}
\put(0,50){\vector(1,0){200}}
\put(-100,580){$\mathds{X}$}

  \put(350,1020){\vector(2,-1){800}}
  \put(650,880){$\lk\mathds{I},\mathds{X}\rk$}

\put(1150,470){\vector(-2,-1){800}}
\put(860,250){$\star$}

    \put(-1200,500){$\bf MF2.1$}
\end{picture}
\end{center}
composes to $\bar{\psi}$ and $\bf MT2$ says that the pasting diagram
\begin{center} \xext=900 \yext=1100
\begin{picture}(\xext,\yext)(\xoff,\yoff)
\putmorphism(250,1050)(0,-1)[1`\cC'`\mathds{I}']{500}{1}r
\putmorphism(250,550)(0,-1)[\phantom{\cC'}`\cX^\cX`r']{500}{1}r
\putmorphism(250,550)(1,0)[\phantom{\cC'}`\cC`F(\mathds{A})]{1000}{-1}l

\put(70,300){$\bar{\varphi}'$}
\put(70,230){$\Ra$}

\put(500,300){$\theta'$}
\put(500,230){$\Ra$}

\put(610,700){$\bar{\phi}$}
\put(600,630){$\Ra$}

\put(200,1050){\line(-1,0){200}}
\put(0,1050){\line(0,-1){1000}}
\put(0,50){\vector(1,0){200}}
\put(-150,580){$\widehat{1_\cX}$}

  \put(350,1020){\vector(2,-1){800}}
  \put(650,880){$\lk\mathds{I},\mathds{X}\rk$}

  \put(1150,470){\vector(-2,-1){800}}
  \put(860,250){$r$}

    \put(-1200,500){$\bf MT2.1$}
\end{picture}
\end{center}
composes to $\bar{\varphi}$. It is easy to see that the pasting diagrams  $\bf MF2.1$ and $\bf MT2.1$ are adjoint and therefore that the pasting diagram $\bf MF2$ commutes iff $\bf MT2$ does. This shows the bijective correspondence of 1-cells of $\Mon_l(\cA)_{/_l\cX^\cX}$ and $\Act_l\Mon_l(\cA,\cX)$.

{\em Correspondence of 2-cells.} Finally, let $(F,\phi,\xi),(F',\phi',\xi'):(\star,\psi)\ra(\star',\psi')$ be two parallel 1-cells in $\Act_l\Mon_l(\cA,\cX)$ and $(F,\phi,\theta),(F',\phi',\theta'):(r,\varphi)\ra(r',\varphi')$ be the corresponding two parallel  1-cells in  $\Mon_l(\cA)_{/_l\cX^\cX}$.

Let $\tau : (F,\phi)\ra (F',\phi')$ be a monoidal transformation. $\tau$ is a transformation of monoidal actions iff $\bf MAT$ commutes, i.e. the pasting diagram
\begin{center} \xext=900 \yext=1100
\begin{picture}(\xext,\yext)(\xoff,\yoff)
\putmorphism(250,1050)(0,-1)[\cC\times\cX`\cC'\times\cX`]{600}{0}r
\putmorphism(140,1050)(0,-1)[\phantom{\cC\times\cX}`\phantom{\cC'\times\cX}`^{\lk F(\mathds{A}),\mathds{X}\rk}]{600}{1}l
\putmorphism(380,1050)(0,-1)[\phantom{\cC\times\cX}`\phantom{\cC'\times\cX}`^{\lk F'(\mathds{A}),\mathds{X}\rk}]{600}{1}r
\putmorphism(250,450)(0,-1)[\phantom{\cC'\times\cX}`\cX`\star']{400}{1}l

\put(140,800){$\lk \tau,1\rk$}
\put(220,700){$\Ra$}

\put(520,490){$\xi'$}
\put(500,400){$\Ra$}

 \put(400,1050){\line(1,0){400}}
 \put(800,1050){\line(0,-1){1000}}
  \put(800,50){\vector(-1,0){400}}
   \put(810,500){$\star$}

    \put(-1000,500){$\bf MAT'$}
\end{picture}
\end{center}
composes to $\xi$. On the other hand, $\tau$ is a 2-cell in the lax slice $\Mon_l(\cA)_{/_l\cX^\cX}$ iff $\bf LST$ for $\tau$
 \begin{center}\xext=800 \yext=500
\begin{picture}(\xext,\yext)(\xoff,\yoff)
 \settriparms[1`1`1;400]
  \putVtriangle(0,50)[r'F(\mathds{A})`r'F'(\mathds{A})`r(\mathds{A});r'(\tau)`\theta`\theta']
\end{picture}
\end{center}
holds, i.e. the pasting diagram
\begin{center} \xext=900 \yext=1100
\begin{picture}(\xext,\yext)(\xoff,\yoff)
\putmorphism(250,1050)(0,-1)[\cC`\cC'`]{600}{0}r
\putmorphism(160,1050)(0,-1)[\phantom{\cC}`\phantom{\cC'}`F(\mathds{A})]{600}{1}l
\putmorphism(340,1050)(0,-1)[\phantom{\cC}`\phantom{\cC'}`F'(\mathds{A})]{600}{1}r
\putmorphism(250,450)(0,-1)[\phantom{\cC'\times\cX}`\cX^\cX`r']{400}{1}l

\put(230,770){$\tau$}
\put(220,700){$\Ra$}

\put(520,460){$\theta'$}
\put(500,400){$\Ra$}

 \put(400,1050){\line(1,0){400}}
 \put(800,1050){\line(0,-1){1000}}
  \put(800,50){\vector(-1,0){400}}
   \put(810,500){$r$}

    \put(-1000,500){$\bf LST'$}
\end{picture}
\end{center}
composes to $\theta$. As $\xi$ is adjoint to $\theta$ and the above two pasting diagrams $\bf MAT'$ and $\bf LST'$ compose to the adjoint 2-cells, the 2-cells
in $\Mon_l(\cA)_{/_l\cX^\cX}$ and $\Act_l\Mon_l(\cA,\cX)$ are in bijective correspondence that respects domains and codomains.

The fact that the above bijective correspondences preserve identities and compositions is left for the reader.
$\Box$

From Propositions \ref{prop-Mon-slice-eq-to-Act} and \ref{prop-lax-slice-2-fib} we get immediately

\begin{corollary}
Let $\cX$ be an exponentiable 0-cell in a 2-category $\cA$ with finite products. Then the forgetful 2-functor $\Act_l\Mon_l(\cA,\cX)\lra \Mon_l(\cA)$ is a 2-fibration.
\end{corollary}

Let $\cC$ be a monoidal object in $\cA$. Clearly, $\cC$ acts on itself (on the left, say). An  {\em ideal} $\cX$ in a monoidal object $(\cC,\otimes)$ is a full and locally full sub-0-cell of $\cX$ of $\cC$ closed under the tensoring with $\cC$ from the left, i.e. we have a 1-cell $\star:\cC\times \cX\ra \cX$. Then $\star$ extends uniquely to an action of $(\cC,\otimes)$ on $\cX$ so that the embedding $\cX\ra \cC$ is a strict morphism of actions. Such an action of $(\cC,\otimes)$ on an ideal $\cX$ in $\cC$ is called a {\em tautologous action}\footnote{The name `tautologous action' is taken from \cite{BD}, but the meaning is slightly extended.}. By the above, if $\cX$ is exponentiable, such a tautologous action $(\star,\psi): (\cC,\otimes,\ldots)\times\cX\ra \cX$ corresponds to a morphism of monoidal categories
\[ (\br,\varphi): (\cC,\otimes)\lra\cX^\cX \]
$\cX$ is a {\em conservative ideal} in $(\cC,\otimes)$ iff  $\br$ is faithful and full on isos. If the strong monoidal 1-cell $(r,\varphi)$ has a right adjoint, then by Theorem \ref{thm-lift-kem}, the $\bkem$-diagram in $\Mon_l(\cA)$ for the induced monad, if it exists, lifts to a $\bkem$-diagram in $\Mon_l(\cA)_{/\cX^\cX}$.

\subsection{Monoids and their actions}\label{subsec-monoids-and-actions}

It was shown in \cite{SiZ} that the formation of monoids for a monoidal object and the formation of actions of monoids along a monoidal object action are weighted limits and therefore they can be considered in any 2-category $\cA$ with finite products. We recall briefly these constructions below.

Fix a monoidal object $(\cC,\otimes,I,\alpha,\lambda,\rho)$ in $\cA$. An {\em object of monoids} for $(\cC,\otimes)$ is
\begin{enumerate}
  \item a 0-cell $\mon(\cC,\otimes)$,
  \item a 1-cell $U^\otimes: \mon(\cC,\otimes)\ra \cC$
  \item two 2-cells $\bm^\otimes: U^\otimes\otimes U^\otimes\ra U^\otimes$, $\be^\otimes : \mathds{I}\ra U^\otimes$
\end{enumerate}
so that $\lk U^\otimes,\be^\otimes,\bm^\otimes\rk$ is a universal monoid in the monoidal category $\cA(\mon(\cC,\otimes),\cC)$, i.e. for any 0-cell $\cY$, it induces by composition
\[ \cA(\cY,\mon(\cC,\otimes))\ra  \mon(\cA(\cY,\cC)) \]
an isomorphism of categories natural in $\cY$, where $\mon(\cA(\cY,\cC))$ is the usual category of monoids in the monoidal category $\cA(\cY,\cC)$ whose monoidal structure is induced `pointwise' from $(\cC,\otimes)$. Thus, for a 0-cell $\cX$, one can identify a 1-cell $\mathds{M}:\cX\ra  \mon(\cC,\otimes)$ with triples $(\mathbf{M},\mathbf{m},\mathbf{e})$ such that $\mathbf{M}:\cX\ra \cC$ is a 1-cell and $\mathbf{m}:\mathbf{M}\otimes\mathbf{M}\ra \mathbf{M}$ and $\mathbf{e}:\mathbf{I}\ra \mathbf{M}$ are two 2-cells in $\cA(\cX,\cC)$ such that the usual equations are satisfied in this context.

\vskip 2mm
{\em Remark.} Note that the 1-cell
\[ \mathbf{M}\otimes\mathbf{M}:\cX\ra \cC\]
is the following composite
\begin{center} \xext=1500 \yext=150
\begin{picture}(\xext,\yext)(\xoff,\yoff)
 \putmorphism(0,50)(1,0)[\cX`\cC`\mathbf{M}]{500}{1}a
 \putmorphism(500,50)(1,0)[\phantom{\cC}`\cC\times\cC`\delta]{500}{1}a
 \putmorphism(1000,50)(1,0)[\phantom{\cC\times\cC}`\cC`\otimes]{500}{1}a
\end{picture}
\end{center}
where $\delta:\cC\ra \cC\times\cC$ is the diagonal morphism or rather in this case the (trivial, unique) comultiplication on the monoidal object $\cC$. Similarly, $\mathbf{I}=I\circ !\circ \mathbf{M}$, where $!:\cC\ra 1$ is the counit. In particular, one needs to have bimonoidal category/object to talk about monoids. Clearly, the comoidal structure always exists when we consider these structures with respect to the products in the ambient category $\cA$, but it is not so for a general monoidal structure on $\cA$.
\vskip 2mm

Let  $(\cC,\otimes,I,\alpha,\lambda,\rho,\cX,\star, \psi,\bar{\psi})$ be a monoidal action. An {\em object of actions} for $(\star,\psi)$ is
\begin{enumerate}
  \item a 0-cell $\act(\star)$,
  \item two 1-cells $V^\star: \act(\star)\ra \mon(\cC,\otimes)$, $\bar{V}^\star: \act(\star)\ra \cX$
  \item a 2-cell $\bd^\star: V^\star\star \bar{V}^\star\ra \bar{V}^\star$,
\end{enumerate}
so that $\lk V^\star,\bar{V}^\star,\bd^\star\rk$ is a universal action of the monoid $\cU^\otimes$ along the action  $\cA(\mon(\cC,\otimes),\cC)\times \cA(\mon(\cC,\otimes),\cX)\ra \cA(\mon(\cC,\otimes),\cX)$ induced by $(\star, \psi)$, i.e. for any 0-cell $\cY$, it induces by composition
\[ \cA(\cY,\mon(\cC,\otimes))\ra  \act\cA(\cY,(\star,\psi)) \]
an isomorphism of categories natural in $\cY$, where $\act\cA(\cY,(\star,\psi))$ is the usual category of actions of  monoids in the monoidal category $\cA(\cY,\cC)$ on objects of the category $\cA(\cY,\cX)$.

\subsection{Monoidal adjunctions}

It is well known that if a lax monoidal functor $(G,\gamma):(\cD,\oplus)\ra(\cC,\otimes)$ has a left adjoint $F\dashv U$ (as a functor), then $\gamma$ induces an oplax monoidal structure on $F$. The same holds true if we replace $\Cat$ by any other 2-category $\cA$ with finite products. It is stated less often that this correspondence extends to monoidal transformations giving rise to a pseudo-natural equivalence of pseudo-2-functors. As we are going to use it in Section \ref{subsec-monoidal-em-obj}, we elaborate on this below.

This shows that in the 2-category $\Mon_l(\cA)$ there is a problem with adjoint 1-cells, a left adjoint (as 1-cell) to a lax monoidal 1-cell is only oplax monoidal and hence it is not in $\Mon_l(\cA)$, in general. But as it will happen often in our applications, the left adjoints of interest will be in fact strong (oplax) monodal, and hence lax monoidal as  well.

Let $\cA$ be a 2-category with finite products. Let $\LAdj(\cA)$, $\RAdj(\cA)$ denote locally full sub-2-categories of $\cA$ with the same 0-cells as $\cA$ with 1-cells  having right, left adjoint 1-cells, respectively. The 2-category $\LAdj(\cA)$ is bi-equivalent to $\RAdj(\cA)^{co,op}$ i.e. with both 1- and 2-cells having exchanged domains and codomains.

Before we extend this correspondence to the monoidal case, we shall describe it in the `pure' case. Let $\cC$ and $\cD$ be two 0-cells in $\cA$. We shall define two functors
\begin{center} \xext=1500 \yext=250
\begin{picture}(\xext,\yext)(\xoff,\yoff)
 \putmorphism(0,100)(1,0)[\LAdj(\cA)(\cC,\cD)`\RAdj(\cA)(\cD,\cC)^{op}`]{1500}{0}a
 \putmorphism(0,150)(1,0)[\phantom{\LAdj(\cA)(\cC,\cD)}`\phantom{\RAdj(\cA)(\cD,\cC)^{op}}`\hat{(-)}]{1500}{1}a
 \putmorphism(0,50)(1,0)[\phantom{\LAdj(\cA)(\cC,\cD)}`\phantom{\RAdj(\cA)(\cD,\cC)^{op}}`\check{(-)}]{1500}{-1}b
\end{picture}
\end{center}
which will establish an adjoint equivalence that (pseudo) respects compositions.

For each left adjoint 1-cell $F:\cC\ra\cD$ we chose an adjunction $(F\dashv \hat{F},\eta,\varepsilon)$. Then for a morphism $\sigma: F\ra F'$ in $\LAdj(\cA)(\cC,\cD)$ we define $\hat{\sigma}: \hat{F'}\ra \hat{F}$ as an adjont 1-cell to
\begin{center} \xext=1200 \yext=150
\begin{picture}(\xext,\yext)(\xoff,\yoff)
 \putmorphism(0,50)(1,0)[F\,\hat{F}'`\phantom{F\,\hat{F}'}`\sigma_{\hat{F}'}]{600}{1}a
 \putmorphism(600,50)(1,0)[F\,\hat{F}'`1`\varepsilon']{600}{1}a
\end{picture}
\end{center}
For each right adjoint 1-cell $G:\cD\ra\cC$ we chose an adjunction $(\check{G}\dashv G,\eta,\varepsilon)$. Then for a morphism $\tau: G'\ra G$ in $\LAdj(\cA)(\cC,\cD)$ we define $\check{\tau}: \check{G}\ra \check{G'}$ as an adjont 1-cell to
\begin{center} \xext=1200 \yext=150
\begin{picture}(\xext,\yext)(\xoff,\yoff)
 \putmorphism(0,50)(1,0)[1` G'\,\check{G}'`\eta']{600}{1}a
 \putmorphism(600,50)(1,0)[\phantom{G'\, \check{G}'}`G \,\check{G}'`\tau_{\check{G}'}]{600}{1}a
\end{picture}
\end{center}
With the chosen adjunctions $(F\dashv \hat{F},\eta,\varepsilon)$ and $(\check{\hat{F}}\dashv \hat{F},\eta',\varepsilon')$, the unit of the adjunction $\hat{(-)}\dashv \check{(-)}$ is
\begin{center} \xext=1200 \yext=150
\begin{picture}(\xext,\yext)(\xoff,\yoff)
 \putmorphism(0,50)(1,0)[F`F\,\hat{F}\, \check{\hat{F}}`F(\eta')]{600}{1}a
 \putmorphism(600,50)(1,0)[\phantom{F\,\hat{F}\, \check{\hat{F}}}`\check{\hat{F}}`\varepsilon_{\check{\hat{F}}}]{600}{1}a
\end{picture}
\end{center}
Moreover, with the chosen adjunctions $(\check{G}\dashv G,\eta,\varepsilon)$ and $(\check{G}\dashv \hat{\check{G}},\eta',\varepsilon')$ the counit of the adjunction $\hat{(-)}\dashv \check{(-)}$ is
\begin{center} \xext=1200 \yext=150
\begin{picture}(\xext,\yext)(\xoff,\yoff)
 \putmorphism(0,50)(1,0)[G`\hat{\check{G}}\,\check{G}\, G`\eta'_G]{600}{1}a
 \putmorphism(600,50)(1,0)[\phantom{\hat{\check{G}}\,\check{G}\, G}`\hat{\check{G}}`\hat{\check{G}}(\varepsilon)]{600}{1}a
\end{picture}
\end{center}
The verification that this gives the adjoint equivalence is well known and it is left for the reader.

Now we shall describe how the above bi-equivalence extends to the monoidal setting.
Let $\LAdj\Mon_{o}(\cA)$ be the 2-category of monoidal objects as 0-cells, oplax monoidal 1-cells that have right adjoints (as 1-cells) and monoidal 2-cells. Similarly, let $\RAdj\Mon_{l}(\cA)$ be the 2-category of monoidal objects as 0-cells, lax monoidal 1-cells that have left adjoints (as 1-cells) and monoidal 2-cells.

Suppose that we have chosen adjoint pairs $(F\dashv G,\eta,\varepsilon)$ and $(F'\dashv G',\eta',\varepsilon')$ either for $F$'s or for $G$'s,
and $\sigma: F\ra F'$ corresponds to $\tau: G'\ra G$ either via $\hat{(-)}$ or via $\check{(-)}$. Moreover, $\sigma : (F,\varphi)\ra (F',\varphi'):(\cC,\otimes)\lra (\cD,\oplus)$ is a monoidal transformation of oplax monoidal 1-cells. Then we define $\bar{\gamma}:I\ra G(\stackrel{_+}{I}): 1\ra \cC$ as the adjoint morphism to $\bar{\varphi}:F(I)\ra \stackrel{_+}{I}: 1\ra \cD$. Moreover, we define $\gamma=\gamma_{\mathds{A},\mathds{B}}$  from $\varphi$ by the following diagram
 \begin{center} \xext=1800 \yext=650
\begin{picture}(\xext,\yext)(\xoff,\yoff)
 \setsqparms[1`-1`1`1;1800`500]
 \putsquare(0,50)[GF(G(\mathds{A})\otimes G(\mathds{B}))`G(FG(\mathds{A})\otimes FG(\mathds{B}))`G(\mathds{A})\otimes G(\mathds{B})`G(\mathds{A}\otimes  \mathds{B});G(\varphi_{\mathds{A},\mathds{B}})`\eta_{G(\mathds{A})\otimes G(\mathds{B})}`G(\varepsilon_\mathds{A}\otimes \varepsilon_\mathds{B})`\gamma_{\mathds{A},\mathds{B}}]
\end{picture}
\end{center}
The definition of $(\varphi,\bar{\varphi})$ from $(\gamma,\bar{\gamma})$ is similar. It is a routine verification that these correspondences are mutually inverse to one another and that  $(F,\varphi,\bar{\varphi})$ is an oplax 1-cell iff $(G,\gamma,\bar{\gamma})$  is a lax 1-cell.

Finally, we need to verify the correspondence at the level of 2-cells. Suppose that we have additionally two 2-cells $\sigma : F\ra F'$ and $\tau : G'\ra G$ so that $\hat{\sigma}=\tau$ (or $\hat{\tau}=\sigma$). As we have a correspondence in the `pure' case, it is enough to verify that if
$\sigma:(F,\varphi)\ra (F',\varphi')$ is a transformation of the oplax functors, then $\tau: (G',\gamma')\ra (G,\varphi)$ is a transformation of lax functors or vice versa. We shall show the first option. First note that the squares
\begin{center} \xext=2000 \yext=500
\begin{picture}(\xext,\yext)(\xoff,\yoff)
 \setsqparms[1`1`1`1;600`400]
 \putsquare(0,50)[1`G\,F`G'\, F'`G\,F';\eta`\eta'`G(\sigma)`\tau_{F'}]

  \setsqparms[1`1`1`1;600`400]
 \putsquare(1400,50)[F\,G'`F\, G`F'\,G'`1;F(\tau)`\sigma_{G'}`\varepsilon`\varepsilon']
\end{picture}
\end{center}
commute, as both composites are adjoint to $\sigma$ and $\tau$, respectively. The inner squares of the following diagram
\begin{center} \xext=3200 \yext=2300
\begin{picture}(\xext,\yext)(\xoff,\yoff)
 \setsqparms[1`1`1`1;1200`600]
 \putsquare(1000,250)[GF(G'(\mathds{A})\otimes G'(\mathds{A}))`G(FG'(\mathds{A})\otimes FG'(\mathds{A}))`
 GF(G(\mathds{A})\otimes G(\mathds{A}))`G(FG(\mathds{A})\otimes FG(\mathds{A}));
 G(\varphi)`GF(\tau\otimes\tau)`G(F(\tau)\otimes F(\tau))`G(\varphi)]

  \setsqparms[1`-1`-1`1;1200`600]
 \putsquare(1000,850)[GF'(G'(\mathds{A})\otimes G'(\mathds{A}))`G(F'G'(\mathds{A})\otimes F'G'(\mathds{A}))`
 \phantom{GF(G'(\mathds{A})\otimes G'(\mathds{A}))}`\phantom{G(FG'(\mathds{A})\otimes FG'(\mathds{A}))};
 G(\varphi')`G(\sigma)`G(\sigma\otimes \sigma)`]

   \setsqparms[1`1`1`1;1200`600]
 \putsquare(1000,1450)[G'F'(G'(\mathds{A})\otimes G'(\mathds{A})`G'(F'G'(\mathds{A})\otimes F'G'(\mathds{A})`
 \phantom{GF'(G'(\mathds{A})\otimes G'(\mathds{A}))}`\phantom{G(F'G'(\mathds{A})\otimes F'G'(\mathds{A})};
 G'(\varphi')`\tau`\tau`]

 \setsqparms[1`1`0`1;1000`1800]
 \putsquare(0,250)[G'(\mathds{A})\otimes G'(\mathds{B})`\phantom{G'(F'G'(\mathds{A})\otimes F'G'(\mathds{A})}`
 G(\mathds{A})\otimes G(\mathds{B})`\phantom{GF(G(\mathds{A})\otimes G(\mathds{A}))};
 \eta'`\tau\otimes\tau``\eta]

  \setsqparms[1`0`1`1;1200`1800]
 \putsquare(2200,250)[\phantom{G'(F'G'(\mathds{A})\otimes F'G'(\mathds{A}))}`G'(\mathds{A}\otimes \mathds{B})`
\phantom{G(FG(\mathds{A})\otimes FG(\mathds{A}))}`G(\mathds{A}\otimes \mathds{B});
 G(\varepsilon'\otimes\varepsilon')``\tau`G(\varepsilon\otimes\varepsilon)]

 \put(100,1950){\line(1,-3){200}}
 \put(300,1350){\vector(1,-1){420}}
 \put(300,1200){$\eta$}

 \put(2600,1400){\line(1,-1){400}}
 \put(3000,1000){\vector(1,-3){220}}
 \put(2800,1210){$G(\varepsilon'\otimes\varepsilon')$}

\put(3400,0){\vector(0,1){200}}
\put(0,0){\line(0,1){200}}
\put(0,0){\line(1,0){3400}}
\put(1900,40){$\gamma$}

\put(3400,2300){\vector(0,-1){200}}
\put(0,2100){\line(0,1){200}}
\put(0,2300){\line(1,0){3400}}
\put(1900,2210){$\gamma'$}

\end{picture}
\end{center}
in $\cA(\cD\times\cD,\cC)$ commute, by the above, the fact that $\eta$, $\tau$, $\varphi$ are natural, and $\sigma$ is a monoidal transformation. Thus $\gamma$'s is compatible with tensor.

The compatibility of $\tau$ with $\bar{\gamma}$'s is expressed by the commutation of the outer triangle in the following diagram
\begin{center} \xext=2000 \yext=1600
\begin{picture}(\xext,\yext)(\xoff,\yoff)
 \setsqparms[1`0`1`1;800`1000]
 \putsquare(850,250)[G'\, F'(I)`G'(\stackrel{_+}{I})`G\,F(I)`G(\stackrel{_+}{I});G'(\bar{\varphi}')``\tau_{\stackrel{_+}{I}}`G(\bar{\varphi})]

 \put(60,800){\vector(3,2){600}}
 \put(60,700){\vector(3,-2){620}}
 \put(-20,720){$I$}
 \put(300,1050){$\eta'_I$}
 \put(300,420){$\eta_I$}

 \putmorphism(850,1250)(0,-1)[\phantom{G'\, F'(I)}`G\, F'(I)`\tau_{F'(I)}]{500}{1}l
  \putmorphism(850,750)(0,1)[\phantom{G\, F(I)}`\phantom{G\, F'(I)}`G(\sigma_I)]{500}{-1}l

\put(1020,670){\vector(3,-2){520}}
\put(1200,580){$G(\bar{\varphi}')$}

\put(1650,1500){\vector(0,-1){150}}
\put(0,1500){\line(1,0){1650}}
\put(0,850){\line(0,1){650}}
\put(500,1410){$\bar{\gamma}$}

\put(1650,0){\vector(0,1){150}}
\put(0,0){\line(1,0){1650}}
\put(0,0){\line(0,1){670}}
\put(500,40){$\bar{\gamma}'$}

\end{picture}
\end{center}
in $\cA(1,\cD)$, where the upper and lower triangles commute by definition of $\bar{\gamma}$ and $\bar{\gamma}'$, respectively, the left square commutes by the above, the right square commutes by naturality of $\tau$ on $\bar{\varphi}'$, and the remaining triangle commutes as $F$ is oplax monoidal.

Thus we have

\begin{proposition}\label{prop-lax-oplax} The above construction describes
bi-equivalence of 2-categories $\LAdj\Mon_{o}(\cA)$ and  $\RAdj\Mon_{l}(\cA)$. $\Box$
 \end{proposition}

\section{Lax monoidal monads}\label{sec-Mon}

In this section we will be working in a 2-category $\cA$ with finite products admitting both Kleisli and Eilenberg-Moore objects such that Kleisli objects commute with finite products. Moreover,  $(\cR,\phi,\eta,\varepsilon)$ is a lax monoidal monad on a monoidal category object $(\cC,\otimes,I,\alpha,\lambda,\rho)$ in $\cA$ i.e. a monad in 2-category $\Mon_l(\cA)$. Clearly, not for everything said below all the assumptions are needed.

\subsection{k-objects}
The monad $(\cC,\otimes,I,\alpha,\lambda,\rho)$  as above admits the standard $\bk$-object $(\cC_\cR,\dot{\otimes},\dot{I},\dot{\alpha},\dot{\lambda},\dot{\rho})$ in $\Mon_l(\cA)$ with $F_\cR$ strict monoidal. All this follows from \cite{Z3}, see also \cite{Mo}, \cite{McC}. We will extract the construction from \cite{Z3} in elementary terms below.

The unit $\dot{I}$ is $F_\cR(I)$. As $\bk$-objects commute with finite products, the monad $\cR\times\cR$ on $\cC\times\cC$ admits the $\bk$-object $\cC_\cR\times\cC_\cR$.  Then one can verify that the morphism
 \begin{center} \xext=2600 \yext=200
\begin{picture}(\xext,\yext)(\xoff,\yoff)
 \putmorphism(0,50)(1,0)[F_\cR(\cR(\mathbf{A})\otimes \cR(\mathbf{B}))`\phantom{F_\cR\cR(\mathbf{A}\otimes \mathbf{B})}` F_\cR(\phi_{\mathbf{\mathbf{A}},\mathbf{\mathbf{B}}})]{1400}{1}a
  \putmorphism(1400,50)(1,0)[F_\cR\cR(\mathbf{A}\otimes \mathbf{B})`F_\cR(\mathbf{A}\otimes \mathbf{B})`\varepsilon_{F_\cR(\mathbf{A}\otimes \mathbf{B})}]{1200}{1}a
\end{picture}
\end{center}
in $\cA(\cC\times\cC,\cC_\cR)$ is a subcoequalizing of $\cR\times\cR$. Thus there is a 1-cell
\[ \dot{\otimes}: \cC_\cR\times\cC_\cR \lra \cC_\cR \]
such that
\[ F_\cR(\mathbf{A})\dot{\otimes} F_\cR(\mathbf{B})= F_\cR(\mathbf{A}\otimes\mathbf{B}), \hskip 5mm \dot{\otimes}(\kappa)=\varepsilon_{F_\cR(X\otimes Y)}\circ F_\cR(\phi). \]
This (will) exhibits $F_\cR$ as a strict monoidal functor, but we still need to define the coherence morphisms for $\dot{I}$ and $\dot{\otimes}$ in $\cC_\cR$. We shall describe $\dot{\alpha}$ leaving the other two for the reader. The 1-cell
\[ F_\cR(\mathbf{A}\otimes(\mathbf{B}\otimes\mathbf{C})) : \cC\times\cC\times\cC \ra \cC_\cR \] 
together with the 2-cell
\begin{center} \xext=2700 \yext=200
\begin{picture}(\xext,\yext)(\xoff,\yoff)
 \putmorphism(0,50)(1,0)[F_\cR(\cR(\mathbf{A})\otimes(\cR(\mathbf{B})\otimes\cR(\mathbf{C}))) `F_\cR\cR(\mathbf{A}\otimes(\mathbf{B}\otimes\mathbf{C}) ` F_\cR(\phi\circ (1\otimes \phi))]{1700}{1}a
  \putmorphism(1700,50)(1,0)[\phantom{F_\cR\cR(\mathbf{A}\otimes(\mathbf{B}\otimes\mathbf{C}))}`F_\cR(\mathbf{A}\otimes(\mathbf{B}\otimes\mathbf{C}))`
  \varepsilon_{F_\cR(\mathbf{A}\otimes(\mathbf{B}\otimes\mathbf{C}))}]{1400}{1}a
\end{picture}
\end{center}
is a subcoequalizing of the monad $\cR\times\cR\times\cR$ on $\cC\times\cC\times\cC$. Similarly
\[  (F_\cR((\mathbf{A}\otimes\mathbf{B})\otimes\mathbf{C}),\varepsilon_{F_\cR((\mathbf{A}\otimes\mathbf{B})\otimes\mathbf{C})}
\circ F_\cR(\phi\circ (\phi \otimes 1))  ) \]
is another subcoequalizing of the monad $\cR\times\cR\times\cR$ on $\cC\times\cC\times\cC$. Then the 2-cell
\[ F_\cR(\alpha):F_\cR(\cR(\mathbf{A})\otimes(\cR(\mathbf{B})\otimes\cR(\mathbf{C}))) \lra  (F_\cR((\mathbf{A}\otimes\mathbf{B})\otimes\mathbf{C})\]
is a morphism of these subcoequalizings and hence, by the couniversal properties of $\cC_\cR\times\cC_\cR\times\cC_\cR$, we get a 2-cell
\[ \dot{\alpha}: \mathds{A}\dot{\otimes}(\mathds{B}\dot{\otimes}\mathds{C})\lra (\mathds{A}\dot{\otimes}\mathds{B})\dot{\otimes}\mathds{C}\]
in $\cA(\cC_\cR\times\cC_\cR\times\cC_\cR,\cC_\cR)$ such that
\[ \dot{\alpha}_{\mathbf{A},\mathbf{B},\mathbf{C}}=F_\cR(\alpha_{\mathbf{A},\mathbf{B},\mathbf{C}}). \]

The functor $U_\cR:\cC_\cR\ra \cC$ is lax monoidal with coherence morphism $\dot{\bu}$ and monoidal adjunction $(F_\cR\dashv (U_\cR,\dot{\bu}),\eta,\varepsilon_\cR)$, with $F_\cR$ a strict monoidal functor. We have $\bar{\dot{\bu}}=\bar{\phi}=\eta_i$ and
$\dot{\bu}$ is defined as the following composition
\begin{center} \xext=1500 \yext=600
\begin{picture}(\xext,\yext)(\xoff,\yoff)
 \setsqparms[1`1`-1`1;1500`400]
  \putsquare(0,50)[U_\cR(\mathds{A})\otimes U_\cR(\mathds{B})`U_\cR(\mathds{A}\dot{\otimes}\mathds{B})`
  U_\cR F_\cR(U_\cR(\mathds{A})\otimes U_\cR(\mathds{B}))`U_\cR (F_\cR U_\cR(\mathds{A})\dot{\otimes}F_\cR U_\cR(\mathds{B}));
  \dot{\bu}_{\mathds{A},\mathds{B}}`\eta_{U_\cR(\mathds{A})\otimes U_\cR(\mathds{B})}`U_\cR(\varepsilon_\mathds{A}\dot{\otimes}\varepsilon_\mathds{B})`=]

\end{picture}
\end{center}


\subsection{em-objects}\label{subsec-monoidal-em-obj}
The existence of $\bem$-objects for monoidal monads is more subtle. It requires reflexive coequalizers. The theorem below describes the situation. The construction of the structure and the proof of this theorem will fill almost the whole subsection.
\begin{theorem}\label{thm-em-Linton}
  Let $\cA$ be a 2-category  with finite products. Let $(\cR,\phi,\eta,\varepsilon)$ be an rc-monad on rc-0-cell $(\cC,\otimes,I,\alpha,\lambda,\rho)$ in $\Mon_l(\cA)$, so that the monad  $(\cR,\eta,\varepsilon)$ in $\cA$ admits $\bem$-object $\cC^\cR$. Then $(\cR,\phi,\eta,\varepsilon)$ admits the standard $\bem$-object $(\cC^\cR,\ddot{\otimes},\ddot{I},\ddot{\alpha},\ddot{\lambda},\ddot{\rho})$ in $\Mon_l(\cA)$. The forgetful functor $U^\cR:\cC^\cR\ra \cC$ is lax monoidal with coherence morphism $\ddot{\bu}$, its left adjoint $(F^\cR,\ddot{\bv})$ is strong monoidal. They give rise to the monoidal adjunction $((F^\cR,\ddot{\bv})\dashv(U^\cR,\ddot{\bu}),\eta,\varepsilon^\cR)$.
\end{theorem}

We fix an rc lax monoidal monad $(\cR,\phi,\eta,\mu)$ on a rc monoidal category $(\cC,\otimes,I,\alpha,\lambda,\rho)$. Moreover, let $(\cC^\cR,U^\cR :\cC^\cR\ra \cC,\beta: \cR U^\cR \ra U^\cR)$ be an $\bem$-object for the monad $(\cR,\eta,\mu)$ in $\cA$.

By an $\cR$-algebra at $\cX$ we mean a 1-cell $\cX\ra \cC^\cR$. By definition of the $\bem$-object $\cC^\cR$ there is a bijective correspondence between subequalizings of $\cR$ at $\cX$, and $\cR$-algebras at $\cX$ that extends to morphisms. Because of this we will sometimes say that we have an $\cR$-algebra when we have in fact a corresponding subequalizing.

Note that $(\cR,\mu)$ is a subequalizing of $\cR$ and hence it gives rise to a 1-cell $F^\cR:\cC\ra \cC^\cR$ such that
\[ U^\cR F^\cR = \cR, \hskip 5mm \beta_{F^\cR}=\mu.\]
One can easily verify that $F^\cR$ is a left adjoint 1-cell to $U^\cR$ with the unit of adjunction $\eta$ and the counit $\varepsilon^\cR: F^\cR U^\cR\ra 1_{\cC^\cR}$ such that $U^\cR(\varepsilon^\cR)=\beta$. Thus we can say that either $F^\cR(\mathbf{A})$ or $\cR(\mathbf{A})$ or even $(\cR(\mathbf{A},\mu_\mathbf{A})$  is a free $\cR$-algebra on $\mathbf{A}=1_\cC:\cC\ra\cC$. All this is in accordance with the usual practice in $\Cat$.

{\em The construction.} First we shall describe the monoidal structure on $\cC^\cR$. The unit is $\ddot{I}=(\cR(I),\mu_I)$.

To define the the tensor $\ddot{\times}$ we consider the diagram in $\cA(\cC^\cR\times\cC^\cR,\cC)$
\begin{center} \xext=2700 \yext=900
\begin{picture}(\xext,\yext)(\xoff,\yoff)
 \setsqparms[0`1`1`0;1500`600]
 \putsquare(0,150)[\cR^2(\cR(\mathbf{A})\otimes\cR(\mathbf{B}))`\cR^2(\mathbf{A}\otimes\mathbf{B})`
 \cR(\cR(\mathbf{A})\otimes\cR(\mathbf{B}))`\cR(\mathbf{A}\otimes\mathbf{B});
 `\mu_{\cR(\mathbf{A})\otimes\cR(\mathbf{A})}`\mu_{\mathbf{A}\otimes\mathbf{B}}`]

 \putmorphism(0,800)(1,0)[\phantom{\cR^2(\cR(\mathbf{A})\otimes\cR(\mathbf{B}))}`\phantom{\cR^2(\mathbf{A}\otimes\mathbf{B})}`
 \cR(\mu\circ\cR(\phi_{\mathbf{A},\mathbf{B}}))]{1500}{1}a
  \putmorphism(0,700)(1,0)[\phantom{\cR^2(\cR(\mathbf{A})\otimes\cR(\mathbf{B}))}`\phantom{\cR^2(\mathbf{A}\otimes\mathbf{B})}`
 \cR^2(\mathbf{a}\otimes\mathbf{b})]{1500}{1}b

 \putmorphism(0,200)(1,0)[\phantom{\cR(\cR(\mathbf{A})\otimes\cR(\mathbf{B}))}`\phantom{\cR(\mathbf{A}\otimes\mathbf{B})}`
 \mu\circ\cR(\phi_{\mathbf{A},\mathbf{B}})]{1500}{1}a
  \putmorphism(0,100)(1,0)[\phantom{\cR(\cR(\mathbf{A})\otimes\cR(\mathbf{B}))}`\phantom{\cR(\mathbf{A}\otimes\mathbf{B})}`
 \cR(\mathbf{a}\otimes\mathbf{b})]{1500}{1}b

 \setsqparms[1`0`1`1;1200`600]
 \putsquare(1500,150)[\phantom{\cR^2(\mathbf{A}\otimes\mathbf{B})}`\cR(\mathbf{A}\ddot{\otimes}\mathbf{B})`
 \phantom{\cR(\mathbf{A}\otimes\mathbf{B})}`\mathbf{A}\ddot{\otimes}\mathbf{B};  \cR(\bq_{\mathbf{A},\mathbf{B}})``\mathbf{a}\bar{\otimes}\mathbf{b}`\bq_{\mathbf{A},\mathbf{B}}]
\end{picture}
\end{center}
In the above diagram $\mathbf{A}=U^\cR(\mathds{A})$, $\mathbf{B}=U^\cR(\mathds{B})$ are the first and the second projections, respectively, composed with the forgetful 1-cell $U^\cR$. $\mathbf{a}=\beta_\mathds{A}$, $\mathbf{b}=\beta_\mathds{B}$ are whiskerings of $\beta$ along projections. Thus $\mathds{A}$ can be thought of as an $\cR$-algebra, $\mathbf{A}$ its universe, and $\mathbf{a}$ its structural map. In the bottom row $\bq$ is defined as a (reflexive) coequalizer. The top row is a coequalizer as well, as $\cR$ preserves reflexive coequalizer. As the left square commutes serially, we have a 2-cell $\mathbf{a}\bar{\otimes}\mathbf{b}$ (a morphism in $\cA(\cC^\cR\times\cC^\cR,\cC)$) making the right square commute. One can easily verify that $(\cC^\cR\times\cC^\cR, \mathbf{A}\ddot{\otimes}\mathbf{B},\mathbf{a}\bar{\otimes}\mathbf{b})$ is a subequalizing of the monad $\cR$, and hence we have a 1-cell
\[ \ddot{\otimes}= \mathds{A}\ddot{\otimes}\mathds{B}:\cC^\cR\times\cC^\cR\lra \cC^\cR\]
such that
\[ U^\cR(\mathds{A}\ddot{\otimes}\mathds{B})=\mathbf{A}\ddot{\otimes}\mathbf{B}, \hskip 6mm  \beta_{\mathds{A}\ddot{\otimes}\mathds{B}}=\mathbf{a}\bar{\otimes}\mathbf{b}. \]

Note that $\bq_{\mathbf{A},\mathbf{B}}$ is a morphism of $\cR$-algebras, in the obvious sense.

In this way the definitions of the unit $\ddot{I}$ and the tensor $\ddot{\otimes}$ look as the usual Linton's definitions in the 2-category $\Cat$  but we interpret these diagrams in an arbitrary 2-category with finite products $\cA$. This allows us to import some statements and even the proofs from \cite{Se} as these proofs, taken as they are and suitably interpreted, in fact are the proofs of these statements in our context. For example, we have statements about the tensor of free $\cR$-algebras

\begin{proposition}\label{prop-coeq-free}
The following diagram  in $\cA(\cC\times\cC,\cC)$
\begin{center} \xext=2900 \yext=650
\begin{picture}(\xext,\yext)(\xoff,\yoff)
 \putmorphism(0,350)(1,0)[\cR(\cR^2(\mathbf{A})\otimes\cR^2(\mathbf{B}))`\cR(\cR(\mathbf{A})\otimes\cR(\mathbf{B}))`]{1400}{0}a
 \putmorphism(0,400)(1,0)[\phantom{\cR(\cR^2(\mathbf{A})\otimes\cR^2(\mathbf{B}))}`\phantom{\cR(\cR(\mathbf{A})\otimes\cR(\mathbf{B}))}`
\mu\circ \cR(\phi_{\cR(\mathbf{A}),\cR(\mathbf{B})})]{1400}{1}a
 \putmorphism(0,300)(1,0)[\phantom{\cR(\cR^1(\mathbf{A})\otimes\cR^2(\mathbf{B}))}`
 \phantom{\cR^2(\cR(\mathbf{A})\otimes\cR(\mathbf{B}))}`\cR(\mu_\mathbf{A}\otimes\mu_\mathbf{B})]{1400}{1}b

  \putmorphism(0,100)(1,0)[\phantom{\cR(\cR^1(\mathbf{A})\otimes\cR^2(\mathbf{B}))}`
 \phantom{\cR^2(\cR(\mathbf{A})\otimes\cR(\mathbf{B}))}`\cR(\eta_{\cR(\mathbf{A})}\otimes\eta_{\cR(\mathbf{B})})]{1400}{-1}b

  \put(1820,400){\vector(3,1){600}}
  \put(2000,560){$\bq_{\mathbf{A},\mathbf{B}}$}
   \put(1820,350){\vector(3,-1){640}}
   \put(2080,300){$\mu\circ\cR(\phi)$}

   \put(2460,40){\vector(-3,1){640}}
      \put(1880,20){$\cR(\eta\otimes\eta)$}
\putmorphism(2750,650)(0,-1)[\cR(\mathbf{A})\ddot{\otimes}\cR(\mathbf{B})`\cR(\mathbf{A}\otimes\mathbf{B})`\bq\circ \cR(\eta\otimes\eta)]{600}{-1}r
\end{picture}
\end{center}
is a split coequalizer with the comparison map $\bq_{\mathbf{A},\mathbf{B}}\circ \cR(\eta_\mathbf{A}\otimes\eta_\mathbf{B})$ an isomorphism. As $(\cC\times\cC,\cR(\mathbf{A}\otimes\mathbf{B}),\mu_{\mathbf{A}\otimes\mathbf{B}})$ is a subequalizing of $\cR$, by universal property of $\cC^\cR$ there is a lift of 1-cell $\cR(\mathbf{A})\ddot{\otimes}\cR(\mathbf{B})$ to $\cR$-algebra that is also denoted
\[ \cR(\mathbf{A})\ddot{\otimes}\cR(\mathbf{B}):\cC\times\cC\lra\cC^\cR  \]
so that
\[  U^\cR(\cR(\mathbf{A})\ddot{\otimes}\cR(\mathbf{B})) = \cR(\mathbf{A})\ddot{\otimes}\cR(\mathbf{B}).\]
\end{proposition}

{\it Proof.}~ To prove the first part of the statement, take the proof of Proposition 2.5.2  from \cite{Se} and interpret it in our context. The second part follows easily. $\Box$

The following proposition describes the presentation of iterated tensor $\ddot{\otimes}$ of $\cR$-algebras
\begin{proposition}\label{prop-coeq-3objects}
The following diagrams  in $\cA(\cC^\cR\times\cC^\cR\times\cC^\cR,\cC)$
\begin{center} \xext=2700 \yext=750
\begin{picture}(\xext,\yext)(\xoff,\yoff)
 \putmorphism(0,550)(1,0)[\cR(\cR(\mathbf{A})\otimes(\cR(\mathbf{B})\otimes\cR(\mathbf{C})))`\cR(\mathbf{A}\otimes(\mathbf{B}\otimes\mathbf{C}))`]{1700}{0}a
 \putmorphism(0,600)(1,0)[\phantom{\cR(\cR(\mathbf{A})\otimes(\cR(\mathbf{B})\otimes\cR(\mathbf{C})))}`
 \phantom{\cR(\mathbf{A}\otimes(\mathbf{B}\otimes\mathbf{C}))}` \mu\circ \cR(\phi\circ(1\otimes\phi))]{1700}{1}a
 \putmorphism(0,500)(1,0)[\phantom{\cR(\cR(\mathbf{A})\otimes(\cR(\mathbf{B})\otimes\cR(\mathbf{C})))}`
 \phantom{\cR(\mathbf{A}\otimes(\mathbf{B}\otimes\mathbf{C}))}`\cR(\mathbf{a}\otimes(\mathbf{b}\otimes\mathbf{c}))]{1700}{1}b

\setsqparms[1`1`-1`1;1200`500]
 \putsquare(1700,50)[\phantom{\cR(\mathbf{A}\otimes(\mathbf{B}\otimes\mathbf{C}))}`\mathbf{A}\ddot{\otimes}(\mathbf{B}\ddot{\otimes}\mathbf{C})`
 \cR(\mathbf{A}\otimes(\cR(\mathbf{B}\otimes\mathbf{C}))`\cR(\mathbf{A}\otimes(\mathbf{B}\ddot{\otimes}\mathbf{C}));  \bq_{\mathbf{A};\mathbf{B},\mathbf{C}}`\cR(1\otimes\eta)`\bq_{\mathbf{A},\mathbf{B}\ddot{\otimes}\mathbf{C}}`\cR(1\otimes\bq)]
\end{picture}
\end{center}
and
\begin{center} \xext=2700 \yext=750
\begin{picture}(\xext,\yext)(\xoff,\yoff)
 \putmorphism(0,550)(1,0)[\cR((\cR(\mathbf{A})\otimes\cR(\mathbf{B}))\otimes\cR(\mathbf{C}))`\cR((\mathbf{A}\otimes\mathbf{B})\otimes\mathbf{C}))`]{1700}{0}a
 \putmorphism(0,600)(1,0)[\phantom{\cR((\cR(\mathbf{A}\otimes(\cR(\mathbf{B}))\otimes\cR(\mathbf{C}))}`
 \phantom{\cR((\mathbf{A}\otimes\mathbf{B})\otimes\mathbf{C}))}` \mu\circ \cR(\phi\circ(\phi\otimes 1))]{1700}{1}a
 \putmorphism(0,500)(1,0)[\phantom{\cR((\cR(\mathbf{A})\otimes\cR(\mathbf{B}))\otimes\cR(\mathbf{C}))}`
 \phantom{\cR((\mathbf{A}\otimes\mathbf{B})\otimes\mathbf{C})}`\cR((\mathbf{a}\otimes\mathbf{b})\otimes\mathbf{c})]{1700}{1}b

\setsqparms[1`1`-1`1;1200`500]
 \putsquare(1700,50)[\phantom{\cR(\mathbf{A}\otimes(\mathbf{B}\otimes\mathbf{C}))}`(\mathbf{A}\ddot{\otimes}\mathbf{B})\ddot{\otimes}\mathbf{C}`
 \cR(\cR(\mathbf{A}\otimes\mathbf{B})\otimes\mathbf{C})`\cR((\mathbf{A}\ddot{\otimes}\mathbf{B})\otimes\mathbf{C});  \bq_{\mathbf{A},\mathbf{B};\mathbf{C}}`\cR(\eta\otimes 1)`\bq_{\mathbf{A}\ddot{\otimes}\mathbf{B},\mathbf{C}}`\cR(\bq\otimes 1)]
\end{picture}
\end{center}
are coequalizers. By universal properties of $\cC^cR$,  we have lifts of 1-cells $\mathbf{A}\ddot{\otimes}(\mathbf{B}\ddot{\otimes}\mathbf{C})$  and   $(\mathbf{A}\ddot{\otimes}\mathbf{B})\ddot{\otimes}\mathbf{C}$ to
\[ \mathds{A}\ddot{\otimes}(\mathds{B}\ddot{\otimes}\mathds{C}),\; (\mathds{A}\ddot{\otimes}\mathds{B})\ddot{\otimes}\mathds{C}\; :\; \cC^\cR\times\cC^\cR\times\cC^\cR\lra \cC^\cR \]
\end{proposition}
{\it Proof.}~ Take the proof of Corollary 2.6.2  from \cite{Se} and interpret it in our context. $\Box$

The coherence morphisms $\ddot{\lambda}=\ddot{\lambda}_\mathds{A}$,  $\ddot{\rho}=\ddot{\rho}_\mathds{A}$ for units are defined in a similar way. We present below the diagram in $\cA(\cC^\cR,\cC^\cR)$ defining $\ddot{\lambda}_\mathds{A}$.
\begin{center} \xext=2600 \yext=2000
\begin{picture}(\xext,\yext)(\xoff,\yoff)
 \setsqparms[1`0`1`1;1200`1600]
 \putsquare(1300,150)[\cR(\cR(\mathds{I})\otimes\mathds{A})`\cR(\mathds{I})\ddot{\otimes}\mathds{A}`
 \cR(\mathds{A})`\mathds{A};\bq_{\cR(\mathds{I}),\mathds{A}}``\ddot{\lambda}_{\mathds{A}}`\mathbf{a}]
  \putmorphism(1300,1750)(0,-1)[\phantom{\cR(\cR(\mathds{I})\otimes\mathds{A})}`
 \cR(\cR(\mathds{I})\otimes\cR(\mathds{A}))` \cR(1\otimes\eta)]{400}{1}r

 \putmorphism(1300,1350)(0,-1)[\phantom{\cR(\cR(\mathds{I})\otimes\cR(\mathds{A}))}`
 \cR^2(\mathds{I}\otimes\mathds{A})` \cR(\phi)]{400}{1}r

  \putmorphism(1300,950)(0,-1)[\phantom{\cR^2(\mathds{I}\otimes\mathds{A})}`
 \cR^2(\mathds{A})` \cR^2(\lambda_\mathds{A})]{400}{1}r

 \putmorphism(1300,550)(0,-1)[\phantom{\cR^2(\mathds{A})}`
 \phantom{\cR(\mathds{A})}` \mu_{\mathds{A}}]{400}{1}r

  \putmorphism(100,1750)(0,-1)[\cR(\cR^2(\mathds{I})\otimes\cR(\mathds{A}))`
 \cR(\cR(\mathds{I})\otimes\cR^2(\mathds{A}))` \cR(1\otimes\eta)]{400}{1}l

 \putmorphism(100,1350)(0,-1)[\phantom{\cR(\cR(\mathds{I})\otimes\cR^2(\mathds{A}))}`
 \cR^2(\mathds{I}\otimes\cR(\mathds{A}))` \cR(\phi)]{400}{1}l

  \putmorphism(100,950)(0,-1)[\phantom{\cR^2(\mathds{I}\otimes\cR(\mathds{A}))}`
 \cR^3(\mathds{A})` \cR^2(\lambda_\mathds{A})]{400}{1}l

 \putmorphism(100,550)(0,-1)[\phantom{\cR^2(\mathds{A})}`
 \cR(\mathds{A})` \mu_{\cR(\mathds{A})}]{400}{1}l

 \putmorphism(180,1800)(1,0)[\phantom{\cR^2(\mathds{I})\otimes\cR(\mathds{A}))}`
 \phantom{\cR(\cR(\mathds{I})\otimes\cR^2(\mathds{A}))}`
 \mu\circ \cR(\phi)]{1200}{1}a
  \putmorphism(180,1700)(1,0)[\phantom{\cR^2(\mathds{I})\otimes\cR(\mathds{A}))}`
  \phantom{\cR(\cR(\mathds{I})\otimes\cR^2(\mathds{A}))}`
 \cR(\mu_\mathds{I}\otimes\mathbf{a})]{1200}{1}b

 \putmorphism(100,200)(1,0)[\phantom{\cR^2(\mathds{A})}`
 \phantom{\cR(\mathds{A})}`
 \mu_{\mathds{A}}]{1200}{1}a
  \putmorphism(100,100)(1,0)[\phantom{\cR^2(\mathds{A})}`
  \phantom{\cR(\mathds{A})}`
 \cR(\mathbf{a})]{1200}{1}b
\end{picture}
\end{center}
Note that there is a small difference ($\mu_\mathds{A}$ replaces $\cR(\ba)$) with respect to the analogous definition in \cite{Se}, but as $\ba\circ \mu =\ba\circ \cR(\ba)$ both definitions come to the same. The one above does not use the algebra map $\ba$ in the (vertical part of) map $\mu\circ \cR^2(\lambda)\circ\cR(\phi)\circ \cR(1\otimes\eta)$ defining $\ddot{\lambda}$.
As  the columns are coequalizers and the left square in these diagrams commutes serially, we get $\ddot{\lambda}$ as the unique map making the right square commute.

The associativity coherence morphism $\ddot{\alpha}=\ddot{\alpha}_{\mathds{A},\mathds{B},\mathds{C}}$ is defined by the diagram in $\cA(\cC^\cR\times\cC^\cR\times\cC^\cR,\cC^\cR)$
\begin{center} \xext=2600 \yext=900
\begin{picture}(\xext,\yext)(\xoff,\yoff)
 \setsqparms[0`1`1`0;1800`600]
 \putsquare(0,150)[\cR(\cR(\mathds{A})\otimes(\cR(\mathds{B})\otimes\cR(\mathds{C})))`\cR(\mathds{A}\otimes(\mathds{B}\otimes\mathds{C}))`
 \cR((\cR(\mathds{A})\otimes\cR(\mathds{B}))\otimes\cR(\mathds{C}))`\cR((\mathds{A}\otimes\mathds{B})\otimes\mathds{C});
 `\cR(\alpha_{\cR(\mathds{A}),\cR(\mathds{B}),\cR(\mathds{C})})`\cR(\alpha_{\mathds{A},\mathds{B},\mathds{C}})`]

 \putmorphism(0,800)(1,0)[\phantom{\cR(\cR(\mathds{A})\otimes(\cR(\mathds{B})\otimes\cR(\mathds{C})))}`
 \phantom{\cR(\mathds{A}\otimes(\mathds{B}\otimes\mathds{C}))}`
 \mu\circ \cR(\phi \circ (1\otimes\phi))]{1800}{1}a
  \putmorphism(0,700)(1,0)[\phantom{\cR((\cR(\mathds{A})\otimes\cR(\mathds{B}))\otimes\cR(\mathds{C}))}`
  \phantom{\cR(\mathds{A}\otimes(\mathds{B}\otimes\mathds{C}))}`
 \cR(\mathbf{a}\otimes(\mathbf{b}\otimes \mathbf{c}))]{1800}{1}b

 \putmorphism(0,200)(1,0)[\phantom{\cR((\cR(\mathds{A})\otimes\cR(\mathds{B}))\otimes\cR(\mathds{C}))}`
 \phantom{\cR((\mathds{A}\otimes\mathds{B})\otimes\mathds{C})}`
 \mu\circ \cR(\phi \circ (\phi\otimes 1))]{1800}{1}a
  \putmorphism(0,100)(1,0)[\phantom{\cR((\cR(\mathds{A})\otimes\cR(\mathds{B}))\otimes\cR(\mathds{C}))}`
  \phantom{\cR((\mathds{A}\otimes\mathds{B})\otimes\mathds{C})}`
 \cR((\mathbf{a}\otimes\mathbf{b})\otimes \mathbf{c})]{1800}{1}b

 \setsqparms[1`0`1`1;1200`600]
 \putsquare(1800,150)[\phantom{\cR(\mathds{A}\otimes(\mathds{B}\otimes\mathds{C}))}`\mathds{A}\ddot{\otimes}(\mathds{B}\ddot{\otimes}\mathds{C})`
 \phantom{\cR((\mathds{A}\otimes\mathds{B})\otimes\mathds{C})}`(\mathds{A}\ddot{\otimes}\mathds{B})\ddot{\otimes}\mathds{C};  \bq_{\mathds{A};\mathds{B},\mathds{C}}``\ddot{\alpha}_{\mathds{A},\mathds{B},\mathds{C}}`\bq_{\mathds{A},\mathds{B};\mathds{C}}]
\end{picture}
\end{center}
as the rows are coequalizers by \ref{prop-coeq-3objects} and the left square commutes serially. Then again the calculations from \cite{Se} suitably interpreted give

\begin{proposition}
The $\bem$-object for the monad $(\cR,\eta,\mu)$ with the structure defined above $(\cC^\cR,\ddot{\otimes},\ddot{I},\ddot{\alpha},\ddot{\lambda},\ddot{\rho})$ is a monoidal object, i.e. a 0-cell in $\Mon_l(\cA)$.
\end{proposition}

We define below the coherence morphisms for 1-cell $U^\cR$. For the unit, we put
\[ \bar{\ddot{\bu}}=\eta_I:I\lra \cR(I)=U^\cR(\ddot{I}) \]
and for the tensor, we put
\begin{center} \xext=2400 \yext=550
\begin{picture}(\xext,\yext)(\xoff,\yoff)
 \putmorphism(0,400)(1,0)[U^\cR(\mathds{A})\otimes U^\cR(\mathds{B})=\mathbf{A}\otimes \mathbf{B}` \cR(\mathbf{A}\otimes\mathbf{B})`\eta_{\mathbf{A}\otimes\mathbf{B}}]{1200}{1}a

 \putmorphism(1200,400)(1,0)[\phantom{\cR(\mathbf{A}\otimes\mathbf{B})}`
 \mathbf{A}\ddot{\otimes}\mathbf{B} = U^\cR(\mathds{A}\ddot{\otimes}\mathds{B})`\bq_{\mathbf{A},\mathbf{B}} ]{1200}{1}a

  \put(450,320){\line(0,-1){200}}
  \put(450,120){\line(1,0){1650}}
 \put(2100,120){\vector(0,1){200}}
\put(1100,30){$\ddot{\bu}_{\mathds{A},\mathds{B}}$}
\end{picture}
\end{center}
in $\cA(\cC^\cR\times\cC^\cR,\cC^\cR)$.

We have

\begin{proposition}\label{prop-UR-monoidal}
The data
\[ (U^\cR,\ddot{\bu}): (\cC^\cR,\ddot{\otimes})\lra (\cC,\otimes)\]
is a monoidal 1-cell.
\end{proposition}

{\it Proof.}~ We need to show that $\ddot{\bu}$ is compatible with  the coherence morphisms $\alpha$, $\lambda$, $\rho$.

The compatibility with $\alpha$'s is the commutation of the outer hexagon in the diagram
 \begin{center} \xext=3000 \yext=2000
\begin{picture}(\xext,\yext)(\xoff,\yoff)

  \put(100,480){\vector(2,-1){700}}
  \put(330,250){$\ddot{\bu}$}

  \put(2900,480){\vector(-2,-1){700}}
  \put(2600,250){$\ddot{\bu}$}

  \put(-200,1550){\line(-1,0){150}}
  \put(-350,1550){\line(0,-1){1000}}
  \put(-350,550){\vector(1,0){150}}
  \put(-550,1000){$1\ddot{\otimes}\ddot{\bu}$}

  \put(3210,1550){\line(1,0){150}}
  \put(3360,1550){\line(0,-1){1000}}
  \put(3360,550){\vector(-1,0){150}}
  \put(3390,1000){$\ddot{\bu}\ddot{\otimes}1$}

  \setsqparms[-1`1`1`0;900`500]
 \putsquare(2000,1050)[\cR((\mathbf{A}\otimes\mathbf{B})\otimes\mathbf{C})`(\mathbf{A}\otimes\mathbf{B})\otimes\mathbf{C}`
 \phantom{\cR(\cR(\mathbf{A}\otimes\mathbf{B})\otimes\mathbf{C})}`\phantom{\cR(\mathbf{A}\otimes\mathbf{B})\otimes\mathbf{C}};
 \eta`\cR(\eta\otimes 1)`\eta\otimes 1`]

  \setsqparms[-1`1`1`0;900`500]
 \putsquare(2000,550)[\cR(\cR(\mathbf{A}\otimes\mathbf{B})\otimes\mathbf{C})`\cR(\mathbf{A}\otimes\mathbf{B})\otimes\mathbf{C}`
 \phantom{\mathbf{A}\otimes(\mathbf{B}\ddot{\otimes}\mathbf{C})}`\phantom{\cR(\mathbf{A}\otimes(\mathbf{B}\ddot{\otimes}\mathbf{C})};
 \eta`\cR(\bq \otimes 1)`\bq \otimes 1`]

 \setsqparms[-1`1`0`0;900`500]
 \putsquare(2000,50)[\cR((\mathbf{A}\ddot{\otimes}\mathbf{B})\otimes\mathbf{C}))`(\mathbf{A}\ddot{\otimes}\mathbf{B})\otimes\mathbf{C}`
 (\mathbf{A}\ddot{\otimes}\mathbf{B})\ddot{\otimes}\mathbf{C}`;\eta`\bq``]

  \setsqparms[1`1`1`0;900`500]
 \putsquare(100,1050)[\mathbf{A}\otimes(\mathbf{B}\otimes\mathbf{C})`\cR(\mathbf{A}\otimes(\mathbf{B}\otimes\mathbf{C}))`
 \phantom{\mathbf{A}\otimes\cR(\mathbf{B}\otimes\mathbf{C})}`\phantom{\cR(\mathbf{A}\otimes\cR(\mathbf{B}\otimes\mathbf{C}))};
 \eta`1\otimes\eta`\cR(1\otimes\eta)`]

  \setsqparms[1`1`1`0;900`500]
 \putsquare(100,550)[\mathbf{A}\otimes\cR(\mathbf{B}\otimes\mathbf{C})`\cR(\mathbf{A}\otimes\cR(\mathbf{B}\otimes\mathbf{C}))`
 \phantom{\mathbf{A}\otimes(\mathbf{B}\ddot{\otimes}\mathbf{C})}`\phantom{\cR(\mathbf{A}\otimes(\mathbf{B}\ddot{\otimes}\mathbf{C}))};
 \eta`1\otimes\bq`\cR(1\otimes\bq)`]

 \setsqparms[1`0`1`0;900`500]
 \putsquare(100,50)[\mathbf{A}\otimes(\mathbf{B}\ddot{\otimes}\mathbf{C})`\cR(\mathbf{A}\otimes(\mathbf{B}\ddot{\otimes}\mathbf{C}))`
 `\mathbf{A}\ddot{\otimes}(\mathbf{B}\ddot{\otimes}\mathbf{C});\eta``\bq`]

\putmorphism(1000,50)(1,0)[\phantom{\mathbf{A}\ddot{\otimes}(\mathbf{B}\ddot{\otimes}\mathbf{C})}`
\phantom{(\mathbf{A}\ddot{\otimes}\mathbf{B})\ddot{\otimes}\mathbf{C})}`
\ddot{\alpha}]{1000}{1}a
\putmorphism(1000,1550)(1,0)[\phantom{\cR(\mathbf{A}\otimes(\mathbf{B}\otimes\mathbf{C}))}`
\phantom{\cR((\mathbf{A}\otimes\mathbf{B})\otimes\mathbf{C})}`
\ddot{\alpha}]{1000}{1}a

  \put(100,1650){\line(0,1){150}}
  \put(100,1800){\line(1,0){2800}}
 \put(2900,1800){\vector(0,-1){150}}
 \put(1500,1840){$\alpha$}
\end{picture}
\end{center}
in $\cA(\cC^\cR\times\cC^\cR\times\cC^\cR,\cC)$, where we write $\mathbf{A}$ for $U^\cR(\mathds{A})$ ($\mathds{A}$ is the first projection), for short, and similar convention we adopt for $\mathbf{B}$ and $\mathbf{C}$.
The inner shapes commute either by naturality of $\eta$ or definitions of either $\ddot{\bu}$ or $\ddot{\alpha}$.

The arguments for the compatibility of $\ddot{\bu}$ with $\lambda$'s and with $\rho$'s are similar. We shall show the compatibility with $\rho$.
Thus with the convention as above we need to show that the outer square in the diagram
\begin{center} \xext=1400 \yext=800
\begin{picture}(\xext,\yext)(\xoff,\yoff)
 \setsqparms[1`1`-1`1;1500`700]
 \putsquare(0,50)[\mathbf{A}\otimes\mathbf{I}`\mathbf{A}`\mathbf{A}\otimes\cR(\mathbf{I})`\mathds{A}\ddot{\otimes}\mathbf{I};
 \rho_{\mathbf{A}}`1\ddot{\otimes}\bar{\ddot{\bu}}`\ddot{\rho}_{\mathbf{A}}`\ddot{\bu}]
 \put(100,150){\vector(2,1){400}}
  \put(900,350){\vector(2,-1){400}}

\put(500,400){$\cR(\mathbf{A}\otimes\mathbf{I})$}
\put(50,330){$\eta_{\mathbf{A}\otimes\cR(\mathbf{I})}$}
\put(1000,330){$\bq_{\mathbf{A},\ddot{\mathbf{I}}}$}
\put(670,600){\framebox(80,80){$\bf \rho$}}
\end{picture}
\end{center}
commutes. As the triangle commutes by definition of $\ddot{\bu}$, it remains to show that the pentagon \framebox(80,80){$\bf \rho$} commutes. To this end we consider the diagram
\begin{center} \xext=2000 \yext=2200
\begin{picture}(\xext,\yext)(\xoff,\yoff)
 \setsqparms[1`0`1`1;2000`2000]
 \putsquare(100,50)[\cR(\mathds{A}\otimes\cR(\mathds{I}))`\mathds{A}\ddot{\otimes}\cR(\mathds{I})`
 \cR(\mathds{A})`\mathds{A};\bq_{\mathds{A},\cR(\mathds{I})}``\ddot{\rho}_{\mathds{A}}`\mathbf{a}]
  \putmorphism(100,2050)(0,-1)[\phantom{\cR(\mathds{A}\otimes\cR(\mathds{I}))}`
 \cR(\cR(\mathds{A})\otimes\cR(\mathds{I}))` \cR(\eta \otimes 1)]{500}{1}l

 \putmorphism(100,1550)(0,-1)[\phantom{\cR(\cR(\mathds{A})\otimes\cR(\mathds{I}))}`
 \cR^2(\mathds{A}\otimes\mathds{I})` \cR(\phi)]{500}{1}l

  \putmorphism(100,1050)(0,-1)[\phantom{\cR^2(\mathds{A}\otimes\mathds{I})}`
 \cR^2(\mathds{A})` \cR^2(\lambda_\mathds{A})]{500}{1}l

 \putmorphism(100,550)(0,-1)[\phantom{\cR^2(\mathds{A})}`
 \phantom{\cR(\mathds{A})}` \mu_{\mathds{A}}]{500}{1}l

 \put(2030,90){\vector(-1,0){1770}}
 \put(1100,130){$\eta_\mathds{A}$}

 \put(1400,1500){$\mathds{A}\otimes\cR(\mathds{I})$}
 \putmorphism(1000,550)(1,0)[\cR(\mathds{A}\otimes\mathds{I})`\mathds{A}\otimes\mathds{I}`\eta_{\mathds{A}\otimes\mathds{I}}]{700}{-1}l

 \put(1000,650){\line(0,1){450}}
 \put(1000,1100){\vector(-2,1){760}}
 \put(550,1340){$\cR(\eta\otimes\eta)$}

  \put(160,960){\vector(2,-1){660}}
  \put(450,740){$\mu$}

 \put(890,650){\vector(-2,1){660}}
 \put(500,860){$\cR(\eta)$}

 \put(890,470){\vector(-2,-1){660}}
 \put(630,260){$\cR(\rho_\mathds{A})$}

  \put(1800,510){\vector(2,-3){260}}
 \put(1800,270){$\rho_\mathds{A}$}

  \put(1650,640){\vector(0,1){800}}
 \put(1680,970){$1\otimes\eta$}

  \put(1450,1600){\vector(-3,1){1100}}
 \put(980,1800){$\eta$}

 \put(1900,1750){\framebox(80,80){$\bf \rho$}}
\end{picture}
\end{center}
where the outer shape commutes by definition of $\ddot{\rho}$ and other shapes than \framebox(80,80){$\bf \rho$} commute by naturality of $\eta$ and $\mu$ and the fact that $\cR$ is a monoidal monad. Thus \framebox(80,80){$\bf \rho$} commutes as well.
$\Box$

We also have
\begin{proposition}\label{prop-beta-monoidal}
The 2-cell
\[ \beta: (\cR U^\cR,\cR(\ddot{\bu})\circ\phi_{U^\cR\times U^\cR},\cR(\bar{\ddot{\bu}})\circ \bar{\phi})\lra    (U^\cR,\ddot{\bu},\bar{\ddot{\bu}}): (\cC^\cR,\ddot{\otimes})\lra (\cC,\otimes)\]
is a monoidal 2-cell.
\end{proposition}

{\it Proof.}~Note that the diagram
\begin{center} \xext=2500 \yext=650
\begin{picture}(\xext,\yext)(\xoff,\yoff)
  \putmorphism(0,500)(1,0)[I`\cR(I)`\bar{\phi}=\eta_I]{1000}{1}a

 \putmorphism(1000,500)(1,0)[\phantom{\cR(I)}` \cR(I)=\cR U^\cR(\ddot{I})` \cR(\ddot{\bu})=\cR(\eta_I)]{1400}{1}a

  \putmorphism(2400,500)(0,-1)[\phantom{ \cR(I)=\cR U^\cR(\ddot{I})}`
 \cR(I)=U^\cR(\ddot{I})` \beta_{\ddot{I}}=\mu_I]{500}{1}r

  \put(50,450){\line(4,-1){1800}}
 \put(1850,00){\vector(1,0){200}}
 \put(400,200){$\bar{\ddot{\bu}}=\eta_I$}
\end{picture}
\end{center}
commutes. This shows compatibility of $\beta$ with units. To show the compatibility of $\beta$ with tensors, as $\ddot{\bu}=\bq\circ \eta$, it is enough to show that the outer heptagon in the diagram below
\begin{center} \xext=2500 \yext=1050
\begin{picture}(\xext,\yext)(\xoff,\yoff)
 \setsqparms[0`1`1`1;1600`800]
 \putsquare(100,50)[\cR(\mathds{A})\otimes\cR(\mathds{B})`\cR^2(\mathds{A}\otimes\mathds{B})`
 \mathds{A}\otimes\mathds{B}`\cR(\mathds{A}\otimes\mathds{B});`\ba\otimes \bb`\mu`\eta]

 \setsqparms[1`0`1`1;800`800]
 \putsquare(1700,50)[\phantom{\cR^2(\mathds{A}\otimes\mathds{B})}`\cR(\mathds{A}\ddot{\otimes}\mathds{B})`\phantom{\cR(\mathds{A}\otimes\mathds{B})}`
 \mathds{A}\ddot{\otimes}\mathds{B};\cR(\bq)``\ba\bar{\otimes} \bb`\bq]

 \putmorphism(100,850)(1,0)[\phantom{\cR(\mathds{A})\otimes\cR(\mathds{B})}`\cR(\mathds{A}\otimes\mathds{B})`\phi]{700}{1}a
 \putmorphism(800,890)(1,0)[\phantom{\cR(\mathds{A}\otimes\mathds{B})}`\phantom{\cR^2(\mathds{A}\otimes\mathds{B})}`\cR(\eta)]{900}{1}a
 \putmorphism(800,810)(1,0)[\phantom{\cR(\mathds{A}\otimes\mathds{B})}`\phantom{\cR^2(\mathds{A}\otimes\mathds{B})}`\eta_\cR]{900}{1}b

 \put(500,400){$\cR(\cR(\mathds{A})\otimes\cR(\mathds{B}))$}
 \put(1300,360){\vector(1,-1){200}}
 \put(990,200){$\cR(\ba\otimes \bb)$}

 \put(1300,520){\vector(1,1){200}}
 \put(1120,600){$\cR(\phi)$}

  \put(250,720){\vector(1,-1){200}}
 \put(400,600){$\eta$}
\end{picture}
\end{center}
commutes, where $\ba=\beta_\mathds{A}$, $\bb=\beta_\mathds{B}$. We have
\[ \bq_{\mathds{A},\mathds{B}}\circ\eta_{\mathds{A}\otimes\mathds{B}}\circ (\ba\otimes \bb) = \]
\[ = \bq_{\mathds{A},\mathds{B}} \circ \cR(\ba\otimes \bb) \circ \eta_{\cR(\mathds{A}\otimes\mathds{B})} = \]
\[ =  \bq_{\mathds{A},\mathds{B}} \circ \mu_{\mathds{A}\otimes\mathds{B}}\circ \cR(\phi_{\mathds{A},\mathds{B}}) \circ \eta_{\cR(\mathds{A})\otimes\cR(\mathds{B})} = \]
\[ =  \bq_{\mathds{A},\mathds{B}} \circ \mu_{\mathds{A}\otimes\mathds{B}}\circ \eta_{\cR(\mathds{A}\otimes\mathds{B})}\circ \phi_{\mathds{A},\mathds{B}} = \]
\[ =  \bq_{\mathds{A},\mathds{B}} \circ \mu_{\mathds{A}\otimes\mathds{B}}\circ \cR(\eta_{\mathds{A}\otimes\mathds{B}})\circ \phi_{\mathds{A},\mathds{B}} = \]
\[ =  (\ba\bar{\otimes}\bb)\circ \cR(\bq_{\mathds{A},\mathds{B}}) \circ \cR(\eta_{\mathds{A}\otimes\mathds{B}})\circ \phi_{\mathds{A},\mathds{B}}. \]
$\Box$

{\it Proof of Theorem \ref{thm-em-Linton}.}~ We need to verify only the universal property of $(\cC^\cR,\ddot{\otimes})$, $(U^\cR,\ddot{\bu})$, and $\beta$.

Assume that $(\cM,\check{\otimes},\check{I},\check{\alpha},\check{\lambda},\check{\rho})$ is a monoidal object, and
$(U,\psi) : (\cM,\check{\otimes})\ra (\cC\otimes)$ is a monoidal 1-cell that together with monoidal 2-cell $\xi: (\cR U,\cR(\psi)\circ\phi_{U\times U}) \ra (U,\psi)$ is a subequalizing of $(\cR,\phi,\eta,\mu)$. By the universal property of the $\bem$-object $\cC^\cR$ in $\cA$, there is a unique 1-cell $L:\cM\ra \cC^\cR$ in $\cA$ such that
\[ U=U^\cR L, \hskip 5mm \xi = \beta_L. \]
We shall define a (unique) coherence structure $\chi$ for 1-cell $L:\cM\ra \cC^\cR$ so that
\begin{equation}\label{eq-chi}
(U^\cR,\ddot{\bu})\circ (L,\chi)= (U,\psi)).
\end{equation}

The morphism $\bar{\chi}$ is defined as adjoint to $\bar{\psi}$, i.e.
\begin{center} \xext=1200 \yext=350
\begin{picture}(\xext,\yext)(\xoff,\yoff)
 \putmorphism(200,250)(1,0)[F^\cR(I)=\ddot{I}` L(\check{I})`\bar{\chi}]{800}{1}a
 \put(0,150){\line(1,0){1200}}
 \putmorphism(200,50)(1,0)[I`U(\check{I})= U^\cR L(\check{I})`\bar{\psi}]{800}{1}b
\end{picture}
\end{center}
Thus $\bar{\chi}$ is a unique such that $\bar{\psi}=U^\cR(\bar{\chi}) \circ \eta_I$.

The 2-cell $\chi$ is defined from the diagram in $\cA(\cM\times\cM,\cC^\cR)$
\begin{center} \xext=2700 \yext=900
\begin{picture}(\xext,\yext)(\xoff,\yoff)
 \setsqparms[0`1`1`0;1500`600]
 \putsquare(0,150)[\cR(\cR U(\mathds{M})\otimes\cR U(\mathds{N}))` \cR(U(\mathds{M})\otimes U(\mathds{N}))`
 \cR^2 U(\mathds{M}\check{\otimes}\mathds{N})`\cR U(\mathds{M}\check{\otimes}\mathds{N});
 `\cR^2(\psi)\circ \cR(\phi)`\cR(\psi)`]

 \putmorphism(0,800)(1,0)[\phantom{\cR(\cR U(\mathds{M})\otimes\cR U(\mathds{N}))}`\phantom{ \cR(U(\mathds{M})\otimes U(\mathds{N}))}`
 \cR(\xi_\mathds{M}\otimes\xi_\mathds{N})]{1500}{1}a
  \putmorphism(0,700)(1,0)[\phantom{\cR(\cR U(\mathds{M})\otimes\cR U(\mathds{N}))}`\phantom{ \cR(U(\mathds{M})\otimes U(\mathds{N}))}`
 \mu\circ \cR(\phi)]{1500}{1}b

 \putmorphism(0,200)(1,0)[\phantom{\cR^2 U(\mathds{M}\check{\otimes}\mathds{N})}`\phantom{\cR U(\mathds{M}\check{\otimes}\mathds{N})}`
 \cR(\xi_{\mathds{M}\check{\otimes}\mathds{N}})]{1500}{1}a
  \putmorphism(0,100)(1,0)[\phantom{\cR^2 U(\mathds{M}\check{\otimes}\mathds{N})}`\phantom{\cR U(\mathds{M}\check{\otimes}\mathds{N})}`
 \mu]{1500}{1}b

 \setsqparms[1`0`1`1;1200`600]
 \putsquare(1500,150)[\phantom{ \cR(U(\mathds{M})\otimes U(\mathds{N}))}`L(\mathds{M})\ddot{\otimes}L(\mathds{N})`
 \phantom{\cR U(\mathds{M}\check{\otimes}\mathds{N})}`L(\mathds{M}\check{\otimes}\mathds{N});  \bq_{L(\mathds{M}),L(\mathds{N})}``\chi_{\mathds{M},\mathds{N}}`\xi_{\mathds{M}\check{\otimes}\mathds{N}}]

\end{picture}
\end{center}
where $\mathds{M}$ and $\mathds{N}$ are the first and the second projections, respectively. The upper square commutes serially as $\xi$ is monoidal and $\mu$ is `natural' i.e. we have internal naturality of  $\mu$ on $\psi$. As the rows are coequalizers, we have a (unique) 2-cell $\chi$ such that $\chi\circ \bq = \xi\circ \cR(\psi)$.

Before we verify that $(L,\chi)$ is a monoidal 1-cell, we shall show that it is unique such that $(U^\cR,\ddot{\bu})\circ (L,\chi)= (U,\psi)$.

As $\eta_I=\bar{\ddot{\bu}}$ by the definition of $\bar{\chi}$, it is unique such that $\bar{\psi}=U^\cR(\bar{\chi}) \circ \bar{\ddot{\bu}}$.
Moreover, in the diagram
\begin{center} \xext=2700 \yext=1040
\begin{picture}(\xext,\yext)(\xoff,\yoff)
 \setsqparms[1`1`1`1;1200`500]
 \putsquare(0,200)[U(\mathds{M})\otimes U(\mathds{N})` \cR(U(\mathds{M})\otimes U(\mathds{N}))`
 U(\mathds{M}\check{\otimes}\mathds{N})`\cR U(\mathds{M}\check{\otimes}\mathds{N});\eta`\psi`\cR(\psi)`\eta]

 \setsqparms[1`0`1`1;1200`500]
 \putsquare(1200,200)[\phantom{\cR(U(\mathds{M})\otimes U(\mathds{N}))}`U^\cR(L(\mathds{M})\ddot{\otimes} L(\mathds{N})`
 \phantom{\cR U(\mathds{M}\check{\otimes}\mathds{N})}` U(\mathds{M}\check{\otimes}\mathds{N});U^\cR(\bq)``U^\cR(\chi_{\mathds{M},\mathds{N}})`\xi]

  \put(0,150){\line(0,-1){150}}
  \put(0,0){\line(1,0){2400}}
 \put(2400,0){\vector(0,1){150}}
\put(1300,20){$1$}

\put(0,750){\line(0,1){150}}
  \put(0,900){\line(1,0){2400}}
 \put(2400,900){\vector(0,-1){150}}
\put(1100,950){$\ddot{\bu}_{L(\mathds{M}),L(\mathds{N})}$}
\end{picture}
\end{center}
the left square commutes by naturality of $\eta$ on $\psi$,  the right square commutes by definition of $\chi$, the top triangle commutes by definition of $\ddot{\bu}$, and  the bottom triangle commutes as $\xi$ is a subequalizing. Hence the outer diagram commutes, i.e. $\psi_{\mathds{M},\mathds{N}}=U^\cR(\chi_{\mathds{M},\mathds{N}})\circ \ddot{\bu}_{L(\mathds{M}),L(\mathds{N})}$ and hence (\ref{eq-chi}) holds.

If $\widetilde{\chi}_{\mathds{M},\mathds{N}} : L(\mathds{M})\ddot{\otimes}L(\mathds{N})\lra L(\mathds{M}\check{\otimes}\mathds{N})$ is another such morphism satisfying (\ref{eq-chi}), then we would have
\[ \psi_{\mathds{M},\mathds{N}}= U^\cR(\widetilde{\chi}_{\mathds{M},\mathds{N}})\circ \ddot{\bu}_{}=\]
\[ = U^\cR(\widetilde{\chi}_{\mathds{M},\mathds{N}})\circ q_{L(\mathds{M}),L(\mathds{N})}\circ \eta_{U(\mathds{M}),U(\mathds{N})}=\]
\[ = U^\cR(\widetilde{\chi}_{\mathds{M},\mathds{N}}\circ q_{L(\mathds{M}),L(\mathds{N})})\circ \eta_{U(\mathds{M}),U(\mathds{N})}\]
Last equality holds as $q_{L(\mathds{M}),L(\mathds{N})}$ is a morphism of algebras.
Thus
\[  U^\cR(\widetilde{\chi}_{\mathds{M},\mathds{N}}\circ q_{L(\mathds{M}),L(\mathds{N})})\circ \eta_{U(\mathds{M}),U(\mathds{N})}= U^\cR(\chi_{\mathds{M},\mathds{N}}\circ q_{L(\mathds{M}),L(\mathds{N})})\circ \eta_{U(\mathds{M}),U(\mathds{N})}\]
and by adjunction $F^\cR\dashv U^\cR$ we have
\[  \widetilde{\chi}_{\mathds{M},\mathds{N}}\circ q_{L(\mathds{M}),L(\mathds{N})}= \chi_{\mathds{M},\mathds{N}}\circ q_{L(\mathds{M}),L(\mathds{N})} \]
As $q_{L(\mathds{M}),L(\mathds{N})}$ is a coequalizer, it is an epi and $\widetilde{\chi}_{\mathds{M},\mathds{N}}= \chi_{\mathds{M},\mathds{N}}$, i.e. $(L,\chi)$ is unique satisfying (\ref{eq-chi}).

Now we verify that $(L,\chi)$ is monoidal, i.e. it is compatible with $\alpha$'s, $\lambda$'s, and $\rho$'s. For compatibility of $\rho$'s we need to show that the square in $\cA(\cM,\cC^\cR)$
\begin{center} \xext=600 \yext=500
\begin{picture}(\xext,\yext)(\xoff,\yoff)
 \setsqparms[1`1`-1`1;800`400]
 \putsquare(0,50)[L(\mathds{M})\ddot{\otimes}\ddot{\mathds{I}}`L(\mathds{M})`L(\mathds{M})\ddot{\otimes}L(\check{\mathds{I}})`L(\mathds{M}\check{\otimes}\check{\mathds{I}});
 \ddot{\rho}_{L(\mathds{M})}`1\ddot{\otimes}\bar{\chi}`L(\check{\rho})`\chi_{\mathds{M},\check{\mathds{I}}}]
\end{picture}
\end{center}
commutes. In the following diagram
\begin{center} \xext=2700 \yext=2000
\begin{picture}(\xext,\yext)(\xoff,\yoff)
  \setsqparms[1`1`1`1;900`400]
 \putsquare(0,300)[\cR(U(\mathds{M})\otimes\cR(\mathds{I}))`\phantom{L(\mathds{M})\ddot{\otimes}\ddot{\mathds{I}}}`\cR(U(\mathds{M})\otimes U(\check{\mathds{I}}))`\phantom{L(\mathds{M}\ddot{\otimes}L(\check{\mathds{I}}))};
 \bq`\cR(1\otimes\bar{\chi})``\bq]
 \setsqparms[1`1`-1`1;900`400]
 \putsquare(900,300)[L(\mathds{M})\ddot{\otimes}\ddot{\mathds{I}}`L(\mathds{M})`L(\mathds{M})\ddot{\otimes}L(\check{\mathds{I}})`
 L(\mathds{M}\check{\otimes}\check{\mathds{I}});\ddot{\rho}_{L(\mathds{M})}`1\ddot{\otimes}\bar{\chi}`L(\check{\rho})`\chi_{\mathds{M},\check{\mathds{I}}}]
\setsqparms[-1`0`-1`-1;900`400]
 \putsquare(1800,300)[\phantom{L(\mathds{M})}`\cR U(\mathds{M})`\phantom{L(\mathds{M}\check{\otimes}\check{\mathds{I}})} `\cR  U(\mathds{M}\check{\otimes}\mathds{I}) ;
 \xi_{\mathds{M}}``\cR U(\check{\rho})`\xi_{\mathds{M}\check{\otimes}\check{\mathds{I}}}]

\putmorphism(500,1200)(1,0)[\cR(\cR U(\mathds{M})\otimes\cR(\mathds{I}))`\phantom{\cR^2( U(\mathds{M})\otimes\mathds{I})}`\cR(\phi)]{900}{1}a
\putmorphism(1400,1200)(1,0)[\cR^2( U(\mathds{M})\otimes\mathds{I})`\cR^2 U(\mathds{M})`\cR^2(\rho)]{900}{1}a

 \put(50,780){\vector(1,1){350}}
\put(310,950){$\cR(\eta\otimes 1)$}

\put(2300,1130){\vector(1,-1){350}}
\put(2180,950){$\mu_{U(\mathds{M})}$}

 \put(0,200){\line(0,-1){200}}
 \put(0,0){\line(1,0){2700}}
\put(2700,0){\vector(0,1){200}}
\put(1300,50){$\cR(\psi)$}

 \put(2700,1800){\vector(0,-1){1000}}
 \put(2500,1850){$\cR(U(\mathds{M})\otimes\mathds{I})$}
\put(2730,1200){$\cR(\rho_{U(\mathds{M})})$}

 \put(0,1800){\vector(0,-1){1000}}
 \put(-200,1850){$\cR(U(\mathds{M})\otimes\mathds{I})$}
\put(-450,1200){$\cR(1\otimes\eta_\mathds{I})$}

  \put(300,1800){\line(3,-1){600}}
 \put(900,1600){\vector(1,-1){300}}
 \put(660,1690){$\cR(\eta)$}

  \put(1600,1300){\line(1,1){300}}
 \put(1900,1600){\vector(3,1){600}}
 \put(2060,1690){$\mu$}

 \put(430,1850){\vector(1,0){2000}}
 \put(1350,1880){$1$}
\end{picture}
\end{center}
the outer shape commutes as $(U,\psi)$ is monoidal. The heptagon commutes by definition of $\ddot{\rho}$, the square on the bottom commutes by definition of $\chi$, the left square commutes by definition of $1\ddot{\otimes}\bar{\chi}$, the right square commutes by naturality of $\xi$ on $\check{\rho}$. Other commutations, except the mid square, are easy. To show that the the mid square commutes, it is enough to show that $\bq \circ \cR(1\otimes\eta_{\mathds{I}}): \cR(U(\mathds{M})\otimes\mathds{I})\ra L(\mathds{M})\ddot{\otimes}\ddot{I}$ is an epi.

In the diagram
\begin{center} \xext=2700 \yext=1650
\begin{picture}(\xext,\yext)(\xoff,\yoff)
\putmorphism(0,1550)(0,-1)[\cR( U(\mathds{M})\otimes\cR(\mathds{I}))`\phantom{\cR(\cR U(\mathds{M})\otimes\cR(\mathds{I}))}`\cR(\eta\otimes 1)]{500}{1}l

  \put(400,1550){\line(1,0){2500}}
 \put(2900,1550){\vector(0,-1){450}}
 \put(1300,1580){$1$}

  \setsqparms[1`0`1`1;1300`1000]
 \putsquare(1600,50)[\phantom{\cR(\cR U(\mathds{M})\otimes\cR^2(\mathds{I}))}` \cR(U(\mathds{M})\otimes \cR(\mathds{I}))`
 \phantom{\cR(U(\mathds{M})\otimes \cR(\mathds{I}))}`L(\mathds{M})\ddot{\otimes}\cR(\mathds{I}));
 \cR(\xi\otimes\mu_\mathds{I})``\bq`\bq]

  \setsqparms[1`1`1`1;1600`500]
 \putsquare(0,550)[\cR(\cR U(\mathds{M})\otimes\cR (\mathds{I}))`\cR(\cR U(\mathds{M})\otimes\cR^2(\mathds{I}))`
 \phantom{\cR^2(U(\mathds{M})\otimes\mathds{I})}`\phantom{\cR^2(U(\mathds{M})\otimes\cR(\mathds{I}))};
 \cR(\cR(1)\otimes\cR(\eta))`\cR(\phi)`\cR(\phi)`]

  \setsqparms[1`1`1`1;1600`500]
 \putsquare(0,50)[\cR^2(U(\mathds{M})\otimes\mathds{I})`\cR^2(U(\mathds{M})\otimes\cR(\mathds{I}))`\cR(U(\mathds{M})\otimes \mathds{I})`\cR(U(\mathds{M})\otimes \cR(\mathds{I}));
 \cR^2(1\otimes\eta)`\mu`\mu`\cR(1\otimes\eta)]
\end{picture}
\end{center}
the right square commutes as $\bq$ is a coequalizer, and other commutations are easy. Thus
\[ \bq = \bq \circ \cR(1\otimes\eta)\circ \mu \circ \cR(\phi)\circ \cR(\eta\otimes 1)\]
and $\bq \circ \cR(1\otimes\eta)$ is epi as $\bq$ is.

The compatibility of $L$ with $\lambda$'s can be shown in a similar way. We shall show that $L$ is compatible with $\alpha$'s. First note that the diagram in \mbox{$\cA(\cM\times\cM\times\cM,\cC)$}
\begin{center} \xext=3200 \yext=1700
\begin{picture}(\xext,\yext)(\xoff,\yoff)
 \setsqparms[1`1`1`0;1500`500]
 \putsquare(800,1050)[U(\mathds{M}_0)\otimes (U(\mathds{M}_1)\otimes U(\mathds{M}_2))`
 \cR(U(\mathds{M}_0)\otimes (U(\mathds{M}_1)\otimes U(\mathds{M}_2)))
 `\phantom{U(\mathds{M}_0)\otimes \cR(U(\mathds{M}_1)\otimes U(\mathds{M}_2))}`
 \phantom{\cR(U(\mathds{M}_0)\otimes \cR(U(\mathds{M}_1)\otimes U(\mathds{M}_2)))};
 \eta`1\otimes\eta`\cR(1\otimes\eta)`]

 \setsqparms[1`1`1`0;1500`500]
 \putsquare(800,550)[U(\mathds{M}_0)\otimes \cR(U(\mathds{M}_1)\otimes U(\mathds{M}_2))`
 \cR(U(\mathds{M}_0)\otimes \cR(U(\mathds{M}_1)\otimes U(\mathds{M}_2)))`
 U(\mathds{M}_0)\otimes (U(\mathds{M}_1)\ddot{\otimes} U(\mathds{M}_2))`
 \cR(U(\mathds{M}_0)\otimes (L(\mathds{M}_1)\ddot{\otimes} L(\mathds{M}_2)));
 \eta`1\otimes \bq`\cR(1\otimes \bq)`]

  \setsqparms[1`0`1`0;1500`500]
  \putsquare(800,50)[\phantom{U(\mathds{M}_0)\otimes (L(\mathds{M}_1)\ddot{\otimes} L(\mathds{M}_2))}`
  \phantom{\cR(U(\mathds{M}_0)\otimes (L(\mathds{M}_1)\ddot{\otimes} L(\mathds{M}_2)))}``
  L(\mathds{M}_0)\ddot{\otimes} (L(\mathds{M}_1)\ddot{\otimes} L(\mathds{M}_2));\eta``\bq`]

 \put(150,1550){\line(-1,0){300}}
 \put(-150,1550){\line(0,-1){1000}}
 \put(-150,550){\vector(1,0){300}}
 \put(-130,800){$1\otimes \ddot{\bu}$}

 \put(900,450){\vector(2,-1){800}}
 \put(1080,250){$\ddot{\bu}$}

  \put(3000,1550){\line(1,0){200}}
  \put(3200,1550){\line(0,-1){1500}}
  \put(3200,50){\vector(-1,0){320}}
  \put(3220,350){$\bq$}
\end{picture}
\end{center}
commutes, i.e. we have
\[ \bq_{U(\mathds{M}_0); U(\mathds{M}_1),U(\mathds{M}_2)}\circ \eta = \ddot{\bu}\circ (1\otimes \ddot{\bu}). \]
Similarly, we can show that
\[ \bq_{U(\mathds{M}_0), U(\mathds{M}_1),U(\mathds{M}_2)}\circ \eta = \ddot{\bu}\circ (\ddot{\bu}\otimes 1). \]
We need to show that in the following diagram in \mbox{$\cA(\cM\times\cM\times\cM,\cC)$}
\begin{center} \xext=3200 \yext=2200
\begin{picture}(\xext,\yext)(\xoff,\yoff)

  \setsqparms[1`1`1`0;1500`500]
 \putsquare(800,1550)[U(\mathds{M}_0)\ddot{\otimes} (U(\mathds{M}_1)\ddot{\otimes} U(\mathds{M}_2))`
 (U(\mathds{M}_0)\ddot{\otimes} U(\mathds{M}_1))\ddot{\otimes} U(\mathds{M}_2)
 `\phantom{\cR(U(\mathds{M}_0)\ddot{\otimes} (U(\mathds{M}_1)\ddot{\otimes} U(\mathds{M}_2)))}`
 \phantom{\cR(U(\mathds{M}_0)\ddot{\otimes} U(\mathds{M}_1))\ddot{\otimes} U(\mathds{M}_2))};\alpha
 `\eta`\eta`]

 \setsqparms[1`1`1`0;1500`500]
 \putsquare(800,1050)[\cR(U(\mathds{M}_0)\ddot{\otimes} (U(\mathds{M}_1)\ddot{\otimes} U(\mathds{M}_2)))`
 \cR(U(\mathds{M}_0)\ddot{\otimes} U(\mathds{M}_1))\ddot{\otimes} U(\mathds{M}_2))
 `\phantom{U^\cR(L(\mathds{M}_0)\ddot{\otimes} (L(\mathds{M}_1)\ddot{\otimes} L(\mathds{M}_2)))}`
 \phantom{(U^\cR(L(\mathds{M}_0)\ddot{\otimes} L(\mathds{M}_1))\ddot{\otimes} L(\mathds{M}_2))};\cR(\alpha)
 `\bq_{U(\mathds{M}_0); U(\mathds{M}_1),U(\mathds{M}_2)}`\bq_{U(\mathds{M}_0), U(\mathds{M}_1);U(\mathds{M}_2)}`]

 \setsqparms[1`1`1`0;1500`500]
 \putsquare(800,550)[U^\cR(L(\mathds{M}_0)\ddot{\otimes} (L(\mathds{M}_1)\ddot{\otimes} L(\mathds{M}_2)))`
 U^\cR((L(\mathds{M}_0)\ddot{\otimes} L(\mathds{M}_1))\ddot{\otimes} L(\mathds{M}_2))`
 U^\cR(L(\mathds{M}_0)\ddot{\otimes} L(\mathds{M}_1\check{\otimes} \mathds{M}_2))`U^\cR(L(\mathds{M}_0\check{\otimes} \mathds{M}_1)\ddot{\otimes} L(\mathds{M}_2));\ddot{\alpha} `1\ddot{\otimes}\chi`\chi\ddot{\otimes} 1`]

  \setsqparms[0`1`1`1;1500`500]
  \putsquare(800,50)[\phantom{U^\cR(L(\mathds{M}_0)\ddot{\otimes} L(\mathds{M}_1\check{\otimes} \mathds{M}_2))}`
  \phantom{U^\cR(L(\mathds{M}_0\check{\otimes} \mathds{M}_1)\ddot{\otimes} L(\mathds{M}_2))}`
  U^\cR(L(\mathds{M}_0\check{\otimes} (\mathds{M}_1\check{\otimes} \mathds{M}_2))`
  U^\cR(L((\mathds{M}_0\check{\otimes} \mathds{M}_1)\check{\otimes} \mathds{M}_2));`\chi`\chi`L(\check{\alpha})]

 \put(200,2050){\line(-1,0){300}}
 \put(-100,2050){\vector(0,-1){1200}}
 \put(-550,750){$U(\mathds{M}_0)\otimes U(\mathds{M}_1\check{\otimes} \mathds{M}_2)$}
 \put(-100,700){\line(0,-1){650}}
 \put(-100,50){\vector(1,0){350}}
 \put(-190,350){$\psi$}
  \put(-360,1350){$1\otimes \psi$}

  \put(2900,2050){\line(1,0){300}}
  \put(3200,2050){\vector(0,-1){1200}}
  \put(2750,750){$U(\mathds{M}_0\otimes \mathds{M}_1)\check{\otimes} U(\mathds{M}_2)$}
  \put(3200,700){\line(0,-1){650}}
  \put(3200,50){\vector(-1,0){350}}
  \put(3220,350){$\psi$}
  \put(3220,1350){$\psi\otimes 1$}
\end{picture}
\end{center}
the hexagon at the bottom commutes. The outer diagram commutes as $(U,\psi)$ is a monoidal 1-cell. The top two squares commute by naturality of $\eta$ on $\alpha$ and the definition of $\ddot{\alpha}$.  The left and right hexagons commute as the above and the fact that $(U,\psi)=(U^\cR,\ddot{\bu})\circ (L,\chi)$ which we have shown earlier. Thus the bottom hexagon commutes when precomposed with $\bq\circ \eta$. But
$\bq$ and all the 2-cells in the bottom hexagon are maps of algebras and $\eta$ is the unit of the adjunction. Thus the bottom hexagon commutes when precomposed with $\bq$ only. But $\bq$ is a coequalizer, so it is an epi and the bottom hexagon commutes as well. This ends the proof that $(L,\chi)$ is a monoidal 1-cell.

To finish the proof we need to verify that for a monoidal morphism of subequalizings
\[ \nu : (U,\psi,\xi)\ra (U',\psi',\xi'):(\cM,\check{\otimes})\ra (\cC,\otimes,\cR,\phi)\]
 the corresponding 2-cell
 \[ n : (L,\chi)\ra (L',\chi'):(\cM,\check{\otimes})\ra (\cC^\cR,\ddot{\otimes})\]
 is monoidal. To see that $n$ respects $\bar{\chi}$, we need to show that
  \begin{center} \xext=600 \yext=450
\begin{picture}(\xext,\yext)(\xoff,\yoff)
 \settriparms[1`1`1;350]
  \putAtriangle(0,50)[\ddot{I}`L(\check{I})`L'(\check{I});\bar{\chi}`\bar{\chi}'`n_{\check{I}}]
  \put(-1400,180){$(*)$}
\end{picture}
\end{center}
commutes. If we apply to the above diagram $U^\cR$ and precompose it with $\eta_I$, we will get
  \begin{center} \xext=600 \yext=1250
\begin{picture}(\xext,\yext)(\xoff,\yoff)
 \settriparms[1`1`1;350]
  \putAtriangle(0,250)[U^\cR(\ddot{I})`U^\cR L(\check{I})`U^\cR L'(\check{I});U^\cR(\bar{\chi})`U^\cR(\bar{\chi}')`U^\cR(n)]

  \putmorphism(350,1250)(0,-1)[I`\cR(I)`\eta_I]{400}{1}l
  \put(340,710){$\|$}

  \put(-560,220){$U(\check{I})=$}
  \put(930,220){$=U'(\check{I})$}
  \put(-460,170){\line(0,-1){170}}
  \put(-460,0){\line(1,0){1620}}
  \put(1160,0){\vector(0,1){170}}
  \put(-30,60){$\nu_{\check{I}}$}

  \put(250,1190){\vector(-2,-3){550}}
  \put(-120,800){$\bar{\psi}$}

  \put(450,1190){\vector(2,-3){550}}
  \put(750,800){$\bar{\psi}'$}
\end{picture}
\end{center}
As $(L,\chi)$, $(L',\chi')$, and $\nu:(U,\psi)\ra (U',\psi')$ are  monoidal, the two side triangles and the outer one commute.  Thus the diagram $(*)$ of $\cR$-algebras commutes when precomposed with $\eta_I$. Hence it also commutes as it is.

Finally, to see that $n$ respects $\chi$ we consider the following diagram
\begin{center} \xext=3500 \yext=1700
\begin{picture}(\xext,\yext)(\xoff,\yoff)
  \setsqparms[1`-1`1`0;3000`500]
 \putsquare(0,1100)[U(\mathds{M}\check{\otimes} \mathds{N})`
 \cR U(\mathds{M}\check{\otimes} \mathds{N})`
 `\phantom{\cR(U(\mathds{M})\check{\otimes} \mathds{N})}`
 \phantom{L(\mathds{M}\check{\otimes} \mathds{N})};\eta`\psi`\xi`]

 \setsqparms[0`1`-1`1;3000`500]
 \putsquare(0,100)[ \phantom{U(\mathds{M})\otimes U(\mathds{N})}`
 \phantom{L'(\mathds{M}\check{\otimes} \mathds{N})}`
 U'(\mathds{M}\check{\otimes} \mathds{N})` \cR U'(\mathds{M}\check{\otimes} \mathds{N});
 `\psi'`\xi'`\eta]

 \put(-300,1600){\line(-1,0){160}}
 \put(-460,1600){\line(0,-1){1500}}
 \put(-460,100){\vector(1,0){160}}
 \put(-750,700){$\nu_{\mathds{M}\check{\otimes}\mathds{N}}$}

  \put(3350,1600){\line(1,0){150}}
 \put(3500,1600){\line(0,-1){1500}}
 \put(3500,100){\vector(-1,0){150}}
 \put(3530,700){$\cR(n_{\mathds{M}\check{\otimes}\mathds{N}})$}

 \put(1200,1200){\vector(4,1){1400}}
 \put(1300,1330){$\cR(\psi)$}

  \put(1200,500){\vector(4,-1){1400}}
 \put(1300,320){$\cR(\psi')$}

  \setsqparms[1`1`1`1;1000`500]
 \putsquare(0,600)[U(\mathds{M})\otimes U(\mathds{N})`\cR(U(\mathds{M})\otimes U(\mathds{N}))`
 U'(\mathds{M})\otimes U'(\mathds{N})`\cR(U(\mathds{M})\otimes U(\mathds{N}));
 \eta`\nu_{\mathds{M}}\otimes\nu_{\mathds{N}}`\cR(\nu_\mathds{M}\otimes \nu_\mathds{N})`\eta]

  \setsqparms[1`0`1`1;1000`500]
  \putsquare(1000,600)[\phantom{\cR(U(\mathds{M})\otimes U(\mathds{N}))}`
  L(\mathds{M})\ddot{\otimes} L(\mathds{N})`\phantom{\cR(U'(\mathds{M})\otimes U'(\mathds{N}))}`
 L(\mathds{M})\ddot{\otimes} L(\mathds{N}); \bq``n_{\mathds{M}}\ddot{\otimes}n_{\mathds{N}}`\bq]

  \setsqparms[1`0`1`1;1000`500]
  \putsquare(2000,600)[\phantom{L(\mathds{M})\ddot{\otimes} L(\mathds{N})}`
  L(\mathds{M}\check{\otimes} \mathds{N})`\phantom{L(\mathds{M})\ddot{\otimes} L(\mathds{N})}`
 L(\mathds{M}\check{\otimes} \mathds{N}); \chi_{\mathds{M},\mathds{N}}``n_{\mathds{M}\check{\otimes}\mathds{N}}`\chi'_{\mathds{M},\mathds{N}}]
\end{picture}
\end{center}
which is partially in $\cA(\cM\times\cM,\cC^\cR)$ and partially in  $\cA(\cM\times\cM,\cC)$. There are four squares (including the outer one) that commute by naturality of $\eta$, the top and bottom squares on the right commute by definitions of $\chi$ and $\chi'$, respectively. The extreme  left square commutes as $\psi$ is a monoidal 2-cell. The extreme right square commutes by naturality of $\xi$ on $U^\cR(n)$. The central square commutes by definition of $n_{\mathds{M}}\ddot{\otimes}n_{\mathds{N}}$. This shows that the  square of $\cR$-algebra morphisms (interesting for us) commutes when composed with $\bq\circ \eta$. Since $\bq$ is an epi and $\eta$ is the unit of adjunction, the square commutes as is, i.e. $n$ is a monoidal 2-cell.
$\Box$

{\em Remark.} From Theorem \ref{thm-em-Linton} just proved, it follows that $F^\cR$ has lax monoidal structure $\ddot{\bv}$ as a left adjoint to $(U^\cR,\ddot{\bu})$. We shall identify here this structure which is in fact strong monoidal. The 2-cell $\bar{\ddot{\bv}}$ is just the identity on $F^\cR(I):1\ra \cC^\cR$. Let $\ddot{\bw}$ denote the oplax monoidal coherence morphism for $F^\cR$, i.e. in the diagram
\begin{center} \xext=3500 \yext=1700
\begin{picture}(\xext,\yext)(\xoff,\yoff)
  \putmorphism(300,550)(0,1)[\phantom{\cR(\cR^2(\mathbf{A})\otimes\cR(\mathbf{B}))}`\cR(\cR^2(\mathbf{A})\otimes\cR(\mathbf{B}))`]{500}{0}l
  \putmorphism(250,550)(0,1)[\phantom{\cR(\cR^2(\mathbf{A})\otimes\cR(\mathbf{B}))}`\phantom{\cR(\cR^2(\mathbf{A})\otimes\cR(\mathbf{B}))}`\mu\circ\phi]{500}{-1}l
  \putmorphism(350,550)(0,1)[\phantom{\cR(\cR^2(\mathbf{A})\otimes\cR(\mathbf{B}))}`\phantom{\cR(\cR^2(\mathbf{A})\otimes\cR(\mathbf{B}))}`\cR(\mu\otimes\mu)]{500}{-1}r

  \setsqparms[1`1`-1`1;2400`1000]
 \putsquare(300,550)[\cR(\mathbf{A}\otimes\mathbf{B})` F^\cR(\mathbf{A})\ddot{\otimes} F^\cR(\mathbf{B})
 `\cR(\cR^2(\mathbf{A})\otimes\cR(\mathbf{B}))`\cR U^\cR(F^\cR(\mathbf{A})\ddot{\otimes} F^\cR(\mathbf{B}));\ddot{\bw}`\cR(\eta\otimes\eta)`
 U^\cR(F^\cR(\mathbf{A})\otimes F^\cR(\mathbf{B})`\cR(\ddot{\bu})=\cR(\bq)\circ \cR(\eta)]

 \put(700,650){\vector(2,1){1650}}
 \put(1750,1100){$\cR(\bq)$}

  \put(400,650){\vector(2,3){220}}
  \put(520,750){$\cR(\phi)$}

   \put(650,1100){\vector(-2,3){220}}
    \put(610,1200){$\mu$}

 \put(500,1000){$\cR^2(\mathbf{A}\otimes\mathbf{B})$}

  \put(1850,700){\framebox(80,80){$\bf q$}}
\end{picture}
\end{center}
the outer shape (starting with $\cR(\cR(\mathbf{A})\otimes\cR(\mathbf{B}))$) commutes. Note that
\[ \mu_{\mathbf{A}\otimes\mathbf{B}}\circ\cR(\phi_{\mathbf{A},\mathbf{B}})\circ \cR(\eta_\mathbf{A}\otimes\eta_\mathbf{B})=1_{\cR(\mathbf{A}\otimes\mathbf{B})} \] and that the triangle \framebox(80,80){$\bf q$} commutes as well. As both $\mu_{\mathbf{A}\otimes\mathbf{B}}\circ\cR(\phi_{\mathbf{A}\otimes\mathbf{B}})$ and $\bq_{F^\cR(\mathbf{A}),F^\cR(\mathbf{B})}$ are coequalizers of the parallel pair $\mu_{\cR(\mathbf{A})\otimes\cR(\mathbf{B})}\circ\phi_{\cR(\mathbf{A})\otimes\cR(\mathbf{B})}$ and $\cR(\mu_{\mathbf{A}}\otimes\mu_{\mathbf{B}})$, it follows that $\bw_{\mathbf{A},\mathbf{B}}$ is an isomorphism, as are morphisms between coequalizers of the same pair of morphisms. Hence its inverse $\bv_{\mathbf{A},\mathbf{B}}$ is the coherence for the strong monoidal 1-cell $(F^\cR,\ddot{\bv})$, the left adjoint to $(U^\cR,\ddot{\bu})$.

\vskip 2mm

By $\cA_{rc}$ we denote the locally full sub-2-category of the 2-category $\cA$ with 0- and 1-cells rc, i.e. 0-cells having reflexive coequalizers and 1-cells preserving them. By an easy application of Johnstone lemma, we get the following.

\begin{lemma}
If a 2-category $\cA$ has finite products, so does $\cA_{rc}$. All 0- and 1-cells in $\cA_{rc}$ are rc.
\end{lemma}

In particular, the rc-monads in $\cA$ are precisely the monads in $\Mon_l(\cA_{rc})$. We also have

\begin{lemma}
 If $(\cC,\cR,\eta,\mu)$ is an rc-monad in $\cA$ and the monad $(\cC,\cR,\eta,\mu)$ admits an $\bem$-object in $\cA$, then
$(\cC,\cR,\eta,\mu)$ admits an $\bem$-object in $\cA_{rc}$ that is preserved by the `forgetful' 2-functor $\cA_{rc} \lra \cA$.
\end{lemma}

{\it Proof.}~ One has to check that if the subequalizing $(U,\psi,\xi)$ of the monad $\cR$, $U:\cM\ra \cC$ is an rc 0-cell, then so is its lift $L:\cM\ra \cC^\cR$. This is an easy adaptation of the proof from $\Cat$. The details are left for the reader.
$\Box$

From Theorem \ref{thm-em-Linton} and the above corollaries we obtain

\begin{corollary} If $(\cC,\otimes,\cR,\phi,\eta,\mu)$ is an rc-monad in $\cA$ and the monad $(\cC,\cR,\eta,\mu)$ admits an $\bem$-object in $\cA$, then
$(\cC,\otimes,\cR,\phi,\eta,\mu)$ admits an $\bem$-object in $\Mon_l(\cA_{rc})$ that is preserved by the `forgetful' 2-functor $\Mon_l(\cA_{rc}) \lra \Mon_l(\cA)$.
\end{corollary}

\subsection{Monoidal monad from a monoidal adjunction}\label{monoidal-mnd from-adj}

Theorem \ref{thm-em-Linton} permits us to refine the description of the process of extension of representation from Subsection \ref{subsec-extensions-of-rep}.
We start with a strong monoidal `representation' 1-cell in $\Mon_l(\cA)$
\[ (\br,\varphi) : (\cC,\otimes) \ra (\cM,\oplus)\]
with an rc-right adjoint $(U,\bu)$ which is necessarily lax monoidal. Hence $(\cR,\phi,\eta,\mu)$ is an rc lax monoidal monad and we get a $\bkem$-diagram in $\Mon_l(\cA)$
\begin{center} \xext=3200 \yext=1050
\begin{picture}(\xext,\yext)(\xoff,\yoff)
 \setsqparms[0`0`0`0;1600`1200]
 \putsquare(0,850)[\phantom{\mon(\cC_\cR,\dot{\otimes})}`\phantom{\mon(\cC,\otimes)}`\cC_\cR`\cC;``\phantom{\cU^{\otimes}}`]

 \putmorphism(0,880)(1,0)[\phantom{\cC_\cR}`\phantom{\cC}`F_\cR]{1600}{-1}a
 \putmorphism(0,820)(1,0)[\phantom{\cC_\cR}`\phantom{\cC}`U_\cR]{1600}{1}b

 \setsqparms[0`0`0`0;1600`1200]
 \putsquare(1600,850)[\phantom{\mon(\cC,{\otimes})}`\phantom{\mon(\cC^{\cR},\ddot{\otimes})}`\phantom{\cC}`\cC^\cR;``\phantom{\cU^{\ddot{\otimes}}}`]

 \putmorphism(1600,880)(1,0)[\phantom{\cC}`\phantom{\cC^\cR}`F^\cR]{1600}{1}a
 \putmorphism(1600,820)(1,0)[\phantom{\cC}`\phantom{\cC^\cR}`U^\cR]{1600}{-1}b

\put(1400,700){$_{(\otimes, I )}$}
\put(-150,700){$_{(\dot{\otimes}, \dot{I} )}$}
\put(3150,700){$_{(\ddot{\otimes}, \ddot{I} )}$}

\put(2000,0){$\cM$}
 \putmorphism(-50,800)(3,-1)[\phantom{\cC}`\phantom{\cC^\cR}`]{2200}{1}b
 \put(400,500){$\dot{\br}$}

 \putmorphism(1560,800)(2,-3)[\phantom{\cC}`\phantom{\cC^\cR}`]{470}{1}b
 \put(1660,500){$\br$}
  \putmorphism(1640,800)(2,-3)[\phantom{\cC}`\phantom{\cC^\cR}`]{470}{-1}b
 \put(1860,500){$U$}

  \putmorphism(3150,800)(-3,-2)[\phantom{\cC}`\phantom{\cC^\cR}`]{1100}{1}b
  \put(2600,500){$\ddot{\br}$}
   \putmorphism(3200,750)(-3,-2)[\phantom{\cC}`\phantom{\cC^\cR}`]{1100}{-1}b
   \put(2760,350){$K$}
 \end{picture}
\end{center}
with the comparison morphisms $(\dot{\br},\dot{\varphi})$ and $(K,\chi)$. $(K,\chi)$ comes from the monoidal subequalizing $(U,U(\varepsilon):U\cR\ra U)$.

The 1-cell $\dot{\br}$ corresponds to the subcoequalizing $(\br,\varepsilon_\br: \br \cR\ra \br)$ of $\cR$, the 2-cell $\bar{\dot{\varphi}}$ is $\bar{\varphi}: \stackrel{+}{I}\ra r(I)=\dot{\br} F_\cR(I)=\dot{\br}(\dot{I})$, and the cell $\dot{\varphi}$ corresponds to the morphism of subcoequalizings
\[ \varphi: (\br(\mathbf{A})\oplus\br(\mathbf{B}),\varepsilon_{\br(\mathbf{A})}\oplus\varepsilon_{\br(\mathbf{B})} ) \lra (\br(\mathbf{A}\otimes\mathbf{B}),\varepsilon_{\br(\mathbf{A}\otimes\mathbf{B})}\circ \br(\phi_{\mathbf{A},\mathbf{B}}) ) \]
of the monad $\cR\times\cR$ on $\cC\times\cC$.

We also have

\begin{proposition}
With the notation as above and $(\cM,\oplus)$ rc 0-cell, the oplax left adjoint $\ddot{\br}$ to the lax monoidal 1-cell $(K,\chi)$ is strong monoidal, i.e. we have a 1-cell $(\ddot{\br},\ddot{\varphi})$ in $\Mon_l(\cA)$ left adjoint to $(K,\chi)$, where $\ddot{\varphi}$ is the inverse of the oplax coherence morphism for $\ddot{\br}$ induced by $(K,\chi)$.
\end{proposition}

{\it Proof.}~ We shall show that the oplax monoidal functor $(\ddot{\br},\ddot{\bz})$, where $\bar{\ddot{\bz}}:\ddot{\br}(\ddot{I})\ra \stackrel{_+}{I}$, is the adjoint morphism to $\bar{\chi}:\ddot{I}\ra K(\stackrel{_+}{I})$ and
\[ \ddot{\bz}_{\mathds{A},\mathds{B}}=\varepsilon_{\ddot{\br}(\mathbb{A})\otimes\ddot{\br}(\mathbb{B})}\circ \ddot{\br}(\chi_{\ddot{\br}(\mathbb{A}),\ddot{\br}(\mathbb{B})}) \circ \ddot{\br}(\ddot{\eta}_\mathbb{A}\ddot{\otimes}\ddot{\eta}_\mathbb{B}):
\ddot{\br}(\mathbb{A}\ddot{\otimes}\mathbb{B}) \lra \ddot{\br}(\mathbb{A})\oplus\ddot{\br}(\mathbb{B}) : \cC^\cR\times\cC^\cR\ra \cM \]
is strong monoidal, i.e. both  $\bar{\ddot{\bz}}$ and $\ddot{\bz}$ are isomorphisms.

By definition of $(K,\chi)$ we have an equality of lax monoidal right adjoint
\[ (U^\cR,\ddot{\bu})\circ(K,\chi) = (U,\bu)\]
1-cells. Thus, by Proposition \ref{prop-lax-oplax}, the oplax left adjoint 1-cells are isomorphic, i.e. we have an oplax monoidal transformation,
\[  \sigma: (\ddot{\br},\bz)\circ (F^\cR,\ddot{\bv})\lra (\br,\varphi) \]
with both $(F^\cR,\ddot{\bv})$, $(r,\varphi)$ strong monoidal. In particular, in the commuting square
\begin{center} \xext=800 \yext=550
\begin{picture}(\xext,\yext)(\xoff,\yoff)
 \setsqparms[1`1`1`1;800`400]
 \putsquare(0,50)[\ddot{\br}F^\cR(I)`\ddot{\br}(\ddot{I})`\br(I)`\stackrel{_+}{I};\bar{\ddot{\bv}}`\sigma_I`\bar{\ddot{\bz}}`\bar{\varphi}]
\end{picture}
\end{center}
$\bar{\ddot{\bz}}$ is an isomorphism as the three remaining morphisms are. As
\begin{center} \xext=1800 \yext=250
\begin{picture}(\xext,\yext)(\xoff,\yoff)
 \putmorphism(0,100)(1,0)[\cR^2(\mathds{A})`\cR(\mathds{A})`]{1000}{0}a
 \putmorphism(0,150)(1,0)[\phantom{\cR^2(\mathds{A})}`\phantom{\cR(\mathds{A})}`\beta_{\cR(\mathds{A})}]{1000}{1}a
 \putmorphism(0,50)(1,0)[\phantom{\cR^2(\mathds{A})}`\phantom{\cR(\mathds{A})}`\mu_\mathds{A}]{1000}{1}b
 \putmorphism(1000,100)(1,0)[\phantom{\cR(\mathds{A})}`\mathds{A}`\beta_\mathds{A}]{800}{1}a
\end{picture}
\end{center}
is a coequalizer in $\cA(\cC^\cR,\cC^\cR)$, and both $\ddot{\otimes}$ and $\ddot{\br}$ preserve this coequalizer, we have a diagram $\cA(\cC^\cR\times\cC^\cR,\cM)$

\begin{center} \xext=2900 \yext=900
\begin{picture}(\xext,\yext)(\xoff,\yoff)

 \setsqparms[0`1`1`0;1500`600]
 \putsquare(0,150)[\ddot{\br}(\cR^2(\mathds{A})\ddot{\otimes}\cR^2(\mathds{B}))`\ddot{\br}(\cR(\mathds{A})\ddot{\otimes}\cR(\mathds{B}))`
 \ddot{\br}\cR^2(\mathds{A})\oplus\ddot{\br}\cR^2(\mathds{B})`\ddot{\br}\cR(\mathds{A})\oplus\ddot{\br}\cR(\mathds{B});
 `\ddot{\bz}_{\cR^2(\mathds{A}),\cR^2(\mathds{B})}`\ddot{\bz}_{\cR(\mathds{A}),\cR(\mathds{B})}`]

 \putmorphism(0,800)(1,0)[\phantom{\ddot{\br}(\cR^2(\mathds{A})\ddot{\otimes}\cR^2(\mathds{B}))}`\phantom{\ddot{\br}(\cR(\mathds{A})\otimes\cR(\mathds{B}))}`
 \ddot{\br}(\beta_{\cR(\mathds{A})}\ddot{\otimes}\beta_{\cR(\mathds{B})})]{1500}{1}a
  \putmorphism(0,700)(1,0)[\phantom{\ddot{\br}(\cR^2(\mathds{A})\ddot{\otimes}\cR^2(\mathds{B}))}`\phantom{\ddot{\br}(\cR(\mathds{A})\otimes\cR(\mathds{B}))}`
 \ddot{\br}(\mu_\mathds{A}\ddot{\otimes}\mu_\mathds{B})]{1500}{1}b

 \putmorphism(0,200)(1,0)[\phantom{\ddot{\br}\cR^2(\mathds{A})\oplus\ddot{\br}\cR^2(\mathds{B})}`\phantom{\ddot{\br}(\cR(\mathds{A})\otimes\cR(\mathds{B}))}`
 \ddot{\br}(\beta_{\cR(\mathds{A})})\oplus\ddot{\br}(\beta_{\cR(\mathds{B})})]{1500}{1}a
  \putmorphism(0,100)(1,0)[\phantom{\ddot{\br}\cR^2(\mathds{A})\oplus\ddot{\br}\cR^2(\mathds{B})}`\phantom{\ddot{\br}(\cR(\mathds{A})\otimes\cR(\mathds{B}))}`
 \ddot{\br}(\mu_\mathds{A})\oplus\ddot{\br}(\mu_\mathds{B})]{1500}{1}b

 \setsqparms[1`0`1`1;1400`600]
 \putsquare(1500,150)[\phantom{\ddot{\br}(\cR(\mathds{A})\otimes\cR(\mathds{B}))}`\ddot{\br}(\mathds{A}\ddot{\otimes}\mathds{B})`
 \phantom{\ddot{\br}(\cR(\mathds{A})\otimes\cR(\mathds{B}))}`\ddot{\br}(\mathds{A})\oplus\ddot{\br}(\mathds{B});  \ddot{\br}(\beta_{\mathds{A}}\ddot{\otimes}\beta_{\mathds{B}})``\ddot{\bz}_{\mathds{A},\mathds{B}}`\ddot{\br}(\beta_{\mathds{A}})\oplus\ddot{\br}(\beta_{\mathds{B}})]
\end{picture}
\end{center}
in which both rows are coequalizers. Thus $\ddot{\bz}=\ddot{\bz}_{\mathds{A},\mathds{B}}$ is an isomorphism if both $\ddot{\bz}_{\cR(\mathds{A}),\cR(\mathds{B})}$ and $\ddot{\bz}_{\cR^2(\mathds{A}),\cR^2(\mathds{B})}$ are, i.e. if $\ddot{\bz}$ is an isomorphism of free algebras. This is true as we have an isomorphism $\sigma$.

The inverse $\ddot{\varphi}$ of $\ddot{\bz}$ is the required lax structure making $(\ddot{\br},\ddot{\varphi})$ a left adjoint $(K,\chi)$ in $\Mon_l(\cA)$. $\Box$
\subsection{Monoids}\label{subsec-monoids}
 Applying the 2-functor $\mon$, cf. \cite{SiZ} and section \ref{subsec-monoids-and-actions}, to the $\bkem$-diagram for $(\cR,\phi,\eta,\varepsilon)$ in $\Mon_l(\cA)$, we get a diagram of six 0-cells and some 1-cells in $\cA$ so that we can form a diagram in $\cA$:
\begin{center} \xext=2000 \yext=1650
\begin{picture}(\xext,\yext)(\xoff,\yoff)
 \setsqparms[0`1`1`0;1000`800]
 \putsquare(0,450)[\mon(\cC_\cR,\dot{\otimes})`\mon(\cC,\otimes)`\cC_\cR`\cC;`\cU^{\dot{\otimes}}`\cU^{\otimes}`]
 \putmorphism(0,1280)(1,0)[\phantom{\mon(\cC_\cR,\dot{\otimes})}`\phantom{\mon(\cC,\otimes)}`\widetilde{F_\cR}]{1000}{-1}a
 \putmorphism(0,1220)(1,0)[\phantom{\mon(\cC_\cR,\dot{\otimes})}`\phantom{\mon(\cC,\otimes)}`\widetilde{U_\cR}]{1000}{1}b

 \putmorphism(0,480)(1,0)[\phantom{\cC_\cR}`\phantom{\cC}`F_\cR]{1000}{-1}a
 \putmorphism(0,420)(1,0)[\phantom{\cC_\cR}`\phantom{\cC}`U_\cR]{1000}{1}b

 \setsqparms[0`0`1`0;1000`800]
 \putsquare(1000,450)[\phantom{\mon(\cC,{\otimes})}`\mon(\cC^{\cR},\ddot{\otimes})`\phantom{\cC}`\cC^\cR;``\cU^{\ddot{\otimes}}`]

 \putmorphism(1000,1280)(1,0)[\phantom{\mon(\cC,\otimes)}`\phantom{\mon(\cC,\otimes)^{\widetilde{\cR}}}`{\widetilde{F^\cR}}]{1000}{1}a
 \putmorphism(1000,1220)(1,0)[\phantom{\mon(\cC,\otimes)}`\phantom{\mon(\cC,\otimes)^{\widetilde{\cR}}}`{\widetilde{U^\cR}}]{1000}{-1}b

 \putmorphism(1000,480)(1,0)[\phantom{\cC}`\phantom{\cC^\cR}`F^\cR]{1000}{1}a
 \putmorphism(1000,420)(1,0)[\phantom{\cC}`\phantom{\cC^\cR}`U^\cR]{1000}{-1}b

  \put(790,1020){\oval(100,100)[bl]}
  \put(740,1020){\line(1,1){185}}
  \put(790,970){\vector(1,1){185}}
  \put(700,865){${\widetilde{\cR}}$}

  \put(790,220){\oval(100,100)[bl]}
  \put(740,220){\line(1,1){185}}
  \put(790,170){\vector(1,1){185}}
   \put(700,65){${\cR}$}

\put(1100,300){$_{(\otimes, I )}$}
\put(-50,300){$_{(\dot{\otimes}, \dot{I} )}$}
\put(1850,300){$_{(\ddot{\otimes}, \ddot{I} )}$}

  \put(50,550){\line(1,1){100}}
  \put(150,650){\line(1,0){800}}
  \put(1050,650){\line(1,0){800}}
  \put(1850,650){\vector(1,-1){100}}
  \put(1300,675){${\bf\Phi}$}

  \put(50,1350){\line(1,1){100}}
  \put(150,1450){\line(1,0){1700}}
  \put(1850,1450){\vector(1,-1){100}}
  \put(1300,1475){${\bf\widetilde{\Phi}}$}

  \put(2800,750){$(1)$}
 \end{picture}
\end{center}
To save space, for 1-cells $K$ in $\cA$, we write $\widetilde{K}$ instead of $\mon(K)$. We describe below the upper row of the above diagram.

The monad $(\widetilde{\cR},\widetilde{\eta},\widetilde{\varepsilon})=\mon(\cR,\eta^\cR,\varepsilon^\cR)$, the lift of the monad $(\cR,\eta^\cR,\varepsilon^\cR)$, is defined on elements as follows. The composition of $\widetilde{\cR}$ with 1-cell
$$\mathds{M}=(\mathbf{M}:\cX\ra \cC,\mathbf{m}:\mathbf{M}\otimes\mathbf{M}\ra\mathbf{M},\mathbf{e}:\mathbf{I}\ra\mathbf{M}):\cX \ra \mon(\cC,\otimes)$$
is equal to
\[ \widetilde{\cR}\circ(\mathbf{M},\mathbf{m},\mathbf{e})=(\widetilde{\cR}(\mathbf{M}),\cR(\mathbf{m})\circ\phi_{\mathbf{M},\mathbf{M}},\,
\cR(\mathbf{e})\circ\bar{\varphi} )\]
Then $\cU^\otimes:(\mon(\cC,\otimes),\widetilde{\cR})\lra (\cC,\cR)$ is a morphism of monads with $\cU^\otimes\circ \widetilde{R}=R\circ \cU^\otimes$.

The composition of 1-cell $\widetilde{F_{\cR}}$ with a 1-cell $\mathds{M}:\cX \ra \mon(\cC,\otimes)$ is equal to
\[  \widetilde{F_{\cR}}(\mathds{M})= (\mathbf{M},\eta_{\mathbf{M}}\circ\mathbf{m}, \eta_{\mathbf{M}}\circ\mathbf{e}).\]

The composition of 1-cell
$\widetilde{U_{\cR}}$ with a 1-cell
$$\dot{\mathds{M}}=(\mathbf{M}:\cX\ra \cC,\mathbf{m}:\mathbf{M}\otimes\mathbf{M}\ra \cR(\mathbf{M}),\,\mathbf{e}:\mathbf{I}\ra \cR(\mathbf{M})): \cX\ra \mon(\cC_\cR,\dot{\otimes})$$
is equal to
\[  \widetilde{U_{\cR}}(\dot{\mathds{M}})= (\cR(\mathbf{M}),\mu_{\mathbf{M}}\circ\cR(\mathbf{m})\circ\phi_{\mathbf{M},\mathbf{M}}, \mathbf{e}).\]

The composition of 1-cell $\widetilde{F^{\cR}}$ with a 1-cell $\mathds{M}:\cX \ra \mon(\cC,\otimes)$ is equal to
\[  \widetilde{F^{\cR}}(\mathds{M})= (\cR(\mathbf{M}),\mu_{\mathbf{M}}, \cR(\mathbf{m})\circ\phi_{\mathbf{M},\mathbf{M}}, \cR(\mathbf{e})).\]

The composition of 1-cell
$\widetilde{U^{\cR}}$ with a 1-cell
$$\ddot{\mathds{M}}=(\mathbf{M}:\cX\ra \cC,\ba:\cR(\mathbf{M})\ra\mathbf{M},\mathbf{m}:\mathbf{M}\otimes\mathbf{M}\ra\mathbf{M},\,\mathbf{e}:\mathbf{I}\ra \cR(\mathbf{M})): \cX\ra \mon(\cC^\cR,\ddot{\otimes})$$
is equal to
\[  \widetilde{U^{\cR}}(\ddot{\mathds{M}})= (\cR(\mathbf{M}),\mu_{\mathbf{M}}\circ\cR(\mathbf{m})\circ\bu_{\mathbf{M},\mathbf{M}}, \mathbf{e}\circ\bar{\bu}).\]

\subsection{Free monoids and distributive law}\label{subsec-free-monoids-dl}

Now assume that the free $\otimes$-monoids exist i.e. $U^\otimes$ has a left adjoint, i.e. we have an adjunction $(F^\otimes\dashv U^\otimes,\eta^\otimes,\varepsilon^\otimes)$. For $\mathbf{A}:\cX\ra\cC$ we denote $F^\otimes(\mathbf{A})=( T^\otimes(\mathbf{A}),\bm_\mathbf{A},\be_\mathbf{A})$.
Then the monad $T^\otimes$ distributes over $\cR$, i.e. we have a distributive law
\begin{center} \xext=1600 \yext=950
\begin{picture}(\xext,\yext)(\xoff,\yoff)
 \settriparms[1`1`-1;800]
  \putVtriangle(0,50)[T^\otimes\cR`\cR T^\otimes=\cU^\otimes \widetilde{\cR}F^\otimes(X)`T^\otimes{\cR}T^\otimes=\cU^\otimes F^\otimes \cU^\otimes \widetilde{\cR}F^\otimes;\tau`T^\otimes\cR(\eta^\otimes)`\cU^\otimes(\varepsilon^\otimes)_{\widetilde{R}F^\otimes}]
\end{picture}
\end{center}
The free monoidal for both $\dot{\otimes}$ and $\ddot{\otimes}$ exists. We shall describe these free monoid 1-cells below.

The composite of the free $\dot{\otimes}$-monoid functor $F^{\dot{\otimes}}$ with $\mathbf{A}:\cX\ra\cC_\cR$ is equal to
\[ F^{\dot{\otimes}}(\mathbf{A})= (T^\otimes(\mathbf{A}), \eta_\mathbf{A}\circ \bm_\mathbf{A}:T^\otimes(\mathbf{A})\otimes T^\otimes(\mathbf{A})\ra \cR T^\otimes(\mathbf{A}), \eta_\mathbf{A}\circ \be_\mathbf{A} : \mathbf{I} \ra \cR T^\otimes(\mathbf{A})):\]
\[\cX\ra \mon(\cC_\cR,\dot{\otimes})  \]

Thus we have an adjunction $(F^{\dot{\otimes}}\dashv U^{\dot{\otimes}}, \eta^{\dot{\otimes}},\varepsilon^{\dot{\otimes}})$ and a monad $$(T^{\dot{\otimes}}=U^{\dot{\otimes}} F^{\dot{\otimes}},\eta^{\dot{\otimes}},\mu^{\dot{\otimes}}=U^{\dot{\otimes}} \varepsilon^{\dot{\otimes}} F^{\dot{\otimes}}).$$

The composite of the free $\ddot{\otimes}$-monoid functor $F^{\ddot{\otimes}}$ with $(\mathbf{A},\ba):\cX\ra\cC^\cR$ is equal to
\[ F^{\ddot{\otimes}}(\mathbf{A},\mathbf{a})=(\ddot{T}(\mathbf{A},\mathbf{a}),\ddot{\mathbf{a}},\ddot{\mathbf{m}}_{(\mathbf{A},\mathbf{a})},
\ddot{\mathbf{e}}_{(\mathbf{A},\mathbf{a})}) \]
Its $\cR$-algebra universe $(\ddot{T}(\mathbf{A},\mathbf{a}),\ddot{\mathbf{a}})$  is defined via the following diagram in $\cA(\cX,\cC)$, where the columns are coequalizers
\begin{center} \xext=2000 \yext=1600
\begin{picture}(\xext,\yext)(\xoff,\yoff)
 \setsqparms[1`1`1`1;1700`600]
 \putsquare(100,50)[\cR^2T^\otimes(\mathbf{A})`\cR T^\otimes(\mathbf{A})`\cR\ddot{T}(\mathbf{A},\mathbf{a})`\ddot{T}(\mathbf{A},\mathbf{a}); \mu^\cR_{T^\otimes(\mathbf{A})}`\cR(\bt_{(\mathbf{A},\mathbf{a})})`\bt_{(\mathbf{A},\mathbf{a})}`\ddot{\mathbf{a}}]
  \setsqparms[1`0`0`0;1700`800]
 \putsquare(100,650)[\cR^2T^\otimes\cR(\mathbf{A})`\cR T^\otimes\cR(\mathbf{A})`\phantom{\cR^2T^\otimes(\mathbf{A})}`\phantom{\cR T^\otimes(\mathbf{A})};\mu^\cR_{T^\otimes\cR(\mathbf{A})}```]

 \putmorphism(150,1450)(0,-1)[\phantom{\cR^2T^\otimes\cR(\mathbf{A})}`\phantom{\cR^2T^\otimes(\mathbf{A})}`\cR(\mu_{T^\otimes(\mathbf{A})}\circ\cR(\tau_{\mathbf{A}}))]{800}{1}r
  \putmorphism(50,1450)(0,-1)[\phantom{\cR^2T^\otimes\cR(\mathbf{A})}`\phantom{\cR^2T^\otimes(\mathbf{A})}`\cR^2T^\otimes(\mathbf{a})]{800}{1}l

 \putmorphism(1850,1450)(0,-1)[\phantom{\cR^2 T^\otimes(\mathbf{A})}`\phantom{\cR T^\otimes(\mathbf{A})}`\mu_{T^\otimes(\mathbf{A})}\circ\cR(\tau_{\mathbf{A}})]{800}{1}r
  \putmorphism(1750,1450)(0,-1)[\phantom{\cR^2 T^\otimes(\mathbf{A})}`\phantom{\cR T^\otimes(\mathbf{A})}`\cR T^\otimes(\mathbf{a})]{800}{1}l
 \end{picture}
\end{center}
The definitions of the multiplication $\ddot{\mathbf{m}}_{(\mathbf{A},\mathbf{a})}$ and of the unit $\ddot{\mathbf{e}}_{(\mathbf{A},\mathbf{a})}$ are left for the reader.
 Again, we have an adjunction $(F^{\ddot{\otimes}}\dashv U^{\ddot{\otimes}}, \eta^{\ddot{\otimes}},\varepsilon^{\ddot{\otimes}})$ and a monad $$(T^{\ddot{\otimes}}=U^{\ddot{\otimes}} F^{\ddot{\otimes}},\eta^{\ddot{\otimes}},\mu^{\ddot{\otimes}}=U^{\ddot{\otimes}} \varepsilon^{\ddot{\otimes}} F^{\ddot{\otimes}}).$$

Now the diagram in $\cA$ with six 0-cells looks as follows
\begin{center} \xext=2000 \yext=1650
\begin{picture}(\xext,\yext)(\xoff,\yoff)
 \setsqparms[0`1`1`0;1000`800]
 \putsquare(0,450)[\mon(\cC_\cR,\dot{\otimes})`\mon(\cC,\otimes)`\cC_\cR`\cC;``\cU^{\otimes}`]
  \put(40,800){$\cU^{\dot{\otimes}}$}
 \putmorphism(0,1280)(1,0)[\phantom{\mon(\cC_\cR,\dot{\otimes})}`\phantom{\mon(\cC,\otimes)}`\widetilde{F_\cR}]{1000}{-1}a
 \putmorphism(0,1220)(1,0)[\phantom{\mon(\cC_\cR,\dot{\otimes})}`\phantom{\mon(\cC,\otimes)}`\widetilde{U_\cR}]{1000}{1}b

 \putmorphism(0,480)(1,0)[\phantom{\cC_\cR}`\phantom{\cC}`F_\cR]{1000}{-1}a
 \putmorphism(0,420)(1,0)[\phantom{\cC_\cR}`\phantom{\cC}`U_\cR]{1000}{1}b

 \setsqparms[0`0`1`0;1000`800]
 \putsquare(1000,450)[\phantom{\mon(\cC,{\otimes})}`\mon(\cC^{\cR},\ddot{\otimes})`\phantom{\cC}`\cC^\cR;``\cU^{\ddot{\otimes}}`]

 \putmorphism(1000,1280)(1,0)[\phantom{\mon(\cC,\otimes)}`\phantom{\mon(\cC,\otimes)^{\widetilde{\cR}}}`{\widetilde{F^\cR}}]{1000}{1}a
 \putmorphism(1000,1220)(1,0)[\phantom{\mon(\cC,\otimes)}`\phantom{\mon(\cC,\otimes)^{\widetilde{\cR}}}`{\widetilde{U^\cR}}]{1000}{-1}b

 \putmorphism(1000,480)(1,0)[\phantom{\cC}`\phantom{\cC^\cR}`F^\cR]{1000}{1}a
 \putmorphism(1000,420)(1,0)[\phantom{\cC}`\phantom{\cC^\cR}`U^\cR]{1000}{-1}b

  \putmorphism(-50,1250)(0,-1)[\phantom{\mon(\cC_\cR,\dot{\otimes})}`\phantom{\cC_\cR}`F_{\dot{\otimes}}]{800}{-1}l
   \putmorphism(950,1150)(0,-1)[\phantom{\mon(\cC,\otimes)}`\phantom{\cC}`F_\otimes]{600}{-1}l
    \putmorphism(1950,1250)(0,-1)[\phantom{\mon(\cC^\cR,\ddot{\otimes})}`\phantom{\cC^\cR}`F_{\ddot{\otimes}}]{800}{-1}l

  \put(790,1020){\oval(100,100)[bl]}
  \put(740,1020){\line(1,1){175}}
  \put(790,970){\vector(1,1){175}}
  \put(630,965){${\widetilde{\cR}}$}

  \put(790,220){\oval(100,100)[bl]}
  \put(740,220){\line(1,1){175}}
  \put(790,170){\vector(1,1){175}}
   \put(700,65){${\cR}$}

 \put(190,220){\oval(100,100)[br]}
  \put(240,220){\line(-1,1){175}}
  \put(190,170){\vector(-1,1){175}}
   \put(200,65){${T^{\dot{\otimes}}}$}

    \put(1190,220){\oval(100,100)[br]}
  \put(1240,220){\line(-1,1){175}}
  \put(1190,170){\vector(-1,1){175}}
   \put(1200,65){${T^\otimes}$}

    \put(2190,220){\oval(100,100)[br]}
  \put(2240,220){\line(-1,1){175}}
  \put(2190,170){\vector(-1,1){175}}
   \put(2200,65){${T^{\ddot{\otimes}}}$}

  \put(50,550){\line(1,1){100}}
  \put(150,650){\line(1,0){750}}
  \put(1050,650){\line(1,0){800}}
  \put(1850,650){\vector(1,-1){100}}
  \put(1300,675){${\bf\Phi}$}

  \put(50,1350){\line(1,1){100}}
  \put(150,1450){\line(1,0){1700}}
  \put(1850,1450){\vector(1,-1){100}}
  \put(1300,1475){${\bf\widetilde{\Phi}}$}

  \put(2800,750){$(1)$}
 \end{picture}
\end{center}
Now $F_\cR:  (\cC,T^\otimes)\ra (\cC_\cR,T^{\dot{\otimes}})$ and $F^\cR:  (\cC,T^{\otimes})\ra (\cC^\cR,T^{\ddot{\otimes}})$
are strict morphisms of monads and $(U_\cR,\tau): (\cC_\cR,T^{\dot{\otimes}})\ra (\cC,T^{\otimes})$ and $(U^\cR,\nu): (\cC^\cR,T^{\ddot{\otimes}})\ra (\cC,T^{\otimes})$ are lax morphisms of monads, where $\nu_{(\mathbf{A},\mathbf{a})}=\bt_{(\mathbf{A},\mathbf{a})}\circ \eta^\cR_{T^\otimes(\mathbf{A})}:T^\otimes(\mathbf{A})\ra \ddot{T}(\mathbf{A},\mathbf{a})$. The functor $\Phi:(\cC_\cR,T^{\dot{\otimes}})\ra (\cC^\cR,T^{\ddot{\otimes}})$ is a strict monad morphism. Moreover, $(\cR,\tau,\eta^\cR,\mu^\cR)$ is a lax monad on the monad $(\cC,\cT^\otimes)$, i.e. a monad in $\Mnd_l(\cA)$.

Thus the bottom part of the above diagram is in fact a $kem$-diagram in  $\Mnd_l(\cA)$. Applying the Eilenberg-Moore 2-functor  $\bem$ to it we get the diagram which is equivalent to the top part of this diagram.

\subsection{Other categories of monoids}\label{subsec-other-monoids}

In this way, if $\bk$-objects for $T$-like monads exist,  we obtain six categories of monoids/algebras with combined $\cR$ and $T^\otimes$ structures. Here, we content ourselves with spelling them out
\[  \cC_{\cR T^\otimes} \ra \mon(\cC,\otimes)_{\widetilde{\cR}}\ra \mon(\cC_\cR,\dot{\otimes})\ra \mon(\cC,\ddot{\otimes})\ra \mon(\cC,\otimes)^{\widetilde{\cR}}\ra \cC^{\cR T^\otimes} \]
The last three are canonically equivalent. This follows from, suitably internalized, considerations in Section 3 of \cite{Beck}. The reader is invited to look at diferent four categories in the examples, at least in case of symmetrization monad on untyped signatures, to see differences. In case of typed signatures the first two 0-cells do not exist as fibrations. We think that this is a reason why they are not considered in the literature so often.

\subsection{Actions of monoidal objects}\label{subsec-Act}
 In this subsection we show how the theory described in previous sections extends when we have not only rc-lax-monoidal monad but also an action of monoidal object.

In this section we assume that we are given an rc-monoidal object  $\cC=(\cC,\otimes,I,\alpha,\lambda,\rho)$ in a 2-category $\cA$ with finite products, an rc-$0$-cells  $\cX$  in $\cA$ and a (strong) action
\[ (\star, \psi): \cC\times \cX \lra \cX \]
of $(\cC,\otimes)$ on $\cX$ in $\cA$.

\subsection*{$\cX$ exponentiable}
If $\cX$ is exponentiable, then by exponential adjunction we get a representation
\[ (\br,\varphi): \cC \lra \cX^\cX \]
that is, a (strong) monoidal morphism into a strict monoidal category object $\cX^\cX$.

\subsection*{Adjunction and $\bkem$-diagram for monoidal monad}
If $\br$ has an rc-right adjoint $U$,
\begin{center}
\xext=1200 \yext=200
\begin{picture}(\xext,\yext)(\xoff,\yoff)
\putmorphism(0,100)(1,0)[\cC`\cX^\cX`]{1200}{0}a
\putmorphism(0,160)(1,0)[\phantom{\cC}`\phantom{\cX^\cX}`U]{1200}{-1}a
\putmorphism(0,40)(1,0)[\phantom{\cC}`\phantom{\cX^\cX}`\br]{1200}{1}b
\end{picture}
\end{center}
then $U$ is again monoidal and the whole adjunction $((\br,\varphi)\dashv (U,\bu), \eta^\otimes, \varepsilon^\otimes)$ is monoidal. By Section \ref{monoidal-mnd from-adj} we get the $\bkem$-diagram in $\Mon_l(\cA)$  with representations $\br$, $\dot{\br}$, $\ddot{\br}$ of algebras

\begin{center} \xext=3200 \yext=1050
\begin{picture}(\xext,\yext)(\xoff,\yoff)
 \setsqparms[0`0`0`0;1600`1200]
 \putsquare(0,850)[\phantom{\mon(\cC_\cR,\dot{\otimes})}`\phantom{\mon(\cC,\otimes)}`\cC_\cR`\cC;``\phantom{\cU^{\otimes}}`]

 \putmorphism(0,880)(1,0)[\phantom{\cC_\cR}`\phantom{\cC}`F_\cR]{1600}{-1}a
 \putmorphism(0,820)(1,0)[\phantom{\cC_\cR}`\phantom{\cC}`U_\cR]{1600}{1}b

 \setsqparms[0`0`0`0;1600`1200]
 \putsquare(1600,850)[\phantom{\mon(\cC,{\otimes})}`\phantom{\mon(\cC^{\cR},\ddot{\otimes})}`\phantom{\cC}`\cC^\cR;``\phantom{\cU^{\ddot{\otimes}}}`]

 \putmorphism(1600,880)(1,0)[\phantom{\cC}`\phantom{\cC^\cR}`F^\cR]{1600}{1}a
 \putmorphism(1600,820)(1,0)[\phantom{\cC}`\phantom{\cC^\cR}`U^\cR]{1600}{-1}b

\put(1400,700){$_{(\otimes, I )}$}
\put(-150,700){$_{(\dot{\otimes}, \dot{I} )}$}
\put(3150,700){$_{(\ddot{\otimes}, \ddot{I} )}$}

\put(2000,0){$\cX^\cX$}
 \putmorphism(-50,800)(3,-1)[\phantom{\cC}`\phantom{\cC^\cR}`]{2200}{1}b
 \put(400,500){$\dot{\br}$}

 \putmorphism(1560,800)(2,-3)[\phantom{\cC}`\phantom{\cC^\cR}`]{470}{1}b
 \put(1660,500){$\br$}
  \putmorphism(1640,800)(2,-3)[\phantom{\cC}`\phantom{\cC^\cR}`]{470}{-1}b
 \put(1860,500){$U$}

  \putmorphism(3150,800)(-3,-2)[\phantom{\cC}`\phantom{\cC^\cR}`]{1100}{1}b
  \put(2600,500){$\ddot{\br}$}
   \putmorphism(3200,750)(-3,-2)[\phantom{\cC}`\phantom{\cC^\cR}`]{1100}{-1}b
   \put(2760,350){$K$}
 \end{picture}
\end{center}

\subsection*{Monoids}

Applying $\mon$ to the diagram in $\Mon_l(\cA)$  constructed above we get a diagram of categories of monoids with forgetful 1-cells to the (universes of) monoidal category objects.

\begin{center} \xext=3200 \yext=2350
\begin{picture}(\xext,\yext)(\xoff,\yoff)
 \setsqparms[0`1`0`0;1600`1200]
 \putsquare(0,850)[\mon(\cC_\cR,\dot{\otimes})`\mon(\cC,\otimes)`\cC_\cR`\cC;`\cU^{\dot{\otimes}}``]
 \put(1600,1850){\line(0,-1){360}}
  \put(1600,1400){\vector(0,-1){450}}
 \put(1440,1300){$\cU^{\otimes}$}
 \putmorphism(0,2080)(1,0)[\phantom{\mon(\cC_\cR,\dot{\otimes})}`\phantom{\mon(\cC,\otimes)}`\widetilde{F_\cR}]{1600}{-1}a
 \putmorphism(0,2020)(1,0)[\phantom{\mon(\cC_\cR,\dot{\otimes})}`\phantom{\mon(\cC,\otimes)}`\widetilde{U_\cR}]{1600}{1}b

 \putmorphism(0,880)(1,0)[\phantom{\cC_\cR}`\phantom{\cC}`F_\cR]{1600}{-1}a
 \putmorphism(0,820)(1,0)[\phantom{\cC_\cR}`\phantom{\cC}`U_\cR]{1600}{1}b

 \setsqparms[0`0`1`0;1600`1200]
 \putsquare(1600,850)[\phantom{\mon(\cC,{\otimes})}`\mon(\cC^{\cR},\ddot{\otimes})`\phantom{\cC}`\cC^\cR;``\cU^{\ddot{\otimes}}`]

 \putmorphism(1600,2080)(1,0)[\phantom{\mon(\cC,\otimes)}`\phantom{\mon(\cC,\otimes)^{\widetilde{\cR}}}`{\widetilde{F^\cR}}]{1600}{1}a
 \putmorphism(1600,2020)(1,0)[\phantom{\mon(\cC,\otimes)}`\phantom{\mon(\cC,\otimes)^{\widetilde{\cR}}}`{\widetilde{U^\cR}}]{1600}{-1}b

  \put(1700,880){\line(1,0){350}} \put(2130,880){\vector(1,0){950}} \put(2400,910){$F^\cR$}
   \put(2130,820){\line(1,0){950}} \put(2050,820){\vector(-1,0){350}} \put(2400,720){$U^\cR$}
\put(1400,700){$_{(\otimes, I )}$}
\put(-150,700){$_{(\dot{\otimes}, \dot{I} )}$}
\put(3150,700){$_{(\ddot{\otimes}, \ddot{I} )}$}

\put(2000,0){$\cX^\cX$}
 \putmorphism(-50,800)(3,-1)[\phantom{\cC}`\phantom{\cC^\cR}`]{2200}{1}b
 \put(400,500){$\dot{\br}$}

 \putmorphism(1560,800)(2,-3)[\phantom{\cC}`\phantom{\cC^\cR}`]{470}{1}b
 \put(1660,500){$\br$}
  \putmorphism(1640,800)(2,-3)[\phantom{\cC}`\phantom{\cC^\cR}`]{470}{-1}b
 \put(1860,500){$U$}

  \putmorphism(3150,800)(-3,-2)[\phantom{\cC}`\phantom{\cC^\cR}`]{1100}{1}b
  \put(2600,500){$\ddot{\br}$}
   \putmorphism(3200,750)(-3,-2)[\phantom{\cC}`\phantom{\cC^\cR}`]{1100}{-1}b
   \put(2760,350){$K$}

\put(2000,1200){$\mon(\cX^\cX)$}
 \putmorphism(-50,2000)(3,-1)[\phantom{\cC}`\phantom{\cC^\cR}`]{2200}{1}b
 \put(400,1700){$\widetilde{\dot{\br}}$}

 \putmorphism(1560,2000)(2,-3)[\phantom{\cC}`\phantom{\cC^\cR}`]{470}{1}b
 \put(1660,1700){$\widetilde{\br}$}
  \putmorphism(1640,2000)(2,-3)[\phantom{\cC}`\phantom{\cC^\cR}`]{470}{-1}b
 \put(1860,500){$U$}

  \putmorphism(3150,2000)(-3,-2)[\phantom{\cC}`\phantom{\cC^\cR}`]{1100}{1}b
  \put(2600,1700){$\widetilde{{\ddot{\br}}}$}
   \putmorphism(3200,1950)(-3,-2)[\phantom{\cC}`\phantom{\cC^\cR}`]{1100}{-1}b
   \put(2760,1550){$\widetilde{K}$}

 \putmorphism(2100,1200)(0,-1)[\phantom{\mon(\cC,\otimes)}`\phantom{\cC}`\cU]{1100}{1}r

 \end{picture}
\end{center}

In the examples, we will be interested in identifying the 0-cells $\cX^\cX$, $\mon(\cX^\cX)$ and the images of 1-cells $\dot{\br}$,  $\ddot{\br}$, $\mon(\dot{\br})$, and  $\mon(\ddot{\br})$ in them.

\subsection*{The lift to the lax slice}
The $\bkem$-diagram for the monoidal monad $(\cR,\phi,\eta,\varepsilon)$ lifts to a $\bkem$-diagram in the lax slice $\Mon_l(\cA)_{/_l\cX^\cX}$ and we get
\begin{center} \xext=3200 \yext=1650
\begin{picture}(\xext,\yext)(\xoff,\yoff)
 \setsqparms[0`0`0`0;1600`1200]
 \putsquare(0,1450)[\phantom{\mon(\cC_\cR,\dot{\otimes})}`\phantom{\mon(\cC,\otimes)}`\cC_\cR`\cC;``\phantom{\cU^{\otimes}}`]

 \putmorphism(0,1480)(1,0)[\phantom{\cC_\cR}`\phantom{\cC}`F_\cR]{1600}{-1}a
 \putmorphism(0,1420)(1,0)[\phantom{\cC_\cR}`\phantom{\cC}`U_\cR]{1600}{1}b

 \setsqparms[0`0`0`0;1600`1200]
 \putsquare(1600,1450)[\phantom{\mon(\cC,{\otimes})}`\phantom{\mon(\cC^{\cR},\ddot{\otimes})}`\phantom{\cC}`\cC^\cR;``\phantom{\cU^{\ddot{\otimes}}}`]

 \putmorphism(1600,1480)(1,0)[\phantom{\cC}`\phantom{\cC^\cR}`F^\cR]{1600}{1}a
 \putmorphism(1600,1420)(1,0)[\phantom{\cC}`\phantom{\cC^\cR}`U^\cR]{1600}{-1}b

\put(1400,1300){$_{(\otimes, I )}$}
\put(-150,1300){$_{(\dot{\otimes}, \dot{I} )}$}
\put(3200,1300){$_{(\ddot{\otimes}, \ddot{I} )}$}

\put(2000,600){$\cX^\cX$}
 \putmorphism(-50,1400)(3,-1)[\phantom{\cC}`\phantom{\cC^\cR}`]{2200}{1}b
 \put(450,1100){$\dot{\br}$}

 \putmorphism(1560,1400)(2,-3)[\phantom{\cC}`\phantom{\cC^\cR}`]{470}{1}b
 \put(1660,1100){$\br$}
  \putmorphism(1640,1400)(2,-3)[\phantom{\cC}`\phantom{\cC^\cR}`]{470}{-1}b
 \put(1860,1100){$U$}

  \putmorphism(3150,1400)(-3,-2)[\phantom{\cC}`\phantom{\cC^\cR}`]{1100}{1}b
  \put(2600,1100){$\ddot{\br}$}
   \putmorphism(3200,1350)(-3,-2)[\phantom{\cC}`\phantom{\cC^\cR}`]{1100}{-1}b
   \put(2760,950){$K$}

  \put(1600,0){$\cX^\cX$}

   \put(1600,800){\vector(0,-1){650}}
  \put(1600,1250){\line(0,-1){350}}
  \put(1530,700){$\br$}

  \putmorphism(-50,1400)(4,-3)[\phantom{\cC}`\phantom{\cC^\cR}`]{1650}{1}b
  \put(600,800){$\dot{\br}$}


  \put(3150,700){\vector(-2,-1){1200}}
  \put(3150,1300){\line(0,-1){600}}
  \put(2600,300){$\ddot{\br}$}

  \putmorphism(2050,600)(-1,-2)[\phantom{\cC}`\phantom{\cC^\cR}`]{250}{1}b
  \put(1950,300){$id$}
  \end{picture}
\end{center}
As a corollary of Theorems \ref{thm-lift-kem} and \ref{thm-em-Linton} we get the following.
\begin{theorem}\label{thm-kem-diag-in-slice}
The $\bkem$-diagram for the rc monoidal monad $(\cR,\phi,\bar{\phi},\eta^\cR,\mu^\cR)$ in $\Mon_l(\cA)$ lifts to a $\bkem$-diagram in the slice $\Mon_l(\cA)_{/_l\cX^\cX}$, for the lifted monad.
\end{theorem}

\subsection*{Back to actions}

Now we can move back the diagram to $\Act_l\Mon_l(\cA,\cX)$ getting the $\bkem$-diagram

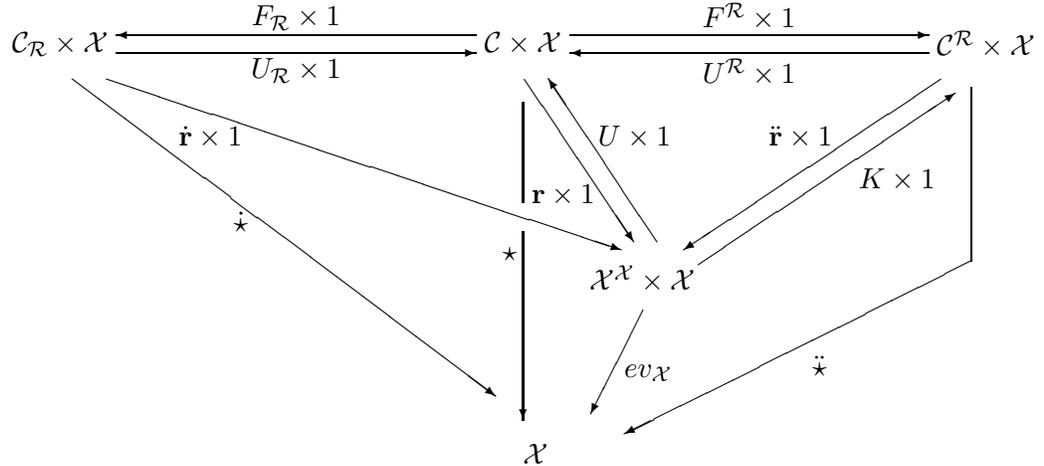
\begin{center} \xext=3200 \yext=1650
\begin{picture}(\xext,\yext)(\xoff,\yoff)
 \setsqparms[0`0`0`0;1600`1200]
 \putsquare(0,1450)[\phantom{\mon(\cC_\cR,\dot{\otimes})}`\phantom{\mon(\cC,\otimes)}`\cC_\cR\times\cX`\cC\times\cX;``\phantom{\cU^{\otimes}}`]

 \putmorphism(0,1480)(1,0)[\phantom{\cC_\cR\times\cX}`\phantom{\cC\times\cX}`F_\cR\times 1]{1600}{-1}a
 \putmorphism(0,1420)(1,0)[\phantom{\cC_\cR\times\cX}`\phantom{\cC\times\cX}`U_\cR\times 1]{1600}{1}b

 \setsqparms[0`0`0`0;1600`1200]
 \putsquare(1600,1450)[\phantom{\mon(\cC,{\otimes})}`\phantom{\mon(\cC^{\cR},\ddot{\otimes})}`\phantom{\cC}`\cC^\cR\times\cX;``\phantom{\cU^{\ddot{\otimes}}}`]

 \putmorphism(1600,1480)(1,0)[\phantom{\cC\times\cX}`\phantom{\cC^\cR\times\cX}`F^\cR\times 1]{1600}{1}a
 \putmorphism(1600,1420)(1,0)[\phantom{\cC\times\cX}`\phantom{\cC^\cR\times\cX}`U^\cR\times 1]{1600}{-1}b

\put(1830,600){$\cX^\cX\times\cX$}
 \putmorphism(-50,1400)(3,-1)[\phantom{\cC}`\phantom{\cC^\cR}`]{2200}{1}b
 \put(410,1100){$\dot{\br}\times 1$}

 \putmorphism(1560,1400)(2,-3)[\phantom{\cC}`\phantom{\cC^\cR}`]{470}{1}b
 \put(1630,900){$\br\times 1$}
  \putmorphism(1640,1400)(2,-3)[\phantom{\cC}`\phantom{\cC^\cR}`]{470}{-1}b
 \put(1860,1100){$U\times 1$}

  \putmorphism(3150,1400)(-3,-2)[\phantom{\cC}`\phantom{\cC^\cR}`]{1100}{1}b
  \put(2450,1100){$\ddot{\br}\times 1$}
   \putmorphism(3200,1350)(-3,-2)[\phantom{\cC}`\phantom{\cC^\cR}`]{1100}{-1}b
   \put(2760,950){$K\times 1$}

  \put(1600,0){$\cX$}

  \put(1600,800){\vector(0,-1){650}}
  \put(1600,1250){\line(0,-1){350}}
  \put(1530,700){$\star$}

  \putmorphism(-50,1400)(4,-3)[\phantom{\cC}`\phantom{\cC^\cR}`]{1650}{1}b
  \put(600,800){$\dot{\star}$}

  \put(3150,700){\vector(-2,-1){1200}}
  \put(3150,1300){\line(0,-1){600}}
  \put(2600,300){$\ddot{\star}$}

  \putmorphism(2050,600)(-1,-2)[\phantom{\cC}`\phantom{\cC^\cR}`]{250}{1}b
  \put(1950,300){$ev_\cX$}
  \end{picture}
\end{center}

As a corollary of Proposition \ref{prop-Mon-slice-eq-to-Act} and Theorem \ref{thm-kem-diag-in-slice}  we get

\begin{theorem}\label{thm-kem-diag-in-actions}
The $\bkem$-diagram for the rc monoidal monad $(\cR,\phi,\bar{\phi},\eta^\cR,\mu^\cR)$ in $\Mon_l(\cA)$ lifts to a $\bkem$-diagram in the slice $\Act_l\Mon_l(\cA,\cX)$, for the lifted monad.
\end{theorem}

\vskip 2mm
{\em Remark.} Recall from Section \ref{sec-mnds-algs} that the left adjoint $\ddot{r}:\cC^\cR\ra \cM$ is given by the (reflexive) coequalizer
\begin{center} \xext=1600 \yext=350
\begin{picture}(\xext,\yext)(\xoff,\yoff)
 \putmorphism(0,150)(1,0)[r \cR U^\cR`r U^\cR`]{800}{0}a
 \putmorphism(0,200)(1,0)[\phantom{r \cR U^\cR}`\phantom{r U^\cR}`r U^\cR(\varepsilon^\cR)]{800}{1}a
 \putmorphism(0,100)(1,0)[\phantom{r \cR U^\cR}`\phantom{r U^\cR}`\varepsilon^\otimes_{r U^\cR}]{800}{1}b
 \putmorphism(800,150)(1,0)[\phantom{r U^\cR}`\ddot{r}`q]{800}{1}a
\end{picture}
\end{center}
in $\cA(\cC^\cR,\cM)$ with the common inverse $r(\eta)_{U^\cR}$. In order to give this diagram a concrete flavour, we indicate below how the $\ddot{r}$ is defined in $Cat$. If $(A,\alpha)$ is an algebra in $\cC^\cR$ and $X$ is an object in $\cX$, then $\ddot{r}(A,\alpha)(X) = (A,\alpha)\ddot{\star} X$ is defined as the coequalizer in $\cX$
\begin{center} \xext=1600 \yext=350
\begin{picture}(\xext,\yext)(\xoff,\yoff)
 \putmorphism(0,150)(1,0)[\cR(A)\star X`A\star X`]{800}{0}a
 \putmorphism(0,200)(1,0)[\phantom{\cR(A)\star X}`\phantom{r U^\cR}`\alpha\star 1]{800}{1}a
 \putmorphism(0,100)(1,0)[\phantom{\cR(A)\star X}`\phantom{r U^\cR}`(\varepsilon^\otimes_A)_X]{800}{1}b
 \putmorphism(800,150)(1,0)[\phantom{r U^\cR}`(A,\alpha)\ddot{\star} X`q]{800}{1}a
\end{picture}
\end{center}

\subsection*{Actions along actions}

Then, applying $\act$ to the previous diagram, we get objects of actions and their representations in $\act(ev)$
\begin{center} \xext=3200 \yext=3050
\begin{picture}(\xext,\yext)(\xoff,\yoff)
 \setsqparms[0`1`0`0;1600`1400]
 \putsquare(0,1450)[\act(\dot{\star})`\act(\star)`\cC_\cR\times\cX`\cC\times\cX;``\cV`]
 \put(40,2200){${\dot{\cV}}$}
 \put(1600,2200){\vector(0,-1){650}}
  \put(1600,2650){\line(0,-1){350}}

 \putmorphism(0,2880)(1,0)[\phantom{\act(\dot{\star})}`\phantom{\act(\star)}`\widetilde{F_\cR}]{1600}{-1}a
 \putmorphism(0,2820)(1,0)[\phantom{\act(\dot{\star})}`\phantom{\act(\star)}`\widetilde{U_\cR}]{1600}{1}b

 \putmorphism(0,1480)(1,0)[\phantom{\cC_\cR\times\cX}`\phantom{\cC\times\cX}`F_\cR\times 1]{1600}{-1}a
 \putmorphism(0,1420)(1,0)[\phantom{\cC_\cR\times\cX}`\phantom{\cC\times\cX}`U_\cR\times 1]{1600}{1}b

 \setsqparms[0`0`1`0;1600`1400]
 \putsquare(1600,1450)[\phantom{\act(\star)}`\act(\ddot{\star})`\phantom{\cC}`\cC^\cR\times\cX;``\ddot{\cV}`]

 \putmorphism(1600,2880)(1,0)[\phantom{\act(\star}`\phantom{\act(\ddot{\star})}`{\widetilde{F^\cR}}]{1600}{1}a
 \putmorphism(1600,2820)(1,0)[\phantom{\act(\star}`\phantom{\act(\ddot{\star})}`{\widetilde{U^\cR}}]{1600}{-1}b

 \putmorphism(1600,1480)(1,0)[\phantom{\cC\times\cX}`\phantom{\cC^\cR\times\cX}`F^\cR\times 1]{1600}{0}a
 \put(2150,1480){\vector(1,0){850}}
  \put(1750,1480){\line(1,0){300}}
 \putmorphism(1600,1420)(1,0)[\phantom{\cC\times\cX}`\phantom{\cC^\cR\times\cX}`U^\cR\times 1]{1600}{0}b
 \put(2150,1420){\vector(1,0){850}}
  \put(1750,1420){\line(1,0){300}}

\put(1400,1300){$_{(\otimes, I )}$}
\put(-150,1300){$_{(\dot{\otimes}, \dot{I} )}$}
\put(3200,1300){$_{(\ddot{\otimes}, \ddot{I} )}$}

\put(2000,2000){$\act(ev)$}
 \putmorphism(-50,2800)(3,-1)[\phantom{\cC}`\phantom{\cC^\cR}`]{2200}{1}b
 \put(400,2500){$\widetilde{\dot{\br}}$}

 \putmorphism(1560,2800)(2,-3)[\phantom{\cC}`\phantom{\cC^\cR}`]{470}{1}b
 \put(1660,2500){$\widetilde{\br}$}
  \putmorphism(1640,2800)(2,-3)[\phantom{\cC}`\phantom{\cC^\cR}`]{470}{-1}b
 \put(1860,1300){$U$}

  \putmorphism(3150,2800)(-3,-2)[\phantom{\cC}`\phantom{\cC^\cR}`]{1100}{1}b
  \put(2600,2500){$\widetilde{\widetilde{\ddot{\br}}}$}
   \putmorphism(3200,2750)(-3,-2)[\phantom{\cC}`\phantom{\cC^\cR}`]{1100}{-1}b
   \put(2760,2350){$\widetilde{K}$}

 \putmorphism(2100,2000)(0,-1)[\phantom{\mon(\cC,\otimes)}`\phantom{\cC}`\phantom{\cU}]{1300}{1}r
  \put(2150,1800){$\cU$}

\put(1830,600){$\cX^\cX\times\cX$}
 \putmorphism(-50,1400)(3,-1)[\phantom{\cC}`\phantom{\cC^\cR}`]{2200}{1}b
 \put(410,1100){$\dot{\br}\times 1$}

 \putmorphism(1560,1400)(2,-3)[\phantom{\cC}`\phantom{\cC^\cR}`]{470}{1}b
 \put(1630,900){$\br\times 1$}
  \putmorphism(1640,1400)(2,-3)[\phantom{\cC}`\phantom{\cC^\cR}`]{470}{-1}b
 \put(1860,1100){$U\times 1$}

  \putmorphism(3150,1400)(-3,-2)[\phantom{\cC}`\phantom{\cC^\cR}`]{1100}{1}b
  \put(2450,1100){$\ddot{\br}\times 1$}
   \putmorphism(3200,1350)(-3,-2)[\phantom{\cC}`\phantom{\cC^\cR}`]{1100}{-1}b
   \put(2760,950){$K\times 1$}

  \put(1600,0){$\cX$}

  \put(1600,800){\vector(0,-1){650}}
  \put(1600,1250){\line(0,-1){350}}
  \put(1530,700){$\star$}

  \putmorphism(-50,1400)(4,-3)[\phantom{\cC}`\phantom{\cC^\cR}`]{1650}{1}b
  \put(600,800){$\dot{\star}$}


  \put(3150,700){\vector(-2,-1){1200}}
  \put(3150,1300){\line(0,-1){600}}
  \put(2600,300){$\ddot{\star}$}

  \putmorphism(2050,600)(-1,-2)[\phantom{\cC}`\phantom{\cC^\cR}`]{250}{1}b
  \put(1950,300){$ev_\cX$}
  \end{picture}
\end{center}
We also have 2-cells
\[ \sigma: \star\circ \cV\ra \pi_2\circ \cV,\hskip 3mm\dot{\sigma}: \dot{\star}\circ \dot{\cV}\ra \dot{\pi}_2\circ \dot{\cV},\hskip 3mm \ddot{\sigma}: \ddot{\star}\circ \ddot{\cV}\ra \ddot{\pi}_2\circ \ddot{\cV},\hskip 3mm \bar{\sigma}: ev_\cX\circ \cU\ra \bar{\pi}_2\circ \cU\]
where $ \pi_2:\cC\times \cX\ra \cX$, $ \dot{\pi}_2:\cC_\cR\times \cX\ra \cX$, $ \ddot{\pi}_2:\cC^\cR\times \cX\ra \cX$, $\bar{\pi}_2: \cX^\cX\times\cX\ra \cX$ are second projections.

\section{Applications}\label{sec-case-study}
\subsection{The general scheme}\label{subsec-general scheme}
In this subsection we list the ingredients that one needs to fix to apply the general considerations from previous sections.
\begin{enumerate}
  \item A {\em good algebraic}\footnote{Here by {\em algebraic} we mean a place (2-category) where one can built algebraic 0-cells of `all sorts of algebras' ($\bk$, $\bem$, $\mon$, $\act$) rather then a 2-category of algebras of some sort.} $2$-category $\cA$. By this we mean a 2-category with all finite products and some (better all) $\bem$- $\bk$- $\mon$- and $\act$-objects. Moreover, we would expect $\bk$-objects to be preserved by finite products and the comparison from $\bk$-objects to the defining adjunctions to be full and faithful.
  \item An rc-monoidal object $(\cC,\otimes,I,\alpha,\lambda,\varrho)$ in $\cA$ such that the free $\otimes$-monoid exists.
  \item A reasonable representation of $(\cC,\otimes)$. This can be provided by a left ideal $\cX$ in $\cC$ such that
  \begin{enumerate}
    \item both $\cX$ and the embedding $\cX\ra \cC$ are rc;
    \item $\cX$ is exponentiable in $\cA$;
    \item the representation $\br:\cC\ra \cX^\cX$ is faithful;
    \item $\br$ has a right adjoint $U$;
    \item the induced monad $\cF=U\br$ is an rc-monad (to have $U$ rc would be enough).
  \end{enumerate}
  \item The induced monad $\cF$ may have interesting submonads that we can identify either directly or via coreflective factorization of $\br\dashv U: \cC\ra \cX^\cX$.
  \item Finally, we can identify all the ingredients described in Section \ref{sec-Mon} for interesting submonads of $\cF$ and the images of the representations of various monoidal categories we obtained on the way.
\end{enumerate}

We apply this scheme to some cases that have been the motivation for the general theory developed in this paper so far.

\subsection{Some 2-categories and their basic properties} Below we discuss the 2-categories of possible interest with respect to the above scheme. As we already said in Section \ref{sec-Alg-2-Cats}, any 2-category $\cA$ determines two sorts of operations
\begin{enumerate}
  \item the global-external operations like: weighted limits and colimits and exponential objects;
  \item the local-internal operations like: internal limits and colimits of 0-cells.
\end{enumerate}
In particular, the local properties of 0-cells depend on the ambient 2-category $\cA$.

We want to apply the abstract considerations from previous sections to 2-categories other than $\Cat$ having 0-cells other than categories, yet we want both global and local operations to be of a very specific kind suitable for us.

The main application we have in mind are to 2-categories  $\cA$  that have (possibly among others) fibration of categories over a fixed base, say $\cB$.
The global operation suitable for us comes from the slice 2-category $\Cat_{/\cB}$ but the local ones come from 2-category $\Fib(\cB)$, \cite{Str1} of fibrations with 1-cells preserving prone morphisms, and 2-cells being vertical natural transformations. For example, the monoidal objects and the actions of monoidal objects that we consider on fibrations are always in those from $\Cat_{/\cB}$ and practically never in $\Fib(\cB)$ (even if the 0-cells on which they are defined are fibrations). The slogan is: the substitution is not cartesian\footnote{i.e. prone.}.  The exponential objects, even if they might not exist in $\Cat_{/\cB}$, when they do, they are much more interesting for us than those from  $\Fib(\cB)$. Moreover, when we perform operations on diagrams of fibrations in $\Cat_{/\cB}$, we do expect the results of our operations to be again fibrations. On the other hand, the local operations that we want to consider `in' fibrations are those from the 2-category $\Fib(\cB)$, i.e. performed in fibers and preserved by reindexing. Note, however, that when the `correct' local operations exist in a fibration, they are `preserved' by the embedding $\Fib(\cB)\ra \Cat_{/\cB}$ (one may say that the limits and colimits in 0-cells are preserved by this embedding).

As the objects in $\Cat$ are often intractable, we consider as a reasonable compromise the 2-category $\Fib_{/\cB}$, the full and locally full sub-2-category of $\Cat_{/\cB}$ whose $0$-cells are fibrations.  We also discus in this section the properties of the `bifibrational' versions of these categories.

All the following categories have finite products and the obvious embedding morphisms preserve them.

The 2-category $\Cat_{/\cB}$  is not cartesian closed with the exponentiable objects called Conduch\'{e} fibrations, cf. \cite{Gi}. It has all weighted limits and colimits. The (co)completeness of 0-cells means the (co)completeness in the fibers. However, overall the 0-cells in $\Cat_{/\cB}$ are not tractable.

The 2-category $\Fib_{/\cB}$ is not cartesian closed either but all bifibration are exponentiable, cf. \cite{Z2}, and the embedding $\Fib_{/\cB} \ra \Cat_{/\cB}$ is cartesian closed (in the sense that it preserves existing exponentials). It has good weighted limits and and some colimits. The thing of main importance to us is that it has $\bk$-objects of monads from $\Fib(\cB)$. The (co)completeness of 0-cells means the (co)completeness in the fibers.  The 0-cells are tractable and 1-cells are sufficiently flexible to describe all the necessary algebraic structures.

The 2-category $\BiFib_{/\cB}$ is the full and locally full sub-2-category of $\Cat_{/\cB}$ whose 0-cells are bifibrations. It is cartesian closed but it is too restrictive: neither limits nor colimits behave well.

The 2-category $\Fib(\cB)$ is cartesian closed but the embedding $\Fib{(\cB)} \ra \Cat_{/\cB}$ is not cartesian closed. The (co)limit in 0-cells are the (co)limits in all fibers that is preserved by reindexing functors.  This is the good notion of (co)limits in fibrations. The 0-cells are tractable but the 1-cells are too restrictive (e.g. monoidal structures on signatures do not live there).

We have the following bijective on object/full and locally full factorization
\begin{center}
\xext=1600 \yext=200
\begin{picture}(\xext,\yext)(\xoff,\yoff)
\putmorphism(0,100)(1,0)[\Fib(\cB)`\Fib_{/\cB}`]{800}{1}a
\putmorphism(800,100)(1,0)[\phantom{\Fib_{/\cB}}`\Cat_{/\cB}`]{800}{1}a
\end{picture}
\end{center}
the first morphism preserves limits and colimits of 0-cells (hence also completeness and cocompleteness of 0-cells), the second preserves the existing global operations (weighted limits, colimits and exponentials). Both morphisms preserve tensors (product $\cE\times \cJ\ra \cE\ra \cB$) and cotensors ($\cE^{(\cJ)}\ra \cB$). This is for the calculation and preservation of limits.

We note for the record.

\begin{proposition}
Let $p: \cE\ra \cB$ be a fibration. We have
\begin{enumerate}
  \item $p$ has (co)limits of type $J$ in $\Fib_{/\cB}$ iff $p$ has (co)limits of type $J$ in the fibers.
  \item $p$ has (co)limits of type $J$ in $\Fib(\cB)$ iff $p$ has (co)limits of type $J$ in the fibers and the reindexing functors preserve them.
\end{enumerate}
\end{proposition}
{\it Proof.}~ Simple check. $\Box$

\begin{proposition}
A 1-cell
\begin{center} \xext=600 \yext=350
\begin{picture}(\xext,\yext)(\xoff,\yoff)
 \settriparms[1`1`1;300]
  \putVtriangle(0,0)[\cE`\cC`\cB ;U`q`p]
\end{picture}
\end{center}
in $\Fib_{/\cB}$ is monadic iff $U$ is morphism of fibrations (i.e. preserves prone morphisms) and it is monadic when restricted to fibers.
\end{proposition}
{\it Proof.}~ Simple check. $\Box$


\subsection{Burroni fibrations and their actions}
Burroni fibrations were introduced in \cite{Bu} and further studied in \cite{Z2}. For more elaborate description the reader may consult these sources.

Let $\cB$ be a category with pullbacks. We shall work in the 2-category $Fib_{/\cB}$. If $T$ is a cartesian monad on $\cB$, then it gives rise to a monoidal object (a lax monoidal fibration) of $T$-graphs $p_T: Gph(T)\ra\cB$ in $Fib_{/\cB}$.

An object $\lk A,O,\gamma,\delta\rk $ of $Gph(T)$ is a span
\begin{center} \xext=800 \yext=430
\begin{picture}(\xext,\yext)(\xoff,\yoff)
\settriparms[1`1`0;400]
  \putAtriangle(0,0)[A`O`T(O);\gamma`\delta`]
\end{picture}
\end{center}
in $\cB$. The morphisms $\gamma$ and $\delta$ are called the {\em codomain map} and the {\em domain map} of the $T$-graph
$\lk A,O,\gamma,\delta\rk$, respectively.
Sometimes we write $A$ instead of $\lk A,O,\gamma,\delta\rk$,  for short,  when it does not lead to a confusion.

A morphism of $T$-graphs $\lk f,u\rk :\lk A,O,\gamma,\delta\rk \lra \lk A',O',\gamma',\delta'\rk$ is a pair of morphisms $f:A\ra A'$ and $u:O\ra O'$ in $\cB$ making the squares
\begin{center} \xext=1900 \yext=550
\begin{picture}(\xext,\yext)(\xoff,\yoff)
 \setsqparms[1`1`1`1;500`400]
  \putsquare(0,50)[A`A'`O`O';f`\gamma`\gamma'`u]
  \setsqparms[1`1`1`1;600`400]
  \putsquare(1300,50)[A`A'`T(O)`T(O');
  f`\delta`\delta'`T(u)]
\end{picture}
\end{center}
commute. The projection functor
\[ p_T:Gph(T)\lra \cB, \]
 sending the morphism $\lk f,u\rk :\lk A,O,\gamma,\delta\rk \lra \lk A',O',\gamma',\delta'\rk $ to the morphism $u:O\ra O'$,
is a fibration, cf. \cite{Bu} p. 235. The monoidal structure in $p_T$ is defined as follows.
Let $\lk A, O, \gamma_A,\delta_A\rk$ and $\lk B, O,
\gamma_B,\delta_B\rk$ be two objects in the fibre over $O$, i.e.
in $Gph(T)_O$. Then the tensor
\[ \lk A, O,\gamma_A,\delta_A\rk\otimes_O\lk B, O, \gamma_B,\delta_B\rk =
\lk A\otimes B, O, \gamma_\otimes, \delta_\otimes \rk\]
 is defined from the following diagram
\begin{center} \xext=1800 \yext=1230
\begin{picture}(\xext,\yext)(\xoff,\yoff)
\settriparms[1`1`0;400]
  \putAtriangle(0,400)[A`O`T(O);\gamma_A``]
    \put(500,530){$\delta_A$}
  \putAtriangle(800,400)[T(B)`\phantom{T(O)}`T^2(O);`T(\delta_B)`]
  \put(1020,530){$T(\gamma_B)$}
  \settriparms[0`1`0;400]
  \putAtriangle(1200,0)[\phantom{T^2(O)}``T(O);`\mu_O`]
  \settriparms[1`1`0;400]
   \putAtriangle(400,800)[A\otimes B`\phantom{A}`\phantom{T(B)};\pi_1`\pi_2`]
\end{picture}
\end{center}
in which the square is a pullback and
\[\gamma_\otimes =  \gamma_A\circ\pi_1,\;\;\; \delta_\otimes = \mu_O\circ T(\delta_B)\circ\pi_2.\]
The unit in the fibre over $O$ is
\begin{center} \xext=800 \yext=430
\begin{picture}(\xext,\yext)(\xoff,\yoff)
\settriparms[1`1`0;400]
  \putAtriangle(0,0)[O`O`T(O);1_O`\eta_O`]
\end{picture}
\end{center}
The coherence morphisms are defined using the universal properties of pullbacks. Note that, as  $T$ is assumed here to be cartesian, the coherence morphisms are isomorphisms. See \cite{Z2} for details. If both $\cB$ and $T$ are rc, then the fibration of graphs $p_T$ is rc in $Fib_{/\cB}$. If, moreover, reflexive coequalizers are pullback stable, for example, if $\cB$ is both rc and Barr-exact, then $p_T$ is rc in $Fib(\cB)$ as well, i.e. reflexive coequalizers not only exist in the fibers of $p_T$ but they are also preserved by reindexing.  Note that even if the fibers of $Gph(T)\ra \cC$ are monoidal categories (coherences are isos), the reindexing functors are automatically lax monoidal but almost never strong monoidal, i.e. coherence morphisms for them are usually non-invertible.

In fact, any endofunctor $N:\cB\ra \cB$ gives rise to a fibration $q_{N}:Gph(N)\lra \cB$ of $N$-graphs over $\cB$ except that such a fibration does not carry the monoidal structure in general. However, if $N$ is a $T$-module, i.e. it is equipped with a natural transformation $\nu:TN\ra N$  such that
\[ \nu\circ T(\nu)=\nu\circ \mu_N ,\hskip 3mm \nu\circ \eta_N =1_N\]
then 
the monoidal fibration $p_T:Gph(T)\ra \cB$ acts on the fibration $q_{N}:Gph(N)\ra \cB$ in a canonical way. The action
\begin{center} \xext=800 \yext=620
\begin{picture}(\xext,\yext)(\xoff,\yoff)
 \settriparms[1`1`1;500]
  \putVtriangle(0,0)[Gph(T)\times_\cB Gph(N)`Gph(N)`\cB;\star_\nu``cod]
\end{picture}
\end{center}
is defined on objects as a pullback
\begin{center} \xext=1800 \yext=1230
\begin{picture}(\xext,\yext)(\xoff,\yoff)
\settriparms[1`1`0;400]
  \putAtriangle(0,400)[A`O`T(O);\gamma_A``]
    \put(500,530){$\delta_A$}
  \putAtriangle(800,400)[T(X)`\phantom{T(O)}`T\,N(O);`T(\delta_X)`]
  \put(1020,530){$T(\gamma_X)$}
  \settriparms[0`1`0;400]
  \putAtriangle(1200,0)[\phantom{T\,N(O)}``N(O);`\nu_O`]
  \settriparms[1`1`0;400]
   \putAtriangle(400,800)[A\star X`\phantom{A}`\phantom{T(X)};\pi_1`\pi_2`]
\end{picture}
\end{center}

By the exponential adjunction, we get a morphism of lax monoidal fibrations
\begin{center} \xext=800 \yext=620
\begin{picture}(\xext,\yext)(\xoff,\yoff)
 \settriparms[1`1`1;500]
  \putVtriangle(0,30)[Gph(T)`Exp(Gph(N))`\cB;rep_T`p_T`p_{exp}]
\end{picture}
\end{center}
that represents $T$-graphs as endofunctors on fibers of $p_{N}:Gph(N)\ra\cB$, cf. \cite{Z2}. Under this representation the $T$-categories
correspond to (some) monads on fibers of $p_{N}:Gph(N)\ra\cB$.

\vskip 2mm
{\em Examples.} A subfunctor $i: T'\ra T$ of a monad $(T,\eta,\mu)$  is {\em left closed (under mutltiplication)} iff
there is a natural transformation $\mu': T\,T'\ra T'$ making the square
\begin{center} \xext=500 \yext=550
\begin{picture}(\xext,\yext)(\xoff,\yoff)
\setsqparms[1`1`1`1;700`400]
 \putsquare(0,50)[T\,T'`T'`T\,T`T;\mu'`T(i)`i`\mu]
 \end{picture}
\end{center}
commute. Any such subfunctor has a canonical structure of a module over $T$ and gives rise to an ideal of $T'$-graphs in the monoidal fibration of $T$-graphs.
If $\cB$ has a terminal object $1$ and $N$ is a constant functor equal $1$, then $p_{N}:Gph(N)\ra\cB$ is basic fibration over $\cB$ and $N$ has a unique $T$-module structure. In this case the induced action is the usual tautologous action.

On the other hand, if $N$ is $T$ and the $T$-module structure is $\mu$, then the induced action is the usual action of a monoidal object on itself.

\subsection{Fibrations of signatures over $Set$}

In this subsection we shall apply the abstract theory developed in previous sections to the Burroni fibration of signatures, i.e. the fibrations of $\bM$-graphs for $\bM$ the monads of monoids in $Set$.  We list below some closed subfunctors of the monad $\bM$
  \begin{enumerate}
    \item the empty word(s!) $\bM_\emptyset\ra \bM$,
    \item the  words of even length  $\bM_2\ra \bM$,
    \item the  words of length divisible by $p$  $\bM_p\ra \bM$ ($p$ positive integer),
  \end{enumerate}
Thus these functors give rise to closed ideals in the fibration of signatures. The ideal $Gph(\bM_\emptyset)\ra Set$ is equivalent to the basic fibration $Set^\ra \ra Set$. This fibration can be identified with the ideal fibration of the signatures that consists of symbols of constants only.

{\em Notation.} An object of the Burroni fibration of signatures $p_\bM:Grph(\bM)\ra Set$ in the fiber over set $O$ will be identified with a pair $(A,\partial:A\ra O^\dagger)$ where $O^\dagger=\coprod_{n\in\o} O^{[n]}$ and $[n]=\{ 0,\ldots,n\}$. We write $\partial_a$ for the value of $\partial$ on $a$ and $A_n$ for the set $\{a\in A\,|\, dom(\partial_a)=[n] \}$, for $n\in \o$. We also use $(n]=\{ 1,\ldots,n\}$ and  $\partial^+_a=\partial_a{\lceil(n]}$.
We often write $A$ when the typing map $\partial$ is understood. For $d: X\ra O$ in $Set_{/O}$ and and $o\in O$, we write $X_o$ for $d^{-1}(\{o\})\subseteq X$
and for $\alpha:(n]\ra O$, $X_\alpha$ denotes $\prod_{i=1}^n X_{\alpha(i)}$. By $\F$, $\E$, $\B$ we denote the skeletal categories of finite sets with objects $(n]$ for $n\in \o$ and all functions, surjections, bijections, respectively.

\vskip 2mm

The fibration $p_\bM$ acts on the basic fibration
\begin{center} \xext=800 \yext=620
\begin{picture}(\xext,\yext)(\xoff,\yoff)
 \settriparms[1`1`1;500]
  \putVtriangle(0,0)[Gph(\bM)\times_{Set}Set^\ra`Set^\ra`Set;\star^{\bM}``cod]
\end{picture}
\end{center}
as described in the previous section, cf. \cite{Z2}, and its exponential adjoint in $Fib_{/Set}$
\begin{center} \xext=800 \yext=620
\begin{picture}(\xext,\yext)(\xoff,\yoff)
 \settriparms[1`1`1;500]
  \putVtriangle(0,30)[Gph(\bM)`Exp(Set^\ra)`Set;\br=rep_{\bM}`p_{\M}`p_{exp}]
\end{picture}
\end{center}
preserves the prone morphisms, i.e. it is a morphism in $Fib(Set)$, cf. \cite{Z2}. $\br$ is also a strong monoidal 1-cell. For $(A,\partial)$ in $Gph(\bM)_O$, the functor
\[ \br(A)=\br(A,\partial): Set_{/O}\lra Set_{/O} \]
is defined for $(X,d:X\ra O)$ as
\[ \br(A)(X,d)= (A\star X, d': A\star X \lra O) \]
so that
\[  A\star X =\{ \lk a,x_1,\ldots, x_n\rk\,|\, n\in \o,\, a\in A_n,\, x_i\in X,\; \partial_a(i)=d(x_i),\; {\rm for}\; i=1,\ldots, n  \} \cong  \coprod_{a\in A} X_{\partial^+(a)} \]
and
\[ d'(\lk a,x_1,\ldots, x_n\rk)=\partial_a(0).\]
The right adjoint $U$ to $\br$ in $Fib_{/Set}$ is defined as follows. For a functor $H:Set_{/O}\ra Set_{/O}$, we have
\[ U_O(H)= \coprod_{\alpha\in O^\dagger} H(\alpha^+)_{\alpha(0)}\]
where for $\alpha\in O^{[n]}$, $\alpha^+=\alpha_{\lceil(n]}$. We put $\partial(a)=\alpha$, for $a\in H(\alpha^+)_{\alpha(0)}$.

Since $\br$ preserves prone morphisms, to see that $U$ is a right adjoint, it is enough to verify this in the fibers only. As arities of operations are finite, $U$ is rc 1-cell in $Fib_{/Set}$. As $U$ is a monoidal 1-cell in $Fib_{/Set}$ and we have a monoidal adjunction $(\br,\varphi)\dashv (U,\bu)$ in $\Mon_l(Fib_{/Set})$.  Thus we have an rc monad $\cF=U\br$ on $p_\bM:Grph(\bM)\ra Set$. The monad $\cF$ is defined as follows. For $(A,\partial)$ in $Grph(\bM)_O$, and $n\in\o$, we have
\[ \cF(A,\partial)_n= \coprod_{m\in \o} \F((m],(n])\times A_m \]
The monad $\cF$ has two obvious submonads related to subcategories $\B$ and $\E$ of $\F$.

The {\em symmetrization monad} $\cS$ is defined as
\[ \cS(A,\partial)_n= \coprod_{m\in \o} \B((m],(n])\times A_m \]
and the monad for the operads with the actions of surjections is defined as
\[ \cR(A,\partial)_n= \coprod_{m\in \o} \E((m],(n])\times A_m. \]

These submonads of $\F$ can be also identified via monoidal subfibratons of $p_{exp}:End(Set^\ra)\ra Set$ as follows.

Recall that the fiber of $p_{exp}$ over the set $O$ consists of endofunctors on $Set_{/O}$. A morphism from $P:Set_{/O}\ra Set_{/O}$ to $Q:Set_{/Q}\ra Set_{/Q}$ over $u:O\ra Q$ in $Set$ is a natural transformation $\tau : Pu^*\ra u^*Q$, i.e. a morphism in $\Cat(Set_{/Q}, Set_{/O})$, see \cite{Z2} for more.

If we restrict endofunctors $P$ to those that preserve pullbacks along monos, and  natural transformation $\tau$ to those whose naturality squares for monos are pullback, we obtain a subfibration of $p_{exp}$ of semi-analytic functors\footnote{It would be reasonable to ask these functors to be finitary but the resulting lax monoidal monad will be the same anyway}. We denote this fibration as  $p_{exp,\,sa}: End_{sa}(Set^\ra)\ra Set$.

Moreover, if we restrict endofunctors $P$ to those that weakly preserve wide pullback and  natural transformation $\tau$ to those that are weakly cartesian, we obtain a subfibration of $p_{exp}$ of analytic functors. We denote it  $p_{exp,\,an}: End_{an}(Set^\ra)\ra Set$.

Thus we have a factorization of $\br$ by embeddings
\[ Gph(\bM) \lra   End_{an}(Set^\ra)\lra End_{sa}(Set^\ra)\lra End(Set^\ra) \]
so that each morphism is strongly monoidal with a (lax monoidal) right adjoint. None of these embeddings is full, but the last two  are full on isomorphisms.

Thus we get monads $\cR$ and $\cS$ on the fibration $Gph(M)\ra Set$  form the embeddings $Gph(\bM)\lra End_{sa}(Set^\ra)$ and $Gph(\bM) \lra   End_{an}(Set^\ra)$, respectively. We shall describe some of the objects that these monads generate that we have studied earlier.

\subsection*{Symmetrization monad $\cS$}
We start with the {\em symmetrization} monad $\cS$. We shall describe the fibrations it generates and some morphisms, i.e. the diagram
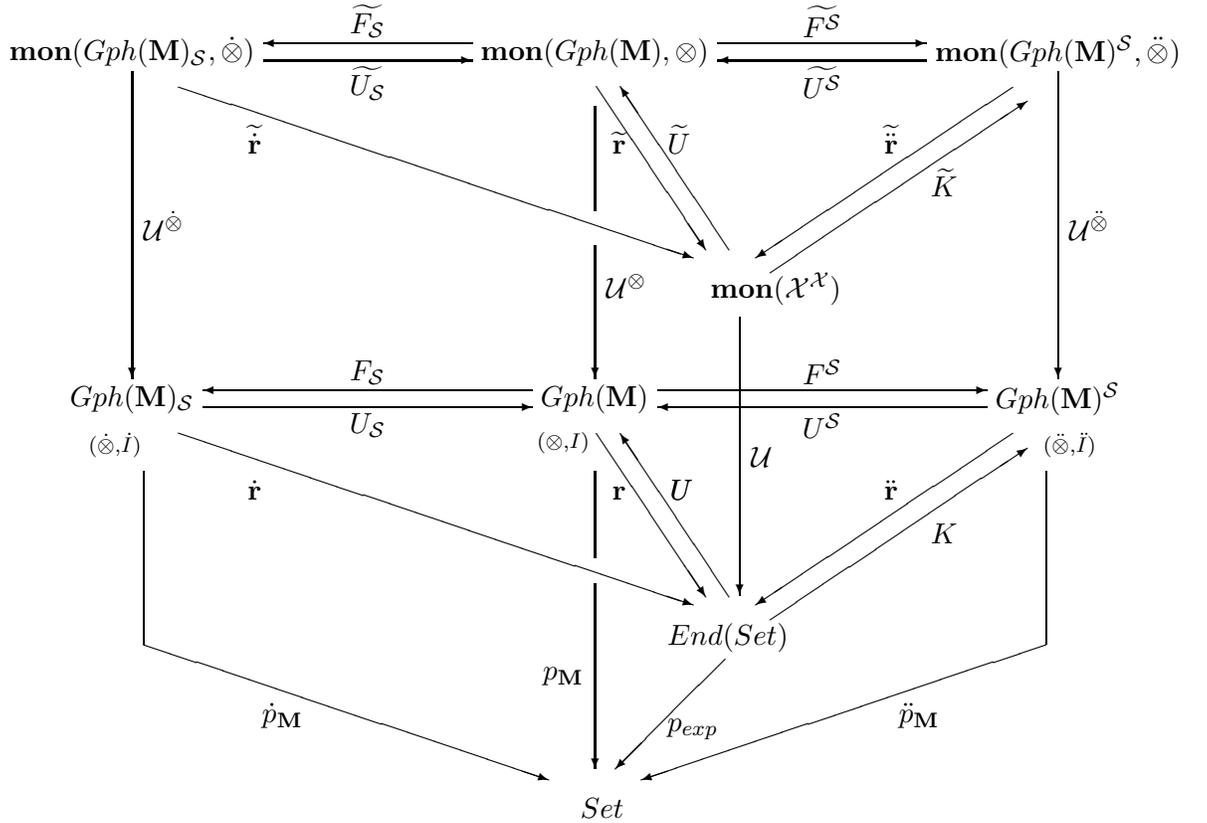
\begin{center} \xext=3200 \yext=2950
\begin{picture}(\xext,\yext)(\xoff,\yoff)
 \setsqparms[0`1`0`0;1600`1200]
 \putsquare(0,1450)[\mon(Gph(\bM)_\cS,\dot{\otimes})`\mon(Gph(\bM),\otimes)`Gph(\bM)_\cS`Gph(\bM);```]
  \put(1640,1800){$\cU^{\otimes}$}
 \putmorphism(1600,2050)(0,-1)[\phantom{\mon(Gph(\bM),\otimes)}`\phantom{(Gph(\bM),\otimes)}`]{600}{1}l
 \put(1600,2460){\line(0,-1){360}}

  \put(40,2000){$\cU^{\dot{\otimes}}$}
 \putmorphism(0,2680)(1,0)[\phantom{\mon(Gph(\bM)_\cS,\dot{\otimes})}`\phantom{\mon(Gph(\bM),\otimes)}`\widetilde{F_\cS}]{1600}{-1}a
 \putmorphism(0,2620)(1,0)[\phantom{\mon(Gph(\bM)_\cS,\dot{\otimes})}`\phantom{\mon(Gph(\bM),\otimes)}`\widetilde{U_\cS}]{1600}{1}b

 \putmorphism(0,1480)(1,0)[\phantom{Gph(\bM)_\cS}`\phantom{Gph(\bM)}`F_\cS]{1600}{-1}a
 \putmorphism(0,1420)(1,0)[\phantom{Gph(\bM)_\cS}`\phantom{Gph(\bM)}`U_\cS]{1600}{1}b

 \setsqparms[0`0`1`0;1600`1200]
 \putsquare(1600,1450)[\phantom{\mon(Gph(\bM),{\otimes})}`\mon(Gph(\bM)^{\cS},\ddot{\otimes})`\phantom{Gph(\bM)}`Gph(\bM)^\cS;``\cU^{\ddot{\otimes}}`]

 \putmorphism(1600,2680)(1,0)[\phantom{\mon(Gph(\bM),\otimes)}`\phantom{\mon(Gph(\bM),\otimes)^{\widetilde{\cS}}}`{\widetilde{F^\cS}}]{1600}{1}a
 \putmorphism(1600,2620)(1,0)[\phantom{\mon(Gph(\bM),\otimes)}`\phantom{\mon(Gph(\bM),\otimes)^{\widetilde{\cS}}}`{\widetilde{U^\cS}}]{1600}{-1}b

 \putmorphism(1600,1480)(1,0)[\phantom{Gph(\bM)}`\phantom{Gph(\bM)^\cS}`F^\cS]{1600}{1}a
 \putmorphism(1600,1420)(1,0)[\phantom{Gph(\bM)}`\phantom{Gph(\bM)^\cS}`U^\cS]{1600}{-1}b

\put(1400,1300){$_{(\otimes, I )}$}
\put(-150,1300){$_{(\dot{\otimes}, \dot{I} )}$}
\put(3150,1300){$_{(\ddot{\otimes}, \ddot{I} )}$}

\put(1850,600){$End(Set)$}
 \putmorphism(-50,1400)(3,-1)[\phantom{Gph(\bM)}`\phantom{Gph(\bM)^\cS}`]{2200}{1}b
 \put(400,1100){$\dot{\br}$}

 \putmorphism(1560,1400)(2,-3)[\phantom{Gph(\bM)}`\phantom{Gph(\bM)^\cS}`]{470}{1}b
 \put(1660,1100){$\br$}
  \putmorphism(1640,1400)(2,-3)[\phantom{Gph(\bM)}`\phantom{Gph(\bM)^\cS}`]{470}{-1}b
 \put(1860,1100){$U$}

  \putmorphism(3150,1400)(-3,-2)[\phantom{Gph(\bM)}`\phantom{Gph(\bM)^\cS}`]{1100}{1}b
  \put(2600,1100){$\ddot{\br}$}
   \putmorphism(3200,1350)(-3,-2)[\phantom{Gph(\bM)}`\phantom{Gph(\bM)^\cS}`]{1100}{-1}b
   \put(2760,950){$K$}

\put(2000,1800){$\mon(\cX^\cX)$}
 \putmorphism(-50,2600)(3,-1)[\phantom{Gph(\bM)}`\phantom{Gph(\bM)^\cS}`]{2200}{1}b
 \put(400,2300){$\widetilde{\dot{\br}}$}

 \putmorphism(1560,2600)(2,-3)[\phantom{Gph(\bM)}`\phantom{Gph(\bM)^\cS}`]{470}{1}b
 \put(1660,2300){$\widetilde{\br}$}
  \putmorphism(1640,2600)(2,-3)[\phantom{Gph(\bM)}`\phantom{Gph(\bM)^\cS}`]{470}{-1}b
 \put(1860,1100){$U$}

  \putmorphism(3150,2600)(-3,-2)[\phantom{Gph(\bM)}`\phantom{Gph(\bM)^\cS}`]{1100}{1}b
  \put(2600,2300){$\widetilde{{\ddot{\br}}}$}
   \putmorphism(3200,2550)(-3,-2)[\phantom{Gph(\bM)}`\phantom{Gph(\bM)^\cS}`]{1100}{-1}b
   \put(2760,2150){$\widetilde{K}$}

 \putmorphism(2100,1800)(0,-1)[\phantom{\mon(Gph(\bM),\otimes)}`\phantom{Gph(\bM)}`\cU]{1100}{1}r

  \put(1550,0){$Set$}
   \putmorphism(1600,880)(0,-1)[\phantom{Gph(\bM)}`\phantom{Set}`p_\bM]{780}{1}l
   \put(1600,1200){\line(0,-1){300}}

    \put(450,330){$\dot{p}_\bM$}
     \put(40,1200){\line(0,-1){600}}
     \put(40,600){\vector(3,-1){1400}}

    \put(2650,330){$\ddot{p}_\bM$}
    \put(3160,1200){\line(0,-1){600}}
     \put(3160,600){\vector(-3,-1){1400}}

    \put(1850,300){$p_{exp}$}
    \put(2050,550){\vector(-1,-1){380}}

    \put(1850,2300){$\widetilde{U}$}
 \end{picture}
\end{center}
As we indicated, $Gph(\bM)_O$ is the category of typed signatures with the set of types $O$ and morphism that preserves typing of the operations (strictly).
The total category of the Kleisli fibration $\dot{p}_\bM$ is the category of signatures (as before) but now the morphisms are a bit richer as they allow for amalgamations i.e. a morphism $(h,\sigma):A\ra B$ in $(Gph(\bM)_\cS)$  over $u:O\ra Q$ is a function $f:A\ra B$, and for $a\in A$ with typing $\partial_a:[n]\ra O$,  $\sigma_a\in S_n$ is a permutation so that
\[  u\circ \partial_a = \partial_{h(a)}\circ \sigma_a  \]
where here, and whenever reasonable, we think of $\sigma_a$ as a bijection from $[n]$ to $[n]$ sending $0$ to $0$.

The total category of the Eilenberg-Moore fibration $\ddot{p}_\bM$ is the category of signatures with actions of symmetric groups $S_n$ on $n$-ary operations.  The morphism  $h:(A,\cdot)\ra (B,\cdot)$ in $(Gph(\bM)^\cS)$  over $u:O\ra Q$ is a function $f:A\ra B$ that preserves both typing i.e.  $u\circ \partial_a=\partial_{h(a)}$ and actions, i.e. for $a\in A_n$ and $\sigma\in S_n$, is a permutation so that
\[  u\circ \partial_a = \partial_{h(a)}\circ \sigma_a.  \]
Thus  $(Gph(\bM)^\cS)_O$ is equivalent to the category of $O$-coloured species. As we said before, $End(Set)_O$ is the category of endofunctors on $Set_{/O}$.

$\ddot{p}_\bM$ is a fibration and, as $\cS$ preserves prone morphisms, $\dot{p}_\bM$ is a fibration, as well. On both $\dot{p}_\bM$ and $\ddot{p}_\bM$ we have substitution tensors making them, by Theorem \ref{thm-em-Linton}, Kleisli and Eilenberg-Moore objects in $Fib_{/Set}$, respectively. Thus we can pass to the fibrations of monoids. The fibration $\mon(Gph(\bM,\otimes))\ra Set$ is the fibration of Lambek's multicategories.  The fibration $\mon(Gph(\bM,\otimes)^\cS)\ra Set$ is the fibration of symmetric multicategories (or coloured symmetric operads), the fiber $\mon(Gph(\bM,\otimes)^\cS)_O$ consists of $O$-coloured symmetric operads.  The fibration $\mon(Gph(\bM,\otimes)_\cS)\ra Set$ is the fibration of symmetric multicategories with the group actions free (or coloured rigid operads) cf. \cite{SZ3} or one-level multicategories with non-standard amalgamation, cf. \cite{HMP}. The fibration $\mon(End(Set))$ is the fibration of monads on slices of $Set$. The embeddings $\br$ and $\tilde{\br}$ are not so interesting as they are not full on isomorphism. The embeddings $\dot{\br}$, $\tilde{\dot{\br}}$, $\ddot{\br}$ and $\tilde{\ddot{\br}}$ are full on isomorphism and their images are polynomial functors and monads and analytic functors and monad on slices of $Set$, respectively. This was proved in \cite{Z2} in Sections 6 and 7.

As free multicategories exist, i.e. $\cU^\otimes$ has a left adjoint, this adjunction induces a monad $\cT^\otimes$ so that $\mon(Gph(\bM,\otimes))\ra Set$ is equivalent to the fibration of $\cT^\otimes$-algebras. As $\cS$ is monoidal, $T^\otimes$ distributes over $\cS$ and this distributive law
\[ \kappa: \cT^\otimes\cS\ra\cS\cT^\otimes\]
is what in \cite{BD} is called combing trees. The signatures in $\cT^\otimes\cS(A)$ consist of `term trees' built from operations from $A$, each decorated with a permutation of its entries. The signatures in $\cS\cT^\otimes(A)$ consist of `term trees' built from operations from $A$
and then the whole term tree is decorated by a suitable (big!) permutation. Thus in passing from  $\cT^\otimes\cS(A)$ to $\cS\cT^\otimes(A)$ we need to push up to the leaves permutation and compose them at the end. All this, but one-step combing hidden in the coherence morphism of the lax monoidal monad $\cS$, is done for us by the abstract theory.

As we noted in Section \ref{subsec-Act}, the above $\bkem$-diagram can be lifted to the lax slice $$\Mon_l(Fib_{/Set}{/_l\, p_{exp}:End(Set^\ra)\ra Set})$$ and hence, by adjunction, we have a $\bkem$-diagram in $$\Act_l\Mon_l(Fib_{/Set},\, p_{exp}:End(Set^\ra)\ra Set).$$ The actions along action $ev: End(Set)\times Set \lra Set$ are algebras for monads on slices of $Set$. The actions along actions of other monoidal categories are algebras for the suitable operads. The representations send algebras for operads to algebras for the corresponding monads.

\subsection*{The monads $\cR$ and $\cF$}

The similar analysis is possible for other submonads of $\cF$, including both $\cR$ and $\cF$. On the signatures and monoids side we `enrich' the above context by allowing actions on operations to be either surjections or all functions between finite sets. Such actions have a very simple interpretation.  If $a$ is an $n$-ary operation and $\sigma:(n]\ra(m]$ is a function (surjection), then $(\sigma\cdot a)(x_1,\ldots,x_m)=a(x_{\sigma(1)},\ldots,x_{\sigma(n)})$. Thus such operations are natural operations on terms and here we just require that our sets of operations are closed with respect to these operations. The categories of monoids are the categories of multicategories  that admit actions of surjections or all finite functions, respectively. The multicategories from $\mon(Gph(\bM)^\cR)$ are represented in $End(Set^\ra)$ as finitary semi-analytic monads on slices of $Set$ and the homomorphism between them as semi-cartesian morphism of monads, i.e. natural transformation such that the naturality squares for monos are pullbacks. This was proved in \cite{SZ4} for multicategories with no colours\footnote{This is what we called `toy model' in the introduction.} (that is in the fiber over 1) but that proof can be easily extended to all fibers the way it was done in case of monad $\cS$ in \cite{Z2}.  The multicategories in $\mon(Gph(\bM)^\cF)$ are represented in $End(Set^\ra)$ as finitary monads on slices of $Set$. We leave for the reader to work out the remaining details of these examples.

%

\subsection{Fibrations over $Cat$}

The considerations from the previous section can be extended to fibrations over $Cat$ or $Gpd$, the categories of small categories or small groupoids, respectively. We concentrate on fibrations over $Cat$ but everything restricts to $Gpd$ and there the theory works even better, cf. \cite{Fi2}.

As $Cat$ has pullbacks, we can consider cartesian monads on $Cat$ and look for representations of graphs over $Cat$ as we did over $Set$. But, as the endofunctors on slices of $Cat$ are not so interesting as the endofunctors on presheaf categories, we will describe a refinement of the notion of Burroni fibration that treats spans over $Cat$ in a more subtle way. For this to work all the components of  both unit and multiplication of the rc cartesian monad $\bT$ on $Cat$ need to be discrete bifibrations. Moreover, we require that $\bT$ preserves discrete two-sided fibrations.

For such a monad $\bT$, we can refine the notion of a $\bT$-graph by asking the span
 \begin{center} \xext=800 \yext=430
\begin{picture}(\xext,\yext)(\xoff,\yoff)
\settriparms[1`1`0;400]
  \putAtriangle(0,0)[A`O`\bT(O);\gamma`\delta`]
\end{picture}
\end{center}
to be a two-sided discrete fibration (with $\gamma$ a fibration and $\delta$ opfibration), cf. \cite{Ri}. We write $(A,\gamma,\delta)$ or just $A$ to denote such a  two-sided discrete fibration. For $o\in O$ and $\vec{o}\in \bT(O)$, we write $A(o;\vec{o})$ for the fiber of $A$  over both $o$ and $\vec{o}$.

The fiber of the fibration $p_\bT:\Gph_{2fi}(\bT)\ra Cat$ over a category $O$ is equivalent to the category of $Set$-valued functors over the category $O^{op}\times \bT(O)$. Such a fibration is still a lax monoidal fibration (with fibers being monoidal categories) but the tensor has to be fitted into this context as follows. Let $A$ and $B$ be two $\bT$-graphs as above.  Applying $\bT$ and $\mu$ to the two-sided fibration $B$, we get a two-sided fibration
 \begin{center} \xext=800 \yext=430
\begin{picture}(\xext,\yext)(\xoff,\yoff)
\settriparms[1`1`0;400]
  \putAtriangle(0,0)[\bT(B)`\bT(O)`\bT(O);\bT(\gamma)`\mu_O\circ \bT(\delta)`]
\end{picture}
\end{center}
denoted $(\bT(B),\bT(\gamma),\mu_O\circ \bT(\delta))$ or just $\bM(B)$. With this notation the fiber over both $o$ and $\vec{o}$ of the tensor product $A\otimes_O B$ is given by the formula
\[ A\otimes B(o;\vec{o}) =\int^{\vec{p}\in \bT(O)} A(o;\vec{p})\times \bT(B)(\vec{p};\vec{o}) \]

In this context the fibration $\pi:\DFib\ra Cat$ of small discrete fibrations over $Cat$, whose morphisms commute over base and preserve prone morphisms, plays the role of the basic fibration over $Cat$. Note that both $p_\bT$ and $\pi$ have reflexive coequalizers in fibers that are preserved by reindexing functor. Thus they are both rc-fibrations even in $\Fib(Cat)$.

The action of $p_\bT$ on $\pi$ \begin{center} \xext=800 \yext=620
\begin{picture}(\xext,\yext)(\xoff,\yoff)
 \settriparms[1`1`1;500]
  \putVtriangle(0,0)[\Gph_{2fi}(\bT)\times_{Cat}\DFib`\DFib`Cat;\star^{\bT}``cod]
\end{picture}
\end{center}
is also given by a coend formula. The fiber of the action of $A$ on a fibration $d:X\ra O$ over $o\in O$ is
\[ A\star X(o) = \int^{\vec{p}\in \bT(O)} A(o;\vec{p})\times X^{\vec{p}}  \]
where
\[ X^{\vec{p}}= \prod_{i=1}^{|\vec{p}|} X(p_i) \]
and $|\vec{p}|$ is the length of the vector $\vec{p}$.

Now we shall specify $\bT$ to be the monad for strict monoidal categories $\bM$. Note that $\bM$ is a lift of the (cartesian) free monoid monad from $Set$ to $Cat$. This is why we denote this monad by the same letter $\bM$. The monad $(\bM,\eta,\mu)$ has all the necessary features i.e. it is rc, cartesian, preserves two-sisded discrete fibrations, and components of both $\eta$ and $\mu$ are bifibrations. This implies that $p_\bM:\Gph_{2fi}(\bM)\ra Cat$ is a monoidal fibration with action defined on the basic fibration
\begin{center} \xext=800 \yext=620
\begin{picture}(\xext,\yext)(\xoff,\yoff)
 \settriparms[1`1`1;500]
  \putVtriangle(0,0)[\Gph_{2fi}(\bM)\times_{Cat}\DFib`\DFib`Cat;\star^{\bM}``cod]
\end{picture}
\end{center}
as above. It has an exponential adjoint in $CAT/Cat$
\begin{center} \xext=800 \yext=680
\begin{picture}(\xext,\yext)(\xoff,\yoff)
 \settriparms[1`1`1;500]
  \putVtriangle(0,30)[\Gph_{2fi}(\bM)`End(\DFib)`Cat;\br=rep_{\bM}`p_{\bM}`\pi_{exp}]
\end{picture}
\end{center}
for which we have

\begin{proposition}\label{prop-rep-cat-fghw1} The functor $\br$ defined above is strong monoidal and is a morphism of bifibrations.  \end{proposition}
{\it Proof.}~  Let $d: X\ra O$ be a (small) discrete fibration. The unit of the tensor in the fiber over $O$ is $(O,1_O,\eta_O)$. Clearly, we have
\[ \br(O)(X)(o) = O\star X (o) \cong X(o).\]
For two $\bM$-graphs $A$ and $B$ we have
\[ r(A\otimes B)(X)(o) = ((A\otimes B)\star X)(o) \cong\]

\[ \cong \int^{\vec{p}\in \bM(O)}(\int^{\vec{q}\in \bM(O)} A(o;\vec{q})\times \bM(B)(\vec{q};\vec{p}) )\times X^{\vec{p}} \cong\]

\[ \cong \int^{\vec{q}\in \bM(O)} A(o;\vec{q})\times (\int^{\vec{p}\in \bM(O)} \bM(B)(\vec{q};\vec{p}) )\times X^{\vec{p}} \cong\]

\[ \cong \int^{\vec{q}\in \bM(O)} A(o;\vec{q})\times \prod_{i=1}^{|\vec{q}|} (\int^{\vec{p^i}\in \bM(O)} B(q_i;\vec{p^i}) )\times X^{\vec{p^i}} \cong\]

\[ \cong \int^{\vec{q}\in \bM(O)} A(o;\vec{q})\times (B\star X)^{\vec{q}} \cong\]

\[ \cong (A\star (B\star X))(o)= r(A)\circ\br(B)(X)(o). \]
Thus $\br$ is strong monoidal.

We shall check that the codomains of supine morphisms are preserved by $\br$. The rest is left for the reader. Let $u:O\ra Q$ be a functor in $Cat$, $A$ an $\bM$-graph over $O$ and $X$ a discrete fibration over $Q$. We have

\[ (u_!(A)\star X)(q) \cong \]

\[ \cong \int^{\vec{q}\in \bM(Q)} u_!(A)(q;\vec{q})\times X^{\vec{q}} \cong \]

\[ \cong \int^{\vec{q}\in \bM(Q)} \int^{p\in O, \vec{p}\in \bM(O)} A(p;\vec{p}) \times Q(q,u(p))\times \bM(Q)(u(\vec{p}),\vec{q})\times X^{\vec{q}}\cong^* \]

\[ \cong \int^{p\in O} \int^{ \vec{p}\in \bM(O)} A(p;\vec{p})\times X^{u(\vec{p})} \times Q(q,u(p))  \cong \]

\[ \cong \int^{p\in O} (A\star u^*(X))(p)\times Q(q,u(p)) \cong \]

\[ \cong u_!(A\star u^*(X))(q). \]

The isomorphism $\cong^*$ uses the fact that $A$ is a two-sided fibration. Going down through $\cong^*$ we send the element represented by a quadruple
\[ \lk a\in A(p;\vec{p}), v\in Q(q,u(p)), \vec{v} \in \bM(O)(\vec{p},\vec{q}), \vec{x}\in X^{\vec{q}} \rk \]
to the element represented by the triple
\[ \lk a\in A(p;\vec{q}), v\in Q(q,u(p)),  \vec{v}^*(\vec{x})\in X^{u(\vec{p})} \rk. \]
Going back, we send the element represented by the triple
\[ \lk a\in A(p;\vec{q}), v\in Q(q,u(p)),  \vec{y}\in X^{u(\vec{p})} \rk \]
to the element represented by a quadruple
\[ \lk a\in A(p;\vec{p}), v\in Q(q,u(p)), 1_{u(\vec{p})} \in \bM(O)(u(\vec{p}),u(\vec{p})), \vec{y}\in X^{u(\vec{p})} \rk. \]

As above, to check that $\br$ preserves prone morphisms, we only check that the domains agree. Let $u:O\ra Q$ be a functor in $Cat$, $B$ an $\bM$-graph over $Q$ and $Y$ a discrete fibration over $Q$. We have

\[ u^*\br(B)u_!u^*(Y)_o \cong \]

\[ \cong  B\star u_!u^*(Y)_{u(o)} \cong \]

\[ \cong \int^{\vec{q}\in \bM(Q)} B(u(o);\vec{q}) \times (u_!u^*(Y)^{\vec{q}} \cong \]

\[ \cong \int^{\vec{q}\in \bM(Q)} B(u(o);\vec{q}) \times \prod_{i=1}^{|\vec{q}|} \int^{o_i\in O} Y(u(o_i))\times Q(q_i,u(o_i)) \cong \]

\[ \cong \int^{\vec{q}\in \bM(Q)} B(u(o);\vec{q}) \times  \int^{\vec{o}\in \bM(O)} Y^{u(\vec{o})}\times \bM(Q)(\vec{q},u(\vec{o}))  \cong \]

\[ \cong \int^{\vec{q}\in \bM(Q)}   \int^{\vec{o}\in \bM(O)} B(u(o);\vec{q}) \times Y^{u(\vec{o})}\times \bM(Q)(\vec{q},u(\vec{o}))  \cong \]

\[ \cong  \int^{\vec{o}\in \bM(O)}  \int^{\vec{q}\in \bM(Q)}  B(u(o);\vec{q})\times \bM(Q)(\vec{q},u(\vec{o}))  \times Y^{u(\vec{o})} \cong^* \]

\[ \cong \int^{\vec{o}\in \bM(O)} B(u(o);u(\vec{o})) \times Y^{u(\vec{o})} \cong \]

\[ \cong u^*(B)\star u^*(Y)(o) \]
where the isomorphism $\cong^*$ uses the fact that $B$ is a two-sided fibration.  Going down through $\cong^*$, we send the element represented by a triple
\[ \lk b\in B(u(o);u(\vec{p})), \vec{v} \in \bM(O)(\vec{q},u(\vec{o})), \vec{y}\in Y^{u(\vec{o})} \rk \]
to the element represented by the pair
\[ \lk v_!(b)\in B(u(o);u(\vec{o})),  \vec{y}\in Y^{u(\vec{o})} \rk. \]
Going back, we send the element represented by the triple
\[ \lk b\in B(u(o);u(\vec{o})), \vec{y}\in Y^{u(\vec{o})}  \rk \]
to the element represented by the triple
\[ \lk b\in B(u(o);u(\vec{o})), 1_{u(\vec{o})} \in \bM(O)(u(\vec{o}),u(\vec{o})), \vec{y}\in Y^{u(\vec{o})} \rk. \]
$\Box$

We have a monoidal embedding
$$\iota=\iota_O:\bM(O)\ra \widehat{O}\simeq\DFib(O)$$
into the category of presheaves over $O$, sending $\vec{o}=\lk o_1,\ldots, o_n\rk$ to the coproduct $Y(o_1)+\ldots Y(o_n)$, where $Y:O^{op}\ra \widehat{O}$ is the Yoneda embedding.
The monoidal structure on $\DFib(O)$ comes from the binary coproducts. Now if we denote by $\sM$ the monad for strict symmetric monoidal categories on $Cat$, then we have a factorization of the monoidal embeding $\iota_O$ into two embeddings
\[ \bM(O)\stackrel{\bu}{\lra} \sM(O) \lra\widehat{O} \]
with the first one being identity  on objects and the second one being full on isomorphisms and symmetric monoidal, as well.

$\br$ has an rc right adjoint $\bU$ such that for a functor $H:\DFib(O)\ra\DFib(O)$, $o\in O$ and $\vec{o}\in \bM(O)$ and the fiber over $(o;\vec{o})$ of the two-sided span $\bU(H)$ is given by

\[ \bU(H)(o;\vec{o}) =  H(\iota(\vec{o}))(o). \]
In a more conceptual way the two-sided discrete fibration
\begin{center} \xext=800 \yext=430
\begin{picture}(\xext,\yext)(\xoff,\yoff)
\settriparms[1`1`0;400]
  \putAtriangle(0,0)[U(H)`O`\bM(O);\gamma`\delta`]
\end{picture}
\end{center}
corresponds to a functor
\[ \overline{H}: O^{op}\times \bM(O)\lra Set \]
which is an adjoint to the composite functor
\[ \bM(O)\stackrel{\iota}{\lra} \widehat{O}\stackrel{H}{\lra} \widehat{O} \]

The unit of the adjunction
\begin{center} \xext=2500 \yext=150
\begin{picture}(\xext,\yext)(\xoff,\yoff)
 \putmorphism(0,50)(1,0)[A(o;\vec{o})`\int^{\vec{q}\in \bM(O)}A(o;\vec{q})\times \widehat{O}(\iota(\vec{q}),\iota(\vec{o}))=\bU\br(A)(o;\vec{o})`(\eta_A)_{o;\vec{o}}]{2000}{1}a
\end{picture}
\end{center}
is given by
\[ a\mapsto [a, 1_{\iota(\vec{o})}] \]
and the counit of the adjunction
\begin{center} \xext=1600 \yext=150
\begin{picture}(\xext,\yext)(\xoff,\yoff)
 \putmorphism(200,50)(1,0)[\br\bU(H)(X)_o =(\bU(H)\star X)(o)=\int^{\vec{p}\in \bM(O)}H(\iota(\vec{p}))_o\times X^{\vec{p}}`H(X)_o`((\varepsilon_H)_X)_{o}]{2200}{1}a
\end{picture}
\end{center}
is given by
\[ [h, \vec{x}:\iota(\vec{p})\ra X] \mapsto H(\vec{x})(h). \]

As $\br$ preserves supine morphisms, we have just shown

\begin{proposition}\label{prop-rep-cat-fghw2} The functor $\br$ defined above has an rc right adjoint $\bU$. $\Box$  \end{proposition}

In this way we get an rc-monoidal monad $\cF$ on $p_\bM:\Gph_{2fi}(\bM)\ra Cat$

\[ \cF(A)(o;\vec{o})= \int^{\vec{q}\in\bM(O)} A(o;\vec{q})\times \widehat{O}(\iota(\vec{q}),\iota(\vec{o})) \]
that has similar submonads as the monad on signatures, coming from the fact that the sets of morphisms $\widehat{O}(\iota(\vec{q}),\iota(\vec{o}))$
can be replaced by some suitable subsets as long as they are closed under compositions with the morphisms coming from $\bM(O)$. If we restrict the hom-sets $\widehat{O}(\iota(\vec{q}),\iota(\vec{o}))$ to $\sM(O)(\vec{q},\vec{o})$, we get the symmetrization monad

\[ \cS(A)(o;\vec{o})= \int^{\vec{q}\in\bM(O)} A(o;\vec{q})\times\sM(O)(\vec{q},\vec{o}) \]
that is giving rise to the notion of an analytic (endo)functor on presheaf categories, cf. \cite{FGHW}, in the sense that

\begin{proposition}\label{prop-rep-cat-fghw3} The objects in the essential image of the representation
\begin{center} \xext=800 \yext=600
\begin{picture}(\xext,\yext)(\xoff,\yoff)
 \settriparms[1`1`1;460]
  \putVtriangle(0,30)[\Gph_{2fi}(\bM)^{\cS}`End(\DFib)`Cat;\ddot{\br}`\ddot{p}_{\bM}`\pi_{exp}]
\end{picture}
\end{center}
consists of analytic endofunctors in the sense of \cite{FGHW}.
\end{proposition}
{\it Proof.}~  Let us fix a small category $O$ and let $\bu:\bM(O)\ra \sM(O)$ denote the embedding functor. We first note that we have an equivalence of categories
\[ Nat(\sM(O),\widehat{O}) \lra (\Gph_{2fi}^\cS )_O\]
sending a functor $F:\sM(O)\ra \widehat{O}$ to an $\cS$-agebra $(F\bu,\phi)$ where
\begin{center} \xext=1600 \yext=150
\begin{picture}(\xext,\yext)(\xoff,\yoff)
 \putmorphism(200,50)(1,0)[ \cS(F\bu)(\vec{o})=\int^{\vec{p}\in \bM(O)} F\bu(\vec{p})\times \sM(O)(\vec{p},\vec{o}))`F\bu(\vec{o})`\phi]{1800}{1}a
\end{picture}
\end{center}
is given by
\[ [a,(\sigma,\vec{v})] \mapsto F(\sigma,\vec{v})(a).\]
We need to show that the triangle
\begin{center} \xext=1200 \yext=650
\begin{picture}(\xext,\yext)(\xoff,\yoff)
  \putmorphism(0,550)(0,-1)[Nat(\sM(O),\widehat{O})`\Gph_{2fi}(\bM)^\cS_O`]{500}{1}l
   \put(400,550){\vector(3,-1){550}}
    \put(600,500){\makebox(100,100){$Lan_\bu$}}
  \put(400,50){\vector(3,1){550}}
  \put(600,00){\makebox(100,100){$\ddot{\br}$}}

  \put(1100,250){\makebox(100,100){$End(\widehat{O})$}}
\end{picture}
\end{center}
commutes up to an isomorphism.

For $X$ in $\widehat{O}$, we have
\[ Lan_\bu(F)(X) = \int^{\vec{p}\in \sM(O)} F(\vec{p})\times X^{\vec{p}}. \]

The representation $\ddot{\br}$ on an $\cS$-agebra $(F\bu,\phi)$ is given by a coequalizer
\begin{center} \xext=2000 \yext=350
\begin{picture}(\xext,\yext)(\xoff,\yoff)
  \putmorphism(0,150)(1,0)[\cS(F\bu)\star X`F\bu\star X`]{1200}{0}l
  \putmorphism(0,100)(1,0)[\phantom{\cS(F\bu)\star X}`\phantom{F\bu\star X}`\phi\star 1_X]{1200}{1}b
   \putmorphism(0,200)(1,0)[\phantom{\cS(F\bu)\star X}`\phantom{F\bu\star X}`\zeta_{F\bu,X}]{1200}{1}a
    \putmorphism(1200,150)(1,0)[\phantom{F\bu\star X}`F\bu\ddot{\star} X`]{800}{1}a
\end{picture}
\end{center}
where for a functor $A:\bM(O)\ra \widehat{O}$, the morphism $\zeta_{A,X}$ is
\begin{center} \xext=2000 \yext=150
\begin{picture}(\xext,\yext)(\xoff,\yoff)
  \putmorphism(0,50)(1,0)[\cS(A)\star X`\cF(A)\star X`]{800}{1}a
  \putmorphism(800,50)(1,0)[\phantom{\cF(A)\star X}`A\star X`\varepsilon_{\br(A)}]{800}{1}b
\end{picture}
\end{center}
such that
\[ [a,(\sigma,\vec{v}),\lk x_i\rk_i] \hskip 5mm \mapsto \hskip 5mm [a, \lk X(v_i)(x_{\sigma(i)})\rk_{i}] \]
and hence the above coequalizer is the coequalizer of the following parallel pair
\begin{center} \xext=2000 \yext=350
\begin{picture}(\xext,\yext)(\xoff,\yoff)
  \putmorphism(200,150)(1,0)[\int^{\vec{p}\in \bM(O)}\int^{\vec{q}\in \bM(O)} F\bu(\vec{q})\times \sM(O)(\vec{q},\vec{p})\times X^{\vec{p}}`\int^{\vec{p}\in \bM(O)}F\bu(\vec{p})\star X`]{2200}{0}a
   \putmorphism(200,200)(1,0)[\phantom{\int^{\vec{p}\in \bM(O)}\int^{\vec{q}\in \bM(O)} F\bu(\vec{q})\times \sM(O)(\vec{q},\vec{p})\times X^{\vec{p}}}`\phantom{\int^{\vec{p}\in \bM(O)}F\bu(\vec{p})\star X}`\zeta_{F\bu,X}]{2200}{1}a
  \putmorphism(200,100)(1,0)[\phantom{\int^{\vec{p}\in \bM(O)}\int^{\vec{q}\in \bM(O)} F\bu(\vec{q})\times \sM(O)(\vec{q},\vec{p})\times X^{\vec{p}}}`\phantom{\int^{\vec{p}\in \bM(O)}F\bu(\vec{p})\star X}`\phi\star 1_X]{2200}{1}b
\end{picture}
\end{center}
i.e. it is the coequalizer formula for the coend defining $Lan_\bu(F)(X)$. The remaining details are left for the reader. $\Box$

{\em Remarks.}

\begin{enumerate}
  \item The essential image of the representation of the Kleisli monoidal fibration for the symmetrization monad $\cS$
\begin{center} \xext=800 \yext=600
\begin{picture}(\xext,\yext)(\xoff,\yoff)
 \settriparms[1`1`1;460]
  \putVtriangle(0,30)[\Gph_{2fi}(\bM)_{\cS}`End(\DFib)`Cat;\dot{\br}`\dot{p}_{\bM}`\pi_{exp}]
\end{picture}
\end{center}
gives rise to a (yet another) notion of polynomial functor on presheaf categories. As these functors are finitary, they are more restrictive then those considered in \cite{GK}. The relation with the discrete generalized polynomial functors of \cite{Fi1} will be studied elsewhere.
  \item The tensor product $\otimes$ on fibration $p_{\bM}$ preserves binary coproducts in left variable and filtered colimits. Hence (cf. \cite{A}, \cite{Ke1}, \cite{BJT}, \cite{Le}), there is a free monoid 1-cell $\cF_{\otimes}$ left adjoint to the forgetful 1-cell $\cU^\otimes: \Mon(\Gph_{2fi}(\bM))\lra \Gph_{2fi}(\bM)$ over $Cat$. The adjunction induces the free monoid  monad $\cT^\otimes$ on the fibration $p_{\bM}$. As the symmetrization is monoidal, we have a distributive law
      \[ \bc : \cT^\otimes \circ \cS \lra \cS \circ T^\otimes \]
      that is `combing trees', cf. \cite{BD}. It pushes through the term tree the permutations of entries of each function symbol to a one big permutation of the level of leaves of the whole tree.
\end{enumerate}

\subsection{The fibrations over $\o$-graphs}

The category of $\o$-categories $\oCat$ is monadic over the category of $\o$-graphs $\oGph$ and let $\cT$ be the resulting monad. The monad $\cT$ is cartesian and we could develop a similar theory in this case for the fibration of $\cT$-graphs $p_\cT: Gph(\cT)\lra \oGph$. However, such a theory is less interesting as the (faithful)  representation morphism
\begin{center} \xext=800 \yext=700
\begin{picture}(\xext,\yext)(\xoff,\yoff)
 \settriparms[1`1`1;500]
  \putVtriangle(0,30)[Gph(\cT)`End(\oGph^\ra)`\oGph;\br=rep_{\cT}`p_{\cT}`p_{exp}]
\end{picture}
\end{center}
is already full on isomorphisms. In fact, a natural transformation $\tau : rep_{\cT}(A)\ra rep_{\cT}(B)$ comes from a morphism in $Gph(\cT)$ iff it preserves pullbacks.  The theory becomes more interesting when we consider the subcategory $\oGph_0$ of those $\o$-graphs that have one 0-cell. The monad $\cT$ restricts to this subcategory and the representation $rep_{\cT_0}$ is not full on isomorphisms anymore. We leave the reader to work out this example in detail.

\section{Appendix}
We finish the paper by making two additional points concerning this story. One is explaining that the Burroni fibrations of interest are very special: they are cartesian bifibrations with the adjoint representation to the tautologous action satisfying Frobenius reciprocity. The other is about lifting actions from monoids to the algebras when we have an action on action along action.

\subsection{Cartesian bifibrations}\label{sec-coslices-}

The Burroni fibrations are more than just fibrations, they are cartesian bi-fibrations. This makes them even more convenient to work with in practice. In this section we briefly discuss this notion.

A functor $p:\cE\ra \cB$ is a {\em cartesian bifibration}, cf. \cite{Z2}, iff it is fibration with fibers having pullbacks, reindexing functors $u^*$ preserving pullbacks, for $u:B\ra B'$ in $\cB$, and having right adjoint $u_!$, so that the unit and counit of the adjunction $u^*\dashv u_!$ are cartesian.

{\em Examples.} If $T$ is a cartesian monad on a category $\bB$ with finite products, then the fibration $p_T:Gph(T)\ra \cB$ is a cartesian bifibration. Moreover, the tautologous actions
\begin{center} \xext=800 \yext=620
\begin{picture}(\xext,\yext)(\xoff,\yoff)
 \settriparms[1`1`1;500]
  \putVtriangle(0,0)[Gph(T)\times_\cB \cB^\ra`\cB^\ra`\cB;\star``cod]
\end{picture}
\end{center}
satisfy kind of Frobenius reciprocity, i.e. for $u:O\ra Q$ and $A$ in $Gph(T)_O$ and  $Y$ in $\cB^\ra_{O}$, we have canonical isomorphisms
\[ u_!(A\star u^*(Y))\lra u_!(A)\star Y. \]
See \cite{Z2} Section 5.2 for more.

The following fact is due to T. Streicher \cite{Str2}.

\begin{theorem}
Let $\cB$ be a category with finite limits. A functor $p:\cE\ra \cB$ is a cartesian bi-fibration with the terminal object iff it is of form $\cC\da F$ for some terminal object preserving functor $F:\cB\ra \cC$ between categories with finite limits.
\end{theorem}

{\it Proof.}~ See \cite{Str2}. $\Box$

\subsection{Actions on actions along actions}

Monoidal endo-1-cells act on actions, i.e. denoting by  $\End_{l,\cA}(\cC,\otimes)$ the category of lax monoidal endo 1-cells on a monoidal object $(\cC,\otimes)$ and by $\Act_{l,\cA}(\cC,\otimes,\cX)$ the category of actions of $(\cC,\otimes)$ on a 0-cell $\cX$ in $\cA$, we have an action
\[ \sharp : \End_{l,\cA}(\cC,\otimes)\times \Act_{l,\cA}(\cC,\otimes,\cX) \lra \Act_{l,\cA}(\cC,\otimes,\cX) \]
\[  \lk \cF,(-)\star(=)\rk \mapsto  \cF(-)\star(=)\]

Let $\br$ be a strong monoidal representation 1-cell in $\cA$ with a lax right adjoint $\bU$, $\cF$ the lax monoidal monad arising from this adjunction,
\begin{center} \xext=1600 \yext=350
\begin{picture}(\xext,\yext)(\xoff,\yoff)
\putmorphism(400,110)(1,0)[\cC`\cX^\cX`]{1200}{0}a
 \putmorphism(400,140)(1,0)[\phantom{\cC}`\phantom{\cX^\cX}`U]{1200}{-1}a
 \putmorphism(400,80)(1,0)[\phantom{\cC}`\phantom{\cX^\cX}`\br]{1200}{1}b

  \put(90,100){\oval(100,100)[l]}
  \put(90,150){\line(1,0){250}}
  \put(90,50){\vector(1,0){250}}
  \put(0,200){${(\cF,\eta,\mu)}$}
 \end{picture}
\end{center}
with $\varepsilon$ the counit of this adjunction. Let $\star:\cC\times\cX\ra \cX$ be an action adjoint to $\br$. Then for any $\bA$ in $\cA$, we have an action
\[ \varepsilon_{\br(\bA)} : \cF(\bA)\star(=)\lra A\star(=) \]
of the action of $\bA$ along the action $\sharp$. If $(\bM,\bm,\be)$ is a monoid in $\cC$ and $\alpha:\bM\star \bX \ra \bX$ is an action along $\star$, then
\begin{center} \xext=1600 \yext=250
\begin{picture}(\xext,\yext)(\xoff,\yoff)
\putmorphism(0,50)(1,0)[\cF(\bM)\star\bX`\bM\star\bX`\varepsilon_{\br(\bM)_\bX}]{1000}{1}a
 \putmorphism(1000,50)(1,0)[\phantom{\bM\star\bX}`\bX`\alpha]{1200}{1}a
 \end{picture}
\end{center}
is an action of a monoid $\cF(\bM)$ on $\bX$ along $\star$.

It is instructive to see how such an action looks like in case $\cC$ is the fibration of signatures and $\cF$ is the symmetrization monad. We leave this to the reader.

\end{document}